\documentclass[numbook]{svjour3}

\smartqed 

\usepackage{amsmath,amssymb,bbold,dsfont,amstext,amsfonts}
\usepackage[numbers]{natbib}
\usepackage[latin1]{inputenc}   
\usepackage[T1]{fontenc} 

\numberwithin{equation}{section}
\newcommand{\Gammabs}{\boldsymbol{\Gamma}}
\newcommand{\Upsilonbs}{\boldsymbol{\Upsilon}}
\newcommand{\Deltabs}{\boldsymbol{\Delta}}

\newcommand{\bs}{\boldsymbol}
\spnewtheorem{assum}{Assumption}[section]{\bf}{\rm}

\journalname{Journal of Theoretical Probability}

\begin{document}

\title{On the almost sure location of the singular values of certain Gaussian block-Hankel large random matrices \thanks{This work was supported by Project ANR-12-MONU-0003 DIONISOS}}

\subtitle{Location of the singular values of block-Hankel random matrices}


\author{Philippe Loubaton
}


\institute{P. Loubaton \at
              Université Paris-Est, 
              Laboratoire d'Informatique Gaspard Monge, UMR CNRS 8049, 
              5 Bd. Descartes, Cité Descartes, Champs sur Marne,
              Marne la Vallée 77454 Cedex 2,
              Tel.: 33-1-60-95-72-93\\
              Fax:  33-1-60-95-72-55\\
              \email{loubaton@univ-mlv.fr}           
}

\date{Received: date / Accepted: date}

\maketitle

\begin{abstract}

This paper studies the almost sure location of the eigenvalues of matrices ${\bf W}_N {\bf W}_N^{*}$
where ${\bf W}_N = ({\bf W}_N^{(1)T}, \ldots, {\bf W}_N^{(M)T})^{T}$ is a $ML \times N$ block-line matrix 
whose block-lines $({\bf W}_N^{(m)})_{m=1, \ldots, M}$ are independent identically distributed $L \times N$ Hankel matrices 
built from i.i.d. standard complex Gaussian sequences. It is shown that if $M \rightarrow +\infty$ and 
$\frac{ML}{N} \rightarrow c_*$ ($c_* \in (0, \infty)$), then the empirical eigenvalue distribution of 
${\bf W}_N {\bf W}_N^{*}$ converges almost surely towards the Marcenko-Pastur distribution. More importantly, 
it is established using the Haagerup-Schultz-Thorbjornsen ideas that if $L = \mathcal{O}(N^{\alpha})$ with $\alpha < 2/3$, then, almost surely, for $N$ large enough, the eigenvalues of ${\bf W}_N {\bf W}_N^{*}$ are located in the neighbourhood of the Marcenko-Pastur distribution. It is 
conjectured that the condition $\alpha < 2/3$ is optimal.

\end{abstract}

\section{Introduction}
\subsection{The addressed problem and the results}
In this paper, we consider independent identically distributed zero mean complex valued Gaussian random variables $(w_{m,n})_{m=1, \ldots, M, n=1, \ldots, N+L-1}$ such that 
$\mathbb{E}| w_{m,n} |^{2} = \frac{\sigma^{2}}{N}$ and $\mathbb{E}(w_{m,n}^{2}) = 0$ where $M,N,L$ are integers. We define the $L \times N$ matrices 
$({\bf W}_N^{(m)})_{m=1, \ldots, M}$ as the Hankel matrices whose entries are given by 
\begin{equation}
\label{eq:def-Wmij}
\left({\bf W}^{(m)}_{N}\right)_{i,j} = w_{m,i+j-1}, \; 1 \leq i \leq L, 1 \leq j \leq N
\end{equation}
and ${\bf W}_N$ represents the  $ML \times N$ matrix 
\begin{equation}
\label{eq:eq:def-W}
{\bf W}_N = \left( \begin{array}{c} {\bf W}_N^{(1)} \\   {\bf W}_N^{(2)} \\ \vdots \\ {\bf W}_N^{(M)} \end{array} \right)
\end{equation}
In this paper, we establish that:
\begin{itemize}
\item the eigenvalue distribution
of $ML \times ML$ matrix ${\bf W}_N {\bf W}_N^{*}$ converges towards the Marcenko-Pastur distribution when $M \rightarrow +\infty$ and 
when $ML$ and $N$ both converge towards $+\infty$ in such a way that $c_N = \frac{ML}{N}$ satisfies  $c_N \rightarrow c_*$ where $0 < c_* < +\infty$
\item more importantly, that if $L = \mathcal{O}(N^{\alpha})$ with $\alpha < 2/3$, then, almost surely, for $N$ large enough, the eigenvalues ${\bf W}_N {\bf W}_N^{*}$ are located 
in the neighbourhood of the support of the Marcenko-Pastur distribution. 
\end{itemize}
\subsection{Motivation}
This work is mainly motivated by detection/estimation problems of certain multivariate time series. Consider a $M$--variate time series 
$({\bf y}_n)_{n \in \mathbb{Z}}$ given by 
\begin{equation}
\label{eq:modele-y}
{\bf y}_n = \sum_{p=0}^{P-1} {\bf a}_p s_{n-p} + {\bf v}_n = {\bf x}_n + {\bf v}_n
\end{equation}
where $(s_n)_{n \in \mathbb{Z}}$ represents a deterministic non observable scalar signal, $({\bf a}_p)_{p=0, \ldots, P-1}$ are deterministic unknown 
$M$--dimensional vectors and $({\bf v}_n)_{n \in \mathbb{Z}}$ represent i.i.d. zero mean complex Gaussian $M$--variate random vectors 
such that $\mathbb{E}({\bf v}_n {\bf v}_n^{*}) = \sigma^{2}  {\bf I}_M$ and $\mathbb{E}({\bf v}_n {\bf v}_n^{T}) = 0$ for each $n$.
The first term of the righthandside of (\ref{eq:modele-y}), that we denote by ${\bf x}_n$,  represents a "useful" non observable signal on which various kinds of informations
have to be infered from the observation of $N$ consecutive samples $({\bf y}_n)_{n=1, \ldots, N}$. Useful informations on $({\bf x}_n)$ may include:
\begin{itemize}
\item Presence versus absence of $({\bf x}_n)$, which is equivalent to a detection problem
\item Estimation of vectors $({\bf a}_p)_{p=0, \ldots, P-1}$
\item Estimation of sequence $(s_n)$ from the observations 
\end{itemize}
The reader may refer e.g. to \cite{vanderveen-talwar-97}, \cite{mou-duh-car-may}, \cite{vanderveen-vanderveen-98}, \cite{abe-mou-lou-1} for more information. 
A number of existing detection/estimation schemes are based on the eigenvalues and eigenvectors of matrix $\frac{{\bf Y}_L {\bf Y}_L^{*}}{N}$ 
where ${\bf Y}_L$ is the block-Hankel $ML \times (N-L+1)$ matrix defined by
$$
{\bf Y}_L = \left( \begin{array}{ccccc} {\bf y}_1 & {\bf y}_2 & \ldots & \ldots & {\bf y}_{N-L+1} \\
                                        {\bf y}_2 & {\bf y}_3 & \ddots &  \ddots & {\bf y}_{N-L+2} \\
                                           \vdots &  \ddots  &  \ddots  & \ddots & \vdots \\
                                        {\bf y}_L & {\bf y}_{L+1} & \ldots  & \ldots & {\bf y}_N \end{array} \right)
$$
and where $L$ is an integer usually chosen greater than $P$. 
We notice that matrix ${\bf Y}_L$ is the sum of deterministic matrix ${\bf X}_L$ and random matrix ${\bf V}_L$ both defined as ${\bf Y}_L$.  
The behaviour of the above mentioned detection/estimation schemes is easy to analyse when $ML$ is fixed and $N \rightarrow +\infty$
because, in this asymptotic regime, it holds that
$$
\| \frac{{\bf Y}_L {\bf Y}_L^{*}}{N} - \left( \frac{{\bf X}_L {\bf X}_L^{*}}{N} + \sigma^{2} {\bf I}_{ML} \right) \| \rightarrow 0
$$
where $\| {\bf A} \|$ represents the spectral norm of matrix ${\bf A}$. However, this asymptotic regime may be unrealistic because $ML$ and $N$ appear sometimes to be of the same order of magnitude. It is therefore 
of crucial interest to evaluate the behaviour of the eigenvalues of matrix $\frac{{\bf Y}_L {\bf Y}_L^{*}}{N}$ when $ML$ and $N$ converge to 
$+\infty$ at the same rate. Matrix ${\bf Y}_L = {\bf X}_L + {\bf V}_L$ can be interpreted as an Information plus Noise model (see \cite{dozier2007empirical})
but in which the noise and the information components are block-Hankel matrices. We believe that in order to understand the behaviour of the eigenvalues of 
$\frac{{\bf Y}_L {\bf Y}_L^{*}}{N}$, it is first quite useful to evaluate the eigenvalue distribution of the noise contribution, 
i.e. $\frac{{\bf V}_L {\bf V}_L^{*}}{N}$, and to check whether its eigenvalues tend to be located in a compact interval. Hopefully, 
the behaviour of the greatest eigenvalues of $\frac{{\bf Y}_L {\bf Y}_L^{*}}{N}$ may be obtained by adapting the approach of
\cite{benaych-rao-2}, at least if the rank of the "Information" component ${\bf X}_L$ is small enough w.r.t. $ML$. 

It is clear that if we replace $N$ by $N+L-1$ in the definition of matrix ${\bf V}_L$, matrix ${\bf W}_N$ is obtained from 
$\frac{{\bf V}_L}{\sqrt{N}}$ by row permutations. Therefore, matrices $\frac{{\bf V}_L {\bf V}_L^{*}}{N}$ and ${\bf W}_N {\bf W}_N^{*}$
have the same eigenvalues. The problem we study in the paper is thus equivalent to the characterization of the eigenvalue distribution
of the noise part of model ${\bf Y}_L$.   
\subsection{On the literature}
Matrix ${\bf W}_N$ can be interpreted as a block-line matrix with i.i.d. $L \times N$ blocks $({\bf W}_N^{m})_{m=1 \ldots, M}$. Such random 
block matrices have been studied in the past e.g. by Girko (\cite{girko-book}, Chapter 16) as well as in \cite{bryc-speicher-2006} in the Gaussian 
case. Using these results, it is easy to check that the eigenvalue distribution of ${\bf W}_N{\bf W}_N^{*}$ converges towards the Marcenko-Pastur 
distribution when $L$ is fixed. However, the case $L \rightarrow +\infty$ and the almost sure location of the eigenvalues of ${\bf W}_N{\bf W}_N^{*}$ around the support 
of the Marcenko-Pastur distribution cannot be addressed using the results of \cite{girko-book} and 
\cite{bryc-speicher-2006}. We note that the $L \times N$ blocks $({\bf W}^{m})_{m=1 \ldots, M}$
are Hankel matrices. We therefore also mention the works \cite{li-liu-wang-2011} and \cite{basu-bose-et-al-2012} 
that are equivalent to the study of the eigenvalue distribution of symetric $ML \times ML$ block matrices, each block being a Toeplitz or a Hankel $L \times L$ matrix built from i.i.d. (possibly non Gaussian) entries. When $L \rightarrow +\infty$ while $M$ remains fixed, it has been shown using the moments method that the eigenvalue distribution of the above matrices converge towards
a non bounded limit distribution. This behaviour generalizes the results of \cite{bryc-dembo-jiang-2006} obtained when $M=1$. When $M$ and $L$ both converge to $+\infty$, it is shown in \cite{basu-bose-et-al-2012} that the eigenvalue 
distribution converges towards the semi-circle law. We however note that the almost sure location of the eigenvalues 
in the neighbourhood of the support of the semi-circle law is not addressed in \cite{basu-bose-et-al-2012}. The behaviour of the singular value distribution of random block Hankel 
matrix (\ref{eq:eq:def-W}) was addressed in \cite{basak-bose-sen-2012} when $M=1$ and $\frac{L}{N} \rightarrow c_*$ but when the $w_{1,n}$ for $N < n <N+L$ are forced to $0$.
The random variables $w_{1,n}$ are also non Gaussian and are possibly dependent in \cite{basak-bose-sen-2012}.
It is shown using the moments method that the singular value distribution converges towards a non bounded limit distribution. The case of block-Hankel matrices where both $M$ and $L$ converge towards 
$\infty$ considered in this paper thus appears simpler because we show that the 
eigenvalue distribution of ${\bf W}_N {\bf W}_N^{*}$ converges towards the Marcenko-Pastur distribution. This behaviour is not surprising in view of the convergence towards the semi-circle law proved in \cite{basu-bose-et-al-2012} when both the number and 
the size of the blocks converge to $\infty$. As mentioned above, the main 
result of the present paper concerns the almost sure location of the eigenvalues of  ${\bf W}_N{\bf W}_N^{*}$ around the support of the Marcenko-Pastur distribution
under the extra-assumption that $L = \mathcal{O}(N^{\alpha})$ with $\alpha < 2/3$. This kind of result is known for a long time for $L=1$
in more general conditions (correlated non Gaussian entries, see e.g. \cite{bai-silverstein-book} and the references therein).
Haagerup and Thorbjornsen introduced in \cite{haagerupnew2005} an efficient approach to address these issues in the context of random matrices built on non 
commutative polynomials of complex Gaussian matrices. The approach of \cite{haagerupnew2005} has been generalized to the real Gaussian case in \cite{schultz-2005}, and used in 
\cite{capitaine2009largest}, \cite{lou-val-2011}, \cite{capitaine2011} to address certain non zero mean random matrix models. We also mention that 
the results of \cite{haagerupnew2005} have been recently generalized in \cite{male-2012} to polynomials of complex Gaussian random matrices and deterministic matrices. \\

To our best knowledge, the existing literature does not allow to prove that the eigenvalues of ${\bf W}_N {\bf W}_N^{*}$ are located 
in the neighbourhood of the bulk of the Marcenko-Pastur distribution. We finally notice that the proof of our main result would have been quite standard if $L$ was assumed fixed, and rather easy if it was assumed that $L \rightarrow +\infty$ and $\frac{L}{M} \rightarrow 0$, a condition very close from $L = \mathcal{O}(N^{\alpha})$ for $\alpha < 1/2$. However, the case 
$1/2 \leq \alpha < 2/3$ needs much more efforts. As explained below, we feel that $2/3$ is the optimal limit. 
\subsection{Overview of the paper}
We first state the main result of this paper. 
\begin{theorem}
\label{theo:main}
When $M \rightarrow +\infty$, and $ML$ and $N$ converge towards $\infty$ in such a way that
$c_N = \frac{ML}{N}$ converges towards $c_* \in (0, +\infty)$, the eigenvalue distribution of 
${\bf W}_N {\bf W}_N^{*}$ converges weakly almost surely towards the Marcenko-Pastur distribution 
with parameters $\sigma^{2}, c_*$. If moreover
\begin{equation}
\label{eq:regime-L}
L = \mathcal{O}(N^{\alpha})
\end{equation}
where $\alpha < 2/3$, then, for each $\epsilon > 0$, almost surely for $N$ large enough, 
all the eigenvalues of ${\bf W}_N {\bf W}_N^{*}$ are located in the interval 
$[ \sigma^{2} (1 - \sqrt{c}_*)^{2} - \epsilon, \sigma^{2} (1 + \sqrt{c}_*)^{2} + \epsilon ]$ 
if $c_* \leq 1$. If $c_* > 1$, almost surely for $N$ large enough, $0$ is 
eigenvalue of ${\bf W}_N {\bf W}_N^{*}$ with multiplicity $ML - N$, and  
the $N$ non zero eigenvalues of ${\bf W}_N {\bf W}_N^{*}$ are located in the interval 
$[ \sigma^{2} (1 - \sqrt{c}_*)^{2} - \epsilon, \sigma^{2} (1 + \sqrt{c}_*)^{2} + \epsilon ]$
\end{theorem}
In order to prove the almost sure location of the eigenvalues of ${\bf W}_N {\bf W}_N^{*}$, 
we follow the approach of \cite{haagerupnew2005} and \cite{schultz-2005}. 
We denote by $t_N(z)$ the Stieltjes transform associated to the Marcenko-Pastur distribution $\mu_{\sigma^{2}, c_N}$ 
with parameters $\sigma^{2}, c_N$, i.e. the unique Stieltjes transform solution of the equation 
\begin{equation}
\label{eq:def-t}
t_N(z) = \frac{1}{-z + \frac{\sigma^{2}}{1 + \sigma^{2} c_N t_N(z)}}
\end{equation}
or equivalently of the system 
\begin{eqnarray}
\label{eq:def-t-symetrique}
t_N(z) & = & \frac{-1}{z \left(1+\sigma^{2} \tilde{t}_N(z) \right)} \\
\label{eq:def-tildet}
\tilde{t}_N(z) & = & \frac{-1}{z \left(1+\sigma^{2} c_N t_N(z) \right)}
\end{eqnarray}
where $\tilde{t}_N(z)$ coincides with 
the Stieltjes transform of $\mu_{\sigma^{2} c_N,1/c_N} =  c_N \mu_{\sigma^{2}, c_N} + (1 - c_N) \delta_0$ where 
$\delta_0$ represents the Dirac distribution at point $0$. 
We denote by $\mathcal{S}^{(0)}_N$ the interval 
\begin{equation}
\label{eq:def-S0}
\mathcal{S}^{(0)}_N = [ \sigma^{2} (1 - \sqrt{c}_N)^{2}, \sigma^{2} (1 + \sqrt{c}_N)^{2} ]
\end{equation}
and by $\mathcal{S}_N$ the support of $\mu_{\sigma^{2}, c_N}$. 
It is well known that $\mathcal{S}_N$ is given by
\begin{eqnarray}
\label{eq:carac-SN-cN<1}
\mathcal{S}_N & = &  \mathcal{S}^{(0)}_N \; \mbox{if $c_N \leq 1$} \\
 \label{eq:carac-SN-cN>1} 
\mathcal{S}_N & = &  \mathcal{S}^{(0)}_N  \cup \{ 0 \} \; \mbox{if $c_N >  1$} 
\end{eqnarray}
Theorem \ref{theo:main} appears to be a consequence of the following identity: 
\begin{equation}
\label{eq:identite-fondamentale}
\mathbb{E} \left[ \frac{1}{ML} \mathrm{Tr} \left( \left({\bf W}_N {\bf W}_N^{*} - z \, {\bf I}_{ML} \right)^{-1} \right) \right] \, - \, t_N(z) = \frac{L}{MN} \, \left( \hat{s}_N(z) \, + \, \frac{L^{3/2}}{MN} \, \hat{r}_N(z) \right)  
\end{equation}
where $\hat{s}_N(z)$ coincides with the Stieltjes transform of a distribution whose support is 
included into $\mathcal{S}^{(0}_N$ and where $\hat{r}_N(z)$ is a function holomorphic in $\mathbb{C}^{+}$ satisfying
\begin{equation}
\label{eq:inegalite-fondamentale-reste}
|\hat{r}_N(z)| \leq P_1(|z|) \, P_2\left(1/\mathrm{Im}(z)\right)
\end{equation}
for $z \in F^{(2)}_N$ where $F^{(2)}_N$ is a subset of $\mathbb{C}^{+}$ defined by 
\begin{equation}
\label{eq:def-FN}
F^{(2)}_N = \{ z \in \mathbb{C}^{+}, \frac{L^{2}}{MN} Q_1(|z|) Q_2(1/\mathrm{Im}(z)) \leq 1 \}
\end{equation}
where $P_1,P_2, Q_1, Q_2$ are polynomials independent of the dimensions $L,M,N$ with positive coefficients. We note that (\ref{eq:regime-L}) is nearly equivalent to $\frac{L^{2}}{MN} \rightarrow 0$ or $\frac{L}{M^{2}} \rightarrow 0$ (in the sense that 
if $\alpha \geq 2/3$, then $\frac{L^{2}}{MN}$ does not converge towards $0$), and that    
$F^{(2)}_N$ appears arbitrary close from $\mathbb{C}^{+}$ when $N$ increases. The present paper is essentially devoted to the proof of (\ref{eq:identite-fondamentale}) under the assumption (\ref{eq:regime-L}). For this, 
we study in the various sections the behaviour of the resolvent ${\bf Q}_N(z)$ of matrix ${\bf W}_N {\bf W}_N^{*}$ defined by
\begin{equation}
\label{eq:def-Q}
{\bf Q}_N(z) = \left( {\bf W}_N {\bf W}_N^{*} - z \, {\bf I}_{ML} \right)^{-1}
\end{equation}
when $z \in \mathbb{C}^{+}$. We use Gaussian tools (integration by parts formula and Poincaré-Nash inequality) as in \cite{pastur-simple} and \cite{pastur-shcherbina-book}
for that purpose. \\

In section \ref{sec:preliminaries}, we present some properties of certain useful operators which map matrices ${\bf A}$ into band Toeplitz matrices whose elements
depend on the sum of the elements of ${\bf A}$ on each diagonal.  Using Poincaré-Nash inequality, we evaluate in section \ref{sec:nash} the variance of 
certain functional of ${\bf Q}_N(z)$ (normalized trace, quadratic forms, and quadratic forms of the $L \times L$ matrix $\hat{{\bf Q}}_N(z)$ obtained as the mean 
of the $M$ $L \times L$ diagonal blocks of ${\bf Q}_N(z)$). In section \ref{sec:novikov}, we use the integration by parts formula in order  to express 
$\mathbb{E}\left({\bf Q}_N(z)\right)$ as
$$
\mathbb{E}\left({\bf Q}_N(z)\right) = {\bf I}_M \otimes {\bf R}_N(z) \, + \, {\bs \Delta}_N(z)
$$
where ${\bf R}_N(z)$ is a certain holomorphic $\mathbb{C}^{L \times L}$ valued function depending on a Toeplitzified version of $\mathbb{E}\left({\bf Q}_N(z)\right)$, 
and where ${\bs \Delta}_N(z)$ is an error term. The goal of section \ref{sec:E(Q)-R} is to control functionals of the error term ${\bs \Delta}_N(z)$. 
We prove that for each $z \in \mathbb{C}^{+}$,
\begin{equation}
\label{eq:eval-delta-intro}
\left| \frac{1}{ML} \mathrm{Tr} \left({\bs \Delta}_N(z)\right) \right| \leq  \frac{L}{MN} \, P_1(|z|) \, P_2\left(1/\mathrm{Im}(z)\right)
\end{equation}
for some polynomials $P_1$ and $P_2$ independent of $L,M,N$ and that, if $\hat{{\bs \Delta}}_N(z)$ represents the
$L \times L$ matrix $\hat{{\bs \Delta}}_N(z) = \frac{1}{M} \sum_{m=1}^{M} {\bs \Delta}_N^{m,m}(z)$, then, it holds that 
\begin{equation}
\label{eq:eval-delta-intro-2}
\left| {\bf b}_{1}^{*} \,  \hat{{\bs \Delta}}_N(z) \, {\bf b}_{2} \right| \leq  \frac{L^{3/2}}{MN} \, P_1(|z|) \, P_2\left(1/\mathrm{Im}(z)\right)
\end{equation}
for deterministic unit norm $L$--dimensional vectors ${\bf b}_1$ and ${\bf b}_2$. In section \ref{sec:marcenko-pastur}, we prove that 
\begin{equation}
\label{eq:convergence-MP-intro}
\mathbb{E} \left[ \frac{1}{ML} \mathrm{Tr}\left( {\bf Q}_N(z) \right) \right] - t_N(z) \rightarrow 0
\end{equation}
for each $z \in \mathbb{C}^{+}$, a property which implies that the eigenvalue distribution of ${\bf W}_N {\bf W}_N^{*}$ converges towards 
the Marcenko-Pastur distribution. We note that (\ref{eq:convergence-MP-intro}) holds as soon as 
$M \rightarrow +\infty$. At this stage, however, the convergence rate of the lefthandside of (\ref{eq:convergence-MP-intro}) is not precised. Under the condition $\frac{L^{3/2}}{MN} \rightarrow 0$ (which 
implies that quadratic forms of $\hat{{\bs \Delta}}_N(z)$ converge towards $0$, see (\ref{eq:eval-delta-intro-2})), we prove in section \ref{sec:improved-R-t} that
\begin{equation}
\label{eq:identite-fondamentale-facile}
\mathbb{E} \left[ \frac{1}{ML} \mathrm{Tr} \left( \left({\bf W}_N {\bf W}_N^{*} - z \, {\bf I}_{ML} \right)^{-1} \right) \right] \, - \, t_N(z) = \frac{L}{MN} \tilde{r}_N(z)
\end{equation}
where $\tilde{r}_N(z)$ is holomorphic in $\mathbb{C}^{+}$ and satisfies 
$$
|\tilde{r}_N(z)| \leq P_1(|z|) P_2(1/\mathrm{Im}z)
$$
for each $z \in F^{(3/2)}_N$, where $F^{(3/2)}_N$ is defined as $F_N^{(2)}$ (see (\ref{eq:def-FN})), but 
when $\frac{L^{2}}{MN}$ is replaced by $\frac{L^{3/2}}{MN}$. In order to establish (\ref{eq:identite-fondamentale-facile}), it is proved in section \ref{sec:convergence-spectral-norm} that the spectral norm of a Toeplitzified version of matrix ${\bf R}_N(z) - t_N(z) \, {\bf I}_L$ is upperbounded 
by a term such as $\frac{L^{3/2}}{MN} \, P_1(|z|) \, P_2\left(1/\mathrm{Im}(z)\right)$. (\ref{eq:identite-fondamentale-facile}) and  
Lemma 5.5.5 of \cite{and-gui-zei-2010} would allow to establish quite easily the almost sure location of the eigenvalues 
of ${\bf W}_N {\bf W}_N^{*}$ under the hypothesis $\frac{L}{M} \rightarrow 0$. However, this condition 
is very restrictive, and, at least intuitively, somewhat similar to $L$ fixed. In section \ref{sec:expre-E(Q)-t}, we establish that under condition (\ref{eq:regime-L}), which is very close from the condition 
$\frac{L^{2}}{MN} \rightarrow 0$, or $\frac{L}{M^{2}} \rightarrow 0$, 
function $\tilde{r}_N(z)$ can be written as $\tilde{r}_N(z) = \hat{s}_N(z) + \frac{L^{3/2}}{MN} \hat{r}_N(z)$  where $\hat{s}_N(z)$ and $\hat{r}_N(z)$ verify the conditions of (\ref{eq:identite-fondamentale}). We first prove that
\begin{equation}
\label{eq:premiere-decomposition}
\mathbb{E} \left[ \frac{1}{ML} \mathrm{Tr}\left( {\bf Q}_N(z) - {\bf I}_M \otimes {\bf R}_N(z) \right) \right] = 
\frac{L}{MN} \left( s_N(z) + \frac{L}{MN} r_N(z) \right)
\end{equation}
where $s_N(z)$ and $r_N(z)$ satisfy the same properties than $\hat{s}_N(z)$ and $\hat{r}_N(z)$.
For this, we compute explicitely $s_N(z)$, and verify that it coincides with the Stieltjes transform 
of a distribution whose support is included into $\mathcal{S}_N^{(0)}$. The most technical part of the paper 
is to establish that 
\begin{equation}
\label{eq:rate-reste}
\mathbb{E} \left[ \frac{1}{ML} \mathrm{Tr}\left( {\bf Q}_N(z) - {\bf I}_M \otimes {\bf R}_N(z) \right) \right] \, - \,  
\frac{L}{MN} \, s_N(z) 
\end{equation}
converges towards $0$ at rate $\left( \frac{L}{MN} \right)^{2}$. For this, the condition 
$\frac{L^{2}}{MN} \rightarrow 0$ appears to be fundamental because it allows, among others, to control the behaviour of the solutions of $L$--dimensional linear systems obtained by inverting the sum of a diagonal matrix with a matrix with 
$\mathcal{O}(\frac{L}{MN})$ entries. Using the results of section \ref{sec:convergence-spectral-norm} concerning  
the spectral norm of a Toeplitzified version of ${\bf R}_N(z) - t_N(z) \, {\bf I}_L$, we obtain easily (\ref{eq:identite-fondamentale}) from(\ref{eq:premiere-decomposition}). Theorem \ref{theo:main} is finally established in section \ref{sec:conclusion}. 
For this, we follow \cite{haagerupnew2005}, \cite{schultz-2005} and \cite{and-gui-zei-2010} (Lemma 5-5-5). We consider a smooth approximation $\phi$ of 
$\mathbb{1}_{[\sigma^{2}(1 - \sqrt{c_{*}})^{2} - \epsilon, \sigma^{2}(1 + \sqrt{c_{*}})^{2} + \epsilon]^{(c)}}$
that vanishes on $\mathcal{S}_N^{(0)}$ for each $N$ large enough, and establish that almost surely, 
\begin{equation}
\label{eq:rate-trace-indicatrice}
\mathrm{Tr} \left( \phi({\bf W}_N {\bf W}_N^{*}) \right) = N \mathcal{O}(\frac{L^{5/2}}{(MN)^{2}}) + [ML - N]_{+} = 
\mathcal{O}\left( (\frac{L}{M^{2}})^{3/2} \right) + [ML - N]_{+}
\end{equation}
(\ref{eq:regime-L}) implies that $\frac{L}{M^{2}} \rightarrow 0$ and that
the righthandside of (\ref{eq:rate-trace-indicatrice}) converges towards $[ML - N]_{+}$ almost surely. 
This, in turn, establishes that 
the number of non zero eigenvalues of ${\bf W}_N {\bf W}_N^{*}$ that are located outside 
$[\sigma^{2}(1 - \sqrt{c_{*}})^{2} - \epsilon, \sigma^{2}(1 + \sqrt{c_{*}})^{2} + \epsilon]$ 
converges towards zero almost surely, and is thus equal to $0$ for $N$ large enough as expected. 
We have not proved that this property does not hold if $L = \mathcal{O}(N^{\alpha})$ with $\alpha \geq 2/3$. We however mention that the hypothesis $\alpha < 2/3$ is used at various crucial independent steps:
\begin{itemize}
\item it is used extensively to establish that (\ref{eq:rate-reste}) converges towards $0$ at rate $\left(\frac{L}{MN}\right)^{2}$
\item it is nearly equivalent to the condition $\frac{L^{2}}{MN} \rightarrow 0$ or 
$\frac{L}{M^{2}} \rightarrow 0$ which implies 
\begin{itemize}
\item that the set $F_N^{(2)}$ defined by (\ref{eq:def-FN}))
is arbitrarily close from $\mathbb{C}^{+}$, a property that appears necessary to generalize 
Lemma 5-5-5 of \cite{and-gui-zei-2010} 
\item that the righthandside of (\ref{eq:rate-trace-indicatrice}) converges towards $[ML - N]_{+}$
\end{itemize}
\end{itemize}
We therefore suspect that the almost sure location of the eigenvalues of ${\bf W}_N {\bf W}_N^{*}$ cannot be established using the approach of \cite{haagerupnew2005} and \cite{schultz-2005}
if $\alpha \geq 2/3$. It would be interesting to study the potential of combinatorial methods in order to be fully convinced that the almost sure location of the eigenvalues of ${\bf W}_N {\bf W}_N^{*}$ does not hold if $\alpha \geq 2/3$. We finally mention that we have performed numerical simulations to check  
whether it is reasonable to conjecture that the almost sure location property of the eigenvalues of  ${\bf W}_N {\bf W}_N^{*}$ holds if and only if $\alpha < 2/3$. For this, we have generated 10.000 independent realizations of the largest eigenvalue $\lambda_{1,N}$ of ${\bf W}_N {\bf W}_N^{*}$ for $\sigma^{2} = 1$,
$N = 2^{14}$, $c_N = ML/N = 1/2$ and for the following values of $(M,L)$ that seem to be in accordance
with the asymptotic regime considered in this paper: $(M,L) = (2^{8}, 2^{5})$, 
$(M,L) = (2^{7}, 2^{6})$, $(M,L) = (2^{6}, 2^{7})$, $(M,L) = (2^{5}, 2^{8})$, corresponding 
to ratios $\frac{L}{M^{2}}$ equal respectively to $2^{-11}$,  $2^{-8}$,  $2^{-5}$, and $1/4$. As condition $\alpha < 2/3$ is nearly equivalent to $\frac{L}{M^{2}} \rightarrow 0$, the first 3 values of $(M,L)$ are in accordance with the asymptotic regime $L = \mathcal{O}(N^{\alpha})$ with $\alpha < 2/3$ while it is of course not the case for the last configuration. The almost sure location property 
of course implies that the largest eigenvalue converges towards $(1 + \sqrt{c_*})^{2}$. 
In order to check this property, we have evaluated the empirical mean $\overline{\lambda}_{1,N}$ of the 10.000 realizations of $\lambda_{1,N}$, and compared $\overline{\lambda}_{1,N}$ with $(1 + \sqrt{1/2})^{2} \simeq 2.91$.

\begin{table}
\centering
\caption{Empirical mean of the largest eigenvalue versus $L/M^{2}$}
\label{table:largest-eigenvalue}      
\begin{tabular}{|c|c|c|c|c|}
\hline\noalign{\smallskip}
$L/M^{2}$ & $2^{-11}$  &  $2^{-8}$  & $2^{-5}$  & $1/4$ \\
\noalign{\smallskip}\hline\noalign{\smallskip}
$\overline{\lambda}_{1,N}$  & 2.91  &  2.92  & 2.94  &  3 \\
\noalign{\smallskip}\hline
\end{tabular}
\end{table}

The values of $\overline{\lambda}_{1,N}$ in terms of $\frac{L}{M^{2}}$ are presented in Table \ref{table:largest-eigenvalue}. It is seen that the difference between $\overline{\lambda}_{1,N}$ and $(1 + \sqrt{1/2})^{2} \simeq 2.91$ 
increases significantly with the ratio $\frac{L}{M^{2}}$, thus suggesting that $\lambda_{1,N}$ does 
not converge almost surely towards $(1 + \sqrt{c_*})^{2}$ when $\frac{L}{M^{2}}$ does not converge towards $0$.

\subsection{General notations and definitions}

{\bf Assumptions on $L,M,N$}
\begin{assum}
\label{as:standard}
\begin{itemize}
\item All along the paper,  we assume that $L, M, N$ satisfy $M \rightarrow +\infty, N \rightarrow +\infty$ in such a way that $c_N = \frac{ML}{N} \rightarrow c_*$, where $0 < c_* < +\infty$. In order to short the notations, $N \rightarrow +\infty$ should be understood as the above asymptotic regime.
\item In sections \ref{sec:convergence-spectral-norm} and \ref{sec:improved-R-t}, $L,M,N$ also satisfy 
$\frac{L^{3/2}}{MN} \rightarrow 0$ or equivalently $\frac{L}{M^{4}} \rightarrow 0$.
\item In sections \ref{sec:expre-E(Q)-t} and \ref{sec:conclusion}, the extra condition $L = \mathcal{O}(N^{\alpha})$
with $\alpha < 2/3$ holds. 
\end{itemize} 
\end{assum}

In the following, we will often drop the index $N$, and will denote ${\bf W}_N, t_N, {\bf Q}_N, \ldots$ by 
${\bf W}, t, {\bf Q}, \ldots$ in order to short the notations. The $N$ columns of matrix ${\bf W}$ are denoted $({\bf w}_j)_{j=1, \ldots, N}$. For $1 \leq l \leq L$, $1 \leq m \leq M$, 
and $1 \leq j \leq N$, ${\bf W}_{i,j}^{m}$ represents the entry $\left(i+(m-1)L, j \right)$ of matrix ${\bf W}$. 

$\mathcal{C}^{\infty}(\mathbb{R})$ (resp. $\mathcal{C}_b^{\infty}(\mathbb{R}), \mathcal{C}_c^{\infty}(\mathbb{R})$) denotes the space of 
all real-valued smooth functions (resp. bounded smooth functions, smooth functions with compact support) 
defined on $\mathbb{R}$.

If ${\bf A}$ is a $ML \times ML$ matrix, 
we denote by ${\bf A}^{m_1,m_2}_{i_1,i_2}$ the entry $(i_1 + (m_1-1)L, i_2 + (m_2-1)L)$ of matrix ${\bf A}$, while 
${\bf A}^{m_1,m_2}$ represents the $L \times L$ matrix $({\bf A}_{i_1,i_2}^{m_1,m_2})_{1 \leq (i_1,i_2) \leq L}$. We also 
denote by $\hat{{\bf A}}$ the $L \times L$ matrix defined by
\begin{equation}
\label{eq:def-hatA}
\hat{{\bf A}} = \frac{1}{M} \sum_{m=1}^{M} {\bf A}^{m,m}
\end{equation}

For each $1 \leq i \leq L$ and 
$1 \leq m \leq M$, ${\bf f}_{i}^{m}$ represents the vector of the canonical basis of $\mathbb{C}^{ML}$ whose non zero component 
is located at index $i+(m-1)L$. If $1 \leq j \leq N$, ${\bf e}_j$ is the $j^{\mathrm{th}}$-vector of the canonical basis of 
$\mathbb{C}^{N}$.

If ${\bf A}$ and ${\bf B}$ are 2 matrices, ${\bf A} \otimes {\bf B}$ represents the Kronecker product of ${\bf A}$ and ${\bf B}$, i.e. the block matrix whose  block $(i,j)$ is
${\bf A}_{i,j} \, {\bf B}$.  $\| {\bf A} \|$ represents the spectral norm of matrix ${\bf A}$.

If $x \in \mathbb{R}$, $[x]_{+}$ represents $\max(x,0)$. $\mathbb{C}^{+}$ denotes the set of complex numbers with strictly positive imaginary parts. The conjuguate of a complex number $z$ is denoted $z^{*}$ or $\overline{z}$ depending on the context. Unless otherwise stated, $z$ represents  an element of $\mathbb{C}^{+}$. If ${\bf A}$ is a square matrix, $\mathrm{Re}({\bf A})$ and $\mathrm{Im}({\bf A})$ represent
the Hermitian matrices $\mathrm{Re}({\bf A}) = \frac{{\bf A} + {\bf A}^{*}}{2}$ and $ \mathrm{Im}({\bf A}) = \frac{{\bf A} - {\bf A}^{*}}{2i}$ respectively. 

If $({\bf A}_N)_{N \geq 1}$ (resp. $({\bf b}_N)_{N \geq 1}$)
is a sequence of matrices (resp. vectors) whose dimensions increase with $N$, $({\bf A}_N)_{N \geq 1}$ (resp. $({\bf b}_N)_{N \geq 1}$) is said to 
be uniformly bounded if $\sup_{N \geq 1} \| {\bf A}_N \| < +\infty$ (resp. $\sup_{N \geq 1} \| {\bf b}_N \| < +\infty$). \\

If $x$ is a complex-valued random variable, the variance of $x$, denoted by $\mathrm{Var}(x)$, is defined by 
$$
\mathrm{Var}(x) = \mathbb{E}(|x|^{2}) - \left| \mathbb{E}(x) \right|^{2}
$$
The zero-mean random variable $x - \mathbb{E}(x)$ is denoted $x^{\circ}$. 

{\bf Nice constants and nice polynomials}. A nice constant is a positive constant independent of the dimensions $L,M,N$ and complex variable $z$. 
A nice polynomial is a polynomial whose degree is independent from $L,M,N$, and whose coefficients are nice constants. In the following, 
$P_1$ and $P_2$ will represent generic nice polynomials whose values may change from one line to another, and $C(z)$ 
is a generic term of the form $C(z) = P_1(|z|) P_2(1/\mathrm{Im}z)$. 

{\bf Properties of matrix ${\bf Q}(z)$.} We recall that ${\bf Q}(z)$ verifies the so-called resolvent identity
\begin{equation}
\label{eq:equation-resolvente}
{\bf Q}(z) = -\frac{{\bf I}_{ML}}{z} + \frac{1}{z} {\bf Q}(z) {\bf W} {\bf W}^{*}
\end{equation}
and that it holds that 
\begin{equation}
\label{eq:borne-QQ*}
{\bf Q}(z) {\bf Q}^{*}(z) \leq \frac{{\bf I}_{ML}}{(\mathrm{Im}z)^{2}}
\end{equation}
and that 
\begin{equation}
\label{eq:norme-resolvente-im}
\| {\bf Q}(z) \|  \leq  \frac{1}{\mathrm{Im}(z)}
\end{equation}
for $z \in \mathbb{C}^{+}$. We also mention that 
\begin{equation}
\label{eq:ImQ>0}
\mathrm{Im}({\bf Q}(z)) > 0, \; \; \mathrm{Im}(z \,{\bf Q}(z)) > 0, \; \mathrm{if} z \in \mathbb{C}^{+}
\end{equation}

{\bf Gaussian tools.} We present the versions of the integration by parts formula (see Eq. (2.1.42) p. 40 in \cite{pastur-shcherbina-book} for the real case and Eq. (17) in 
\cite{hac-kho-lou-naj-pas-08} for the present complex case) and of the 
Poincaré-Nash (see Proposition 2.1.6 in \cite{pastur-shcherbina-book} for the real case
and Eq. (18) in \cite{hac-kho-lou-naj-pas-08} for the complex case) that we use in this paper. 
\begin{proposition}
   {\bf Integration by parts formula.} Let $\bs \xi=[\xi_1, \ldots, \xi_K]^T$ be a complex Gaussian random vector such that
   $\mathbb{E}[\bs \xi]=\bs 0$,
   $\mathbb{E}[\bs \xi \bs \xi^T]=\bs 0$
   and
   $\mathbb{E}[\bs \xi \bs \xi^*]=\bs \Omega$.
   If $\Gammabs : (\bs \xi) \mapsto \Gammabs(\bs \xi, \overline{\bs \xi})$
   is a $\mathcal{C}^1$ complex function polynomially bounded together with its derivatives, then
   \begin{equation}
     \mathbb{E}[\xi_p\Gammabs (\bs \xi)]=\sum_{m=1}^K \bs\Omega_{pm} \mathbb{E}\left[\frac{\partial\Gammabs (\bs \xi)}{\partial \overline{\xi}_m}\right].
     \label{eq:IPP-th}
   \end{equation}
  \end{proposition}

\begin{proposition}
   {\bf Poincaré-Nash inequality.} Let $\bs \xi=[\xi_1, \ldots, \xi_K]^T$ be a complex Gaussian random vector such that
   $\mathbb{E}[\bs \xi]=\bs 0$,
   $\mathbb{E}[\bs \xi \bs \xi^T]=\bs 0$
   and
   $\mathbb{E}[\bs \xi \bs \xi^*]=\bs \Omega$.
   If $\Gammabs : (\bs \xi) \mapsto \Gammabs (\bs \xi, \overline{\bs \xi})$
   is a $\mathcal{C}^1$ complex function polynomially bounded together with its derivatives, then,
   noting $\nabla_{\bs \xi} {\bs \Gamma}=[\frac{\partial\Gammabs}{\partial \xi_1},\ldots,\frac{\partial\Gammabs}{\partial \xi_K}]^T$
   and $\nabla_{\overline{\bs \xi}} \Gammabs =[\frac{\partial\Gammabs}{\partial \overline\xi_1},\ldots,\frac{\partial\Gammabs}{\partial \overline\xi_K}]^T$,
   \begin{equation}
     \mathrm{Var}(\Gammabs(\bs\xi))\leq
     \mathbb{E}\left[ \nabla_{\bs \xi} \Gammabs(\bs \xi)^T \ \bs\Omega \ \overline{\nabla_{\bs \xi} \Gammabs(\bs\xi)} \right]
     + \mathbb{E}\left[ \nabla_{\overline{\bs \xi}} \Gammabs(\bs \xi)^* \ \bs\Omega \ \nabla_{\overline{\bs \xi}} \Gammabs(\bs\xi) \right]
     \label{eq:NPinq-th}
   \end{equation}
  \end{proposition}

The above two propositions are used below in the case where ${\bs \xi}$ coincides with the $LMN$--dimensional vector 
$\mathrm{vec}({\bf W}_N)$. In the following, the particular structure ${\bf W}^{m}_{i,j} = w_{m,i+j-1}$ of ${\bf W}_N$ is encoded by the correlation structure of the entries of ${\bf W}_N$:
\begin{equation}
\label{eq:correlation-W}
\mathbb{E}\left( {\bf W}_{i_1,j_1}^{m_1} \overline{{\bf W}}_{i_2,j_2}^{m_2} \right) = 
\frac{\sigma^{2}}{N} \, \delta(i_1 - i_2 = j_2 - j_1) \, \delta(m_1=m_2)
\end{equation}

{\bf A useful property of the Stieltjes transform $t_N(z)$ of the Marcenko-Pastur 
$\mu_{\sigma^{2},c_N}$.} 

The following lemma is more or less known. A proof is provided in the Appendix of \cite{loubaton-arxiv} for the reader's convenience. 
\begin{lemma}
\label{le:proprietes-t-tilde-t}
It holds that 
\begin{equation}
\label{eq:zttildet}
\sigma^{4} c_N |z t_N(z) \tilde{t}_N(z)|^{2} < 1
\end{equation}
for each $z \in \mathbb{C}^{+}$. Moreover, for each $N$ and for each $z \in \mathbb{C}^{+}$, it holds that 
\begin{equation}
\label{eq:borneinf-zttildet}
1 - \sigma^{4} c_N |z t_N(z) \tilde{t}_N(z)|^{2} > C \, \frac{(\mathrm{Im}z)^{4}}{(\eta^{2} + |z|^{2})^{2}}
\end{equation}
for some nice constants $C$ and $\eta$. Finally, for each $N$, it holds that 
\begin{equation}
\label{eq:borne-sup-inverse-zttildet-distribution}
\left(1 - \sigma^{4} c_N |z t(z) \tilde{t}(z)|^{2}\right)^{-1} \leq C \, \max \left(1, \frac{1}{(\mathrm{dist}(z, \mathcal{S}^{(0)}_N))^{2}} \right)
\end{equation}
for some nice constant $C$ and for each $z \in \mathbb{C} - \mathcal{S}^{(0)}_N$. 
\end{lemma}

\section{Preliminaries}
\label{sec:preliminaries}
In this section, we introduce certain Toeplitzification operators, and establish some useful related properties.

\begin{definition}
\begin{itemize}
\item If ${\bf A}$ is a $K \times K$ Toeplitz matrix, we denote by
$({\bf A}(k))_{k=-(K-1), \ldots, K-1}$ the sequence such that 
${\bf A}_{k,l} = {\bf A}(k-l)$.
\item For any integer $K$, $J_K$ is the $K \times K$ ``shift'' matrix defined by $(J_K)_{i,j} = \delta(j-i=1)$. In order to 
short the notations, matrix $J_K^{*}$ is denoted $J_K^{-1}$, although $J_K$ is of course not invertible. 
\item For any $PK \times PK$ block matrix ${\bf A}$ with $K \times K$ blocks $({\bf A}^{p_1,p_2})_{1 \leq (p_1,p_2) \leq P}$, we define $(\tau^{(P)}({\bf A})(k))_{k=-(K-1), \ldots, K-1}$ as the sequence
\begin{equation}
\label{eq:def-calt(A)}
 \tau^{(P)}({\bf A})(k) = \frac{1}{PK} \mathrm{Tr} \left[ {\bf A} ({\bf I}_P \otimes {\bf J}_K^{k}) \right] = \frac{1}{PK} \sum_{i-j=k} \sum_{p=1}^{P} {\bf A}^{(p,p)}_{i,j} = \frac{1}{PK} \sum_{p=1}^{P} \sum_{u=1}^{K} {\bf A}^{p,p}_{k+u,u} \, \mathbb{1}_{1 \leq k+u\leq K}
\end{equation}
\item For any $PK \times PK$ block matrix ${\bf A}$ and for 2 integers $R$ and $Q$ such that $R \geq Q$ and $Q \leq K$, matrix $\mathcal{T}^{(P)}_{R,Q}({\bf A})$ represents the $R \times R$ Toeplitz matrix
given by 
\begin{equation}
\label{eq:def-calT(A)} 
\mathcal{T}^{(P)}_{R,Q}({\bf A}) = \sum_{q=-(Q-1)}^{Q-1} \tau^{(P)}({\bf A})(q) \;  {\bf J}_R^{*q}
\end{equation}
In other words, for $(i,j) \in \{ 1, 2, \ldots, R \}$, it holds that
\begin{equation}
\label{eq:def-calT(A)-bis} 
\left(\mathcal{T}^{(P)}_{R,Q}({\bf A})\right)_{i,j} = \tau^{(P)}({\bf A})(i-j) \, \mathbb{1}_{|i-j| \leq Q-1}
\end{equation}
When $P=1$, sequence  $(\tau^{(1)}({\bf A})(k))_{k=-(K-1), \ldots, K-1}$ and matrix $\mathcal{T}^{(1)}_{R,Q}({\bf A})$ are denoted  $(\tau({\bf A})(k))_{k=-(K-1), \ldots, K-1}$ and matrix $\mathcal{T}_{R,Q}({\bf A})$ in order to simplify the notations. We note that if ${\bf A}$ is a $PK \times PK$ 
block matrix, then, sequence $(\tau^{(P)}({\bf A})(k))_{k=-(K-1), \ldots, K-1}$ coincides with sequence 
$\left(\tau\left(\hat{{\bf A}} \right)(k)\right)_{k=-(K-1), \ldots, K-1}$ 
where we recall that $\hat{{\bf A}} = \frac{1}{P} \sum_{p=1}^{P} {\bf A}^{p,p}$; matrix $\mathcal{T}^{(P)}_{R,Q}({\bf A})$ is equal to $\mathcal{T}_{R,Q}(\hat{{\bf A}})$. 
\end{itemize}
\end{definition}
The reader may check that the following straightforward identities hold: 
\begin{itemize}
\item If ${\bf A}$ is a $R \times R$ Toeplitz matrix, for any 
$R \times R$ matrix ${\bf B}$, it holds that
\begin{equation}
\label{eq:Tr(AB)}
\frac{1}{R} \mathrm{Tr} ({\bf A} {\bf B}) = \sum_{k=-(R-1)}^{R-1} 
{\bf A}(-k) \tau({\bf B})(k) = \frac{1}{R} \mathrm{Tr} \left( {\bf A} \mathcal{T}_{R,R}({\bf B}) \right)
\end{equation}
\item If ${\bf A}$ and ${\bf B}$ are both $R \times R$ matrices, and if $Q \leq R$, then,  
\begin{equation}
\label{eq:Tr(T(A)B)}
\frac{1}{R} \mathrm{Tr} \left(\mathcal{T}_{R,Q}({\bf A}) {\bf B}\right) = \sum_{q=-(Q-1)}^{Q-1} 
\tau({\bf A})(-q) \; \tau({\bf B})(q) = \frac{1}{R} \mathrm{Tr} \left( {\bf A} \mathcal{T}_{R,Q}({\bf B}) \right)
\end{equation}
\item 
If ${\bf A}$ is a $PK \times PK$ matrix, if ${\bf B}$ is a $R \times R$ matrix, and if  $R \geq Q$ and $Q \leq K$, 
then it holds that
\begin{equation}
\label{eq:Tr(BT(A)}
\frac{1}{R} \mathrm{Tr}\left( {\bf B} \mathcal{T}^{(P)}_{R,Q}({\bf A}) \right) = \sum_{k=-(Q-1)}^{Q-1} \tau({\bf B})(k) \; \tau^{(P)}({\bf A})(-k) = \frac{1}{PK} \mathrm{Tr} \left( \left({\bf I}_M \otimes \mathcal{T}_{K,Q}({\bf B}) \right) {\bf A} \right)
\end{equation}
\item If ${\bf C}$ is a $PK \times PK$ matrix, ${\bf B}$ is a $K \times K$ matrix
and ${\bf D}, {\bf E}$  $R \times R$ matrices with $K \leq R$, then, it holds that 
\begin{equation}
\label{eq:utile}
\frac{1}{K} \mathrm{Tr} \left[ {\bf B} \mathcal{T}_{K,K}\left( {\bf D}  \mathcal{T}_{R,K}^{(P)}({\bf C}) {\bf E} \right) \right] = \frac{1}{PK} \mathrm{Tr} \left[ {\bf C} \left( {\bf I}_{P} \otimes
\mathcal{T}_{K,K}[ {\bf E} \mathcal{T}_{R,K}({\bf B}) {\bf D} ] \right) \right]
\end{equation}
\end{itemize}

We now establish useful properties of matrix $\mathcal{T}^{(P)}_{R,Q}({\bf A})$. 
\begin{proposition} 
\label{prop:contractant}
If ${\bf A}$ is a $PK \times PK$ matrix, then, for each integer $R \geq K$, it holds that
\begin{equation}
\label{eq:Tcontractant}
\left\|  \mathcal{T}^{(P)}_{R,K}({\bf A}) \right\| \leq \sup_{\nu \in [0,1]} \left| {\bf a}_K(\nu)^{*} 
\left( \frac{1}{P} \sum_{p=1}^{P} {\bf A}^{p,p} \right) {\bf a}_K(\nu) \right| \leq \| {\bf A} \|
\end{equation}
where ${\bf a}_K(\nu)$ represents the $K$--dimensional vector defined by 
\begin{equation}
\label{eq:def-anu}
{\bf a}_K(\nu) = \frac{1}{\sqrt{K}} \, \left(1, e^{2i \pi \nu}, \ldots, e^{2i \pi (K-1) \nu} \right)^{T}
\end{equation}
If ${\bf A}$ is a $K \times K$ matrix and if $R \leq K$, then, it holds that 
\begin{equation}
\label{eq:Tcontractant-bis}
\left\|  \mathcal{T}_{R,R}({\bf A}) \right\| \leq \sup_{\nu \in [0,1]} \left| {\bf a}_K(\nu)^{*}  {\bf A} \, {\bf a}_K(\nu) \right| \leq \| {\bf A} \|
\end{equation}
\end{proposition}
{\bf Proof.} We first establish (\ref{eq:Tcontractant}). As $R \geq K$, matrix 
$\mathcal{T}^{(P)}_{R,K}({\bf A})$ is a submatrix of the infinite band Toeplitz matrix 
with $(i,j)$ elements $\tau^{(P)}({\bf A})(i-j) \mathbb{1}_{|i-j| \leq K-1}$. 
The norm of this matrix is known to be equal to the $\mathbb{L}_{\infty}$ norm 
of the corresponding symbol (see \cite{bottcher-silbermann}, Eq. (1-14), p. 10). Therefore, it 
holds that  
$$
\| \mathcal{T}^{(P)}_{R,K}({\bf A}) \| \leq \sup_{\nu \in [0,1]} \left| \sum_{k=-(K-1)}^{K-1} \tau^{(P)}({\bf A})(k) e^{-2 i \pi k \nu} \right|
$$
We now verify the following useful identity:
\begin{equation}
\label{eq:expre-symbole}
\sum_{k=-(K-1)}^{K-1} \tau^{(P)}({\bf A})(k) e^{-2 i \pi k \nu}  = {\bf a}_K(\nu)^{*} \, \left (\frac{1}{P} \sum_{p=1}^{P} {\bf A}^{(p,p)}\right) \, {\bf a}_K(\nu)
\end{equation}
Using the definition (\ref{eq:def-calt(A)}) of $\tau^{(P)}({\bf A})(k)$, the term $\sum_{k=-(K-1)}^{K-1} \tau^{(P)}({\bf A})(k) e^{-2 i \pi k \nu}$ can also be written as
$$
\sum_{k=-(K-1)}^{K-1} \tau^{(P)}({\bf A})(k) e^{-2 i \pi k \nu} = \frac{1}{K} \sum_{k=-(K-1)}^{K-1} \mathrm{Tr}\left( \left (\frac{1}{P} \sum_{p=1}^{P} {\bf A}^{(p,p)}\right) \, e^{-2 i \pi k \nu} {\bf J}_K^{k} \right)  
$$
or equivalently as 
$$
\mathrm{Tr}\left( \left(\frac{1}{P} \sum_{p=1}^{P} {\bf A}^{(p,p)} \right) \, \frac{1}{K} \left(\sum_{k=-(K-1)}^{K-1} e^{-2 i \pi k \nu} {\bf J}_K^{k} \right) \right)  
$$
It is easily seen that 
$$
\frac{1}{K} \left (\sum_{k=-(K-1)}^{K-1} e^{-2 i \pi k \nu} {\bf J}_K^{k} \right) = {\bf a}_K(\nu) {\bf a}_K(\nu)^{*}
$$
from which (\ref{eq:expre-symbole}) and (\ref{eq:Tcontractant}) follow immediately. 

In order to justify (\ref{eq:Tcontractant-bis}), we remark that $R \leq K$ implies that $\mathcal{T}_{R,R}({\bf A})$ is a submatrix of 
$\mathcal{T}_{K,K}({\bf A})$ whose norm is bounded by $\sup_{\nu} \left|  {\bf a}_K(\nu)^{*}  {\bf A} \, {\bf a}_K(\nu) \right|$ by (\ref{eq:Tcontractant}). \\

We also prove that the operators $\mathcal{T}$ preserve the positivity of matrices. 
\begin{proposition} 
\label{prop:positivite}
If ${\bf A}$ is a $PK \times PK$  positive definite matrix, then, for each integer $R \geq K$, it holds that
\begin{equation}
\label{eq:positivite}
\mathcal{T}_{R,K}^{(P)}({\bf A}) > 0
\end{equation}
If ${\bf A}$ is a $K \times K$ positive definite matrix and if $R \leq K$, then, it holds that 
\begin{equation}
\label{eq:positivite-bis}
\mathcal{T}_{R,R}({\bf A}) > 0
\end{equation}
\end{proposition}
{\bf Proof.} We first prove (\ref{eq:positivite}). (\ref{eq:expre-symbole}) implies that 
$$
\sum_{k=-(K-1)}^{K-1} \tau^{(P)}({\bf A})(k) e^{-2 i \pi k \nu} > 0
$$
for each $\nu$. $(\tau^{(P)}({\bf A})(k))_{k=-(K-1), \ldots, K-1}$ thus coincide the Fourier coefficients of a 
positive function. Elementary results related the trigonometric moment problem (see e.g. \cite{grenander-szego}, 1.11 (a)) imply that for each $R \geq K$, matrix $\mathcal{T}_{R,K}^{(P)}({\bf A})$ is positive definite. We finally justify 
(\ref{eq:positivite-bis}). As $R \leq K$, matrix $\mathcal{T}_{R,R}({\bf A})$ is a submatrix of $\mathcal{T}_{K,K}({\bf A})$ which is positive definite
by (\ref{eq:positivite}). \\

We finally give the following useful result proved in the Appendix.  
\begin{proposition}
\label{prop:T(H)T(H)*}
If ${\bf A}$ is a $K \times K$ matrix and if $R \geq K$, then, it holds that 
\begin{equation}
\label{eq:TNL(H)T(H)*}
\mathcal{T}_{R,K}({\bf A}) \left(\mathcal{T}_{R,K}({\bf A}) \right)^{*} \leq  \mathcal{T}_{R,K}({\bf A} {\bf A}^{*})
\end{equation}
If ${\bf A}$ is a $K \times K$ matrix and if $R \leq K$, then 
\begin{equation}
\label{eq:TLL(H)T(H)*}
\mathcal{T}_{R,R}({\bf A}) \left(\mathcal{T}_{R,R}({\bf A}) \right)^{*} \leq  \mathcal{T}_{R,R}({\bf A} {\bf A}^{*})
\end{equation}
\end{proposition}

\section{Poincaré-Nash variance evaluations}
\label{sec:nash}
In this section, we take benefit of the Poincaré-Nash inequality to evaluate the 
variance of certain important terms. In particular, we prove the following useful result. 
\begin{proposition}
\label{prop-nash-poincare}
Let ${\bf A}$ be a deterministic $ML \times ML$ matrix for which $\sup_{N} \| {\bf A} \| \leq \kappa$, and consider 
2 $ML$--dimensional deterministic vectors ${\bf a}_1, {\bf a}_2$ such that $\sup_{N} \|{\bf a}_i\| \leq \kappa$ for $i=1,2$ as well as
2 $L$--dimensional deterministic vectors ${\bf b}_1, {\bf b}_2$ such that $\sup_{N} \|{\bf b}_i\| \leq \kappa$ for $i=1,2$. Then, 
for each $z \in \mathbb{C}^{+}$, it holds that 
\begin{align}
\label{eq:vartraceQ}
\mathrm{Var}\left( \frac{1}{ML} \mathrm{Tr}\left( {\bf A} {\bf Q}(z) \right) \right)  \leq  C(z) \, \kappa^{2} \, \frac{1}{MN} \\
\label{eq:varentreesQ}
\mathrm{Var}\left( {\bf a}_1^{*} {\bf Q}(z) {\bf a}_2 \right) \leq  C(z) \, \kappa^{4} \,  \frac{L}{N} \\
\label{eq:vartraceentreesQ}
\mathrm{Var}\left(  {\bf b}_1^{*} \left[ \frac{1}{M} \sum_{m=1}^{M} \left({\bf Q}(z) \right)^{m,m} \right] {\bf b}_2  \right)   \leq  C(z) \, \kappa^{4} \,  \frac{L}{MN} 
\end{align}
where $C(z)$ can be written as $C(z) = P_1(|z|) P_2\left(\frac{1}{\mathrm{Im}(z)} \right)$ for some nice polynomials $P_1$ and $P_2$. 
Moreover, if ${\bf G}$ is a $N \times N$ deterministic
matrix verifying $\sup_{N} \| {\bf G} \| \leq \kappa$,  
the following evaluations hold: 
\begin{align}
\label{eq:vartraceQW}
\mathrm{Var}\left( \frac{1}{ML} \mathrm{Tr}\left( {\bf A} {\bf Q}(z) {\bf W} {\bf G} {\bf W}^{*} \right) \right)  \leq  C(z) \, \kappa^{4} \, \frac{1}{MN} \\
\label{eq:varentreesQW} 
\mathrm{Var}\left(  {\bf a}_1^{*} {\bf Q}(z) {\bf W} {\bf G} {\bf W}^{*} {\bf a}_2 \right)  \leq  C(z) \, \kappa^{6} \,  \frac{L}{N} \\
\label{eq:vartraceentreesQW} 
\mathrm{Var}\left(  {\bf b}_1^{*} \left[ \frac{1}{M} \sum_{m=1}^{M} \left({\bf Q}(z) {\bf W} {\bf G} {\bf W}^{*} \right)^{m,m} \right] {\bf b}_2  \right)   \leq  C(z) \, \kappa^{6} \,\frac{L}{MN} 
\end{align}
where $C(z)$ can be written as above. 
\end{proposition}
{\bf Proof.} We first establish (\ref{eq:vartraceQ}) and denote by $\xi$ the random variable 
$\xi = \frac{1}{ML} \mathrm{Tr}\left( {\bf A} {\bf Q}(z) \right)$. As the various entries of 2 different blocks ${\bf W}^{m_1}, {\bf W}^{m_2}$ 
are independent, the Poincaré-Nash inequality can be written as  
\begin{eqnarray}
\label{eq:forme-nash}
\mathrm{Var} \, \xi \; & \leq & \sum_{m,i_1,i_2,j_1,j_2} \mathrm{E} \left[ \left( \frac{\partial \xi}{\partial \overline{{\bf W}}_{i_1,j_1}^{m}} \right)^{*} 
\mathbb{E}\left({\bf W}_{i_1,j_1}^{m} \overline{{\bf W}}_{i_2,j_2}^{m} \right)  \frac{\partial \xi}{\partial \overline{{\bf W}}_{i_2,j_2}^{m}} \right] +\\
   &   &  \sum_{m,i_1,i_2,j_1,j_2} \mathrm{E} \left[  \frac{\partial \xi}{\partial {\bf W}_{i_1,j_1}^{m}}
\mathbb{E}\left({\bf W}_{i_1,j_1}^{m} \overline{{\bf W}}_{i_2,j_2}^{m} \right) \left( \frac{\partial \xi}{\partial 
{\bf W}_{i_2,j_2}^{m}}  \right)^{*} \right] 
\label{eq:terme2-nash}
\end{eqnarray}
In the following, we just evaluate the right hand side of (\ref{eq:forme-nash}), denoted by $\beta$, because 
the behaviour of the term defined by (\ref{eq:terme2-nash}) can be established similarly. 
It is easy to check that 
$$
\frac{\partial {\bf Q}}{\partial \overline{{\bf W}}_{i,j}^{m}} = - {\bf Q} {\bf W} {\bf e}_j ({\bf f}_{i}^{m})^{T} {\bf Q}
$$
so that 
$$
\frac{\partial \xi}{\partial \overline{{\bf W}}_{i,j}^{m}} = -\frac{1}{ML} \mathrm{Tr} \left( {\bf A} {\bf Q} {\bf W} {\bf e}_j ({\bf f}_{i}^{m})^{T} {\bf Q} \right)
$$
which can also be written $-\frac{1}{ML} ({\bf f}^{m}_i)^{T} {\bf Q} {\bf A} {\bf Q} {\bf W} {\bf e}_j$. 
We recall that $\mathbb{E}\left( {\bf W}_{i_1,j_1}^{m} \overline{{\bf W}}_{i_2,j_2}^{m} \right) = 
\frac{\sigma^{2}}{N}\delta(i_1 - i_2 = j_2 - j_1)$ (see (\ref{eq:correlation-W})). Therefore, $\beta$ is equal to  
the mathematical expectation of the term 
$$
 \frac{1}{(ML)^{2}} \frac{\sigma^{2}}{N} \sum_{m,i_1,i_2,j_1,j_2} 
\delta(j_2-j_1 = i_1 - i_2) {\bf e}_{j_1}^{T} {\bf W}^{*} {\bf Q}^{*} {\bf A}^{*} {\bf Q}^{*} {\bf f}^{m}_{i_1}  
({\bf f}^{m}_{i_2})^{T} {\bf Q} {\bf A} {\bf Q} {\bf W} {\bf e}_{j_2}
$$
We put $u=i_1-i_2$ and remark that  $\sum_{m, i_1-i_2=u} {\bf f}^{m}_{i_1}  ({\bf f}^{m}_{i_2})^{T} = {\bf I}_M \otimes {\bf J}_L^{*u}$. 
We thus obtain that 
$$
\beta =  \frac{1}{(ML)^{2}} \frac{\sigma^{2}}{N} \, \mathbb{E} \left[\sum_{u=-(L-1)}^{L-1} \sum_{j_2-j_1=u}  {\bf e}_{j_1}^{T} {\bf W}^{*} {\bf Q}^{*} {\bf A}^{*}
{\bf Q}^{*} ({\bf I}_M \otimes {\bf J}_L^{*u}) {\bf Q} {\bf A} {\bf Q} {\bf W} {\bf e}_{j_2} \right]
$$
Using that $\sum_{j_2-j_1=u} {\bf e}_{j_2} {\bf e}_{j_1}^{T} = {\bf J}_N^{*u}$, we get that
$$
\beta = \frac{1}{ML} \frac{\sigma^{2}}{N} \, \mathbb{E} \left[\sum_{u=-(L-1)}^{L-1} \frac{1}{ML} \mathrm{Tr} \left( {\bf Q} {\bf A} {\bf Q} {\bf W} {\bf J}_N^{*u}
  {\bf W}^{*} {\bf Q}^{*} {\bf A}^{*} {\bf Q}^{*}
({\bf I}_M \otimes {\bf J}_L^{*u}) \right) \right]
$$
If ${\bf B}$ is a $ML \times N$ matrix, the Schwartz inequality as well as the inequality $(xy)^{1/2} \leq 1/2(x+y)$ lead to 
$$
\left| \frac{1}{ML} \mathrm{Tr} \left( {\bf B} {\bf J}_N^{*u} {\bf B}^{*} ({\bf I}_M \otimes {\bf J}_L^{*u}) \right) \right| \leq 
\frac{1}{2ML} \mathrm{Tr} \left( {\bf B} {\bf J}_N^{*u} {\bf J}_N^{u} {\bf B}^{*} \right) + \frac{1}{2ML} \mathrm{Tr} \left( {\bf B}^{*} ({\bf I}_M \otimes {\bf J}_L^{*u} {\bf J}_L^{u} )
 {\bf B} \right)
$$
It is clear that matrices ${\bf J}_N^{*u} {\bf J}_N^{u}$ and ${\bf J}_L^{*u} {\bf J}_L^{u}$ are less than ${\bf I}_N$ and ${\bf I}_L$ respectively. Therefore, 
\begin{equation}
\label{eq:schwartz-B}
\left| \frac{1}{ML} \mathrm{Tr} \left( {\bf B} {\bf J}_N^{*u} {\bf B}^{*} ({\bf I}_M \otimes {\bf J}_L^{*u}) \right) \right|  \leq \frac{1}{ML} \mathrm{Tr} \left( {\bf B} {\bf B}^{*} \right)
\end{equation}
Using (\ref{eq:schwartz-B}) for ${\bf B} =  {\bf Q} {\bf A} {\bf Q} {\bf W}$ for each $u$ leads to 
$$
\beta \leq \frac{\sigma^{2}}{MN} \mathbb{E} \left[ \frac{1}{ML} \mathrm{Tr} \left( {\bf Q} {\bf A} {\bf Q} {\bf W} {\bf W}^{*} {\bf Q}^{*} {\bf A}^{*} {\bf Q}^{*} \right)
\right]
$$
The resolvent identity (\ref{eq:equation-resolvente}) can also be written as 
${\bf Q} {\bf W} {\bf W}^{*} = {\bf I} + z {\bf Q}$. This implies that the greatest eigenvalue of
$ {\bf Q} {\bf W} {\bf W}^{*} {\bf Q}^{*}$ coincides with the greatest eigenvalue of $({\bf I} + z {\bf Q}) {\bf Q}^{*}$ which is itself less 
than $\| {\bf Q} \| + |z| \| {\bf Q} \|^{2}$. As $\| {\bf Q} \| \leq \frac{1}{\mathrm{Im}z}$, we obtain 
that 
\begin{equation}
\label{eq:inegalite-QWW*Q*}
 {\bf Q} {\bf W} {\bf W}^{*} {\bf Q}^{*} \leq \frac{1}{\mathrm{Im}z} \left( 1 + \frac{|z|}{\mathrm{Im}z} \right) \, {\bf I}.
\end{equation}
Therefore, it holds that 
\begin{equation}
\label{eq:inequation-varxi}
\beta \leq \frac{1}{\mathrm{Im}z} \left( 1 + \frac{|z|}{\mathrm{Im}z} \right) \, \frac{1}{MN} \mathbb{E} \left[ \frac{1}{ML} \mathrm{Tr} \left( {\bf Q} {\bf A} {\bf A}^{*} {\bf Q}^{*} \right) \right]
\end{equation}
We eventually obtain that 
$$
\beta  \leq   \, \kappa^{2} \, \frac{1}{MN} C(z) \, \frac{1}{(\mathrm{Im}z)^{3}} \left( 1 + \frac{|z|}{\mathrm{Im}z} \right)
$$
The conclusion follows from the observation that 
$$
\frac{1}{(\mathrm{Im}z)^{3}} \left( 1 + \frac{|z|}{\mathrm{Im}z} \right) \leq \left[ \frac{1}{(\mathrm{Im}z)^{3}} + \frac{1}{(\mathrm{Im}z)^{4}} \right] \, (|z|+1)
$$ 
In order to prove (\ref{eq:varentreesQ}) and (\ref{eq:vartraceentreesQ}), we remark that 
\begin{eqnarray*}
{\bf a}_1^{*} \, {\bf Q} \, {\bf a}_2 & = & ML \, \frac{1}{ML} \mathrm{Tr}\left( {\bf Q} {\bf a}_2 {\bf a}_1^{*} \right) \\  
{\bf b}_1^{*} \left[ \frac{1}{M} \sum_{m=1}^{M} \left({\bf Q}(z) \right)^{m,m} \right] {\bf b}_2  & =  & L \, 
\frac{1}{ML} \mathrm{Tr} \left( {\bf Q}({\bf I}_M \otimes {\bf b}_2 {\bf b}_1^{*}) \right)
\end{eqnarray*}
(\ref{eq:varentreesQ}) and (\ref{eq:vartraceentreesQ}) follow immediately from this and inequality (\ref{eq:inequation-varxi})
used in the case ${\bf A} =  {\bf a}_2 {\bf a}_1^{*}$ and ${\bf A} = {\bf I}_M \otimes {\bf b}_2 {\bf b}_1^{*}$ respectively. \\

We finally provide a sketch of proof of (\ref{eq:vartraceQW}), and omit the proof 
of (\ref{eq:vartraceentreesQW}) and (\ref{eq:varentreesQW}) which can be obtained as above. We still denote by $\xi$ the random variable
$\xi = \frac{1}{ML} \mathrm{Tr} \left( {\bf Q}(z) {\bf W} {\bf G} {\bf W}^{*} \right)$, and only evaluate the behaviour of the
right hand side $\beta$ of (\ref{eq:forme-nash}). After easy calculations using  tricks similar to those used in the course 
of the proof of (\ref{eq:vartraceQ}), we obtain that 
\begin{eqnarray}
\label{eq:terme1-vartraceQW}
\beta &  \leq & \frac{2 \sigma^{2}}{MN} \mathbb{E} \left[ \frac{1}{ML} \mathrm{Tr} \left( {\bf Q} {\bf W} {\bf G} {\bf W}^{*} {\bf A} {\bf Q} {\bf W} {\bf W}^{*} {\bf Q}^{*} {\bf A}^{*} {\bf W} 
{\bf G}^{*} {\bf W}^{*} {\bf Q}^{*} \right) \right] + \\
\label{eq:terme2-vartraceQW} 
   &   &  \frac{2 \sigma^{2}}{MN} \mathbb{E} \left[  \frac{1}{ML} \mathrm{Tr} \left( {\bf G}^{*} {\bf W}^{*} {\bf Q}^{*} {\bf A}^{*} {\bf A} {\bf Q} {\bf W} {\bf G} \right) \right] 
\end{eqnarray} 
The term defined by (\ref{eq:terme2-vartraceQW}) is easy to handle because 
${\bf Q}^{*} {\bf A}^{*} {\bf A} {\bf Q} \leq \frac{\kappa^{2}}{(\mathrm{Im}(z))^{2}} \, {\bf I}$. Therefore, (\ref{eq:terme2-vartraceQW})
is less than \\ $\frac{2 \sigma^{2} \kappa^{2}}{(\mathrm{Im}(z))^{2}} \, \frac{1}{MN} \, \mathbb{E} \left[ \frac{1}{ML} \mathrm{Tr} \left( {\bf W} {\bf G} {\bf G}^{*} {\bf W}^{*} \right) \right]$
which is itself lower bounded by $ \frac{1}{MN}  \, \frac{2 \sigma^{4} \kappa^{4}}{(\mathrm{Im}(z))^{2}}$
because $\mathbb{E}\left( \frac{1}{ML} \mathrm{Tr}({\bf W}{\bf W}^{*}) \right) = \sigma^{2}$. 
To evaluate the righthandside of (\ref{eq:terme1-vartraceQW}),
we use (\ref{eq:inegalite-QWW*Q*}) twice, and obtain immediately that is less than $\frac{C(z) \kappa^{4}}{MN}$.

\section{Expression of matrix $\mathbb{E}({\bf Q})$ obtained using the integration by parts formula}
\label{sec:novikov}
In this section, we use the integration by parts formula in order to express 
$\mathbb{E}\left({\bf Q}(z)\right)$ as a term which will appear to be close from $t(z) {\bf I}_{ML}$ 
where we recall that $t(z)$ represents the Stieltjes transform of the Marcenko-Pastur distribution 
$\mu_{\sigma^{2}, c_N}$. For this, 
we have first to introduce useful matrix valued functions of the complex variable 
$z$ and to study their properties. 
\begin{lemma}
\label{le:proprietes-H-R}
For each $z \in \mathbb{C}^{+}$, matrix ${\bf I}_N + \sigma^{2} c_N \mathcal{T}_{N,L}^{(M)}\left(\mathbb{E}({\bf Q}(z))\right)$
is invertible. We denote by ${\bf H}(z)$ its inverse, i.e. 
\begin{equation}
\label{eq:def-H}
{\bf H}(z) = \left[ {\bf I}_N + \sigma^{2}  c_N  \mathcal{T}^{(M)}_{N,L}\left( \mathbb{E}({\bf Q}(z)) \right) \right]^{-1}
\end{equation}
Then, function $z \rightarrow {\bf H}(z)$ is holomorphic in $\mathbb{C}^{+}$ and verifies
\begin{equation}
\label{eq:borne-HH*}
{\bf H}(z) {\bf H}(z)^{*} \leq \left( \frac{|z|}{\mathrm{Im}z} \right)^{2} \, {\bf I}_N
\end{equation}
Moreover, for each $z \in \mathbb{C}^{+}$, matrix $ -z \, {\bf I} + \sigma^{2} \, \mathcal{T}_{L,L}\left({\bf H}(z)\right)$
is invertible. We denote by ${\bf R}(z)$ its inverse, i.e. 
\begin{equation}
\label{eq:def-R}
{\bf R}(z) = \left[- z {\bf I}_L + \sigma^{2} \mathcal{T}_{L,L}({\bf H}(z)) \right]^{-1}
\end{equation}
Then, function $z \rightarrow {\bf R}(z)$ is holomorphic in $\mathbb{C}^{+}$, and it exists a positive matrix valued 
measure ${\bs \mu}_{{\bf R}}$ carried by $\mathbb{R}^{+}$, satisfying ${\bs \mu}_{{\bf R}}(\mathbb{R}^{+}) = {\bf I}_L$, and for which 
$$
{\bf R}(z) = \int_{\mathbb{R}^{+}} \frac{d \, {\bs \mu}_{{\bf R}}(\lambda)}{\lambda - z}
$$
Finally, it holds that 
\begin{equation}
\label{eq:borne-RR*}
{\bf R}(z) {\bf R}(z)^{*} \leq \left( \frac{1}{\mathrm{Im}z} \right)^{2} \, {\bf I}_L
\end{equation}
\end{lemma}
{\bf Proof.} The proof is sketched in the appendix. \\

In order to be able the integration by parts formula, we use the identity (\ref{eq:equation-resolvente})
which implies that 
\begin{equation}
\label{eq:equation-resolvente-bis}
\mathbb{E}\left[ {\bf Q}^{m_1,m_2}_{i_1,i_2} \right]= -\frac{1}{z} \delta(i_1-i_2) \delta(m_1-m_2) + \frac{1}{z} \mathbb{E} \left[\left({\bf Q} {\bf W} {\bf W}^{*}\right)^{m_1,m_2}_{i_1,i_2} \right]
\end{equation}
We express $\left({\bf Q} {\bf W} {\bf W}^{*}\right)^{m_1,m_2}_{i_1,i_2}$ as 
$$
\left({\bf Q} {\bf W} {\bf W}^{*}\right)^{m_1,m_2}_{i_1,i_2} = \sum_{j=1}^{N} \left({\bf Q} {\bf w}_j {\bf w}_j^{*}\right)^{m_1,m_2}_{i_1,i_2} = \sum_{j=1}^{N} \left({\bf Q} {\bf w}_j\right)_{i_1}^{m_1} \overline{{\bf W}}_{i_2,j}^{m_2}
$$
where we recall that $({\bf w}_j)_{j=1, \ldots, N}$ represent the columns of ${\bf W}$. In order to be able to evaluate $\mathbb{E}\left[ \left( {\bf Q} {\bf w}_j {\bf w}_j^{*} \right)_{i_1,i_2}^{m_1,m_2} \right]$, 
it is necessary to express  $\mathbb{E}\left[ \left( {\bf Q} {\bf w}_k {\bf w}_j^{*} \right)_{i_1,i_2}^{m_1,m_2} \right]= 
\mathbb{E}\left[ \left( {\bf Q} {\bf w}_k \right)_{i_1}^{m_1} \left( {\bf w}_j^{*} \right)_{i_2}^{m_2} \right]$ for each 
pair $(k,j)$. For this, we use the identity
$$
\mathbb{E}\left[ \left( {\bf Q} {\bf w}_k \right)_{i_1}^{m_1} \left( {\bf w}_j^{*} \right)_{i_2}^{m_2} \right] =
\sum_{i_3,m_3} \mathbb{E} \left( {\bf Q}_{i_1,i_3}^{m_1,m_3} {\bf W}_{i_3,k}^{m_3} \overline{{\bf W}}_{i_2,j}^{m_2} \right)
$$
and use the integration by parts formula
$$
\mathbb{E} \left( {\bf Q}_{i_1,i_3}^{m_1,m_3} {\bf W}_{i_3,k}^{m_3} \overline{{\bf W}}_{i_2,j}^{m_2} \right) = 
\sum_{i^{'}, j^{'}} \mathbb{E}\left( {\bf W}_{i_3,k}^{m_3} \overline{{\bf W}}_{i^{'},j^{'}}^{m_3} \right) 
\mathbb{E} \left[ \frac{ \partial \left(  {\bf Q}_{i_1,i_3}^{m_1,m_3}  \overline{{\bf W}}_{i_2,j}^{m_2} \right)}{\partial \overline{{\bf W}}_{i^{'},j^{'}}^{m_3}} \right]
$$
It is easy to check that 
$$
\frac{ \partial \left(  {\bf Q}_{i_1,i_3}^{m_1,m_3}  \overline{{\bf W}}_{i_2,j}^{m_2} \right)}{\partial \overline{{\bf W}}_{i^{'},j^{'}}^{m_3}} =
 {\bf Q}_{i_1,i_3}^{m_1,m_3} \delta(m_2=m_3) \delta(i^{'}=i_2) \delta(j=j^{'}) - 
 \left({\bf Q} {\bf w}_{j^{'}}\right)_{i_1}^{m_1}  {\bf Q}_{i^{'},i_3}^{m_3,m_3} \overline{{\bf W}}_{i_2,j}^{m_2}
$$
(\ref{eq:def-Wmij}) implies that  $\mathbb{E}\left( {\bf W}_{i_3,k}^{m_3} \overline{{\bf W}}_{i^{'},j^{'}}^{m_3} \right) = 
\frac{\sigma^{2}}{N}\delta(i_3-i^{'} = j^{'} - k)$. Therefore, we obtain that 
\begin{align*}
\mathbb{E} \left( {\bf Q}_{i_1,i_3}^{m_1,m_3} {\bf W}_{i_3,k}^{m_3} \overline{{\bf W}}_{i_2,j}^{m_2} \right) =
\frac{\sigma^{2}}{N} \delta(i_3-i_2=j-k) \delta(m_2=m_3) \mathbb{E} \left( {\bf Q}_{i_1,i_3}^{m_1,m_3} \right)  \\ - \; 
\frac{\sigma^{2}}{N} \sum_{i^{'}, j^{'}}   \delta(i_3-i^{'}=j^{'}-k) \, \mathbb{E} \left[ \left({\bf Q} {\bf w}_{j^{'}}\right)_{i_1}^{m_1}
\left( {\bf w}_j^{*} \right)_{i_2}^{m_2}  {\bf Q}_{i^{'},i_3}^{m_3,m_3} \right]
\end{align*}
and that 
\begin{align*}
\mathbb{E}\left[ \left( {\bf Q} {\bf w}_k \right)_{i_1}^{m_1} \left( {\bf w}_j^{*} \right)_{i_2}^{m_2} \right] =
\frac{\sigma^{2}}{N} \sum_{i_3,m_3} \delta(i_3-i_2=j-k) \delta(m_2=m_3) \mathbb{E} \left( {\bf Q}_{i_1,i_3}^{m_1,m_3} \right)  \\ - \; 
 \frac{\sigma^{2}}{N} \sum_{i_3,m_3}  \sum_{i^{'}, j^{'}}  \delta(i_3-i^{'} = j^{'} - k) \, \mathbb{E} \left[ \left({\bf Q} {\bf w}_{j^{'}}\right)_{i_1}^{m_1}
\left( {\bf w}_j^{*} \right)_{i_2}^{m_2}  {\bf Q}_{i^{'},i_3}^{m_3,m_3} \right]
\end{align*}
We put $i=i^{'} - i_3$ in the above sum, and get that
\begin{align*}
\mathbb{E}\left[ \left( {\bf Q} {\bf w}_k \right)_{i_1}^{m_1} \left( {\bf w}_j^{*} \right)_{i_2}^{m_2} \right] =
\frac{\sigma^{2}}{N}  \mathbb{E} \left( {\bf Q}_{i_1,i_2-(k-j)}^{m_1,m_2} \right) \mathbb{1}_{1 \leq i_2-(k-j) \leq L}  
 \\
- \; \sigma^{2} c_N \sum_{i=-(L-1)}^{L-1} \mathbb{1}_{1 \leq k-i \leq N} \, 
\mathbb{E} \left[ \left( {\bf Q} {\bf w}_{k-i} \right)_{i_1}^{m_1} \left( {\bf w}_j^{*} \right)_{i_2}^{m_2} 
\frac{1}{ML} \sum_{i^{'}-i_3=i} \sum_{m_3}  {\bf Q}_{i^{'},i_3}^{m_3,m_3} \right]
\end{align*}
or, using the definition (\ref{eq:def-calt(A)}), 
\begin{align}
\label{eq:intermediaire-1}
\mathbb{E}\left[ \left( {\bf Q} {\bf w}_k \right)_{i_1}^{m_1} \left( {\bf w}_j^{*} \right)_{i_2}^{m_2} \right] =
\frac{\sigma^{2}}{N}  \mathbb{E} \left( {\bf Q}_{i_1,i_2-(k-j)}^{m_1,m_2} \right) \mathbb{1}_{1 \leq i_2-(k-j) \leq L} 
 \\
\nonumber
- \; \sigma^{2} c_N \sum_{i=-(L-1)}^{L-1} \mathbb{1}_{1 \leq k-i \leq N} \, 
\mathbb{E} \left[ \tau^{(M)}({\bf Q})(i)\left( {\bf Q} {\bf w}_{k-i} \right)_{i_1}^{m_1} \left( {\bf w}_j^{*} \right)_{i_2}^{m_2} \right]
\end{align}
Setting $u=k-i$, the second term of the righthandside of the above equation can also be written as
$$
- \;  \sigma^{2} c_N \mathbb{E} \left[ \sum_{u=1}^{N} \tau^{(M)}({\bf Q})(k-u) \mathbb{1}_{-(L-1) \leq k-u \leq L-1} \left({\bf Q} {\bf w}_u \right)_{i_1}^{m_1} \left( {\bf w}_j^{*} \right)_{i_2}^{m_2} \right]
$$
or, using the observation that $\tau^{(M)}({\bf Q})(k-u) \mathbb{1}_{-(L-1) \leq k-u \leq L-1} = \left(\mathcal{T}^{(M)}_{N,L}({\bf Q}) \right)_{k,u}$
(see Eq. (\ref{eq:def-calT(A)-bis})), as 
$$
-\; \sigma^{2} c_N \; \mathbb{E} \left[ {\bf e}_k^{T} \; \mathcal{T}^{(M)}_{N,L}({\bf Q}) \; \left( \begin{array}{c} \left({\bf Q} {\bf w}_1 \right)_{i_1}^{m_1} \left( {\bf w}_j^{*} \right)_{i_2}^{m_2} \\
                                                                                            \left({\bf Q} {\bf w}_2 \right)_{i_1}^{m_1} \left( {\bf w}_j^{*} \right)_{i_2}^{m_2} \\ 
                                                                                            \vdots \\
                                                                                            \left({\bf Q} {\bf w}_N \right)_{i_1}^{m_1} \left( {\bf w}_j^{*} \right)_{i_2}^{m_2} \end{array} \right) \right]
$$
We express matrix ${\bf Q}$ as ${\bf Q} = \mathbb{E}({\bf Q}) + {\bf Q}^{\circ}$ and define the following $N \times N$ matrices $ {\bf A}_{i_1,i_2}^{m_1,m_2},   {\bf B}_{i_1,i_2}^{m_1,m_2}, \Upsilonbs_{i_1,i_2}^{m_1,m_2}$
\begin{align*}
\left( {\bf A}_{i_1,i_2}^{m_1,m_2} \right)_{k,j} =   \mathbb{E} \left[ \left({\bf Q} {\bf w}_k \right)_{i_1}^{m_1} \left( {\bf w}_j^{*} \right)_{i_2}^{m_2} \right] \\
\left( {\bf B}_{i_1,i_2}^{m_1,m_2} \right)_{k,j} =  \mathbb{E} \left[ {\bf Q}^{m_1,m_2}_{i_1,i_2-(k-j)} \mathbb{1}_{1 \leq i_2 -(k-j) \leq L} \right]
\end{align*}
$$
 \Upsilonbs_{i_1,i_2}^{m_1,m_2}   = -\sigma^{2} c_N \,   \mathbb{E} \left[
\mathcal{T}^{(M)}_{N,L}({\bf Q^{\circ}}) \; \left( \begin{array}{c} \left({\bf Q} {\bf w}_1 \right)_{i_1}^{m_1}  \\ \left({\bf Q} {\bf w}_2 \right)_{i_1}^{m_1} \\  \vdots  \\
 \left({\bf Q} {\bf w}_N \right)_{i_1}^{m_1}  \end{array} \right) \; \left( \begin{array}{cccc} \left( {\bf w}_1^{*} \right)_{i_2}^{m_2} & \left( {\bf w}_2^{*} \right)_{i_2}^{m_2} & \ldots & 
 \left( {\bf w}_N^{*} \right)_{i_2}^{m_2}  \end{array} \right) \right]
$$
 We notice that matrix
$$
 \left[
 \left( \begin{array}{c} \left({\bf Q} {\bf w}_1 \right)_{i_1}^{m_1}  \\ \left({\bf Q} {\bf w}_2 \right)_{i_1}^{m_1} \\  \vdots  \\
 \left({\bf Q} {\bf w}_N \right)_{i_1}^{m_1}  \end{array} \right) \; \left( \begin{array}{cccc} \left( {\bf w}_1^{H} \right)_{i_2}^{m_2} & \left( {\bf w}_2^{H} \right)_{i_2}^{m_2} & \ldots & 
 \left( {\bf w}_N^{H} \right)_{i_2}^{m_2}  \end{array} \right) \right]
$$
can also be written as
$$
\left( \begin{array}{c} {\bf w}_1^{T} {\bf Q}^{T} \\ \vdots \\ {\bf w}_N^{T} {\bf Q}^{T} \end{array} \right) \left( {\bf f}_{i_1}^{m_1} \right)
\left( {\bf f}_{i_2}^{m_2} \right)^{T} \left( \overline{{\bf w}}_1, \ldots, \overline{{\bf w}}_N \right)
$$
or as 
$$
{\bf W}^{T} {\bf Q}^{T} \left( {\bf f}_{i_1}^{m_1} \right) \left( {\bf f}_{i_2}^{m_2} \right)^{T} \overline{{\bf W}}
$$
Therefore, 
\begin{equation}
\label{eq:expre-upsilon}
 \Upsilonbs_{i_1,i_2}^{m_1,m_2}   = -\sigma^{2} c_N \, \mathbb{E} \left[ \mathcal{T}^{(M)}_{N,L}({\bf Q}^{\circ}) \; {\bf W}^{T} {\bf Q}^{T} \left( {\bf f}_{i_1}^{m_1} \right) \left( {\bf f}_{i_2}^{m_2} \right)^{T} \overline{{\bf W}} \right]
\end{equation}
It is useful to notice that matrix $ {\bf B}_{i_1,i_2}^{m_1,m_2}$ is a band Toeplitz matrix whose $(k,l)$ element is zero if $|k-l| \geq L$. It is clear that Eq. (\ref{eq:intermediaire-1}) is equivalent to
$$
\left[ {\bf I}_N + \sigma^{2}  c_N  \mathcal{T}^{(M)}_{N,L}\left( \mathbb{E}({\bf Q}) \right) \right] \;  {\bf A}_{i_1,i_2}^{m_1,m_2}  = 
\frac{\sigma^{2}}{N} {\bf B}_{i_1,i_2}^{m_1,m_2}  + \; \Upsilonbs_{i_1,i_2}^{m_1,m_2} 
$$
Lemma  \ref{le:proprietes-H-R} implies that matrix $\left[ {\bf I}_N + \sigma^{2}  c_N  \mathcal{T}^{(M)}_{N,L}\left( \mathbb{E}({\bf Q}(z)) \right) \right]$ is invertible for each $z \in \mathbb{C}^{+}$, and we recall that its inverse is denoted ${\bf H}(z)$.  
We obtain that 
\begin{equation}
\label{eq:expre-A}
{\bf A}_{i_1,i_2}^{m_1,m_2} = \frac{\sigma^{2}}{N} \; {\bf H}  \;  {\bf B}_{i_1,i_2}^{m_1,m_2} + \; {\bf H} \Upsilonbs_{i_1,i_2}^{m_1,m_2}
\end{equation}
The term $\mathbb{E} \left( {\bf Q} {\bf W} {\bf W}^{*} \right)_{i_1,i_2}^{m_1,m_2}$ coincides with $\mathrm{Tr} \left( {\bf A}_{i_1,i_2}^{m_1,m_2} \right)$, so that 
\begin{equation}
\label{eq:intermediaire-3}
\mathbb{E} \left( {\bf Q} {\bf W} {\bf W}^{*} \right)_{i_1,i_2}^{m_1,m_2} = \sigma^{2} \; \frac{1}{N} \mathrm{Tr} \left( {\bf H}  \;  {\bf B}_{i_1,i_2}^{m_1,m_2} \right) + \; 
\mathrm{Tr} \left( {\bf H}  \; \Upsilonbs_{i_1,i_2}^{m_1,m_2} \right)
\end{equation}
As matrix ${\bf B}_{i_1,i_2}^{m_1,m_2}$ is Toeplitz, it holds that (see Eq. (\ref{eq:Tr(AB)}))
$$
 \frac{1}{N} \mathrm{Tr} \left( {\bf H}  \;  {\bf B}_{i_1,i_2}^{m_1,m_2} \right) = \sum_{u=-(N-1)}^{N-1} \tau({\bf H})(u) \mathbb{E}\left( {\bf Q}_{i_1, i_2+u}^{m_1,m_2} \right) \mathbb{1}_{1 \leq i_2+u \leq L}
$$
which also coincides with
$$
 \frac{1}{N} \mathrm{Tr} \left( {\bf H}  \;  {\bf B}_{i_1,i_2}^{m_1,m_2} \right) = \sum_{u=-(L-1)}^{L-1} \tau({\bf H})(u) \mathbb{E}\left( {\bf Q}_{i_1, i_2+u}^{m_1,m_2} \right) \mathbb{1}_{1 \leq i_2+u \leq L}
$$
because $\mathbb{1}_{1 \leq i_2+u \leq L} = 0$ if $|u| \geq L$. Setting $v=i_2+u$, this term can be written as
$$
 \frac{1}{N} \mathrm{Tr} \left( {\bf H}  \;  {\bf B}_{i_1,i_2}^{m_1,m_2} \right) = \sum_{v=1}^{L}  \mathbb{E}\left( {\bf Q}_{i_1, v}^{m_1,m_2} \right) \tau({\bf H})(v-i_2) 
$$
or, using definition (\ref{eq:def-calT(A)-bis}), as
\begin{eqnarray*}
\frac{1}{N} \mathrm{Tr} \left( {\bf H}  \;  {\bf B}_{i_1,i_2}^{m_1,m_2} \right) & =  & \sum_{v=1}^{L}  \mathbb{E}\left( {\bf Q}_{i_1, v}^{m_1,m_2} \right) \; \left( \mathcal{T}_{L,L}({\bf H}) \right)_{v,i_2} \\
                                                                                &  =  &  \left( \mathbb{E}({\bf Q}^{m_1,m_2})  \mathcal{T}_{L,L}({\bf H})  \right)_{i_1,i_2}
\end{eqnarray*}      
 Eq. (\ref{eq:intermediaire-3}) eventually leads to 
\begin{equation}
\label{eq:intermediaire-4}
\mathbb{E} \left[ \left( {\bf Q} {\bf W} {\bf W}^{*} \right)^{m_1,m_2} \right] =  \sigma^{2} \mathbb{E}({\bf Q}^{m_1,m_2})  \mathcal{T}_{L,L}({\bf H})  + \; \Upsilonbs({\bf H})^{m_1,m_2}
\end{equation}
where, for each $N \times N$ matrix ${\bf F}$, $\Upsilonbs({\bf F})$ represents the $ML \times ML$ matrix defined by  
\begin{equation}
\label{eq:def-Gamma}
 \Upsilonbs({\bf F})^{m_1,m_2}_{i_1,i_2} = \mathrm{Tr} \left( {\bf F}  \; \Upsilonbs_{i_1,i_2}^{m_1,m_2} \right)
\end{equation}
(\ref{eq:expre-upsilon}) implies that matrix $\Upsilonbs({\bf F})$ can be written as
\begin{equation}
\label{eq:expre-upsilonF}
\Upsilonbs({\bf F}) = -\sigma^{2} c_N \mathbb{E} \left[ {\bf Q} {\bf W} \left(\mathcal{T}_{N,L}^{(M)}({\bf Q}^{\circ})\right)^{T} {\bf F}^{T} {\bf W}^{*} \right]
\end{equation} 
By (\ref{eq:equation-resolvente}), it holds that $\left( {\bf Q} {\bf W} {\bf W}^{*} \right)^{m_1,m_2} = \delta(m_1=m_2) \, {\bf I}_L  + z {\bf Q}^{m_1,m_2}$. Therefore, 
we deduce from (\ref{eq:intermediaire-4}) that 
\begin{equation}
\label{eq:intermediaire-5}
 \mathbb{E}({\bf Q}^{m_1,m_2}) \; \left( - z {\bf I}_L + \sigma^{2}  \mathcal{T}_{L,L}({\bf H}) \right) = {\bf I}_L \delta(m_1=m_2) \, - \, \Upsilonbs({\bf H})^{m_1,m_2}
\end{equation}
By Lemma \ref{le:proprietes-H-R}, $- z {\bf I}_L + \sigma^{2}  \mathcal{T}_{L,L}({\bf H}(z))$ is invertible for $z \in \mathbb{C}^+$ and we recall 
that its inverse is denoted by ${\bf R}$. We thus obtain that 
\begin{equation}
\label{eq:expre-E(Q)}
\mathbb{E}({\bf Q}) = {\bf I}_M \otimes  {\bf R} + \Deltabs
\end{equation}
where $ \Deltabs$ is the $ML \times ML$ matrix defined by 
\begin{equation}
\label{eq:def-Delta}
 \Deltabs= -  \Upsilonbs({\bf H}) \, \left( {\bf I}_M \otimes  {\bf R} \right)
\end{equation}
The above evaluations also allow to obtain a similar expression of matrix $\mathbb{E}( {\bf Q} {\bf W} {\bf G} {\bf W}^{*})$ 
where ${\bf G}$ is a $N \times N$ matrix. For this, we express $\mathbb{E}\left[ ( {\bf Q} {\bf W} {\bf G} {\bf W}^{*})^{m_1,m_2}_{i_1,i_2} \right]$
as 
$$
\mathbb{E}\left[ ( {\bf Q} {\bf W} {\bf G} {\bf W}^{*})^{m_1,m_2}_{i_1,i_2} \right] = \sum_{(k,j)=1}^{N} {\bf G}_{k,j} 
\mathbb{E}\left[ ( {\bf Q} {\bf w}_k)_{i_1}^{m_1} ({\bf w}_j^{*})^{m_2}_{i_2} \right] 
$$
or equivalently as 
$$
\mathbb{E}\left[ ( {\bf Q} {\bf W} {\bf G} {\bf W}^{*})^{m_1,m_2}_{i_1,i_2} \right] = \mathrm{Tr} \left( {\bf G}^{T} {\bf A}^{m_1,m_2}_{i_1,i_2} \right)
$$
Therefore, using (\ref{eq:expre-A}), it holds that 
$$
\mathbb{E}\left[  \left({\bf Q} {\bf W} {\bf G} {\bf W}^{*}\right)_{i_1,i_2}^{m_1,m_2} \right] = \frac{\sigma^{2}}{N} \, \mathrm{Tr} \left( {\bf G}^{T} {\bf H} 
{\bf B}^{m_1,m_2}_{i_1,i_2} \right) + \;  \mathrm{Tr} \left( {\bf G}^{T} {\bf H} {\bs \Upsilon}^{m_1,m_2}_{i_1,i_2} \right)
$$
Replacing matrix ${\bf H}$ by matrix ${\bf G}^{T} {\bf H}$ in the above calculations, we obtain that 
$$
\mathbb{E}\left[  {\bf Q} {\bf W} {\bf G} {\bf W}^{*} \right] = \sigma^{2} \mathbb{E}({\bf Q})  \; \left( {\bf I}_M \otimes
\mathcal{T}_{L,L}({\bf G}^{T} {\bf H}) \right) \;  + \;   {\bs \Upsilon}({\bf G}^{T} {\bf H})
$$
Using (\ref{eq:expre-E(Q)}), we eventually get that
\begin{equation}
\label{eq:expre-E(QWGW*)}
\mathbb{E} \left( {\bf Q} {\bf W} {\bf G} {\bf W}^{*} \right) = \sigma^{2} \left( {\bf I}_M \otimes {\bf R} \, \mathcal{T}_{L,L}({\bf G}^{T} {\bf H}) \right) + \sigma^{2} {\bs \Delta} \, \left( {\bf I}_M \otimes \mathcal{T}_{L,L}({\bf G}^{T} {\bf H}) \right)  + \; {\bs \Upsilon}({\bf G}^{T} {\bf H})
\end{equation}

\section{Controls of the error term ${\bs \Delta}$}
\label{sec:E(Q)-R}
In this section, we evaluate the behaviour of various terms depending on ${\bs \Delta}$, 
i.e. normalized traces $\frac{1}{ML} \mathrm{Tr} \Deltabs {\bf A}$, quadratic forms ${\bf a}_1^* {\bs \Delta} {\bf a}_2$, 
and quadratic forms of matrix $\hat{{\bs \Delta}} = \frac{1}{M} \sum_{m=1}^{M} {\bs \Delta}^{m,m}$. 
Using rough estimates based on the results of section \ref{sec:nash} and the Schwartz inequality, we establish that the normalized traces are $\mathcal{O}(\frac{L}{MN})$, and that two other terms are $\mathcal{O}(\sqrt{\frac{L}{M}} \frac{L}{N})$ and $\mathcal{O}(\frac{L^{3/2}}{MN})$
respectively.  We first establish the following proposition. 
\begin{proposition}
\label{prop:controle-Delta}
Let ${\bf A}$ be a $ML \times ML$ matrix satisfying
$\sup_{N} \| {\bf A} \| \leq \kappa$. Then, it holds 
that  
\begin{equation}                   
\label{eq:controle-Delta}
\left| \frac{1}{ML} \mathrm{Tr} \Deltabs {\bf A} \right| \leq  \kappa  \, \frac{L}{MN} \, C(z)
\end{equation}
where $C(z)$ can be written as $C(z) = P_1(|z|) \, P_2\left( (\mathrm{Im}z)^{-1} \right)$ for some 
nice polynomials $P_1$ and $P_2$. 
\end{proposition}
{\bf Proof.}
 As matrix ${\bf R}$ verifies $\| {\bf R} \| \leq \left( \mathrm{Im}z \right)^{-1}$, it is sufficient to establish (\ref{eq:controle-Delta}) 
when $\Deltabs$ 
 is replaced by $\Upsilonbs({\bf H})$. In order to simplify the notations, matrix  $\Upsilonbs({\bf H})$ is denoted by  $\Upsilonbs$ in this section. We denote by $\gamma$ the term $\gamma = \frac{1}{ML} \mathrm{Tr} \Upsilonbs {\bf A}$
 which is given by
 $$
 \gamma = \frac{1}{M} \sum_{m_1,m_2} \frac{1}{L} \sum_{i_1,i_2} \Upsilonbs_{i_1,i_2}^{m_1,m_2} {\bf A}_{i_2,i_1}^{m_2,m_1}
 $$
Using the expression (\ref{eq:expre-upsilonF}) of matrix $\Upsilon$, we obtain that $\gamma$ can be written as
$$
\gamma =  -\sigma^{2} \mathbb{E} \left[ \frac{1}{N} \mathrm{Tr} \left(  \left(\mathcal{T}_{N,L}^{(M)}({\bf Q}^{\circ})\right)^{T}  {\bf H}^{T} {\bf W}^{*} {\bf A} {\bf Q} {\bf W} 
 \right) \right]
$$
Using Eq. (\ref{eq:Tr(BT(A)}) and the identity $\tau^{(M)}\left(({\bf Q}^{\circ})^{T} \right)(-u) = \tau^{(M)}\left({\bf Q}^{\circ} \right)(u)$,  we get that 
\begin{equation}
\label{eq:expre-gamma}
\gamma = -\sigma^{2} c_N \mathbb{E} \left[ \sum_{u=-(L-1)}^{L-1} \tau^{(M)}({\bf Q}^{\circ})(u) \; \frac{1}{ML} \mathrm{Tr} \left( {\bf Q} {\bf W} {\bf J}_N^{u} {\bf H}^{T} {\bf W}^{*} {\bf A} \right) \right]
\end{equation}
(\ref{eq:vartraceQ}, \ref{eq:vartraceQW}) imply that 
$\mathbb{E} \left| \tau^{(M)}({\bf Q}^{\circ})(-u) \right|^{2}$ and $\mathrm{Var}\left( \frac{1}{ML} \mathrm{Tr} \left( {\bf Q} {\bf W} {\bf J}_N^{u} {\bf H}^{T} {\bf W} {\bf A} \right) \right)$ 
are upperbounded by terms of the form $\frac{C(z)}{MN}$ and $\kappa^{2} \frac{C(z)}{MN}$ respectively. The Cauchy-Schwartz inequality thus implies immediately (\ref{eq:controle-Delta}). \\

We now evaluate the behaviour of quadratic forms of matrix ${\bs \Delta}$ and of matrix $\hat{{\bs \Delta}}$. 
\begin{proposition}
\label{prop:entrees-Delta}
Let ${\bf a}_1$ and ${\bf a}_2$ 2 $ML$--dimensional vectors such that $\sup_{N} \| {\bf a}_i \| \leq \kappa$  
for $i=1,2$. Then, it holds that 
\begin{equation}
\label{eq:entrees-Delta}
{\bf a}_1^{*} {\bs \Delta} {\bf a}_2 \leq \kappa^{2} \, C(z) \, \sqrt{\frac{L}{M}} \frac{L}{N} 
\end{equation}
for each $z \in \mathbb{C}^{+}$, where $C(z)$ is as in Proposition \ref{prop:controle-Delta}. Let ${\bf b}_i$, $i=1,2$  be 2 deterministic $L$--dimensional vectors
such that $\sup_{N} \|{\bf b}_i\| < \kappa$. Then, it holds that 
\begin{equation}
\label{eq:entrees-trace-Delta}
\left| {\bf b}_1^{*} \left( \frac{1}{M} \sum_{m=1}^{M} {\bs \Delta}^{m,m} \right) {\bf b}_2 \right| \leq \kappa^{2} \, C(z) \frac{L^{3/2}}{MN} 
\end{equation}
\end{proposition}
{\bf Proof.} As above, it is sufficient to establish the proposition when ${\bs \Delta}$ is replaced 
by ${\bs \Upsilon}$. We first establish (\ref{eq:entrees-Delta}). We remark that ${\bf a}_1^{*} {\bs \Upsilon} {\bf a}_2 = ML \, \frac{1}{ML} \mathrm{Tr}({\bs \Upsilon} {\bf a}_2 {\bf a}_1^{*})$. Using Eq. (\ref{eq:expre-gamma}) in the case ${\bf A}= {\bf a}_2 {\bf a}_1^{*}$, we obtain that 
$$
{\bf a}_1^{*} {\bs \Upsilon} {\bf a}_2 = -\sigma^{2} \mathbb{E} \left[ \sum_{u=-(L-1)}^{L-1}  \tau^{(M)}({\bf Q}^{\circ})(u) \, {\bf a}_1^{*} {\bf Q} {\bf W} {\bf J}_N^{u} {\bf H}^{T} {\bf W}^{*} {\bf a}_2 \right]
$$
(\ref{eq:varentreesQW}, \ref{eq:vartraceQ}) and the Schwartz inequality lead immediately to 
$$
\left| {\bf a}_1^{*} {\bs \Upsilon} {\bf a}_2 \right| \leq \kappa^{2} \, C(z)  L \, \frac{1}{\sqrt{MN}} \sqrt{\frac{L}{N}}  = \kappa^{2}\,  C(z)  \sqrt{\frac{L}{M}} \frac{L}{N}.
$$
We now establish (\ref{eq:entrees-trace-Delta}). We remark that
$$
{\bf b}_1^* \left( \frac{1}{M} \sum_{m=1}^{M} {\bs \Upsilon}^{m,m} \right) {\bf b}_2 = L \, 
\frac{1}{ML} \mathrm{Tr} \left( {\bs \Upsilon} ({\bf I}_M \otimes {\bf b}_2 {\bf b}_1^*) \right) 
$$
Using Eq. (\ref{eq:expre-gamma}) in the case ${\bf A}= {\bf I}_M \otimes {\bf b}_2 {\bf b}_1^{*}$, we obtain immediately that 
\begin{equation}
\label{eq:expre-quadratique-Gammahat}
\begin{array}{c} 
{\bf b}_1^* \left( \frac{1}{M} \sum_{m=1}^{M} {\bs \Upsilon}^{m,m} \right) {\bf b}_2  =  \\
\sum_{u=-(L-1)}^{L-1} \mathbb{E} \left[ \tau^{(M)}({\bf Q}^{\circ})(u) \, {\bf b_1}^{*} \left( \frac{1}{M} \sum_{m=1}^{M} ({\bf Q} {\bf W} {\bf J}_N^{u} {\bf H}^{T} {\bf W}^{*})^{m,m} \right) {\bf b}_2 \right] 
\end{array}
\end{equation}
(\ref{eq:entrees-trace-Delta}) thus appears as a direct consequence of (\ref{eq:vartraceQ}), \ref{eq:vartraceentreesQW}) and of the Schwartz inequality.

We finally mention a useful corollary of (\ref{eq:entrees-trace-Delta}). 
\begin{corollary}
\label{co:converge-toeplitzifie-Q-R}
 It holds that 
\begin{equation}
\label{eq:converge-toeplitzifie-Q-R}
\| \mathcal{T}^{(M)}_{N,L}\left( \mathbb{E}({\bf Q}) - ({\bf I}_M \otimes {\bf R}) \right) \| \leq C(z) \frac{L^{3/2}}{MN}
\end{equation}
for each $z \in \mathbb{C}^{+}$ where $C(z)$ can be written as $C(z) = P_1(|z|) \, P_2\left( (\mathrm{Im}z)^{-1} \right)$ for some 
nice polynomials $P_1$ and $P_2$. 
\end{corollary}
Taking into account Proposition \ref{prop:contractant}, (\ref{eq:converge-toeplitzifie-Q-R}) follows immediately 
from (\ref{eq:entrees-trace-Delta}) by considering the unit norm vector ${\bf b} = {\bf a}_L(\nu)$.

\section{Convergence towards the Marcenko-Pastur distribution}
\label{sec:marcenko-pastur}

In the following, we establish that 
\begin{equation}
\label{eq:convergence-MP}
\frac{1}{ML} \mathrm{Tr}\left( \mathbb{E}({\bf Q}(z)) - t(z) {\bf I}_{ML} \right) \rightarrow 0
\end{equation}
for each $z \in \mathbb{C}^{+}$. (\ref{eq:vartraceQ}) does not imply in general that 
$\frac{1}{ML} \mathrm{Tr}\left( {\bf Q}(z) - \mathbb{E}({\bf Q}(z)) \right)$ converges towards $0$ almost surely 
(this would be the case if $M$ was of the same order of magnitude than $N^{\kappa}$ for some $\kappa > 0$). 
However, the reader may check using the Poincaré-Nash inequality that the variance of $\left[\frac{1}{ML} \mathrm{Tr}( {\bf Q}^{\circ}(z))\right]^{2}$ 
is a $\mathcal{O}(\frac{1}{(MN)^{2}})$ term. As 
$$
\mathbb{E} \left| \frac{1}{ML} \mathrm{Tr} ({\bf Q}^{\circ}(z)) \right|^{4} =  \left| \mathbb{E} \left[ \frac{1}{ML} \mathrm{Tr} ({\bf Q}^{\circ}(z)) \right] ^{2}  \right|^{2} + 
\mathrm{Var}\left[\frac{1}{ML} \mathrm{Tr}( {\bf Q}^{\circ}(z))\right]^{2}
$$ 
(\ref{eq:vartraceQ}) implies that the fourth-order moment of $\frac{1}{ML} \mathrm{Tr}\left( {\bf Q}^{\circ}(z)\right)$ 
is also a $\mathcal{O}(\frac{1}{(MN)^{2}})$ term, and that $\frac{1}{ML} \mathrm{Tr}\left( {\bf Q}(z) - \mathbb{E}({\bf Q}(z)) \right)$ converges towards $0$ almost surely. 
Consequently, (\ref{eq:convergence-MP}) allows to prove that the eigenvalue value distribution 
of ${\bf W} {\bf W}^{*}$ has almost surely the same behaviour than the Marcenko-Pastur distribution $\mu_{\sigma^{2}, c_N}$. As $c_N \rightarrow c_*$, this of course establishes the almost sure convergence of the eigenvalue distribution 
of ${\bf W}_N {\bf W}_N^{*}$ towards the Marcenko-Pastur $\mu_{\sigma^{2},c_*}$.\\

In the following, we thus prove (\ref{eq:convergence-MP}). 
(\ref{eq:expre-E(Q)}) and Proposition \ref{prop:controle-Delta} imply that 
for each uniformly bounded $L \times L$ matrix ${\bf A}$, then, it holds that 
\begin{equation}
\label{eq:E(Q)-R}
\frac{1}{ML} \mathrm{Tr} \left[ \left(\mathbb{E}({\bf Q}(z)) - {\bf I}_M \otimes {\bf R}(z) \right) \left({\bf I}_M \otimes {\bf A} \right)\right] = \mathcal{O}(\frac{L}{MN})
\end{equation}
for each $z \in \mathbb{C}^{+}$. We now establish that
$$
\frac{1}{ML} \mathrm{Tr} \left[ \left({\bf I}_M \otimes {\bf R}(z) - t(z) {\bf I}_{ML} \right) \left( {\bf I}_M \otimes {\bf A} \right)\right] \rightarrow 0 
$$
or equivalently that 
\begin{equation}
\label{eq:R-t}
\frac{1}{L} \mathrm{Tr} \left[ \left({\bf R}(z) - t(z) {\bf I}_{L} \right)  {\bf A} \right] \rightarrow 0
\end{equation}
for each $z \in \mathbb{C}^{+}$.
For this, we first mention that straighforward computations lead to 
\begin{equation}
\label{eq:expre-R-t}
{\bf R} - t {\bf I} = -\sigma^{4} c_N \, z t(z) \tilde{t}(z)  \; 
{\bf R} \; \mathcal{T}_{L,L}\left( {\bf H}  \mathcal{T}^{(M)}_{N,L}\left[ \mathbb{E}({\bf Q}) - t {\bf I}_{ML} \right] 
\right) 
\end{equation}
Therefore, 
$$
\frac{1}{L} \mathrm{Tr} \left[ ({\bf R} - t {\bf I}_L ) {\bf A} \right] = -\sigma^{4} c_N \, z t(z) \tilde{t}(z)
  \; 
\frac{1}{L} \mathrm{Tr} {\bf A} {\bf R}  \mathcal{T}_{L,L}\left( {\bf H}  \mathcal{T}^{(M)}_{N,L}\left[\mathbb{E}({\bf Q}) - t {\bf I}_{ML} \right] 
 \right)
$$
Direct application of (\ref{eq:utile}) to the case $P=M, K=L, R=L, {\bf C} = \mathbb{E}({\bf Q}) - t {\bf I}_{ML}, {\bf B} = {\bf A} {\bf R}$
and ${\bf D} = {\bf H}$ implies that 
$$
\frac{1}{L} \mathrm{Tr} \left( ({\bf R} - t {\bf I}_L)  {\bf A} \right) = -\sigma^{4} c_N \, z t(z) \tilde{t}(z) 
\; \frac{1}{ML} \mathrm{Tr} \left[ \left( \mathbb{E}({\bf Q}) - t {\bf I}_{ML} \right) \left( {\bf I}_M \otimes \mathcal{T}_{L,L}( \mathcal{T}_{N,L}({\bf A}{\bf R}) {\bf H} \right) \right]
$$
In the following, we denote by ${\bf G}({\bf A})$ the $L \times L$ matrix defined by
\begin{equation}
\label{eq:def-G(A)}
{\bf G}({\bf A}) = \mathcal{T}_{L,L}\left(\mathcal{T}_{N,L}({\bf A}{\bf R}) {\bf H} \right)
\end{equation}
Writing that $ {\bf E}({\bf Q}) - t {\bf I}_{ML} =  {\bf E}({\bf Q}) - {\bf I}_M \otimes {\bf R} +  {\bf I}_M \otimes {\bf R} - t {\bf I}_{ML}$, we obtain that
\begin{eqnarray}
\label{eq:decomposition}
\frac{1}{L} \mathrm{Tr} \left[ ({\bf R} - t {\bf I}_L)  {\bf A} \right] = & -\sigma^{4} c_N \, z t(z) \tilde{t}(z) \; \frac{1}{ML} \mathrm{Tr} \left[ \left( {\bf E}({\bf Q}) - {\bf I}_M \otimes {\bf R} \right) \left( {\bf I}_M \otimes {\bf G}({\bf A}) \right) \right] - \\
                                                                          &  \sigma^{4} c_N \, z t(z) \tilde{t}(z)\; \frac{1}{L} 
\mathrm{Tr} \left[  ({\bf R} - t {\bf I}_L)  {\bf G}({\bf A}) \right] \nonumber
\end{eqnarray}
We now prove that 
\begin{equation}
\label{eq:supremum}
\sup_{\| {\bf B} \| \leq 1} \left| \frac{1}{L} \mathrm{Tr} \left( ({\bf R} - t {\bf I}_L)  {\bf B} \right) \right| = \mathcal{O}(\frac{L}{MN}) 
\end{equation}
when $z$ belongs to a certain domain. For this, we first remark that (\ref{eq:Tcontractant}) implies that $\|  {\bf G}({\bf A}) \| \leq \| {\bf H} \| \| {\bf R} \| \| {\bf A} \|$. 
By Lemma \ref{le:proprietes-H-R}, it holds that $\| {\bf H} \| \| {\bf R} \| \leq \frac{|z|}{(\mathrm{Im}(z))^{2}}$. Consequently, we obtain that
\begin{equation}
\label{eq:inegalite-norme-G(A)}
\|  {\bf G}({\bf A}) \| <  \frac{|z|}{(\mathrm{Im}(z))^{2}} \, \| {\bf A} \|
\end{equation}
This implies that for each $L \times L$ matrix ${\bf A}$ such that $\| {\bf A} \| \leq 1$, 
then, it holds that 
\begin{align*}
\left|  \frac{1}{ML} \mathrm{Tr} \left[ \left( {\bf E}({\bf Q}) - {\bf I}_M \otimes {\bf R} \right) \left( {\bf I}_M \otimes {\bf G}({\bf A}) \right) \right] \right|  & \leq   
\frac{|z|}{(\mathrm{Im}(z))^{2}} \; \sup_{\| {\bf B} \| \leq 1} \left|  \frac{1}{ML} \mathrm{Tr} \left[ \left( {\bf E}({\bf Q}) - {\bf I}_M \otimes {\bf R}) \right) {\bf B} \right] \right| \; , \\
 \left| \frac{1}{L} \mathrm{Tr} \left[ ({\bf R} - t {\bf I}_L)  {\bf G}({\bf A}) \right] \right| & \leq  
\frac{|z|}{(\mathrm{Im}(z))^{2}} \;  \sup_{\| {\bf B} \| \leq 1}  \left| \frac{1}{L} \mathrm{Tr} \left[ ({\bf R} - t {\bf I}_L)  {\bf B} \right] \right|
\end{align*}
Proposition \ref{prop:controle-Delta} implies that 
$$
\sup_{\| {\bf B} \| \leq 1} \left|  \frac{1}{ML} \mathrm{Tr} \left[ \left( {\bf E}({\bf Q}) - {\bf I}_M \otimes {\bf R}) \right) {\bf B} \right] \right| =  \mathcal{O}(\frac{L}{MN})
$$
This and Eq. (\ref{eq:decomposition}) eventually imply that 
$$
\sup_{\| {\bf B} \| \leq 1}  \left| \frac{1}{L} \mathrm{Tr} \left( ({\bf R} - t {\bf I}_L)  {\bf B} \right) \right| \leq \mathcal{O}(\frac{L}{MN}) + 
\sigma^{4} c_N \, |z t(z) \tilde{t}(z)| \;  \frac{|z|}{(\mathrm{Im}(z))^{2}} \; \sup_{\| {\bf B} \| \leq 1}  \left| \frac{1}{L} \mathrm{Tr} \left( ({\bf R} - t {\bf I}_L)  {\bf B} \right) \right|
$$ 
It also holds that $|z t(z) \tilde{t}(z)| \leq  \frac{|z|}{(\mathrm{Im}(z))^{2}}$. Therefore, if $z$ belongs to the domain $\sigma^{4} c_N \frac{|z|^{2}}{(\mathrm{Im}(z))^{4}} < \frac{1}{2}$,  we obtain that 
\begin{equation}
\label{eq:comportement-sup}
\sup_{\| {\bf B} \| \leq 1}  \left| \frac{1}{L} \mathrm{Tr} \left[ ({\bf R} - t {\bf I}_L)  {\bf B} \right] \right| =  \mathcal{O}(\frac{L}{MN})
\end{equation}
This establishes (\ref{eq:R-t}) for each uniformly bounded $L \times L$ matrix ${\bf A}$ whenever $z$ is well chosen. 
Moreover, for these values of $z$,  
$\frac{1}{L} \mathrm{Tr}\left(({\bf R} - t \, {\bf I}) {\bf A} \right)$, and thus $\frac{1}{ML} \mathrm{Tr} \left( \mathbb{E}({\bf Q}(z) - t(z) \, {\bf I}_{ML}) {\bf A} \right)$,  are $\mathcal{O}(\frac{L}{MN})$ terms. 
A standard application 
of Montel's theorem implies that (\ref{eq:R-t})  holds on $\mathbb{C}^{+}$.  This, in turn, 
establishes (\ref{eq:convergence-MP}). 

\begin{remark}
\label{re:matrices-test}
We have proved that for each uniformely bounded $L \times L$ matrix ${\bf A}$, 
then it holds that 
$$
\frac{1}{ML} \mathrm{Tr}\left[ \left( \mathbb{E}({\bf Q}(z) - t(z) {\bf I}_{ML}) \right) \left({\bf I}_M \otimes {\bf A} 
\right) \right] \rightarrow 0
$$
for each $z \in \mathbb{C}^{+}$. It is easy to verify that matrix $ {\bf I}_M \otimes {\bf A}$ 
can be replaced by any uniformly bounded $ML \times ML$ matrix ${\bf B}$. In effect, Proposition  
\ref{prop:controle-Delta} implies that it is sufficient to establish that 
$$
\frac{1}{ML} \mathrm{Tr}\left[ \left( {\bf I}_M \otimes {\bf R}(z) - t(z) {\bf I}_{ML} \right) {\bf B} \right] \rightarrow 0
$$
The above term can also be written as
$$
\frac{1}{L} \mathrm{Tr} \left[ \left( {\bf R}(z) - t(z) \, {\bf I}_L \right) \left( \frac{1}{M} \sum_{m=1}^{M} 
{\bf B}^{m,m} \right) \right] 
$$
and converges towards 0 because matrix $\frac{1}{M} \sum_{m=1}^{M} {\bf B}^{m,m}$ is uniformly bounded.
\end{remark}

\section{Convergence of the spectral norm of $\mathcal{T}_{N,L}({\bf R}(z) - t(z) {\bf I}_N)$}
\label{sec:convergence-spectral-norm}
From now on, we assume that $L,M,N$ satisfy the following extra-assumption: 
\begin{assum}
\label{as:L32/MN}
$\frac{L^{3/2}}{MN} \rightarrow 0$ or equivalently, $\frac{L}{M^{4}} \rightarrow 0$. 
\end{assum}
The goal of this section is prove Theorem \ref{theo:convergence-norme-tau(R-t)} which will be used extensively
in the following.  
\begin{theorem}
\label{theo:convergence-norme-tau(R-t)}
Under assumption \ref{as:L32/MN}, it exists 2 nice polynomials $P_1$ and 
$P_2$ for which
\begin{equation}
\label{eq:haagerup-norme-1}
\| \mathcal{T}_{N,L}({\bf R}(z) - t(z) {\bf I}_N) \| \leq \sup_{\nu \in [0,1]} \left|  {\bf a}_L(\nu)^{*} \left( {\bf R}(z) - t(z) {\bf I}_L \right) {\bf a}_L(\nu) \right|  \leq \frac{L^{3/2}}{MN} \, P_1(|z|) P_2(\frac{1}{\mathrm{Im}(z)})
\end{equation}
for each $z \in \mathbb{C}^{+}$. 
\end{theorem}
{\bf Proof.} \\

{\bf First step. } The first step consists in showing that 
\begin{equation}
\label{eq:norme-T(R-t)}
\sup_{\nu \in [0,1]} \left|  {\bf a}_L(\nu)^{*} \left( {\bf R}(z) - t(z) {\bf I}_L \right) {\bf a}_L(\nu) \right| 
\rightarrow 0
\end{equation}
for each $z \in \mathbb{C}^{+}$, which implies that $\| \mathcal{T}_{N,L}\left( {\bf R} - t {\bf I}_L \right) \| \rightarrow 0$ for each $z \in \mathbb{C}^{+}$ (see (\ref{eq:Tcontractant})). We first establish that (\ref{eq:norme-T(R-t)}) holds for certain values of $z$, and extend the property
to $\mathbb{C}^{+}$ using Montel's theorem. We take (\ref{eq:expre-R-t}) as a starting point, and write $\mathbb{E}({\bf Q} - t \, {\bf I}_{ML})$ as
$$
\mathbb{E}({\bf Q} - t \, {\bf I}_{ML}) = \mathbb{E}({\bf Q}) - {\bf I}_M \otimes {\bf R} + ({\bf I}_M \otimes {\bf R} - t \, {\bf I}_{ML})
$$
(\ref{eq:expre-R-t}) can thus be rewritten as
\begin{eqnarray}
\label{eq:equation-R-t}
{\bf R} - t \, {\bf I}_L & = & - \sigma^{4} c_N z \, t(z) \, \tilde{t}(z) \, {\bf R} \mathcal{T}_{L,L} \left( {\bf H} \; \mathcal{T}^{(M)}_{N,L} \left[\mathbb{E}({\bf Q}) - {\bf R}_M \right]  \right) -  \\ \nonumber &   & \sigma^{4} c_N z \, t(z) \, \tilde{t}(z) \, {\bf R} \mathcal{T}_{L,L} \left( {\bf H} \; \mathcal{T}_{N,L} \left[{\bf R} - t \, {\bf I}_L \right] \right)
\end{eqnarray} 
Therefore, for each deterministic uniformly bounded $L$--dimensional vector ${\bf b}$, then, it holds that 
\begin{eqnarray}
\label{eq:expre-entrees-R-t-1}
{\bf b}^{*} \left( {\bf R} - t \, {\bf I} \right) {\bf b} & = & 
- z t(z) \tilde{t}(z) \sigma^{4} c_N {\bf b}^{*} {\bf R} \mathcal{T}_{L,L} \left( {\bf H} \, \mathcal{T}^{(M)}_{N,L} \left[\mathbb{E}({\bf Q}) - {\bf I}_M \otimes {\bf R} \right] \, \right) {\bf b}  - \\
\label{eq:expre-entrees-R-t-2}
  &   &   z t(z) \tilde{t}(z) \sigma^{4} c_N {\bf b}^{*} {\bf R} \mathcal{T}_{L,L}\left( {\bf H}  \, \mathcal{T}_{N,L} \left[ {\bf R} - t \, {\bf I} \right] \right) {\bf b}  
\end{eqnarray}
Proposition \ref{prop:contractant} implies that 
$$
\|  \mathcal{T}_{L,L}\left( \mathcal{T}_{N,L} \left[ {\bf R} - t \, {\bf I} \right] \, {\bf H} \right) \| 
\leq \| {\bf H} \| \, \|  \mathcal{T}_{N,L} \left[ {\bf R} - t \, {\bf I} \right] \| \leq 
\| {\bf H} \| \sup_{\nu} \left| {\bf a}_L(\nu)^{*} \left( {\bf R} - t \, {\bf I} \right) {\bf a}_L(\nu) \right|
$$
and that 
$$
\|  \mathcal{T}_{L,L}\left(  {\bf H} \; \mathcal{T}^{(M)}_{N,L} \left[ \mathbb{E}({\bf Q}) - {\bf I}_M \otimes {\bf R}  \right] \right) \| 
\leq \| {\bf H} \| \, \|  \mathcal{T}^{(M)}_{N,L} \left[ \mathbb{E}({\bf Q}) - {\bf I}_M \otimes {\bf R}   \right] \| \leq 
\| {\bf H} \| \sup_{\nu} \left| {\bf a}_L(\nu)^{*} \, \hat{{\bs \Delta}}  {\bf a}_L(\nu) \right|
$$
where we recall that ${\bs \Delta} = \mathbb{E}({\bf Q}) - {\bf I}_M \otimes {\bf R}$ and that $\hat{{\bs \Delta}} = \frac{1}{M} \sum_{m=1}^{M} \bs{\Delta}^{(m,m)}$. We denote by $\beta$ and $\delta$ the terms \- 
$\beta = \sup_{\nu} \left| {\bf a}_L(\nu)^{*} \left( {\bf R} - t \, {\bf I} \right) {\bf a}_L(\nu) \right|$ and 
 $\delta = \sup_{\nu} \left| {\bf a}_L(\nu)^{*} \, \hat{{\bs \Delta}} {\bf a}_L(\nu) \right|$. We remark that $\delta = \mathcal{O}\left( \frac{L^{3/2}}{MN} \right)$ (see (\ref{eq:entrees-trace-Delta})).
We choose ${\bf b} = {\bf a}_L(\mu)$ in (\ref{eq:expre-entrees-R-t-1}), evaluate the modulus of the left handside of
(\ref{eq:expre-entrees-R-t-1}),  and take the supremum over $\mu$. This immediately leads to 
\begin{equation}
\label{eq:inegalite-entrees-R-t}
\beta \leq |z t(z) \tilde{t}(z)| \sigma^{4} c_N \| {\bf R} \| \| {\bf H} \| \delta + 
 |z t(z) \tilde{t}(z)| \sigma^{4} c_N \| {\bf R} \|  \| {\bf H} \| \beta 
\end{equation}
Moreover, (see Lemma (\ref{le:proprietes-H-R})), it holds that 
$$
|z t(z) \tilde{t}(z)| \sigma^{4} c_N \| {\bf R} \|  \| {\bf H} \| \leq \sigma^{4} c_N  
\frac{|z|^{2}}{(\mathrm{Im}(z))^{4}}
$$
(\ref{eq:inegalite-entrees-R-t}) implies that if $z$ satisfies 
\begin{equation}
\label{eq:condition-z-gn}
\sigma^{4} c_N \frac{|z|^{2}}{(\mathrm{Im}(z))^{4}} \leq \frac{1}{2}, 
\end{equation}
then $\beta = \mathcal{O}\left( \frac{L^{3/2}}{MN} \right)$ and therefore, converges towards $0$. We now extend this 
property on $\mathbb{C}^{+}$ using Montel's theorem. For this, we consider an integer sequence 
$K(N)$ for which $\frac{L(N)}{K(N)} \rightarrow 0$, and denote for each $N$ and $0 \leq k \leq K(N)$ by 
$\nu_k^{(N)}$ the element of $[0,1]$ defined by $\nu_k^{(N)} = \frac{k}{K(N)}$. We denote by $\phi(k,N)$ the 
one-to-one correspondance between the set of integer couples $(k,N)$, $k \leq K(N)$ and the set of integers $\mathbb{N}$ 
defined by $\phi(0,0) = 0$, $\phi(k+1,N) = \phi(k,N) + 1$ for $k < K(N)$ and $\phi(0,N+1) = \phi(K(N),N) + 1$. Each integer $n$ can therefore be written in a unique way as $n = \phi(k,N)$ for a certain couple $(k,N)$, $0 \leq k \leq K(N)$. We define 
a sequence of analytic functions $(g_n(z))_{n \in \mathbb{N}}$ defined on $\mathbb{C}^{+}$ by 
\begin{equation}
\label{eq:def-g_n}
g_{\phi(k,N)}(z) = {\bf a}_L(\nu_k^{(N)})^{*} \left( {\bf R}(z) - t(z) \, {\bf I}_L \right) {\bf a}_L(\nu_k^{(N)})
\end{equation}
If $z$ satisfies (\ref{eq:condition-z-gn}), the sequence $g_n(z)$ converges towards 0. Moreover, 
$(g_n(z))_{n \in \mathbb{N}}$ is a normal family of $\mathbb{C}^{+}$. Consider a subsequence extracted 
from $(g_n)_{n \in \mathbb{Z}}$ converging uniformly on compact subsets of  $\mathbb{C}^{+}$
towards an analytic function $g_{*}$. As $g_*(z) = 0$ if $z$ satifies (\ref{eq:condition-z-gn}), function 
$g_*$ is zero. This shows that all convergent subsequences extracted from $(g_n)_{n \in \mathbb{N}}$ 
converges towards 0, so that the whole sequence $(g_n)_{n \in \mathbb{N}}$ converges towards 0. This immediately
implies that 
\begin{equation}
\label{eq:sup-discret}
\lim_{N \rightarrow +\infty} \sup_{0 \leq k \leq K(N)} |g_{\phi(k,N)}(z)| = 0
\end{equation}
for each $z \in \mathbb{C}^{+}$. For each $\nu \in [0,1]$, it exists an index $k$, $0 \leq k \leq K(N)$ such that 
$|\nu - \nu_k^{(N)}| \leq \frac{1}{2 K(N)}$. It is easily checked that 
$$
\| {\bf a}_L(\nu) - {\bf a}_L(\nu_k^{(N)}) \| = \mathcal{O}\left( L(N) |\nu - \nu_k^{(N)}|  \right)= \mathcal{O}\left(\frac{L(N)}{K(N)}\right)) = o(1)
$$
and that 
$$
\left| {\bf a}_L(\nu)^{*} \left( {\bf R}(z) - t(z) \, {\bf I}_L \right) {\bf a}_L(\nu) - {\bf a}_L(\nu_k^{(N)})^{*} \left( {\bf R}(z) - t(z) \, {\bf I}_L \right) {\bf a}_L(\nu_k^{(N)}) \right| \rightarrow 0
$$
for each $z \in \mathbb{C}^{+}$. We deduce from (\ref{eq:sup-discret}) that (\ref{eq:norme-T(R-t)}) 
holds for each $z \in \mathbb{C}^{+}$ as expected. \\

{\bf Second step.} The most difficult part of the proof consists in evaluating the 
rate of convergence of \\ $\sup_{\nu} \left| {\bf a}_L(\nu)^{*} ({\bf R}(z) - t(z) {\bf I}_N) {\bf a}_L(\nu) \right|$.

By (\ref{eq:expre-symbole}), the quadratic form $ {\bf a}_L(\nu)^{*} ({\bf R}(z) - t(z) {\bf I}_N) {\bf a}_L(\nu)$ can also be written as 
$$
 {\bf a}_L(\nu)^{*} ({\bf R}(z) - t(z) {\bf I}_N) {\bf a}_L(\nu) =  \sum_{l=-(L-1)}^{L-1} \tau({\bf R} - t \, {\bf I})(l) e^{-2 i \pi l \nu} 
$$
where we recall that $\tau({\bf R} - t \, {\bf I})(l)  = \frac{1}{L} \mathrm{Tr} \left( ({\bf R} - t \, {\bf I}) {\bf J}_L^{l} \right)$. 
In order to study more thoroughly \\ $\sup_{\nu} \left| {\bf a}_L(\nu)^{*} ({\bf R}(z) - t(z) {\bf I}_N) {\bf a}_L(\nu) \right|$, it is 
thus possible to evaluate the coefficients $(\tau({\bf R} - t \, {\bf I})(l))_{l=-(L-1), \ldots, L-1}$. In the following, for a $L \times L$ matrix ${\bf X}$, we denote by 
${\bs \tau}({\bf X})$ the  $2L-1$--dimensional vector defined by
$$
{\bs \tau}({\bf X}) = \left(\tau({\bf X})(-(L-1)), \ldots, \tau({\bf X})(L-1) \right)^{T}
$$
(\ref{eq:equation-R-t}) can be associated to a linear equation whose unknown 
is vector $\tau({\bf R} - t \, {\bf I})$. Writing $\mathcal{T}_{N,L} \left[ {\bf R} - t \, {\bf I} \right]$ as $\sum_{l=-(L-1)}^{L-1} \tau({\bf R} - t \, {\bf I})(l) {\bf J}_N^{*l}$, multiplying (\ref{eq:equation-R-t}) from both sides by ${\bf J}_L^{k}$, and taking the normalized trace, we obtain that 
\begin{equation}
\label{eq:matrice-equation-R-t}
\tau({\bf R} - t \, {\bf I}) = {\bs \tau}({\bs \Gamma}) + {\bf D}^{(0)} \, \tau({\bf R} - t \, {\bf I})
\end{equation}
where ${\bf D}^{(0)}$ is the $(2L - 1) \times (2L - 1)$ matrix whose entries ${\bf D}^{(0)}_{k,l}$, $(k,l) \in \{ -(L-1), \ldots, L-1) \}$ are defined by 
$$
{\bf D}^{(0)}_{k,l} = - \sigma^{4} c_N z \, t(z) \, \tilde{t}(z) \frac{1}{L} \mathrm{Tr}  \left[ {\bf R} \mathcal{T}_{L,L} \left( {\bf H} {\bf J}_N^{*l} \right) {\bf J}_L^{k} \right] 
$$
and where matrix ${\bs \Gamma}$ represents the first term of the righthanside of (\ref{eq:equation-R-t}), i.e. 
\begin{equation}
\label{eq:def-Upsilon-D0}
{\bs \Gamma} = - \sigma^{4} c_N z \, t(z) \, \tilde{t}(z) \, {\bf R} \mathcal{T}_{L,L} \left( {\bf H} \; \mathcal{T}^{(M)}_{N,L} \left[\mathbb{E}({\bf Q}) - {\bf I}_M \otimes {\bf R} \right]  \right)
\end{equation}
Equation (\ref{eq:matrice-equation-R-t}) should be inverted, and the effect of the inversion on vector ${\bs \tau}({\bs \Gamma})$ should be analysed in 
order to evaluate the behaviour of $\| \mathcal{T}_{N,L}({\bf R}(z) - t(z) {\bf I}_N) \|$. The invertibility of matrix ${\bf I} - {\bf D}^{(0)}$ and
the control of its inverse are however non trivial, and need some efforts. \\

In the following, we denote by ${\bs \Phi}^{(0)}$ the operator defined on $\mathbb{C}^{L \times L}$ 
by 
\begin{equation}
\label{eq:def-Phi0}
{\bs \Phi}^{(0)}({\bf X}) = - \sigma^{4} c_N z \, t(z) \, \tilde{t}(z) \, {\bf R} \mathcal{T}_{L,L} \left( {\bf H} \; \mathcal{T}_{N,L} \left[ {\bf X} \right] \right)
\end{equation}
for each $L \times L$ matrix ${\bf X}$. Eq. (\ref{eq:equation-R-t}) can thus be written as
$$
{\bf R} - t \, {\bf I}_L = {\bs \Gamma} +  {\bs \Phi}^{(0)}({\bf R} - t \, {\bf I}_L)
$$ 
We also remark that matrix ${\bs \Gamma}$ is given by
\begin{equation}
\label{eq:expre-Upsilon-D0}
{\bs \Gamma} = {\bs \Phi}^{(0)} \left( \mathbb{E}(\hat{{\bf Q}}) - {\bf R} \right)
\end{equation}

Moreover, it is clear that vector ${\bs \tau}\left( {\bs \Phi}^{(0)}({\bf X}) \right)$ can be written as
\begin{equation}
\label{eq:representation-D0}
{\bs \tau}\left( {\bs \Phi}^{(0)}({\bf X}) \right) = {\bf D}^{(0)} \,  {\bs \tau}({\bf X})
\end{equation}
In order to study the properties of operator ${\bs \Phi}^{(0)}$ and of matrix ${\bf D}^{(0)}$, we introduce the operator ${\bs \Phi}$ and the corresponding 
$(2L-1) \times (2L-1)$ matrix ${\bf D}$ defined respectively by 
\begin{equation}
\label{eq:def-Phi}
{\bs \Phi}({\bf X}) = \sigma^{4} c_N  {\bf R} \mathcal{T}_{L,L} \left( {\bf H} \; \mathcal{T}_{N,L} \left[ {\bf X} \right] {\bf H}^{*} \right)
{\bf R}^{*}
\end{equation}
and
\begin{equation}
\label{eq:def-D}
{\bf D}_{k,l} = \sigma^{4} c_N  \frac{1}{L} \mathrm{Tr}  \left[ {\bf R} \mathcal{T}_{L,L} \left( {\bf H} {\bf J}_N^{*l} {\bf H}^{*} \right) 
{\bf R}^{*} {\bf J}_L^{k} \right] 
\end{equation}
for $(k,l) \in  \{ -(L-1), \ldots, L-1) \}$. Matrix ${\bf D}$ of course satisfies 
\begin{equation}
\label{eq:relation-Phi-D}
{\bs \tau} \left({\bs \Phi}({\bf X})\right) = {\bf D} {\bs \tau}({\bf X})
\end{equation}
Before establishing the relationships between $({\bs \Phi}_0, {\bf D}^{(0)})$ and $({\bs \Phi}, {\bf D})$, we prove the following proposition. 
\begin{proposition}
\label{prop:proprietes-Phi}
\begin{itemize}
\item 
If ${\bf X}$ is positive definite, then matrix ${\bs \Phi}({\bf X})$ is also positive definite. Moreover, if ${\bf X}_1 \geq {\bf X}_2$, then 
${\bs \Phi}({\bf X}_1) \geq {\bs \Phi}({\bf X}_2)$.
\item 
It exists 2 nice polynomials $P_1$ and $P_2$ and an integer $N_1$ such that 
the spectral radius $\rho({\bf D})$ of matrix ${\bf D}$ verifies $\rho({\bf D}) < 1$ for $N \geq N_1$ and for each $z \in E_N$ where $E_N$ is the subset of $\mathbb{C}^{+}$ defined by
\begin{equation}
\label{eq:def-EN}
E_N = \{ z \in \mathbb{C}^{+}, \frac{L^{3/2}}{MN} P_1(|z|) P_2(1/\mathrm{Im}z) \leq 1 \}. 
\end{equation}
\item for $N \geq N_1$,  matrix ${\bf I} - {\bf D}$ is invertible for $z \in E_N$. If
we denote by ${\bf f} = (f_{-(L-1)}, \ldots, f_0, \ldots, f_{L-1})^{T}$ the 
$(2L-1)$--dimensional vector defined by 
\begin{equation}
\label{eq:def-f}
{\bf f} = ({\bf I} - {\bf D})^{-1} {\bs \tau}({\bf I}) = ({\bf I} - {\bf D})^{-1} {\bf e}_0
\end{equation} 
where ${\bf e}_0 = (0, \ldots, 0, 1, 0, \ldots, 0)^{T}$,  
then, for each $\nu \in [0,1]$, the term $\sum_{l=-(L-1)}^{L-1}  {\bf f}_l \, e^{-2 i \pi l \nu}$ is real and positive, and 
\begin{equation}
\label{eq:controle-TF-(I-D)-1}
\sup_{\nu \in [0,1]} \sum_{l=-(L-1)}^{L-1}  {\bf f}_l  \, e^{-2 i \pi l \nu} \leq C \, \frac{(|\eta_1|^{2} + |z|^{2})^{2}}{(\mathrm{Im} z)^{4}}
\end{equation}
for some nice constants $C$ and $\eta_1$. 
\end{itemize}
\end{proposition}
{\bf Proof.} The first item follows immediately from the basic properties of operators $\mathcal{T}$. 
The starting point of the proof of item 2 consists in writing matrix $\mathbb{E} ( \hat{{\bf Q}} )  = \frac{1}{M} \sum_{m=1}^{M} \mathbb{E}({\bf Q}^{m,m})$
as  $\mathbb{E} ( \hat{{\bf Q}} ) = {\bf R} + \hat{{\bs \Delta}}$, and in expressing the imaginary part of $\mathbb{E}(\hat{{\bf Q}})$ as 
$\mathrm{Im}\left(\mathbb{E}(\hat{{\bf Q}})\right) =  \mathrm{Im}\left( \mathbb{E}(\hat{{\bs \Delta}})\right) +  \mathrm{Im}({\bf R})$. Writing 
$ \mathrm{Im}({\bf R})$ as
$$
 \mathrm{Im}({\bf R}) = \frac{{\bf R} - {\bf R}^{*}}{2i} = \frac{1}{2i} \; {\bf R} \left( {\bf R}^{-*} - {\bf R}^{-1} \right) {\bf R}^{*}
$$
and expressing ${\bf R}^{-1}$ in terms of ${\bf H}$, and using the same tricks for ${\bf H}$, we eventually obtain that
\begin{equation}
\label{eq:expre-ImhatQ}
\mathrm{Im}\left(\mathbb{E}(\hat{{\bf Q}})\right) = \mathrm{Im}\left( \mathbb{E}(\hat{{\bs \Delta}})\right) + \mathrm{Im} z \, {\bf R} {\bf R}^{*} 
+ \sigma^{4} c_N \; {\bf R} \mathcal{T}_{L,L} \left[ {\bf H} \, \mathcal{T}_{N,L}\left(\mathrm{Im}\left(\mathbb{E}(\hat{{\bf Q}})\right)\right) \, {\bf H}^{*} \right] {\bf R}^{*}
\end{equation}
In order to simplify the notations, we denote by ${\bf X}$ and ${\bf Y}$ the matrices $\mathrm{Im}\left(\mathbb{E}(\hat{{\bf Q}})\right)$ 
and $\mathrm{Im}\left( \mathbb{E}(\hat{{\bs \Delta}})\right) + \mathrm{Im} z \, {\bf R} {\bf R}^{*}$ respectively. (\ref{eq:expre-ImhatQ})
implies that for each $z \in \mathbb{C}^{+}$, then the positive definite matrix ${\bf X}$ satisfies 
\begin{equation}
\label{eq:expre-ImhatQ-Phi}
{\bf X} = {\bf Y} + {\bs \Phi}({\bf X})
\end{equation}
Iterating this relation, we obtain that for each $n \geq 1$
\begin{equation}
\label{eq:expre-ImhatQ-Phi-itere}
{\bf X} = {\bf Y} + \sum_{k=1}^{n} {\bs \Phi}^{k}({\bf Y}) + {\bs \Phi}^{n+1}({\bf X})
\end{equation}
The general idea of the proof is to recognize that matrix $\mathcal{T}_{N,L}({\bf Y})$ 
is positive definite if $z$ belongs to a set $E_N$ defined by (\ref{eq:def-EN}). This implies that for $z \in E_N$, then ${\bs \Phi}^{k}({\bf Y}) > 0$ for each $k \geq 1$. Therefore, (\ref{eq:expre-ImhatQ-Phi-itere})
and ${\bs \Phi}^{n+1}({\bf X}) > 0$ imply that for each $n$, the positive definite matrix $\sum_{k=1}^{n} {\bs \Phi}^{k}({\bf Y})$ satisfies
\begin{equation}
\label{eq:inegalite-itere}
\sum_{k=1}^{n} {\bs \Phi}^{k}({\bf Y}) \leq {\bf X} - {\bf Y}
\end{equation}
so that the series $\sum_{k=1}^{+\infty} {\bs \Phi}^{k}({\bf Y})$ appears to be convergent for $z \in E_N$. As shown below, this implies that 
$\rho({\bf D}) < 1$. We begin to prove that $\mathcal{T}_{N,L}({\bf Y})$ is positive definite on a set $E_N$. 
\begin{lemma}
\label{le:T(Y)positif}
It exists 2 nice polynomials $P_1$ and $P_2$, a nice constant $\eta_1$ and an integer $N_1$ such that 
\begin{equation}
\label{eq:borne-inf-T}
\mathcal{T}_{N,L}({\bf Y}) >  \frac{(\mathrm{Im} z)^{3}}{32(\eta_1^{2} + |z|^{2})^{2}} \; {\bf I}
\end{equation}
for $N \geq N_1$ and $z \in E_N$ where $E_N$ is defined by (\ref{eq:def-EN}). 
\end{lemma}
{\bf Proof.} We show that it exist a nice constant $\eta_1> 0$ and 2 nice polynomials $P_1$ and $P_2$ such that for each $\nu \in [0,1]$, 
\begin{equation}
\label{eq:borneinf-TNL(Y)}
{\bf a}_L(\nu)^{*} \, {\bf Y} \, {\bf a}_L(\nu) > \frac{(\mathrm{Im} z)^{3}}{16(\eta_1^{2} + |z|^{2})^{2}} - \frac{L^{3/2}}{MN} \, P_1(|z|) \, P_2(1/\mathrm{Im}z)
\end{equation}
For this, we first note that 
$$ 
{\bf a}_L(\nu)^{*} {\bf R} {\bf R}^{*} {\bf a}_L(\nu) \geq \left| {\bf a}_L(\nu)^{*} {\bf R} {\bf a}_L(\nu) \right|^{2} \geq \left( {\bf a}_L(\nu)^{*} \mathrm{Im}({\bf R}) {\bf a}_L(\nu) \right)^{2}
$$
As ${\bf R}(z)$ is the Stieltjes transform of a positive matrix-valued measure ${\bs \mu}_{{\bf R}}$ (see Lemma \ref{le:proprietes-H-R}), it holds that 
$$
{\bf a}_L(\nu)^{*} \mathrm{Im}({\bf R}) {\bf a}_L(\nu) = \mathrm{Im} z \; \int_{\mathbb{R}^{+}} \frac{ {\bf a}_L(\nu)^{*} \, d{\bs \mu}_{{\bf R}}(\lambda) \, {\bf a}_L(\nu)}{|\lambda - z|^{2}}
$$
We claim that it exists $\eta_1 > 0$ and an integer $N_0$ such that 
\begin{equation}
\label{eq:tension}
{\bf a}_L(\nu)^{*} \, {\bs \mu}_{{\bf R}}\left([0,\eta_1] \right) \, {\bf a}_L(\nu) > \frac{1}{2}
\end{equation}
for each $\nu \in [0,1]$ and for each $N > N_0$.  In effect, as $c_N \rightarrow c_*$, 
it exists a nice constant $\eta_1$ for which $\mu_{\sigma^{2}, c_N}([0,\eta_1]) > \frac{3}{4}$ for each $N$. We consider the sequence of analytic functions $(g_n(z))_{n \in \mathbb{N}}$ defined by (\ref{eq:def-g_n}). 
If $n= \phi(k,N)$, $g_n(z)$ is the Stieltjes transform of measure $\mu_n$ defined by  $\mu_n = {\bf a}_L(\nu_k^{(N)})^{*} \, {\bs \mu}_{{\bf R}} \, {\bf a}_L(\nu_k^{(N)}) - \mu_{\sigma^{2}, c_N}$. 
Therefore, (\ref{eq:sup-discret}) implies that sequence $(\mu_n)_{n \in \mathbb{N}}$ converges weakly towards $0$. As the Marcenko-Pastur distribution is absolutely continuous, this leads to
$$
\lim_{N \rightarrow +\infty} \sup_{0 \leq k \leq K(N)} \left| {\bf a}_L(\nu_k^{(N)})^{*} \, {\bs \mu}_{{\bf R}}\left([0,\eta_1] \right) \, {\bf a}_L(\nu_k^{(N)})
- \mu_{\sigma^{2}, c_N}([0,\eta_1]) \right| = 0
$$
This implies the existence of $N_0^{'} \in \mathbb{N}$ such that 
$$
\sup_{0 \leq k \leq K(N)} {\bf a}_L(\nu_k^{(N)})^{*} \, {\bs \mu}_{{\bf R}}\left([0,\eta_1] \right) \, {\bf a}_L(\nu_k^{(N)}) > \frac{5}{8}
$$
for each $N \geq N_0^{'}$. As mentioned above, for each $\nu \in [0,1]$, it exists an index $k$, $0 \leq k \leq K(N)$ such that 
$|\nu - \nu_k^{(N)}| \leq \frac{1}{2 K(N)}$. As 
$$
\| {\bf a}_L(\nu) - {\bf a}_L(\nu_k^{(N)}) \| = \mathcal{O}\left( L(N) |\nu - \nu_k^{(N)}|  \right)= o(1)
$$
it is easy to check that
$$
{\bf a}_L(\nu)^{*} {\bs \mu}_{{\bf R}}\left([0,\eta_1] \right) {\bf a}_L(\nu) - {\bf a}_L(\nu_k^{(N)})^{*} {\bs \mu}_{{\bf R}}\left([0,\eta_1] \right) {\bf a}_L(\nu_k^{(N)})
\rightarrow 0
$$
which implies the existence of an integer $N_0 \geq N_0^{'}$ for which 
$$
\sup_{\nu \in [0,1]} {\bf a}_L(\nu)^{*} \, {\bs \mu}_{{\bf R}}\left([0,\eta_1] \right) \, {\bf a}_L(\nu) > \frac{1}{2}
$$
for each $N \geq N_0$, as expected. 

It is clear that 
$$
{\bf a}_L(\nu)^{*} \mathrm{Im}({\bf R}) {\bf a}_L(\nu) \geq  \mathrm{Im} z \; \int_{0}^{\eta_1} \frac{ {\bf a}_L(\nu)^{*} \, d{\bs \mu}_{{\bf R}}(\lambda) \, {\bf a}_L(\nu)}{|\lambda - z|^{2}}
$$
As $|\lambda - z|^{2} \leq 2(\lambda^{2}+|z|^{2}) \leq 2(\eta_1^{2} + |z|^{2})$ if $\lambda \in [0,\eta_1]$, it holds that
$$
{\bf a}_L(\nu)^{*} \mathrm{Im}({\bf R}) {\bf a}_L(\nu) \geq  \frac{\mathrm{Im}z}{4(\eta_1^{2} + |z|^{2})}
$$
and that 
$$ 
{\bf a}_L(\nu)^{*} {\bf R} {\bf R}^{*} {\bf a}_L(\nu) \geq \frac{(\mathrm{Im}z)^{2}}{16(\eta_1^{2} + |z|^{2})^{2}}
$$
for each $\nu \in [0,1]$. (\ref{eq:entrees-trace-Delta}) implies that for each $\nu$, 
\begin{equation}
\label{eq:inegalite-ImDelta}
\left| {\bf a}_L(\nu)^{*} \, \mathrm{Im} \hat{{\bs \Delta}} \,  {\bf a}_L(\nu) \right| \leq \frac{L^{3/2}}{MN} \, P_1(|z|) P_2(\frac{1}{\mathrm{Im} z})
\end{equation}
for some nice polynomials $P_1$ and $P_2$, which, in turn, leads to (\ref{eq:borneinf-TNL(Y)}). If we denote by $E_{N}$ the subset of $\mathbb{C}^{+}$ defined 
by $\frac{L^{3/2}}{MN} \, P_1(|z|) P_2(\frac{1}{\mathrm{Im} z}) < \frac{1}{2} \frac{(\mathrm{Im}z)^{3}}{16(\eta_1^{2} + |z|^{2})^{2}}$, then,
${\bf Y} = \mathrm{Im}(\hat{{\bs \Delta}}) + \mathrm{Im}z \, {\bf R} {\bf R}^{*}$ verifies
\begin{equation}
\label{eq:pourn=0}
\inf_{\nu \in [0,1]} {\bf a}_L(\nu)^{*} \, {\bf Y} \, {\bf a}_L(\nu) \, >  \frac{(\mathrm{Im}z)^{3}}{32(\eta_1^{2} + |z|^{2})^{2}}
\end{equation}
for each $z \in E_{N}$. As 
$$
{\bf a}_L(\nu)^{*} \, {\bf Y} \, {\bf a}_L(\nu) = \sum_{l=-(L-1)}^{L-1} \tau({\bf Y})(l) e^{-2 i \pi l \nu} 
$$
we obtain that 
$$
\inf_{\nu \in [0,1]} \sum_{l=-(L-1)}^{L-1}  \tau({\bf Y})(l) e^{-2 i \pi l \nu} >  \frac{(\mathrm{Im}z)^{3}}{32(\eta_1^{2} + |z|^{2})^{2}}
$$
for $z \in E_N$. If we denote $\alpha(z) =  \frac{(\mathrm{Im}z)^{3}}{32(\eta_1^{2} + |z|^{2})^{2}}$, 
this implies that $\left(  \tau({\bf Y})(l) - \alpha \, \delta(l=0) \right)_{l=-(L-1)}^{L-1}$
coincide with Fourier coefficients of a positive function. Therefore, matrix 
$\mathcal{T}_{N,L}({\bf Y}) - \alpha {\bf I}$ is positive definite (see \cite{grenander-szego}, 1.11 (a)), which implies  that 
(\ref{eq:borne-inf-T}) holds. Lemma \ref{le:T(Y)positif} follows from the observation that the set $E_{N}$ can be written as
(\ref{eq:def-EN}) for some other pair of nice polynomials $P_1, P_2$. \\

We now complete the proof of item 2 of Proposition (\ref{prop:proprietes-Phi}). We establish that for $N$ fixed and large enough and $z \in E_N$, then for each $L$--dimensional 
vector ${\bf b}$, 
${\bf D}^{n} {\bf b} \rightarrow 0$ when $n \rightarrow +\infty$, a property equivalent to $\rho({\bf D}) < 1$. We emphasize that in the forthcoming analysis, $N$, and therefore $L$, are assumed to be fixed parameters. As matrix $\mathcal{T}_{N,L}({\bf Y}) > \alpha(z) {\bf I}_N > 0$ on the set $E_N$ for $N$ large enough, (\ref{eq:inegalite-itere}) is valid there. This implies that the positive definite matrix-valued series $\sum_{n=1}^{+\infty} {\bs \Phi}^{n}({\bf Y})$ is convergent, in the sense that for 
each unit norm $L$--dimensional vector ${\bf u}$, then $\sum_{n=1}^{+\infty} {\bf u}^{*} {\bs \Phi}^{n}({\bf Y}) {\bf u} < +\infty$. Using the 
polarization identity, we obtain that the series $\sum_{n=1}^{+\infty} {\bf u}_1^{*} {\bs \Phi}^{n}({\bf Y}) {\bf u}_2$ is convergent for each pair of 
unit norm vectors $({\bf u}_1, {\bf u}_2)$. This implies that each entry of $ {\bs \Phi}^{n}({\bf Y})$ converges towards $0$ when $n \rightarrow +\infty$, and that the same property holds true
for each component of vector ${\bs \tau}\left( {\bs \Phi}^{n}({\bf Y}) \right)$. This vector of course coincides with ${\bf D}^{n} {\bs \tau}({\bf Y})$. 
We have thus shown that ${\bf D}^{n} {\bs \tau}({\bf Y}) \rightarrow 0$ when $n \rightarrow +\infty$. We now establish that this property holds, not only for 
vector ${\bs \tau}({\bf Y})$, but also for each $(2L-1)$--dimensional vector. We consider any positive hermitian $L \times L$ matrix ${\bf Z}$
such that $\mathcal{T}_{N,L}({\bf Y}) - \mathcal{T}_{N,L}({\bf Z}) \geq 0$. Then, it is clear that for each $n \geq 1$, 
$0 \leq {\bs \Phi}^{n}({\bf Z}) \leq {\bs \Phi}^{n}({\bf Y})$, and that the series $\sum_{n=1}^{\infty} {\bs \Phi}^{n}({\bf Z})$ is convergent. As above, 
this implies that ${\bf D}^{n} {\bs \tau}({\bf Z}) \rightarrow 0$ when $n \rightarrow +\infty$. If now ${\bf Z}$ is any positive hermitian matrix, 
it holds that $0 \leq  \mathcal{T}_{N,L} \left( \frac{\alpha(z)}{\| {\bf Z} \|} {\bf Z} \right) \leq \mathcal{T}_{N,L}({\bf Y})$
because $\mathcal{T}_{N,L}({\bf Z}) \leq \| \mathcal{T}_{N,L}({\bf Z}) \| \, {\bf I} \leq \| {\bf Z} \| \, {\bf I}$. This implies that 
${\bf D}^{n}\left( \frac{\alpha(z)}{\| {\bf Z} \|} \, {\bs \tau}({\bf Z})  \right) \rightarrow 0$, or equivalently that ${\bf D}^{n}{\bs \tau}({\bf Z}) \rightarrow 0$
for each positive hermitian matrix ${\bf Z}$. This property holds in particular for positive rank one matrices ${\bf h} {\bf h}^{*}$, and thus for 
linear combination (with complex coefficients) of such matrices, and in particular for hermitian (non necessarily positive) matrices. We now consider any $L \times L$ matrix 
${\bf B}$. It can be written as ${\bf B} = \mathrm{Re}({\bf B}) + i \, \mathrm{Im}({\bf B})$, i.e. as a linear combination of hermitian matrices. 
Therefore, it holds that ${\bf D}^{n}{\bs \tau}({\bf B}) \rightarrow 0$ for any $L \times L$ matrix. The conclusion follows from the obvious 
observation that any $(2L-1)$--dimensional vector ${\bf b}$ can be written as ${\bf b} = {\bs \tau}({\bf B})$ for some $L \times L$ matrix ${\bf B}$. This completes the proof 
of item 2 of Proposition (\ref{prop:proprietes-Phi}). \\

We finally establish item 3. We assume that $z \in E_N$ and that $N$ is large enough. We first remark that, as $\mathcal{T}_{N,L}({\bf Y}) \geq \alpha(z) {\bf I}_N$, then, for each $n \geq 1$, it holds that 
${\bs \Phi}^{n}({\bf Y}) \geq \alpha(z) \, {\bs \Phi}^{n}({\bf I})$. We also note that ${\bs \Phi}^{n}({\bf I}) > 0$
for each $n$ which implies that 
$$
{\bf a}_L(\nu)^{*} \, {\bs \Phi}^{n}({\bf Y}) \, {\bf a}_L(\nu) \geq \alpha(z) \, {\bf a}_L(\nu)^{*} \, {\bs \Phi}^{n}({\bf I}) \, {\bf a}_L(\nu) > 0
$$
for each $\nu$. We also remark that this inequality also holds for $n=0$ (see (\ref{eq:pourn=0})). We recall that for each $L \times L$ matrix ${\bf B}$, then 
\begin{equation}
\label{eq:expre-symbole-TB}
{\bf a}_L(\nu)^{*} \, {\bf B} \, {\bf a}_L(\nu) = \sum_{l=-(L-1)}^{L-1} {\bs \tau}({\bf B})(l) e^{-2i \pi l \nu}
\end{equation}
Using this identity for ${\bf B} = {\bs \Phi}^{n}({\bf Y})$ and ${\bf B} = {\bs \Phi}^{n}({\bf I})$ and using that ${\bs \tau}({\bf I}) = {\bf e}_0$, we obtain that
$$
\sum_{l=-(L-1)}^{L-1} \left({\bf D}^{n} {\bs \tau}({\bf Y}) \right)(l) e^{-2i \pi l \nu} \geq \alpha(z) \, 
\sum_{l=-(L-1)}^{L-1} \left({\bf D}^{n} {\bf e}_0 \right)(l) e^{-2i \pi l \nu} > 0
$$
As $({\bf I} - {\bf D})^{-1} = \sum_{n=0}^{+\infty} {\bf D}^{n}$, 
we finally obtain that  
$$
0 < \sum_{l=-(L-1)}^{L-1} {\bf f}_l e^{-2i \pi l \nu} \leq \frac{1}{\alpha(z)} \sum_{l=-(L-1)}^{L-1} \left(({\bf I} - {\bf D})^{-1} {\bs \tau}({\bf Y})\right)(l) e^{-2 i \pi \l \nu}
$$
The conclusion follows from the observation that ${\bs \tau}({\bf X}) =  {\bs \tau}({\bf Y}) + {\bf D} \, {\bs \tau}({\bf X})$ and that
${\bs \tau}({\bf X}) = \left({\bf I} - {\bf D}\right)^{-1} {\bs \tau}({\bf Y})$. Therefore, 
$$
\sum_{l=-(L-1)}^{L-1} \left(({\bf I} - {\bf D})^{-1} {\bs \tau}({\bf Y})\right)(l) e^{-2 i \pi \l \nu}
$$ 
coincides with ${\bf a}_L(\nu)^{*} \, {\bf X} \, {\bf a}_L(\nu)$, a term which is upperbounded by $\frac{1}{\mathrm{Im}z}$ on $\mathbb{C}^{+}$. \\

We now make the appropriate connections between $({\bs \Phi}_0, {\bf D}^{(0)})$ and $({\bs \Phi}, {\bf D})$, and establish the following Proposition.
\begin{proposition}
\label{prop:lien-D0-D}
If $N$ is large enough and if $z$ belongs to the set $E_N$ defined by (\ref{eq:def-EN}), matrix ${\bf I} - {\bf D}^{(0)}$ is invertible, and
for each matrix $L \times L$ matrix ${\bf X}$, it holds that 
\begin{equation}
\label{eq:controle-(I-D0)-1}
\sup_{\nu \in [0,1]} \left| \sum_{l=-(L-1)}^{L-1} \left( ({\bf I} - {\bf D}^{(0)})^{-1} {\bs \tau}({\bf X}) \right)(l)  e^{-2 i \pi l \nu} \right| 
\, \leq \, \frac{\| \mathcal{T}_{N,L}({\bf X}) \|}{2} \, \left( \frac{1}{1 - \sigma^{4}c_N |z t(z) \tilde{t}(z)|^{2}} \; + 
\sum_{l=-(L-1)}^{L-1} {\bf f}_l e^{-2 i \pi l \nu} \right)
\end{equation}
\end{proposition}
{\bf Proof.} We first establish by induction that
\begin{equation}
\label{eq:inegalite-Phi0-Phi}
({\bs \Phi}^{(0)})^{n}({\bf X}) \left({\bs \Phi}^{(0)})^{n}({\bf X}) \right)^{*} \leq \| \mathcal{T}_{N,L}({\bf X}) \|^{2} \; \left(\sigma^{4}c_N |z t(z) \tilde{t}(z)|^{2}\right)^{n} \;  {\bs \Phi}^{n}({\bf I})
\end{equation}
for each $n \geq 1$. We first verify that (\ref{eq:inegalite-Phi0-Phi}) holds for $n=1$. Using Proposition (\ref{prop:T(H)T(H)*}), we obtain that 
$$
\mathcal{T}_{L,L}\left( {\bf H} \mathcal{T}_{N,L}({\bf X}) \right) \left[ \mathcal{T}_{L,L}\left( {\bf H} \mathcal{T}_{N,L}({\bf X}) \right) \right]^{*}
\leq  \mathcal{T}_{L,L}\left( {\bf H} \mathcal{T}_{N,L}({\bf X}) \mathcal{T}_{N,L}({\bf X})^{*} {\bf H}^{*} \right) 
$$
Remarking that $\mathcal{T}_{N,L}({\bf X}) \mathcal{T}_{N,L}({\bf X})^{*}  \leq \| \mathcal{T}_{N,L}({\bf X}) \|^{2} \, {\bf I}$, we get that
$$
\mathcal{T}_{L,L}\left( {\bf H} \mathcal{T}_{N,L}({\bf X}) \right) \left[ \mathcal{T}_{L,L}\left( {\bf H} \mathcal{T}_{N,L}({\bf X}) \right) \right]^{*}
\leq 
\| \mathcal{T}_{N,L}({\bf X}) \|^{2} \, \mathcal{T}_{L,L}\left( {\bf H} {\bf H}^{*} \right)
$$
This and the identity ${\bs \Phi}({\bf I}) = \sigma^{4} c_N {\bf R} \mathcal{T}_{L,L}({\bf H} {\bf H}^{*}) {\bf R}^{*}$ imply immediately  (\ref{eq:inegalite-Phi0-Phi}) for $n=1$. We assume that (\ref{eq:inegalite-Phi0-Phi}) holds until integer $n-1$. By
Proposition \ref{prop:T(H)T(H)*}, we get that 
\begin{multline}
({\bs \Phi}^{(0)})^{n}({\bf X}) \left(({\bs \Phi}^{(0)})^{n}({\bf X}) \right)^{*} \leq  \\ \left|\sigma^{4}c_N z t(z) \tilde{t}(z) \right|^{2}  \, 
{\bf R} \mathcal{T}_{L,L}\left[ {\bf H} \mathcal{T}_{N,L}\left(({\bs \Phi}^{(0)})^{n-1}({\bf X})\right) \left( \mathcal{T}_{N,L} \left(({\bs \Phi}^{(0)})^{n-1}({\bf X})\right)\right)^{*} {\bf H}^* \right] {\bf R}^{*}
\end{multline}
Using again Proposition (\ref{prop:T(H)T(H)*}), we obtain that 
$$
\mathcal{T}_{N,L}\left(({\bs \Phi}^{(0)})^{n-1}({\bf X})\right) \left( \mathcal{T}_{N,L}\left(({\bs \Phi}^{(0)})^{n-1}({\bf X})\right)\right)^{*} \leq \mathcal{T}_{N,L}\left(({\bs \Phi}^{(0)})^{n-1}({\bf X})\left[({\bs \Phi}^{(0)})^{n-1}({\bf X})\right]^{*} \right)
$$
 (\ref{eq:inegalite-Phi0-Phi}) for integer $n-1$ yields to 
$$
({\bs \Phi}^{(0)})^{n}({\bf X}) \left(({\bs \Phi}^{(0)})^{n}({\bf X}) \right)^{*} \leq \| \mathcal{T}_{N,L}({\bf X}) \|^{2} \,  (\sigma^{4}c_N)^{n+1} \,  |z t(z) \tilde{t}(z)|^{2n} \, 
{\bf R} \mathcal{T}_{N,L}\left( {\bf H} {\bs \Phi}^{n-1}({\bf I}) {\bf H}^{*} \right) {\bf R}^{*}
$$
(\ref{eq:inegalite-Phi0-Phi}) for integer $n$  directly follows from ${\bs \Phi}^{n}({\bf I}) = \sigma^{4} c_N \, {\bf R} \mathcal{T}_{N,L}\left( {\bf H} {\bs \Phi}^{n-1}({\bf I}) {\bf H}^{*} \right) {\bf R}^{*}$. 

We now prove that if $z \in E_N$  defined by (\ref{eq:def-EN}) and if $N$ is large enough, then, for each $(2L-1)$--dimensional vector ${\bf x}$, it holds that $\left({\bf D}^{(0)}\right)^{n} {\bf x} \rightarrow 0$, a condition which is equivalent to 
$\rho({\bf D}^{(0)}) < 1$. For this, we observe that each vector ${\bf x}$ can be written as ${\bf x} = {\bs \tau}({\bf X})$
for some $L \times L$ matrix ${\bf X}$. The entries of Toeplitz matrix 
$\mathcal{T}_{L,L}\left( ({\bs \Phi}^{(0)})^{n}({\bf X}) \right)$ are the components of 
vector $\left({\bf D}^{(0)}\right)^{n} {\bs \tau}({\bf X})$. Therefore, condition 
$\left({\bf D}^{(0)}\right)^{n} {\bf x} \rightarrow 0$ is equivalent to 
$\| \mathcal{T}_{L,L}\left( ({\bs \Phi}^{(0)})^{n}({\bf X}) \right) \| \rightarrow 0$.
We now prove that 
$$
\sup_{\nu \in [0,1]} \left| {\bf a}_L(\nu)^{*} ({\bs \Phi}^{(0)})^{n}({\bf X}) {\bf a}_L(\nu) \right| \rightarrow 0
$$
a condition which implies $\| \mathcal{T}_{L,L}\left( ({\bs \Phi}^{(0)})^{n}({\bf X}) \right) \| \rightarrow 0$ by Proposition \ref{prop:contractant}, and thus that $\rho({\bf D}^{(0)}) < 1$. 
It is clear that 
\begin{equation}
\label{eq:schwartz}
\left| {\bf a}_L(\nu)^{*} ({\bs \Phi}^{(0)})^{n}({\bf X}) {\bf a}_L(\nu) \right|^{2} \leq  {\bf a}_L(\nu)^{*} ({\bs \Phi}^{(0)})^{n}({\bf X}) \left(({\bs \Phi}^{(0)})^{n}({\bf X})\right)^{*} {\bf a}_L(\nu) 
\end{equation}
Inequality (\ref{eq:inegalite-Phi0-Phi}) implies that 
\begin{equation}
\label{eq:consequence-inegalite-Phi0-Phi}
 {\bf a}_L(\nu)^{*} ({\bs \Phi}^{(0)})^{n}({\bf X}) \left(({\bs \Phi}^{(0)})^{n}({\bf X})\right)^{*} {\bf a}_L(\nu) \; \leq \; 
\| \mathcal{T}_{N,L}({\bf X}) \|^{2} \,  \left( \sigma^{4}c_N \,  |z t(z) \tilde{t}(z)|^{2} \right)^{n} \, {\bf a}_L(\nu)^{*} {\bs \Phi}^{n}({\bf I}) {\bf a}_L(\nu)
\end{equation}
By (\ref{eq:borneinf-zttildet}), it exists 2 nice constants $C$ and $\eta > 0$ such that 
\begin{equation}
\label{eq:borne-inf-MP}
\sigma^{4}c_N \,  |z t(z) \tilde{t}(z)|^{2}  \leq 1 - C \frac{(\eta^{2} + |z|^{2})^{2}}{(\mathrm{Im}(z))^{4}}
\end{equation}
for $N$ large enough. Moreover, it has been shown before that each entry of matrix ${\bs \Phi}^{n}({\bf I})$ converges towards 0, which 
implies that $\sup_{\nu \in [0,1]} \, {\bf a}_L(\nu)^{*} {\bs \Phi}^{n}({\bf I}) {\bf a}_L(\nu) \rightarrow 0$ (we recall 
that $L$ is assumed fixed in the present analysis). Therefore, 
$$
\sup_{\nu \in [0,1]} \, {\bf a}_L(\nu)^{*} ({\bs \Phi}^{(0)})^{n}({\bf X}) \left(({\bs \Phi}^{(0)})^{n}({\bf X})\right)^{*} {\bf a}_L(\nu) \rightarrow 0
$$
which implies that $\| \mathcal{T}_{L,L}\left( ({\bs \Phi}^{(0)})^{n}({\bf X}) \right) \|$ and $\left({\bf D}^{(0)}\right)^{n} {\bs \tau}({\bf X})$
converge towards $0$. We have thus established that $\rho({\bf D}^{(0)}) < 1$, and that matrix ${\bf I} - {\bf D}^{(0)}$ is invertible. 

We finally establish Eq. (\ref{eq:controle-(I-D0)-1}). 
Using $({\bf I} - {\bf D}^{(0)})^{-1} = \sum_{n=0}^{+\infty} \; \left({\bf D}^{(0)}\right)^{n}$ and 
$$ 
\sum_{l=-(L-1)}^{L-1} \left( ({\bf D}^{(0)})^{n} {\bs \tau}({\bf X}) \right)(l) \, e^{-2 i \pi l \nu} = {\bf a}_L(\nu)^{*} ({\bs \Phi}^{(0)})^{n}({\bf X}) {\bf a}_L(\nu)
$$ 
we first remark that 
$$
\left| \sum_{l=-(L-1)}^{L-1} \left( ({\bf I} - {\bf D}^{(0)})^{-1} {\bs \tau}({\bf X}) \right)(l) \, e^{-2 i \pi l \nu} \right|
\leq \sum_{n=0}^{+\infty}  \left| {\bf a}_L(\nu)^{*} ({\bs \Phi}^{(0)})^{n}({\bf X}) {\bf a}_L(\nu) \right|
$$
Inequalities (\ref{eq:schwartz}, \ref{eq:consequence-inegalite-Phi0-Phi}) 
imply that 
$$
\left| {\bf a}_L(\nu)^{*} ({\bs \Phi}^{(0)})^{n}({\bf X}) {\bf a}_L(\nu) \right| = \left| \sum_{l=-(L-1)}^{L-1} \left( ({\bf D}^{(0)})^{n} {\bs \tau}({\bf X}) \right)(l) \, e^{-2 i \pi l \nu} \right| 
$$
is less than $\| \mathcal{T}_{N,L}({\bf X}) \| \,  \left( \sigma^{4}c_N \,  |z t(z) \tilde{t}(z)|^{2} \right)^{n/2} \,  \left( {\bf a}_L(\nu)^{*} {\bs \Phi}^{n}({\bf I}) {\bf a}_L(\nu) \right)^{1/2}$. Using the inequality $|ab| \leq \frac{(a^{2} + b^{2})}{2}$, we obtain that
$$
\left| {\bf a}_L(\nu)^{*} ({\bs \Phi}^{(0)})^{n}({\bf X}) {\bf a}_L(\nu) \right| \leq \frac{\| \mathcal{T}_{N,L}({\bf X}) \|}{2}  \,  
\left[ \left( \sigma^{4}c_N \,  |z t(z) \tilde{t}(z)|^{2} \right)^{n} +  {\bf a}_L(\nu)^{*} {\bs \Phi}^{n}({\bf I}) {\bf a}_L(\nu) \right]
$$
Summing over $n$ eventually leads to (\ref{eq:inegalite-Phi0-Phi}). \\

We are now in position to establish the main result of this section, which, eventually, implies (\ref{eq:haagerup-norme-1}).  
\begin{proposition}
\label{prop:norme-R-t-provisoire}
It exists 2 nice polynomials $P_1$ and 
$P_2$ for which
\begin{equation}
\label{eq:haagerup-norme-1-bis}
\sup_{\nu \in [0,1]} \left| {\bf a}_L(\nu)^{*} \left( {\bf R}(z) - t(z) \, {\bf I}_L \right) \, {\bf a}_L(\nu) \right| \leq \frac{L^{3/2}}{MN} \, P_1(|z|) P_2(\frac{1}{\mathrm{Im}(z)})
\end{equation}
for $N$ large enough and for each $z \in \mathbb{C}^{+}$
\end{proposition}
{\bf Proof.} 
We recall that ${\bf a}_L(\nu)^{*} \left( {\bf R}(z) - t(z) \, {\bf I}_L \right) \, {\bf a}_L(\nu)$ coincides with
$\sum_{l=-(L-1)}^{L-1} \tau({\bf R} - t {\bf I})(l) e^{-2 i \pi l \nu}$ (see (\ref{eq:expre-symbole})), and recall that by Eq. (\ref{eq:equation-R-t}), 
vector ${\bs \tau}({\bf R} - t {\bf I})$ satisfies the equation 
$$
{\bs \tau}({\bf R} - t {\bf I})  = {\bs \tau}({\bs \Gamma}) + {\bf D}^{(0)} {\bs \tau}({\bf R} - t {\bf I})
$$
where matrix ${\bs \Gamma}$ is defined by (\ref{eq:def-Upsilon-D0}). 
Proposition \ref{prop:proprietes-Phi}, Proposition \ref{prop:lien-D0-D} used in the case ${\bf X} = {\bs \Gamma}$ as well as (\ref{eq:borne-inf-MP}) imply that 
for $N$ large and $z \in E_N$, it holds that
\begin{equation}
\label{eq:controle-norme-T-R-t}
\left| \sum_{l=-(L-1)}^{L-1} \tau({\bf R} - t {\bf I})(l) e^{-2 i \pi l \nu} \right| \leq C \, \frac{(|z|^{2}+\eta_2^{2})^{2}}{(\mathrm{Im}(z))^{4}} \, 
\| \mathcal{T}_{N,L}({\bs \Gamma}) \|
\end{equation}
for some nice constant $C$ and for $\eta_2 = \max(\eta, \eta_1)$. It is clear that 
\begin{equation}
\label{eq:controle-norme-T-Upsilon}
\| \mathcal{T}_{N,L}({\bs \Gamma}) \| \leq  P_1(|z|) P_2(\frac{1}{\mathrm{Im}(z)}) \, \| \mathcal{T}_{N,L}^{(M)}\left( \mathbb{E}({\bf Q}) - {\bf R}_M \right) \|
\end{equation}
Corollary \ref{co:converge-toeplitzifie-Q-R} thus implies that (\ref{eq:haagerup-norme-1-bis}) holds for $N$ large enough and $z \in E_N$. It remains to establish that (\ref{eq:haagerup-norme-1-bis}) also holds on the complementary 
$E_N^{c}$ of $E_N$. For this, we remark that on $E_N^{c}$, $1 < \frac{L^{3/2}}{MN} \, P_1(|z|) P_2(\frac{1}{\mathrm{Im}(z)})$. 
As $\sup_{\nu \in [0,1]} \left| {\bf a}_L(\nu)^{*} \left( {\bf R}(z) - t(z) \, {\bf I}_L \right) \, {\bf a}_L(\nu) \right| \leq \frac{2}{\mathrm{Im}(z)}$ on $\mathbb{C}^{+}$, we obtain that 
$$
\sup_{\nu \in [0,1]} \left| {\bf a}_L(\nu)^{*} \left( {\bf R}(z) - t(z) \, {\bf I}_L \right) \, {\bf a}_L(\nu) \right| \leq  \frac{1}{\mathrm{Im}(z)} \; \frac{L^{3/2}}{MN} \, P_1(|z|) P_2(\frac{1}{\mathrm{Im}(z)})
$$
for $z \in E_N^{c}$. This, in turn, shows that (\ref{eq:haagerup-norme-1-bis}) holds for $N$ large enough and for each $z \in \mathbb{C}^{+}$. \\

\begin{remark}
\label{re:rate-b*(R-tI)b}
We note that this property also implies that any quadratic form of ${\bf R} - t \, {\bf I}$ converges towards $0$ at rate 
$\frac{L^{3/2}}{MN}$. Using the polarization identity, it is sufficient to prove that ${\bf b}^{*} \left( {\bf R} - t \, {\bf I} \right)  {\bf b}$ 
is a $\mathcal{O}(\frac{L^{3/2}}{MN})$ term for each uniformly bounded deterministic vector ${\bf b}$. We consider 
Eqs. (\ref{eq:expre-entrees-R-t-1}, \ref{eq:expre-entrees-R-t-2}), and note that the righthandside of (\ref{eq:expre-entrees-R-t-1}) and (\ref{eq:expre-entrees-R-t-2})
 are bounded, up to constant terms depending on $z$ (and not on the dimensions $L,M,N$) by 
$\| \mathcal{T}^{(M)}_{N,L}\left[ \mathbb{E}({\bf Q}) - {\bf R}_M \right] \|$ and $\| \mathcal{T}_{N,L}\left( {\bf R} - t \, {\bf I} \right) \|$
respectively. 
\end{remark}

\section{Proof of (\ref{eq:identite-fondamentale-facile})}
\label{sec:improved-R-t}
The purpose of this section is to establish the identity (\ref{eq:identite-fondamentale-facile}). For this, we have essentially to 
control the term $\frac{1}{L} \mathrm{Tr}\left( {\bf R} - t \, {\bf I} \right)$. More precisely, we prove the following proposition. 
\begin{proposition}
\label{prop:Trace-R-t}
It exists nice polynomials $P_1$ and $P_2$ such that 
\begin{equation}
\label{eq:sup-(R-t)A}
\sup_{\| {\bf A} \| \leq 1} \left| \frac{1}{L} \mathrm{Tr} \left[ ({\bf R} - t {\bf I}_L)  {\bf A} \right] \right| \, \leq \frac{L}{MN} \, P_1(z) P_2(1/\mathrm{Im}z)
\end{equation}
for each $z \in F_N^{(3/2)}$ where $F_N^{(3/2)}$ is a subset of $\mathbb{C}^{+}$ defined by
\begin{equation}
\label{eq:def-FN-bis}
F_N^{(3/2)} = \{ z \in \mathbb{C}^{+}, \frac{L^{3/2}}{MN} \, Q_1(z) Q_2(1/\mathrm{Im}z) \leq 1 \}
\end{equation}
for some nice polynomials $Q_1$ and $Q_2$. 
\end{proposition}
{\bf Proof.} In the following, we denote by $\beta({\bf A})$ the term $\frac{1}{L} \mathrm{Tr} \left[ ({\bf R} - t {\bf I}_L)  {\bf A} \right]$. 
We write (\ref{eq:decomposition}) as
\begin{eqnarray}
\label{eq:decomposition-bis}
\frac{1}{L} \mathrm{Tr} \left[ ({\bf R} - t {\bf I}_L)  {\bf A} \right] = & -  \sigma^{4}  c_N z t(z) \tilde{t}(z)\; \frac{1}{ML} \mathrm{Tr} \left[ \left( {\bf E}({\bf Q}) - {\bf I}_M \otimes {\bf R} \right) \left( {\bf I}_M \otimes {\bf G}({\bf A}) \right) \right] - \\  &   \sigma^{4}  c_N z t(z) \tilde{t}(z) \; \frac{1}{L} 
\mathrm{Tr} \left( {\bf R} - \, t {\bf I} ) \mathcal{T}_{L,L} \left[ \left(\mathcal{T}_{N,L}({\bf A} {\bf R})\right) {\bf H}  \right]  \right) \nonumber
\end{eqnarray}
We denote by $\epsilon({\bf A})$ the first term of the righthandside of (\ref{eq:decomposition-bis}). (\ref{eq:inegalite-norme-G(A)}) and Proposition \ref{prop:controle-Delta} imply that
$\sup_{\| {\bf A} \| \leq 1} \left| \epsilon({\bf A} \right| \leq \frac{L}{MN} P_1(|z|) P_2(1/\mathrm{Im}z)$ for some nice polynomials $P_1$ and $P_2$. 
In order to evaluate the contribution of the second term of the righthandside of  (\ref{eq:decomposition-bis}), 
we remark that matrices ${\bf R}(z)$ and ${\bf H}(z)$ should be ``close'' from 
$t(z) {\bf I}_L$ and $- z \tilde{t}(z) \, {\bf I}_N$ respectively. It is thus appropriate to rewrite (\ref{eq:decomposition-bis}) as
\begin{align}
\label{eq:expre-Tr-R-t}
\frac{1}{L} \mathrm{Tr}\left(({\bf R} - t \, {\bf I}) {\bf A}\right) = -z t(z) \tilde{t}(z) \sigma^{4} c_N\;  \frac{1}{ML} \mathrm{Tr} \left[ \left( \mathbb{E}({\bf Q} - {\bf I}_M \otimes {\bf R} \right) {\bf I}_M \otimes {\bf G}({\bf A}) \right] + \\
\nonumber
(z t(z) \tilde{t}(z))^{2} \sigma^{4} c_N\; \frac{1}{L} \mathrm{Tr} \left[ ({\bf R} - t \, {\bf I}) \mathcal{T}_{L,L}\left(  \mathcal{T}_{N,L}({\bf A}) \right) \right] + \\
\nonumber
(z  \tilde{t}(z))^{2} t(z) \sigma^{4} c_N\; \frac{1}{L} \mathrm{Tr} \left[ ({\bf R} - t \, {\bf I}) \mathcal{T}_{L,L}\left(  \mathcal{T}_{N,L} \left[ {\bf A}  ({\bf R} - t \, {\bf I}) \right] 
\right) \,\right] - \\
\nonumber
z (t(z))^{2} \tilde{t}(z) \sigma^{4} c_N\;  \frac{1}{L} \mathrm{Tr}  \left[ \left({\bf R} - t \, {\bf I} \right)  \mathcal{T}_{L,L}\left(   \mathcal{T}_{N,L}({\bf A}) ({\bf H} + z \tilde{t}(z) \, {\bf I}) \, \right) \right] - \\
\nonumber
z t(z) \tilde{t}(z) \sigma^{4} c_N\;  \frac{1}{L} \mathrm{Tr} \left[ \left( {\bf R} - t \, {\bf I} \right)  \mathcal{T}_{L,L}\left(  \mathcal{T}_{N,L}\left[ {\bf A} ({\bf R} - t \, {\bf I}) \right]({\bf H} + z \tilde{t}(z) \, {\bf I}) \,  \right)  \right] 
\end{align}
We denote by $\alpha_1({\bf A}), \alpha_2({\bf A}), \alpha_3({\bf A})$, and $\alpha_4({\bf A})$ the second, third, fourth and fifth terms of the righthandside 
of the above equation respectively. 

We first study the term $\alpha_1({\bf A})$. We first recall that for each $z \in \mathbb{C}^{+}$ and $N$ large enough, it holds that
$$
\sigma^{4} c_N | z t(z) \tilde{t}(z) |^{2} < 1 - C \, \frac{(\mathrm{Im}z)^{4}}{(\eta^{2} + |z|^{2})^{2}}
$$
where $C$ and $\eta$ are nice constants (see Eq. (\ref{eq:borneinf-zttildet})). Moreover, for each ${\bf A}, \|{\bf A} \| \leq 1$, it is clear that
$$
\left| \frac{1}{L} \mathrm{Tr} \left[ ({\bf R} - t \, {\bf I}) \mathcal{T}_{L,L}\left(  \mathcal{T}_{N,L}({\bf A}) \right) \right] \right| \leq 
\sup_{\| {\bf B} \| \leq 1} |\beta({\bf B})| \; \|  \mathcal{T}_{L,L}\left(  \mathcal{T}_{N,L}({\bf A}) \right) \| \leq \sup_{\| {\bf B} \| \leq 1} |\beta({\bf B})|
$$
because $\|  \mathcal{T}_{L,L}\left(  \mathcal{T}_{N,L}({\bf A}) \right) \| \leq \| {\bf A} \| \leq 1$ (see Proposition \ref{prop:contractant}). This shows that
$$
\sup_{\| {\bf A} \| \leq 1} |\alpha_1({\bf A})| \leq \left( 1 - C \, \frac{(\mathrm{Im}z)^{4}}{(\eta^{2} + |z|^{2})^{2}} \right)  \, \sup_{\| {\bf A} \| \leq 1} |\beta({\bf A}|
$$ 
We now evaluate the behaviour of $\alpha_2({\bf A})$. We first use (\ref{eq:utile}) to obtain that
$$
\alpha_2({\bf A}) = (z  \tilde{t}(z))^{2} t(z) \sigma^{4} c_N\; \frac{1}{L} \mathrm{Tr} \left[  {\bf A} \,  ({\bf R} - t \, {\bf I}) \mathcal{T}_{L,L}\left(  \mathcal{T}_{N,L}({\bf R} - t \, {\bf I}) \right) \right] 
$$
We remark that for each matrix ${\bf A}, \| {\bf A} \| \leq 1$, it holds that
$$
\left| \frac{1}{L} \mathrm{Tr} \left[ ({\bf R} - t \, {\bf I}) \mathcal{T}_{L,L}\left(  \mathcal{T}_{N,L}({\bf R} - t \, {\bf I}) \right) \, {\bf A} \right] \right|
\leq \sup_{\| {\bf B} \| \leq 1} \beta({\bf B}) \, \| {\bf A} \| \| \mathcal{T}_{N,L}({\bf R} -t \, {\bf I}) \|
$$
(\ref{eq:haagerup-norme-1}) implies that 
$$
 \sup_{\| {\bf A} \| \leq 1}  |\alpha_2({\bf A})| <  \sup_{\| {\bf A} \| \leq 1} \beta({\bf A}) \; \frac{L^{3/2}}{MN} \, P_1(|z|) P_2(1/\mathrm{Im}z)
$$ 
for each $z \in \mathbb{C}^{+}$. 
The terms $\alpha_3({\bf A})$ and $\alpha_4({\bf A})$ can be handled similarly by writing ${\bf H} + z \tilde{t}(z) {\bf I}$ as
$$
{\bf H} + z \tilde{t}(z) {\bf I} = \sigma^{2} c_N z \tilde{t}(z) \; {\bf H} \, \mathcal{T}_{N,L}^{(M)}\left( \mathbb{E}({\bf Q}) - {\bf I}_M \otimes {\bf R} \right) +
\sigma^{2} c_N z \tilde{t}(z) \; {\bf H} \, \mathcal{T}_{N,L}\left( {\bf R} - t \, {\bf I} \right)
$$
In particular, it can be shown that for $i=3,4$ and $N$ large enough, it holds that 
$$
 \sup_{\| {\bf A} \| \leq 1}  |\alpha_i({\bf A})| <   \sup_{\| {\bf A} \| \leq 1} \beta({\bf A}) \;  \frac{L^{3/2}}{MN} \, P_1(|z|) P_2(1/\mathrm{Im}z)
$$ 
Therefore, it holds that 
$$
\sup_{\| {\bf A} \| \leq 1} |\beta({\bf A})|  \leq  \sup_{\| {\bf A} \| \leq 1} |\epsilon({\bf A})| + 
                                                 \sup_{\| {\bf A} \| \leq 1} \beta({\bf A}) \;  \left[ \left( 1 - C \, \frac{(\mathrm{Im}z)^{4}}{(\eta + |z|^{2})^{2}} \right) +  
\frac{L^{3/2}}{MN} \, P_1(|z|) P_2(1/\mathrm{Im}z) \right]
$$
We define the set $F_N^{(3/2)}$ as 
$$
F_N^{(3/2)} = \{ z \in \mathbb{C}^{+}, \frac{L^{3/2}}{MN} \, P_1(|z|) P_2(1/\mathrm{Im}z) \leq C/2 \; \frac{(\mathrm{Im}z)^{4}}{(\eta^{2} + |z|^{2})^{2}} \} 
$$
which can also be written as 
$$
F_N^{(3/2)} = \{ z \in \mathbb{C}^{+}, \frac{L^{3/2}}{MN} \, Q_1(|z|) Q_2(1/\mathrm{Im}z) \leq 1 \} 
$$
for some nice polynomials $Q_1$ and $Q_2$. Then, it is clear that for each $z \in F_N^{(3/2)}$, then it holds that 
$$
\sup_{\| {\bf A} \| \leq 1} |\beta({\bf A})| \leq 2/C \; \frac{(\eta^{2} + |z|^{2})^{2}}{(\mathrm{Im}z)^{4}} \; \sup_{\| {\bf A} \| \leq 1} |\epsilon({\bf A})| 
\leq \frac{L}{MN} \, P_1(|z|) P_2(1/\mathrm{Im}z) 
$$ 
for some nice polynomials $P_1$ and $P_2$. This completes the proof of Proposition \ref{prop:Trace-R-t}. \\

We conclude this section by the corollary:
\begin{corollary}
\label{co-final}
The mathematical expectation of the Stieltjes transform $\frac{1}{ML} \mathrm{Tr}({\bf Q}(z))$ of the empirical eigenvalue distribution of ${\bf W} {\bf W}^{*}$ can be written for $z \in \mathbb{C}^{+}$ as
\begin{equation}
\label{eq:ecriture-haagerup-bis}
\mathbb{E} \left[ \frac{1}{ML} \mathrm{Tr}\left({\bf Q}(z)\right) \right]= \, t(z) \, + \, \frac{L}{MN} \, \tilde{r}(z) 
\end{equation}
where $\tilde{r}(z)$ is holomorphic in $\mathbb{C}^{+}$ and satisfies 
\begin{equation}
\label{eq:borne-polynomiale}
|\tilde{r}(z)| \leq P_1(|z|) P_2(\frac{1}{\mathrm{Im}(z)})
\end{equation}
for each $z \in F_N^{(3/2)}$ defined by (\ref{eq:def-FN-bis}).  
\end{corollary}
{\bf Proof.} In order to establish (\ref{eq:ecriture-haagerup-bis}), we have to prove that 
$$
\left| \frac{1}{ML} \mathrm{Tr} \left( \mathbb{E}({\bf Q}(z)) \right) \, - \, t(z) \right| \leq P_1(|z|) P_2(\frac{1}{\mathrm{Im}(z)}) \, \frac{L}{MN}
$$
for $z \in F_N^{(3/2)}$. 
$\mathbb{E}({\bf Q}(z)) \, - \, t(z) {\bf I}$ can be written as 
$$
\mathbb{E}({\bf Q}(z))  \, - \, t(z) {\bf I}_{ML} = {\bs \Delta}(z) + {\bf I}_M \otimes {\bf R}(z) \,  - \, t(z) \, {\bf I}_{ML}
$$
Therefore, Proposition \ref{prop:controle-Delta} implies that we have just to verify that
$$
\left| \frac{1}{L} \mathrm{Tr}({\bf R} - t \, {\bf I}_L) \right|  \leq P_1(|z|) P_2(\frac{1}{\mathrm{Im}(z)}) \, \frac{L}{MN}
$$
for $z \in F_N^{(3/2)}$, a consequence of Proposition \ref{prop:Trace-R-t}. 

\section{Expansion of $\frac{1}{ML} \mathrm{Tr} \left( \mathbb{E}({\bf Q}_N(z)) \right) \, - \, t_N(z)$.}
\label{sec:expre-E(Q)-t}
{\bf Notations and definitions used in section \ref{sec:expre-E(Q)-t}.} In order to simplify the exposition of
the results presented in this section, we define the following simplified notations: 
\begin{itemize}
\item  Let $(\beta_N)_{N \geq 1}$ be a sequence depending on $N$. A term $\phi_N(z)$ depending on $N$ defined for $z \in \mathbb{C}^{+}$ will be said to be a $\mathcal{O}(\beta_N)$ term 
if it exists 2 nice polynomials $P_1$ and $P_2$ such that 
$$
|\phi_N(z)| \leq \beta_N P_1(|z|) P_2(1/\mathrm{Im} z)
$$
for $N$ large enough and for each $z$ belonging to a set defined as $F_N^{(2)}$, but possibly with other nice polynomials.
\item $C_N(z,u_1, \ldots, u_k)$ will represent a generic term depending on $N$, $z$, and on 
indices $u_1, \ldots, u_k \in \{-(L-1), \ldots, L-1 \}$, and satisfying 
$\sup_{u_1, \ldots,u_k} |C_N((z,u_1, \ldots, u_k)| = \mathcal{O}(1)$ in the sense of the above definition of 
operator $\mathcal{O}(.)$. Very often, we will not mention the dependency of $C_N(z,u_1, \ldots, u_k)$ w.r.t. 
$N$ and $z$, and use the notation $C(u_1, \ldots, u_k)$. 
\item By a real distribution, we mean a real valued continuous (in an appropriate sense) linear 
form $D$ defined on the space $\mathbb{C}^{\infty}_c(\mathbb{R})$ of all real valued 
compactly supported smooth functions defined on $\mathbb{R}$. Such a distribution can of course be extended to complex valued smooth functions defined on $\mathbb{R}$ by setting $<D,\phi_1 + i \phi_2> = 
<D, \phi_1> + i <D,\phi_2>$ for $\phi_1,\phi_2 \in \mathcal{C}^{\infty}_c(\mathbb{R})$. We also 
recall that a compactly supported distribution $D$ can be extended to a continuous linear form to 
the space $\mathbb{C}^{\infty}_b(\mathbb{R})$ of all bounded smooth functions. In particular, 
$<D, \mathbb{1}>$ represents $<D, \phi>$ where $\phi$ is any function of $\mathcal{C}^{\infty}_c(\mathbb{R})$
that is equal to $1$ on the support of $D$. 
\end{itemize}

From now on, we assume that $L$ satisfies the condition
\begin{equation}
\label{eq:stronger-condition}
L = \mathcal{O}(N^{\alpha}), \; \mbox{where $\alpha < \frac{2}{3}$}
\end{equation}
which implies that
\begin{equation}
\label{eq:consequence-stronger-condition}
\frac{L^{2}}{MN} \rightarrow 0, \, \mbox{i.e.} \,  \frac{L}{M^{2}} \rightarrow 0
\end{equation}
The goal of this section is to establish the following theorem. 
\begin{theorem}
\label{theo:developpement-r}
Under (\ref{eq:stronger-condition}),  
$\frac{1}{ML} \mathrm{Tr} \left( \mathbb{E}({\bf Q}_N(z)) \right) \, - \, t_N(z)$ can be 
expanded as 
\begin{equation}
\label{eq:developpement-r}
\frac{1}{ML} \mathrm{Tr} \left( \mathbb{E}({\bf Q}_N(z)) \right) \, - \, t_N(z) = \frac{L}{MN} \left(\hat{s}_N(z) + \frac{L^{3/2}}{MN} \hat{r}_N(z) \right) 
\end{equation}
where $\hat{s}_N(z)$ coincides with the Stieltjes transform of a distribution $\hat{D}_N$
whose support is included into $\mathcal{S}_N^{(0)} = [\sigma^{2}(1 - \sqrt{c_N})^{2}, \sigma^{2}(1 + \sqrt{c_N})^{2}]$
and which verifies $<\hat{D}_N, \mathbb{1} > = 0$, 
and where $|\hat{r}_N(z)| \leq P_1(|z|) P_2(\frac{1}{\mathrm{Im}z})$ when $z$ belongs to a set 
$F_N^{(2)}$ defined by 
\begin{equation}
\label{eq:def-EN2}
F_N^{(2)} = \{ z \in \mathbb{C}^{+}, \frac{L^{2}}{MN} Q_1(|z|) Q_2(1/\mathrm{Im} z) \leq 1 \}
\end{equation}
for some nice polynomials $Q_1$ and $Q_2$.
\end{theorem}
As shown below in section \ref{sec:conclusion}, (\ref{eq:developpement-r}) provides the desired 
almost sure location of the eigenvalues of ${\bf W}_N{\bf W}_N^{*}$.  In order to establish (\ref{eq:developpement-r}), we express 
$\frac{1}{ML} \mathrm{Tr} \left( \mathbb{E}({\bf Q}_N(z)) \right) \, - \, t_N(z)$ as 
$$
\frac{1}{ML} \mathrm{Tr} \left( \mathbb{E}({\bf Q}_N(z)) \right) \, - \, t_N(z) = \frac{1}{ML} \mathrm{Tr} {\bs \Delta}_N(z) \, + \, \frac{1}{L} \mathrm{Tr}\left( {\bf R}_N(z) \, - \, t_N(z) \, {\bf I} \right)
$$
and study the 2 terms separately. We first establish that if (\ref{eq:stronger-condition}) holds, then 
\begin{equation}
\label{eq:resultat-trace-Delta}
\frac{1}{ML} \mathrm{Tr} {\bs \Delta}_N(z) = \frac{L}{MN} s_N(z) + \left( \frac{L}{MN} \right)^{2} r_N(z)
\end{equation}
where $s_N(z)$ is the Stieltjes transform of a distribution whose support is included in 
$\mathcal{S}_N^{(0)}$, and where 
$$
|r_N(z)| \leq P_1(|z|) P_2(1/\mathrm{Im} z)
$$
for some nice polynomials $P_1$ and $P_2$ and for $z \in F_N^{(2)}$. Using Theorem 
\ref{theo:convergence-norme-tau(R-t)}, (\ref{eq:developpement-r}) will follow 
easily from (\ref{eq:resultat-trace-Delta}). 

The proof of (\ref{eq:resultat-trace-Delta}) 
is quite demanding. It needs to establish a number of intermediate results that are presented in subsection 
\ref{subsec:preuve-difficile}, and used in subsection \ref{subsec:preuve-trace-delta}. \\

\subsection{Useful results concerning the Stieltjes transforms of compactly supported distributions.}
\label{subsec:distributions}
Before establishing (\ref{eq:resultat-trace-Delta}), we need to recall some results concerning
the Stieltjes transform of compactly 
supported real distributions, and to establish that the so-called Hellfer-Sjöstrand formula, 
valid for probability measures,  
can be generalized to compactly supported distributions. \\

The following 
useful result was used in \cite{schultz-2005}, Theorem 5.4 and Lemma 5.6 (see also Theorem 4.3 in \cite{capitaine2009largest}). 
\begin{lemma}
\label{le:characterization-stieljes-support}
If $D$ is a real distribution with compact support $\mathrm{Supp}(D)$, its Stieltjes transform 
$s(z)$ is defined for each $z \in \mathbb{C} - \mathrm{Supp}(D)$ by 
$$
s(z) = <D, \frac{1}{\lambda - z} >.
$$
Then, $s$ is analytic on $\mathbb{C} - \mathrm{Supp}(D)$ and verifies the following properties:
\begin{itemize}
\item $(a)$ $s(z) \rightarrow 0$ if $|z| \rightarrow +\infty$
\item It exists a compact $\mathcal{K} \subset \mathbb{R}$ containing $\mathrm{Supp}(D)$ such that
\begin{itemize}
\item $(b)$ $s(z^{*}) = (s(z))^{*}$ for each $z \in \mathbb{C}- \mathcal{K}$
\item $(c)$ It exists an integer $n_0$ and a constant $C$ such that for each $z \in \mathbb{C} - \mathcal{K}$,
\begin{equation}
\label{eq:condition-support}
|s(z)| \leq  C \, \mathrm{Max}\left( \frac{1}{\left(\mathrm{Dist}(z,\mathcal{K})\right)^{n_0}}, 1 \right)
\end{equation}
\end{itemize}
\item If $\phi$ is an element of $\mathcal{C}^{\infty}_c(\mathbb{R})$, then
the following inversion formula holds 
\begin{equation}
\label{eq:inversion-stieltjes-distributions}
\frac{1}{\pi} \, \lim_{y \rightarrow 0^{+}} \int \phi(\lambda) \, \mathrm{Im}(s(\lambda + i y)) \, d\lambda = <D, \phi>
\end{equation}
\item If  $\lim_{|z| \rightarrow +\infty} |z s(z)| = 0$, then, it holds that 
\begin{equation}
\label{eq:masse-D-0}
< D, \mathbb{1} > = 0
\end{equation}
\end{itemize}
Conversely, if $\mathcal{K}$ is a compact subset of $\mathbb{R}$, and if $s(z)$ is a function 
analytic on $\mathbb{C} - \mathcal{K}$ satisfying $(a), (b), (c)$, then $s(z)$ is the Stieltjes transform 
of a compactly supported real distribution $D$ such that $\mathrm{Supp}(D) \subset \mathcal{K}$. In this case, 
$\mathrm{Supp}(D)$ is the set of singular points of $s(z)$. 
\end{lemma}
\begin{remark}
\begin{itemize}
\item We note that (\ref{eq:condition-support}) of course implies that 
\begin{equation}
\label{eq:stieltjes-distribution-borne-imz}
|s(z)| \leq C \, \mathrm{Max}\left( \frac{1}{(\mathrm{Im}z)^{n_0}}, 1 \right) \leq C \, \left( 1 +  \frac{1}{(\mathrm{Im}z)^{n_0}} \right)
\end{equation}
for each $z \in \mathbb{C} - \mathbb{R}$.
\item We have chosen to present Lemma \ref{le:characterization-stieljes-support} as it is stated in 
\cite{schultz-2005}. However, we mention that (b) and (c) hold for each compact subset $\mathcal{K}$ of $\mathbb{R}$ 
containing $\mathrm{Supp}(D)$. $n_0$ does not depend on the compact $\mathcal{K}$ and is related to 
the order of $D$. However, the constant $C$ does depend on $\mathcal{K}$. . 
\end{itemize}
\end{remark}

We now provide a useful example of such functions $s(z)$. 
\begin{lemma}
\label{le:exemple}
If $p \geq 1$, then function $s_N(z)$ 
defined by  
$$
s_N(z) = (t_N(z))^{p} (z \tilde{t}_N(z))^{q} \frac{1}{\left( 1 - a_N \, \sigma^{4} c_N (z \, t_N(z) \,  \tilde{t}_N(z))^{2} \right)^{n}}
$$ 
for $|a_N| \leq 1$ coincides with the Stieltjes transform of a real bounded distribution $D_N$ whose support is included in 
$\mathcal{S}_N$ for each integers $q \geq 0$ and $n \geq 0$. Moreover, $D_N$ satisfies 
(\ref{eq:masse-D-0})  as soon as $p \geq 2$. 
\end{lemma}
{\bf Proof.}  It is clear that $s_N(z^*) = (s_N(z))^{*}$ and that $s_N(z) \rightarrow 0$ if $|z| \rightarrow +\infty$ 
because $p \geq 1$ and that $z \tilde{t}(z) \rightarrow -1$. We use  
Lemma \ref{le:proprietes-t-tilde-t} to manage the term 
$$
\frac{1}{\left( 1 - a_N \, \sigma^{4} c_N (z \, t_N(z) \,  \tilde{t}_N(z))^{2} \right)^{n}}
$$
and use that $|t_N(z)| \leq \frac{1}{\mathrm{dist}(z,\mathcal{S}_N)}$ for $z \in \mathbb{C} - \mathcal{S}_N$
We also remark that 
$$
z \tilde{t}_N(z) = c_N \, \int_{\mathcal{S}_N} \frac{z}{\lambda-z} \, d\mu_{\sigma^{2},c_N}(\lambda) - (1-c_N)
$$
or equivalently that 
$$
z \tilde{t}_N(z) = c_N \, \int_{\mathcal{S}_N} \frac{\lambda}{\lambda-z} \, d\mu_{\sigma^{2},c_N}(\lambda) \, - \, 1
$$
Therefore, 
$$
|z\tilde{t}_N(z)| \leq C \, (1 + \frac{1}{\mathrm{dist}(z,\mathcal{S}_N)}) \leq C \, \max\left( 1, \frac{1}{\mathrm{dist}(z,\mathcal{S}_N)} \right)
$$
for each $z \in \mathbb{C} - \mathcal{S}_N$. Moreover, it holds that $zs(z) \rightarrow 0$ if 
$|z| \rightarrow +\infty$ as soon as $p \geq 2$.   \\

We now briefly justify that the Hellfer-Sjöstrand formula  
can be generalized to compactly supported distributions. 
In order to introduce this formula, used in the context of large random matrices in \cite{and-gui-zei-2010}, 
\cite{and-2013} 
and \cite{naj-yao-2013}, we have to define some notations. $\chi$ is a function of $\mathcal{C}^{\infty}_c(\mathbb{R})$ 
with support $[-1,1]$, and which is equal to $1$ in a neighborhood of $0$. If 
$\phi(x) \in \mathcal{C}^{\infty}_c(\mathbb{R})$, we denote by $\overline{\phi}_k$ the function of $\mathcal{C}^{\infty}_c(\mathbb{R}^{2},\mathbb{C})$
defined for $z=x+iy$ by
$$
\overline{\phi}_k(z) = \sum_{l=0}^{k} \phi^{(l)}(x) \, \frac{(iy)^{l}}{l!} \, \chi(y)
$$
Function $\partial \overline{\phi}_k$ is the "derivative"
$$
\partial \overline{\phi}_k(z) = \frac{\partial  \overline{\phi}_k(z)}{\partial x} + i  \frac{\partial  \overline{\phi}_k(z)}{\partial y}
$$
and is given by 
\begin{equation}
\label{eq:expre-derivee-phibar}
\partial \overline{\phi}_k(z) = \phi^{(k+1)}(x) \frac{(iy)^{k}}{k!}
\end{equation}
in the neighborhood of $0$ in which $\chi(y) = 1$. If $s(z)$ is the Stieltjes transform of a probability measure $\mu$, 
$s(z)$ verifies $|s(z)| \leq \frac{1}{\mathrm{Im}z}$ on $\mathbb{C}^{+}$. Therefore, (\ref{eq:expre-derivee-phibar}) implies that if $k \geq 1$, then function $\partial \overline{\phi}_k(z) \, s(z)$ is well defined near the real axis. 
The Hellfer-Sjöstrand allows to reconstruct $\int \phi(\lambda) \, d\mu(\lambda)$ as:
\begin{equation}
\label{eq:hellfer-mesure}
\int \phi(\lambda)  \, d\mu(\lambda) = \frac{1}{\pi} \, \mathrm{Re} \left( \int_{\mathbb{C}^{+}} \partial \overline{\phi}_k(z) \, s(z)  \, dxdy \right)
\end{equation}
The following Lemma extends formula (\ref{eq:hellfer-mesure}) to real compactly supported distributions.
\begin{lemma}
\label{le:hellfer-distribution}
We consider a compactly supported distribution $D$ and $s(z)$ is Stieljes transform. Then, if 
$k$ is greater than the index $n_0$ defined by (\ref{eq:stieltjes-distribution-borne-imz}), 
then $ \partial \overline{\phi}_k(z) \, s(z)$ is well defined near the real axis, and 
\begin{equation}
\label{eq:hellfer-distribution}
<D, \phi>  = \frac{1}{\pi} \, \mathrm{Re} \left( \int_{\mathbb{C}^{+}} \partial \overline{\phi}_k(z) \, s(z)  \, dxdy \right)
\end{equation}
\end{lemma}
{\bf Sketch of proof.} It is clear that $ \partial \overline{\phi}_k(z) \, s(z)$ is well defined near the real axis. Therefore, the integral at the righthandside of (\ref{eq:hellfer-distribution}) exists. By linearity, it is sufficient to establish (\ref{eq:hellfer-distribution}) if $D$ coincides with a derivative of a Dirac distribution $D = \delta_{\lambda_0}^{(p)}$ for $p \leq n_0 - 1$, i.e. $s(z) = \frac{1}{(\lambda_0 - z)^{p+1}}$. Using the integration by parts formula and the analyticity of $s(z)$ on $\mathbb{C}^{+}$, we obtain that  
$$
\frac{1}{\pi} \, \mathrm{Re} \left( \int_{\mathbb{C}^{+}} \partial \overline{\phi}_k(z) \, s(z)  \, dxdy \right) = \lim_{\epsilon \rightarrow 0} \frac{1}{\pi} \, \mathrm{Re} \left( - i  \int_{\mathbb{R}} \overline{\phi}_k(x + i \epsilon) s(x + i \epsilon) dx \right)
$$
$<D, \phi>$ is of course equal to 
$$
<D, \phi> = (-1)^{p} < \delta_{\lambda_0}, \phi^{(p)} >
$$
As the Hellfer-Sjöstrand formula is valid for measure $\delta_{\lambda_0}$ and that the 
Stieltjes transform of $\delta_{\lambda_0}$ is $\frac{1}{\lambda_0 - z}$, it holds that 
$$
< \delta_{\lambda_0}, \phi^{(p)} > =  \lim_{\epsilon \rightarrow 0} \frac{1}{\pi} \, \mathrm{Re} \left( - i  \int_{\mathbb{R}} \left(\overline{\phi^{(p)}}\right)_k(x + i \epsilon) \, \frac{1}{\lambda_0 - (x + i \epsilon)} \, dx
\right)
$$
It is clear that $ \left(\overline{\phi^{(p)}}\right)_k(x + i \epsilon) = \frac{d^{p}}{d x^{p}} \overline{\phi}_k(x + i \epsilon)$. Therefore, the integration by parts leads to 
$$
\int_{\mathbb{R}} \left(\overline{\phi^{(p)}}\right)_k(x + i \epsilon) \, \frac{1}{\lambda_0 - (x + i \epsilon)} \, dx = 
(-1)^{p} \int_{\mathbb{R}} \overline{\phi}_k(x + i \epsilon) \, \frac{1}{(\lambda_0 - (x + i \epsilon))^{p+1}} \, dx  
$$
from which (\ref{eq:hellfer-distribution}) follows immediately.

\subsection{Some useful evaluations.}
\label{subsec:preuve-difficile}
(\ref{eq:def-Delta}) and (\ref{eq:expre-gamma}) imply that $ \frac{1}{ML} \mathrm{Tr}\left( {\boldsymbol \Delta}(z) \right)$ is given by
$$
\frac{1}{ML} \mathrm{Tr}\left( {\boldsymbol \Delta}(z) \right) = \sigma^{2} c_N \sum_{l_1=-(L-1)}^{L-1} \mathbb{E} \left( \tau^{(M)}({\bf Q}^{\circ})(l_1) \, \frac{1}{ML} \mathrm{Tr}\left( {\bf Q} {\bf W} {\bf J}_N^{l_1} {\bf H}^{T} {\bf W}^{*}({\bf I}_M \otimes {\bf R}) \right)^{\circ} \right)
$$
In order to establish (\ref{eq:resultat-trace-Delta}), it is necessary to evaluate the righthandside of 
the above equation up to $\mathcal{O}(\frac{L}{MN})^{2}$ terms using the integration by parts formula. 
If we denote by $\kappa^{(2)}(l_1, l_2)$ the term defined by 
$\kappa^{(2)}(l_1, l_2) = \mathbb{E}\left(\tau^{(M)}({\bf Q}^{\circ})(l_1) \tau^{(M)}({\bf Q}^{\circ})(l_2)\right)$, 
then, we establish in the following that
\begin{multline}
\label{eq:programme-rejouissances}
\frac{1}{ML} \mathrm{Tr}\left( {\boldsymbol \Delta}(z) \right) = (\sigma^{2} c_N)^{2} 
\sum_{l_1,l_2=-(L-1)}^{L-1} \kappa^{(2)}(l_1,l_2) \, 
\mathbb{E} \left[ \frac{1}{ML} \mathrm{Tr}\left( {\bf Q} {\bf W} {\bf J}_N^{l_2} {\bf H}^{T} {\bf W}^{*}
\left({\bf I}_M \otimes \sigma^{2} {\bf R} \mathcal{T}_{L,L}({\bf H} {\bf J}_N^{*l_1} {\bf H}) {\bf R}\right) \right) \right] \\
-  (\sigma^{2} c_N)^{2} 
\sum_{l_1,l_2=-(L-1)}^{L-1} \kappa^{(2)}(l_1,l_2) \, 
\mathbb{E} \left[ \frac{1}{ML} \mathrm{Tr}\left( {\bf Q} {\bf W} {\bf J}_N^{l_2} {\bf H}^{T}  {\bf J}_N^{l_1} {\bf H}^{T}{\bf W}^{*} ({\bf I}_M \otimes {\bf R}) \right) \right]  \\ 
+ \frac{\sigma^{4} c_N}{MLN} \sum_{l_1,i=-(L-1)}^{L-1} \mathbb{E} \left[ \frac{1}{ML}  \mathrm{Tr}\left({\bf Q} ({\bf I}_M \otimes {\bf J}_L^{i}) {\bf Q} ({\bf I}_M \otimes {\bf J}_L^{l_1}) {\bf Q} {\bf W} {\bf J}_N^{i} {\bf H}^{T} {\bf W}^{*} ({\bf I}_M \otimes \sigma^{2} {\bf R} \mathcal{T}_{L,L}({\bf H} {\bf J}_N^{*l_1} {\bf H}) {\bf R}  \right) \right]  \\
- \frac{\sigma^{4} c_N}{MLN} \sum_{l_1,i=-(L-1)}^{L-1} \mathbb{E} \left[ \frac{1}{ML}  \mathrm{Tr}\left({\bf Q} ({\bf I}_M \otimes {\bf J}_L^{i}) {\bf Q} ({\bf I}_M \otimes {\bf J}_L^{l_1}) {\bf Q} {\bf W} {\bf J}_N^{i} {\bf H}^{T} {\bf J}_N^{l_1} {\bf H}^{T} {\bf W}^{*} ({\bf I}_M \otimes {\bf R} \right) \right] \\
+ (\sigma^{2} c_N)^{2} 
\sum_{l_1,l_2=-(L-1)}^{L-1)} \mathbb{E} \left[ \tau^{(M)}({\bf Q}^{\circ})(l_1) \tau^{(M)}({\bf Q}^{\circ})(l_2) \, 
\frac{1}{ML} \mathrm{Tr}\left( {\bf Q} {\bf W} {\bf J}_N^{l_2} {\bf H}^{T} {\bf W}^{*}
\left({\bf I}_M \otimes \sigma^{2} {\bf R} \mathcal{T}_{L,L}({\bf H} {\bf J}_N^{*l_1} {\bf H}) {\bf R}\right) \right)^{\circ} \right]  \\
 - (\sigma^{2} c_N)^{2} 
\sum_{l_1,l_2=-(L-1)}^{L-1} \mathbb{E} \left[ \tau^{(M)}({\bf Q}^{\circ})(l_1) \tau^{(M)}({\bf Q}^{\circ})(l_2) \, 
\frac{1}{ML} \mathrm{Tr}\left( {\bf Q} {\bf W} {\bf J}_N^{l_2} {\bf H}^{T}  {\bf J}_N^{l_1} {\bf H}^{T}{\bf W}^{*} ({\bf I}_M \otimes {\bf R}) \right)^{\circ} \right]
\end{multline}
We evaluate in closed form the third and the fourth term of the righthandside of (\ref{eq:programme-rejouissances})
up to $\mathcal{O}(\frac{L}{MN})^{2}$, prove that 
$\kappa^{(2)}(u_1,u_2) = \frac{1}{MN} C(z,u_1)  \delta(u_1+u_2=0) \, + \, \mathcal{O}(\frac{L}{(MN)^{2}})$, 
and establish that the 2 last terms of (\ref{eq:programme-rejouissances}) are $\mathcal{O}(\frac{L}{MN})^{2}$. 
In Paragraph \ref{subsub:TrQQQW}, we calculate useful quantities similar to the  
third and the fourth term of the righthandside of (\ref{eq:programme-rejouissances}), and in Paragraph 
\ref{subsub:eval-kappa}, we 
evaluate $\kappa^{(2)}(u_1,u_2)$.

\subsubsection{Evaluation of the third and fourth terms of the righthandside of (\ref{eq:programme-rejouissances}).}
\label{subsub:TrQQQW}
We first state 2 technical Lemmas.
\begin{lemma}
\label{le:generalisation-NP}
We consider uniformy bounded $ML \times ML$ matrices $({\bf C}^{s})_{s=1, \ldots, r}$ and ${\bf A}$, 
and a uniformly bounded $N \times N$ matrix ${\bf G}$. Then, for each $p \geq 2$, it holds that 
\begin{eqnarray}
\label{eq:generalisation-NP-1}
\mathbb{E} \left( \frac{1}{ML} \mathrm{Tr}\left( \Pi_{s=1}^{r} {\bf Q} {\bf C}^{s} \right)^{\circ} \right)^{p}
& = & \mathcal{O}(\frac{1}{(MN)^{p/2}}) \\
\label{eq:generalisation-NP-2}
\mathbb{E} \left[ \frac{1}{ML} \mathrm{Tr}\left( (\Pi_{s=1}^{r} {\bf Q} {\bf C}^{s}) {\bf W} {\bf G} {\bf W}^{*} {\bf A} \right)^{\circ} \right]^{p}
& = & \mathcal{O}(\frac{1}{(MN)^{p/2}})
\end{eqnarray}
\end{lemma}
{\bf Proof.} We just provide a sketch of proof. We first establish (\ref{eq:generalisation-NP-1}) and 
(\ref{eq:generalisation-NP-2}) by induction for even integers $p=2q$. For $q=1$, we use the Poincaré-Nash inequality, 
and for $q \geq 1$, we take benefit of the identity 
$$
\mathbb{E}|x|^{2q} = \left| \mathbb{E}(x^{q}) \right|^{2} + \mathrm{Var}(x^{q})
$$
and of the Poincaré-Nash inequality. We obtain (\ref{eq:generalisation-NP-1}) and 
(\ref{eq:generalisation-NP-2}) for odd integers using the Schwartz inequality. \\

We now evaluate the expectation of normalized traces of matrices such as 
$\Pi_{s=1}^{r} {\bf Q} {\bf C}^{s}$. Proposition \ref{coro:expre-traces-Q-Q-Q-W} is used in the sequel 
in the case $r=2$ and $r=3$. 
\begin{proposition}
\label{coro:expre-traces-Q-Q-Q-W}
For each $ML \times ML$ deterministic uniformly bounded matrices $({\bf C}^{s})_{s=1, \ldots, r+1}$
and ${\bf A}$, it holds that 
\begin{align}
\label{eq:expre-trace-Q-Q}
\mathbb{E} \left( \frac{1}{ML} \mathrm{Tr}\left( \Pi_{s=1}^{r+1} {\bf Q} {\bf C}^{s} \right) \right)
= \mathbb{E} \left( \frac{1}{ML} \mathrm{Tr}\left[ (\Pi_{s=1}^{r} {\bf Q} {\bf C}^{s})   ({\bf I}_M \otimes {\bf R}) {\bf C}^{r+1}) \right] \right) + \mathcal{O}(\frac{L}{MN}) \\
+ \sigma^{2} c_N \sum_{s=1}^{r} \sum_{i=-(L-1)}^{L-1} \mathbb{E} \left[ \frac{1}{ML} 
\mathrm{Tr}  \left( (\Pi_{t=s}^{r} {\bf Q} {\bf C}^{s}) {\bf Q} ({\bf I}_M \otimes {\bf J}_L^{i}) \right) \right] \; 
 \mathbb{E} \left[ \frac{1}{ML} \mathrm{Tr} \left( (\Pi_{t=1}^{s-1} {\bf Q} {\bf C}^{s}) 
{\bf Q} {\bf W} {\bf J}_N^{i} {\bf H}^{T} 
{\bf W}^* ({\bf I}_M \otimes {\bf R}) {\bf C}^{r+1} \right) \right] 
\nonumber
\end{align}
and that
\begin{multline}
\label{eq:expre-trace-Q-QW}
\mathbb{E} \left[ \frac{1}{ML} \mathrm{Tr}\left( (\Pi_{s=1}^{r} {\bf Q} {\bf C}^{s}) {\bf Q} 
{\bf W} {\bf G} {\bf W}^{*} {\bf A}  \right) \right]  
= \mathbb{E} \left( \frac{1}{ML} \mathrm{Tr}\left[ \Pi_{s=1}^{r} {\bf Q} {\bf C}^{s} ({\bf I}_M \otimes \sigma^{2} {\bf R} \mathcal{T}_{L,L}({\bf G}^{T} {\bf H})) {\bf A} ) \right] \right)  + 
\mathcal{O}(\frac{L}{MN}) + \\
\sigma^{2} c_N \sum_{s=1}^{r} \sum_{i=-(L-1)}^{L-1} \mathbb{E} \left[ \frac{1}{ML} 
\mathrm{Tr}  \left( (\Pi_{t=s}^{r} {\bf Q} {\bf C}^{s}) {\bf Q} ({\bf I}_M \otimes {\bf J}_L^{i}) \right) \right] 
 \mathbb{E} \left[ \frac{1}{ML} \mathrm{Tr}  (\Pi_{t=1}^{s-1} {\bf Q} {\bf C}^{s}) 
{\bf Q} {\bf W} {\bf J}_N^{i} {\bf H}^{T} {\bf W}^{*} ({\bf I}_M \otimes \sigma^{2} {\bf R} \mathcal{T}_{L,L}({\bf G}^{T} {\bf H})) {\bf A}  \right] \\
- \sigma^{2} c_N \sum_{s=1}^{r} \sum_{i=-(L-1)}^{L-1} \mathbb{E} \left[ \frac{1}{ML} 
\mathrm{Tr}  \left( (\Pi_{t=s}^{r} {\bf Q} {\bf C}^{s}) ({\bf Q} ({\bf I}_M \otimes {\bf J}_L^{i}) \right) \right] \; 
 \mathbb{E} \left[ \frac{1}{ML} \mathrm{Tr} \left( (\Pi_{t=1}^{s-1} {\bf Q} {\bf C}^{s}) 
{\bf Q} {\bf W} {\bf J}_N^{i} {\bf H}^{T} {\bf G} {\bf W}^{*} {\bf A} \right) \right] 
\end{multline}
\end{proposition}
The proof of this result is similar to the proof of (\ref{eq:expre-E(Q)}) and (\ref{eq:expre-E(QWGW*)}), 
but is of course more tedious. To establish (\ref{eq:expre-trace-Q-Q}) and (\ref{eq:expre-trace-Q-QW}), 
it is sufficient to evaluate matrix $\mathbb{E} \left[ \Pi_{s=1}^{r} {\bf Q}_{l_s, l_s^{'}}^{n_s, n_s^{'}} \;  {\bf Q} {\bf W} {\bf G} {\bf W}^{*} \right]$
using the integration by parts formula for each multi-indices $(l_1^{'}, \ldots, l_r^{'})$ and 
$(n_1^{'}, \ldots, n_r^{'})$. A proof is provided in \cite{loubaton-arxiv}. \\

We now use Proposition \ref{coro:expre-traces-Q-Q-Q-W} to study the behaviour of certain useful terms. For this, 
it is first necessary to give the following lemma. If ${\bf A}$ is a matrix, $||| {\bf A} |||_{\infty}$ 
is defined as 
$$
||| {\bf A} |||_{\infty} = \sup_{i} \sum_{j} |{\bf A}_{i,j}|
$$
\begin{lemma}
\label{le:inversion-D}
We consider the $(2L-1) \times (2L-1)$ diagonal matrix ${\bf D}(z) = \mathrm{Diag}(d(-(L-1),z), \ldots, d(0), \ldots, 
d(L-1,z)$ where for each $l \in \mathbb{Z}$, $d(l,z)$ is defined as
\begin{equation}
\label{eq:def-dl}
d(l,z) = \sigma^{4} c_N \,  (z \, t(z) \, \tilde{t}(z))^{2} \, (1 - |l|/L)_{+} \, (1 - |l|/N)_{+} 
\end{equation}
We consider a $(2L-1) \times (2L-1)$ deterministic matrix $\bs{\Upsilon}$ whose entries 
$(\epsilon_{k,l})_{-(L-1) \leq k,l \leq L-1}$ depend on $z,L,M,N$ and satisfy 
\begin{equation}
\label{eq:hypothese-perturbation}
|\epsilon_{k,l}| \leq \frac{L}{MN} P_1(|z|) P_2(\frac{1}{\mathrm{Im}(z)})
\end{equation}
for some nice polynomials $P_1$ and $P_2$ for each $z \in \mathbb{C}^{+}$. Then, for each $z$ belonging to a set $E_N$ defined by 
\begin{equation}
\label{eq:justif-EN-nonvoid}
E_N = \{ z \in \mathbb{C}^{+}, \frac{L^{2}}{MN} Q_1(|z|) Q_2(\frac{1}{\mathrm{Im}(z)}) < 1 \}
\end{equation}
for some nice polynomials $Q_1, Q_2$, matrix $\left({\bf I} - ({\bf D} + \bs{\Upsilon}) \right)$ is invertible 
and for each $L,M,N$, and for each $z \in E_N$, it holds that  
\begin{equation}
\label{eq:norme-infinie-inverse}
\sup_{L,M,N} ||| \left({\bf I} - ({\bf D} + \bs{\Upsilon}) \right)^{-1} |||_{\infty} < C \frac{(\eta^{2} + |z|^{2})^{2}}{(\mathrm{Im}(z))^{4}} 
\end{equation}
for some nice constants $\eta$ and $C$. 
\end{lemma}
{\bf Proof.}
It is well known (see e.g. \cite{horn-johnson}, Corollary 6.1.6 p. 390) that 
$$
\rho({\bf D} + \bs{\Upsilon}) \leq ||| {\bf D} + \bs{\Upsilon} |||_{\infty}
$$
Therefore, we obtain that 
$$
\rho({\bf D} + \bs{\Upsilon}) \leq \sigma^{4} c_N  |z \, t(z) \, \tilde{t}(z)|^{2} + 
\frac{L^{2}}{MN} P_1(|z|) P_2(\frac{1}{\mathrm{Im}(z)})
$$
As $\sigma^{4} c_N  |z \, t(z) \, \tilde{t}(z)|^{2} \leq 1 - C \frac{(\mathrm{Im}(z))^{4}}{(\eta^{2} + |z|^{2})^{2}}$ for some nice constants $C$ and $\eta$ (see Eq. \ref{eq:borneinf-zttildet})), we get that
$$
\rho({\bf D} + \bs{\Upsilon}) < 1 - \frac{C}{2} \frac{(\mathrm{Im}(z))^{4}}{(\eta^{2} + |z|^{2})^{2}}
$$
if $z$ satisfies 
$$
 C \frac{(\mathrm{Im}(z))^{4}}{(\eta^{2} + |z|^{2})^{2}} - \frac{L^{2}}{MN} P_1(|z|) P_2(\frac{1}{\mathrm{Im}(z)}) > \frac{C}{2} \frac{(\mathrm{Im}(z))^{4}}{(\eta^{2} + |z|^{2})^{2}}
$$
a condition that can be written as $z \in E_N$ for well chosen nice polynomials $Q_1, Q_2$. We note that 
a similar result holds for $\rho(|{\bf D}| + |\bs{\Upsilon}|)$ where for any matrix ${\bf A}$, 
$|{\bf A}|$ is the matrix defined by $(|{\bf A}|)_{i,j} = |{\bf A}|_{i,j}$. 
This implies that for $z \in E_N$, matrices ${\bf I} -{\bf D} - \bs{\Upsilon}$ and 
${\bf I} -|{\bf D}| - |\bs{\Upsilon}|$ are invertible, and that 
$({\bf I} -{\bf D} - \bs{\Upsilon})^{-1} = \sum_{n=0}^{+\infty} ({\bf D} + \bs{\Upsilon})^{n}$
and $({\bf I} -|{\bf D}| - |\bs{\Upsilon}|)^{-1} = \sum_{n=0}^{+\infty} (|{\bf D}| + |\bs{\Upsilon}|)^{n}$. We note that for each $k,l$, $|\left(({\bf D} + \bs{\Upsilon})^{n}\right)_{k,l}| \leq 
\left( (|{\bf D}| + | \bs{\Upsilon}|)^{n} \right)_{k,l}$. Therefore, 
\begin{equation}
\label{eq:inegalite-inverse}
\left| \left(({\bf I} -{\bf D} - \bs{\Upsilon})^{-1}\right)_{k,l} \right| \leq 
\left(({\bf I} -|{\bf D}| - |\bs{\Upsilon}|)^{-1}\right)_{k,l} 
\end{equation}
We denote by ${\bf 1}$ the 
$2L-1$ dimensional vector with all components equal to 1, and by ${\bf b}$ the vector 
${\bf b} = \left( {\bf I} -|{\bf D}| - |\bs{\Upsilon}| \right) \, {\bf 1}$. 
It is clear that for each $l \in \{ -(L-1), \ldots, L-1 \}$, ${\bf b}_l$ is equal to 
$$
{\bf b}_l = 1 - \sigma^{4} c_N  |z \, t(z) \, \tilde{t}(z)|^{2} \, (1 - |l|/L)(1 - |l|/N) \, - \, 
\sum_{k} |\epsilon_{l,k}|
$$
which is greater than $\frac{C}{2} \frac{(\mathrm{Im}(z))^{4}}{(\eta^{2} + |z|^{2})^{2}}$ if $z \in E_N$. Therefore, 
for each $l$, for $z \in E_N$, it holds that 
$$
1 = \sum_{k} \left( {\bf I} -|{\bf D}| - |\bs{\Upsilon}| \right)^{-1}_{l,k} \, {\bf b}_k 
> \frac{C}{2} \frac{(\mathrm{Im}(z))^{4}}{(\eta^{2} + |z|^{2})^{2}} \sum_{k} \left( {\bf I} -|{\bf D}| - |\bs{\Upsilon}| \right)^{-1}_{l,k}
$$
which implies that 
$$
 ||| \left({\bf I} - (|{\bf D}| + |\bs{\Upsilon})| \right)^{-1} |||_{\infty} < \frac{2}{C} \, \frac{(\eta^{2} + |z|^{2})^{2}}{(\mathrm{Im}(z))^{4}} 
$$
(\ref{eq:norme-infinie-inverse}) follows immediately from (\ref{eq:inegalite-inverse}). \\

We now introduce $\omega(u_1,u_2,z)$ defined for $-(L-1) \leq u_i \leq (L-1)$ for $i=1,2$
by 
\begin{equation}
\label{eq:def-omega-2}
\omega(u_1,u_2,z) =  \frac{1}{ML} \mathrm{Tr} \left( {\bf Q} ({\bf I}_M \otimes {\bf J}_L^{u_1})
{\bf Q} ({\bf I}_M \otimes {\bf J}_L^{u_2}) \right) 
\end{equation}
and prove the following result. 
\begin{proposition}
\label{prop:expre-omega-2}
$\mathbb{E}(\omega(u_1,u_2,z))$ can be expressed as 
\begin{equation}
\label{eq:expre-omega-2}
\mathbb{E} \left(\omega(u_1,u_2,z) \right) = \delta(u_1 + u_2=0) \, \overline{\omega}(u_1,z) 
+ \mathcal{O}(\frac{L}{MN}) 
\end{equation}
for each $z \in E_N$ where $E_N$ is defined by (\ref{eq:justif-EN-nonvoid}) and where $\overline{\omega}(u_1,z)$
is defined by 
$$
\overline{\omega}(u_1,z) =  \frac{(1- |u_1|/L) \; t^{2}(z)}{1 - \sigma^{4} c_N  (z \, t(z) \, \tilde{t}(z))^{2}(1- |u_1|/L)(1- |u_1|/N)}
$$
\end{proposition}
{\bf Proof.} We use (\ref{eq:expre-trace-Q-Q}) for $r=1, {\bf C}^{1} = ({\bf I}_M \otimes {\bf J}_L^{u_1}), 
{\bf C}^{2} = ({\bf I}_M \otimes {\bf J}_L^{u_2})$. Using that 
$$
\mathbb{E} \left( \frac{1}{ML} \mathrm{Tr}({\bf Q} {\bf C}^{1} ({\bf I}_M \otimes {\bf R}) {\bf C}^{2}) \right)= 
\frac{1}{ML} \mathrm{Tr}\left( ({\bf I}_M \otimes {\bf R}) {\bf C}^{1} ({\bf I}_M \otimes {\bf R}) {\bf C}^{2}) \right) + \mathcal{O}(\frac{L}{MN})
$$
we obtain that 
\begin{equation}
\label{eq:systeme-omega-2}
\mathbb{E}(\omega(u_1,u_2)) = \frac{1}{L} \mathrm{Tr}\left({\bf R} {\bf J}_L^{u_1} {\bf R} {\bf J}_L^{u_2}\right) \; + \; \sigma^{2} c_N \sum_{i=-(L-1)}^{L-1} \mathbb{E} \left( \frac{1}{ML} \mathrm{Tr}({\bf Q} {\bf W} {\bf J}_N^{i} {\bf H}^{T} 
{\bf W}^{*} ({\bf I}_M \otimes {\bf R} {\bf J}_L^{u_2}) \right) \, \mathbb{E}(\omega(u_1, i)) + \mathcal{O}(\frac{L}{MN}) 
\end{equation}
For each $u_1$ fixed, this equation can be interpreted as a linear system whose 
unknowns are the \\ $\left(\mathbb{E}(\omega(u_1, u_2))\right)_{u_2 = -(L-1), \ldots, L-1}$. (\ref{eq:expre-E(QWGW*)}) implies that
$$
\mathbb{E} \left( \frac{1}{ML} \mathrm{Tr}({\bf Q} {\bf W} {\bf J}_N^{i} {\bf H}^{T} 
{\bf W}^{*} ({\bf I}_M \otimes {\bf R} {\bf J}_L^{u_2}) \right) = \frac{\sigma^{2}}{L} \mathrm{Tr} {\bf R} \mathcal{T}_{L,L}({\bf H} {\bf J}_L^{*i} {\bf H}) {\bf R} {\bf J}_L^{u_2} +  \mathcal{O}(\frac{L}{MN})
$$ 
Moreover, we check that, up to a $\mathcal{O}(\frac{L}{MN})$ term, matrices ${\bf R}$ and ${\bf H}$ can be replaced into the righthandside of the above equation 
by $t(z) {\bf I}_L$ and $- z \tilde{t}(z) {\bf I}_L$ respectively. In other words, 
\begin{eqnarray*}
\mathbb{E} \left( \frac{1}{ML} \mathrm{Tr}({\bf Q} {\bf W} {\bf J}_N^{i} {\bf H}^{T} 
{\bf W}^{*} ({\bf I}_M \otimes {\bf R} {\bf J}_L^{u_2}) \right) & = & \sigma^{2}  (z t(z) \, \tilde{t}(z))^{2} \frac{1}{L} \mathrm{Tr} \left( \mathcal{T}_{L,L}({\bf J}_L^{*i}) {\bf J}_L^{u_2} \right) + \mathcal{O}(\frac{L}{MN})  \\ 
         & = & \delta(i-u_2) \, \sigma^{2}  (z t(z) \, \tilde{t}(z))^{2}  (1- |u_2|/L)(1- |u_2|/N) + \mathcal{O}(\frac{L}{MN})
\end{eqnarray*}
We write ${\bf R} \mathcal{T}_{L,L}\left( {\bf H} {\bf J}_N^{*i} {\bf H} \right) {\bf R} {\bf J}_L^{u_2}$ as
\begin{align*}
{\bf R} \mathcal{T}_{L,L}\left( {\bf H} {\bf J}_N^{*i} {\bf H} \right) {\bf R} {\bf J}_L^{u_2} = 
({\bf R} - t {\bf I}) \mathcal{T}_{L,L}\left( {\bf H} {\bf J}_N^{*i} {\bf H}  \right) {\bf R} {\bf J}_L^{u_2}+ \\
t \, \mathcal{T}_{L,L}\left( ({\bf H} + z \tilde{t} {\bf I}) {\bf J}_N^{*i}  {\bf H} \right) {\bf R} {\bf J}_L^{u_2}
- z t \tilde{t}  \mathcal{T}_{L,L}\left(  {\bf J}_N^{*i} ({\bf H} + z \tilde{t} {\bf I})  \right) {\bf R} {\bf J}_L^{u_2} +\\
 t (z \tilde{t})^{2}  \mathcal{T}_{L,L}\left(  {\bf J}_N^{*i}  \right) ({\bf R} - t \, {\bf I}) {\bf J}_L^{u_2} 
+t^{2}   (z \tilde{t})^{3}   \mathcal{T}_{L,L}\left(  {\bf J}_N^{*u} \right) {\bf J}_L^{u_2}
\end{align*}
The terms 
$ \frac{1}{L} \mathrm{Tr} \left( ({\bf R} - t \, {\bf I}) \mathcal{T}_{L,L}\left( {\bf H} {\bf J}_N^{*i} {\bf H} \right) {\bf R} {\bf J}_L^{u_2}\right)$  and 
$\frac{1}{L} \mathrm{Tr} \left( \mathcal{T}_{L,L}\left(  {\bf J}_N^{*i}  \right) ({\bf R} - t \, {\bf I}) {\bf J}_L^{u_2}  \right)$ are $\mathcal{O}(\frac{L}{MN})$ by Proposition \ref{prop:Trace-R-t}. We just study the term 
$\frac{1}{L} \mathrm{Tr} \left( t \,   \mathcal{T}_{L,L}\left( ({\bf H} + z \tilde{t} {\bf I}) {\bf J}_N^{*i} {\bf H}  \right) {\bf R}  {\bf J}_L^{u_2} \right)$
and omit $\frac{1}{L} \mathrm{Tr} \left( \mathcal{T}_{L,L}\left(  {\bf J}_N^{*i} ({\bf H} + z \tilde{t} {\bf I})  \right) {\bf R} {\bf J}_L^{u_2} \right)$  because it can be handled similarly. 
We express ${\bf H} + z \tilde{t} {\bf I}$ as
\begin{eqnarray*}
{\bf H} + z \tilde{t} {\bf I} & = & \sigma^{2} c_N \, z \tilde{t} \, {\bf H} \mathcal{T}_{N,L}^{(M)}\left( \mathbb{E}({\bf Q}) - t \, {\bf I} \right) \\
& = & \sigma^{2} c_N \, z \tilde{t} \, {\bf H} \mathcal{T}_{N,L}^{(M)}  \left( \mathbb{E}({\bf Q}) - {\bf I}_M \otimes {\bf R} \right) + 
 \sigma^{2} c_N \, z \tilde{t} \, {\bf H} \mathcal{T}_{N,L}\left( {\bf R} - t \, {\bf I} \right)
\end{eqnarray*}
Property (\ref{eq:utile}) and Proposition \ref{prop:Trace-R-t} imply that $\frac{1}{L} \mathrm{Tr} \left( t \,   \mathcal{T}_{L,L}\left( ({\bf H} + z \tilde{t} {\bf I}) {\bf J}_N^{*i} {\bf H}  \right) {\bf R} {\bf J}_L^{u_2} \right)$ is a $\mathcal{O}(\frac{L}{MN})$. 
We have thus shown that for $i,u_2 \in -(L-1), \ldots, L-1$, then, it holds that
\begin{equation}
\label{eq:D-perturbee}
\sigma^{2} c_N \, \mathbb{E} \left( \frac{1}{ML} \mathrm{Tr}({\bf Q} {\bf W} {\bf J}_N^{i} {\bf H}^{T} 
{\bf W}^{*} ({\bf I}_M \otimes {\bf R} {\bf J}_L^{u_2}) \right) = \delta(i+u_2=0) \, d(i,z) + \mathcal{O}(\frac{L}{MN})
\end{equation}
Similarly, it holds that 
\begin{eqnarray*}
\frac{1}{L} \mathrm{Tr}\left({\bf R} {\bf J}_L^{u_1} {\bf R} {\bf J}_L^{u_2}\right) & = & t(z)^{2} \, \frac{1}{L} \mathrm{Tr}({\bf J}_L^{u_1}{\bf J}_L^{u_2}) + \mathcal{O}(\frac{L}{MN}) \\
    &  =  & \delta(u_1 + u_2=0) \, (t(z))^{2} \, (1- |u_1|/L) + \mathcal{O}(\frac{L}{MN})
\end{eqnarray*}

We denote by $\bs{\omega}(u_1)$ the $(2L-1)$ dimension vector $(\omega(u_1,u_2))_{u_2=-(L-1), \ldots, L-1}$, and by $\overline{\bs{\gamma}}(u_1)$ 
the vector such that 
$$
\overline{\bs{\gamma}}(u_1)_{u_2} = \delta(u_1 + u_2=0) \, (t(z))^{2} \, (1- |u_1|/L)
$$
The linear system (\ref{eq:systeme-omega-2}) can be written as
$$
\mathbb{E}(\bs{\omega}(u_1)) = ({\bf D} + \bs{\Upsilon}) \, \mathbb{E}(\bs{\omega}(u_1)) + \overline{\bs{\gamma}}(u_1) + \bs{\epsilon}
$$
where the elements of matrix $\bs{\Upsilon}$ and the components of vector $\bs{\epsilon}$ are $\mathcal{O}(\frac{L}{MN})$ terms. Matrices 
${\bf D}$ and $\bs{\Upsilon}$ verify the assumptions of Lemma \ref{le:inversion-D}. Therefore, it holds that 
$$
\mathbb{E}(\bs{\omega}(u_1)) = \left( {\bf I} - {\bf D} - \bs{\Upsilon} \right)^{-1} \, (\overline{\bs{\gamma}}(u_1) + \bs{\epsilon})
$$
when $z$ belongs to a set $E_N$ defined as in (\ref{eq:justif-EN-nonvoid}). 
Writing matrix $\left( {\bf I} - {\bf D} - \bs{\Upsilon} \right)^{-1}$ as
$$
\left( {\bf I} - {\bf D} - \bs{\Upsilon} \right)^{-1} = \left( {\bf I} - {\bf D} \right)^{-1} + \left( {\bf I} - {\bf D} - \bs{\Upsilon} \right)^{-1} \, \bs{\Upsilon} \, \left( {\bf I} - {\bf D} \right)^{-1}
$$
we obtain that 
$$
\mathbb{E}(\bs{\omega}(u_1)) = \left( {\bf I} - {\bf D} \right)^{-1} \, \overline{\bs{\gamma}}(u_1) + \left( {\bf I} - {\bf D} - \bs{\Upsilon} \right)^{-1} \, \bs{\Upsilon} \, \left( {\bf I} - {\bf D} \right)^{-1} \, \overline{\bs{\gamma}}(u_1) \; + \; \left( {\bf I} - {\bf D} - \bs{\Upsilon} \right)^{-1} \, \bs{\epsilon}
$$
(\ref{eq:norme-infinie-inverse}) implies that for each $u_2$, 
$$
\left( ( {\bf I} - {\bf D} - \bs{\Upsilon} )^{-1} \, \bs{\epsilon} \right)_{u_2} = \mathcal{O}(\frac{L}{MN})
$$
Moreover, as vector $\overline{\bs{\gamma}}(u_1)$ has only 1 non zero component, it is clear that each component of vector 
$\bs{\Upsilon} \, \left( {\bf I} - {\bf D} \right)^{-1} \, \overline{\bs{\gamma}}(u_1)$ is a $\mathcal{O}(\frac{L}{MN})$ term. Hence, 
(\ref{eq:norme-infinie-inverse}) leads to 
$$
\left( ({\bf I} - {\bf D} - \bs{\Upsilon})^{-1} \, \bs{\Upsilon} \, ( {\bf I} - {\bf D} )^{-1} \, \overline{\bs{\gamma}}(u_1) \right)_{u_2} = \mathcal{O}(\frac{L}{MN})
$$
This establishes (\ref{eq:expre-omega-2}). We notice that Lemma \ref{le:inversion-D} plays an important role in the above calculations. The control 
of $||| \left({\bf I} - (|{\bf D}| + |\bs{\Upsilon})| \right)^{-1} |||_{\infty}$ allows in particular to show that 
$\mathbb{E}(\omega(u_1,u_2)) = \mathcal{O}(\frac{L}{MN})$ if $u_1+u_2 \neq 0$, instead of $\mathcal{O}(\frac{L^{2}}{MN})$
in the absence of control on $||| \left({\bf I} - (|{\bf D}| + |\bs{\Upsilon})| \right)^{-1} |||_{\infty}$. As Lemma 
\ref{le:inversion-D} is a consequence of $\frac{L^{2}}{MN} \rightarrow 0$, this discussion confirms the importance of condition (\ref{eq:stronger-condition}), and strongly suggests that it is a necessary condition to obtain positive results. \\

It is also necessary to evaluate $\mathbb{E}(\omega(u_1,u_2,u_3,z))$ where $\omega(u_1,u_2,u_3,z)$ is defined by
\begin{equation}
\label{eq:def-omega-3}
\omega(u_1,u_2,u_3,z) = \mathbb{E} \left[ \frac{1}{ML} \mathrm{Tr} \left( {\bf Q} ({\bf I}_M \otimes {\bf J}_L^{u_1})
{\bf Q} ({\bf I}_M \otimes {\bf J}_L^{u_2}) {\bf Q} ({\bf I}_M \otimes {\bf J}_L^{u_3})\right) \right]
\end{equation}
It holds that for $z \in E_N$ defined as in (\ref{eq:justif-EN-nonvoid})
\begin{proposition}
\label{prop:expre-omega-3}
$\mathbb{E}(\omega(u_1,u_2, u_3,z))$ can be expressed as 
\begin{equation}
\label{eq:expre-omega-3}
\mathbb{E}(\omega(u_1,u_2, u_3,z)) = \delta(u_1 + u_2 + u_3=0) \, \overline{\omega}(u_1,u_2,z) + \mathcal{O}(\frac{L}{MN}) 
\end{equation}
where $\overline{\omega}(u_1,u_2,z)$ is given by
\begin{equation}
\label{eq:expre-C(u_1,u_2,z)}
(t(z))^{3} \; \frac{\frac{1}{L} \mathrm{Tr}({\bf J}_L^{u_2} {\bf J}_L^{u_1} {\bf J}_L^{*(u_1+u_2)}) + \sigma^{6} c_N^{2} 
(z t(z) \, \tilde{t}(z))^{3} \, (1- |u_1|/L) (1- |u_2|/L) (1- |u_1+u_2|/L)_{+} \,  \frac{1}{N} \mathrm{Tr}({\bf J}_N^{u_1} {\bf J}_N^{u_2} {\bf J}_N^{*(u_1+u_2)})}{(1 - d(u_1,z)) \, (1 - d(u_2,z)) \, (1 - d(u_1+u_2,z))}
\end{equation}
\end{proposition}
{\bf Proof.} The proof is somewhat similar to the proof of Proposition \ref{prop:expre-omega-2}, but it 
needs rather tedious calculations. We just provide the main steps and omit the straightforward details. 
 We use again (\ref{eq:expre-trace-Q-Q}), but for $r=2$, and ${\bf C}^{s} = ({\bf I}_M \otimes {\bf J}_L^{u_s})$ 
for $s=1, 2, 3$. We obtain immediately that 
\begin{align}
\label{eq:systeme-omega-3}
\mathbb{E}(\omega(u_1,u_2,u_3)) = \frac{1}{ML} \mathbb{E} \left[ \mathrm{Tr}\left( {\bf Q} ({\bf I}_M \otimes {\bf J}_L^{u_1})
{\bf Q} ({\bf I}_M \otimes {\bf J}_L^{u_2} {\bf R} {\bf J}_L^{u_3})  \right) \right] + \\
\nonumber
\sigma^{2} c_N \sum_{i=-(L-1)}^{L-1} \mathbb{E} \left[ \frac{1}{ML} \left( \mathrm{Tr}({\bf Q} {\bf W} {\bf J}_N^{i} {\bf H}^{T} 
{\bf W}^{*} ({\bf I}_M \otimes {\bf R} {\bf J}_L^{u_3}) \right) \right] \, \mathbb{E}(\omega(u_1, u_2, i)) + \\
\nonumber
\sigma^{2} c_N \sum_{i=-(L-1)}^{L-1} \frac{1}{ML} \mathbb{E} \left[ \mathrm{Tr}\left( {\bf Q} ({\bf I}_M \otimes {\bf J}_L^{u_1}) {\bf Q} {\bf W} {\bf J}_N^{i} {\bf H}^{T} {\bf W}^{*} ({\bf I}_M \otimes  {\bf R} {\bf J}_L^{u_3}) \right) \right] \; \mathbb{E}(\omega(u_2,i)) + \mathcal{O}(\frac{L}{MN})
\end{align}
(\ref{eq:systeme-omega-3}) can still be interpreted as a linear system whose unknown are the 
$\left(\mathbb{E}(\omega(u_1,u_2,u_3))\right)_{u_3 \in \{ -(L-1), \ldots, L-1 \}}$. The matrix governing the system
is the same matrix ${\bf D} + \bs{\Upsilon}$ as in the proof of Proposition \ref{prop:expre-omega-2} (but for a different matrix $\bs{\Upsilon}$). In order 
to use the same arguments, it is sufficient to establish that 
\begin{equation}
\label{eq:verif-1}
\frac{1}{ML} \mathbb{E} \left[ \mathrm{Tr}\left( {\bf Q} ({\bf I}_M \otimes {\bf J}_L^{u_1})
{\bf Q} ({\bf I}_M \otimes {\bf J}_L^{u_2} {\bf R} {\bf J}_L^{u_3})  \right) \right] = C(u_1,u_2,z) \delta(u_1+u_2+u_3=0) + \mathcal{O}(\frac{L}{MN}) 
\end{equation}
and 
\begin{multline}
\label{eq:verif-2}
\sum_{i=-(L-1)}^{L-1} \frac{1}{ML} \mathbb{E} \left[ \mathrm{Tr}\left( {\bf Q} ({\bf I}_M \otimes {\bf J}_L^{u_1}) {\bf Q} {\bf W} {\bf J}_N^{i} {\bf H}^{T} {\bf W}^{*} ({\bf I}_M \otimes  {\bf R} {\bf J}_L^{u_3}) \right) \right] \; \mathbb{E}(\omega(u_2,i)) = \\ 
 C(u_1,u_2,z) \delta(u_1+u_2+u_3=0) + \mathcal{O}(\frac{L}{MN})
\end{multline}
To check (\ref{eq:verif-1}), we use (\ref{eq:expre-trace-Q-Q}) for $r=1, {\bf C}^{1} = {\bf I}_M \otimes {\bf J}_L^{u_1}, {\bf C}^{2} = {\bf I}_M \otimes {\bf J}_L^{u_2} {\bf R} {\bf J}_L^{u_3}$. This leads to 
\begin{align*}
\frac{1}{ML} \mathbb{E} \left[ \mathrm{Tr}\left( {\bf Q} ({\bf I}_M \otimes {\bf J}_L^{u_1})
{\bf Q} ({\bf I}_M \otimes {\bf J}_L^{u_2} {\bf R} {\bf J}_L^{u_3})  \right) \right] = 
\frac{1}{L} \mathrm{Tr} \left( {\bf R} {\bf J}_L^{u_1}  {\bf R} {\bf J}_L^{u_2}   {\bf R} {\bf J}_L^{u_3} \right) + \\
 \sigma^{2} c_N \sum_{i=-(L-1)}^{L-1} \mathbb{E}(\omega(u_1,i)) \,  \mathbb{E} \left[ \frac{1}{ML} \left( \mathrm{Tr}({\bf Q} {\bf W} {\bf J}_N^{i} {\bf H}^{T} 
{\bf W}^{*} ({\bf I}_M \otimes {\bf R}^{2} {\bf J}_L^{u_2} {\bf R}  {\bf J}_L^{u_3}) \right) \right]
+ \mathcal{O}(\frac{L}{MN})
\end{align*}
Up to a $\mathcal{O}(\frac{L}{MN})$ term, it is possible to replace ${\bf R}(z)$ by $t(z) {\bf I}$ into the
first term of the righthandside of the above equation. This leads to 
\begin{eqnarray*}
\frac{1}{L} \mathrm{Tr} \left( {\bf R} {\bf J}_L^{u_1}  {\bf R} {\bf J}_L^{u_2}   {\bf R} {\bf J}_L^{u_3} \right) & = &
(t(z))^{3} \, \frac{1}{L} \mathrm{Tr} {\bf J}_L^{u_1}  {\bf J}_L^{u_2}  {\bf J}_L^{u_3} +  \mathcal{O}(\frac{L}{MN}) \\
 & = &  (t(z))^{3} \, \frac{1}{L} \mathrm{Tr} {\bf J}_L^{u_1}  {\bf J}_L^{u_2}  {\bf J}_L^{*(u_1+u_2)} \delta(u_1+u_2+u_3=0)
+  \mathcal{O}(\frac{L}{MN}) 
\end{eqnarray*}
Similarly, it is easy to check that 
$$
\mathbb{E} \left[ \frac{1}{ML} \left( \mathrm{Tr}({\bf Q} {\bf W} {\bf J}_N^{i} {\bf H}^{T} 
{\bf W}^{*} ({\bf I}_M \otimes {\bf R}^{2} {\bf J}_L^{u_2} {\bf R}  {\bf J}_L^{u_3}) \right) \right] = 
C(u_2,u_3,z) \delta(i=u_2+u_3) + \mathcal{O}(\frac{L}{MN})
$$
As $\mathbb{E}(\omega(u_1,i,z)) = \overline{\omega}(u_1,z) \delta(i+u_1=0) + \mathcal{O}(\frac{L}{MN})$, we get immediately that if 
$u_1 + u_2 + u_3 \neq 0$, then, 
$$
 \sigma^{2} c_N \sum_{i=-(L-1)}^{L-1} \mathbb{E}(\omega(u_1,i)) \,  \mathbb{E} \left[ \frac{1}{ML} \left( \mathrm{Tr}({\bf Q} {\bf W} {\bf J}_N^{i} {\bf H}^{T} 
{\bf W}^{*} ({\bf I}_M \otimes {\bf R}^{2} {\bf J}_L^{u_2} {\bf R}  {\bf J}_L^{u_3}) \right) \right] = 
\mathcal{O}(\frac{L}{MN}) + L \, \mathcal{O}\left( (\frac{L}{MN})^{2} \right)
$$
(\ref{eq:verif-1}) follows from the observation that, as $\frac{L^{2}}{MN} \rightarrow 0$, then 
$L (\frac{L}{MN})^{2} = \frac{L^{2}}{MN} \frac{L}{MN} = o(\frac{L}{MN})$. 

Finally, (\ref{eq:verif-2}) holds because, using (\ref{eq:expre-trace-Q-QW}) for $r=1, {\bf C}^{1} = 
{\bf I}_M \otimes {\bf J}_L^{u_1}$, ${\bf G} = {\bf J}_N^{i} {\bf H}^{T}$, ${\bf A} = {\bf I}_M \otimes {\bf R} {\bf J}_L^{u_3}$, it can be shown that 
$$
\frac{1}{ML} \mathbb{E} \left[ \mathrm{Tr}\left( {\bf Q} ({\bf I}_M \otimes {\bf J}_L^{u_1}) {\bf Q} {\bf W} {\bf J}_N^{i} {\bf H}^{T} {\bf W}^{*} ({\bf I}_M \otimes  {\bf R} {\bf J}_L^{u_3}) \right) \right] = 
C(u_1,u_3,z) \, \delta(i=u_1+u_3) + \mathcal{O}(\frac{L}{MN})
$$
As $\mathbb{E}(\omega_2(i,u_2,z)) = \delta(i+u_2=0) \, \overline{\omega}(u_2,z) + \mathcal{O}(\frac{L}{MN})$, $\frac{L^{2}}{MN} \rightarrow 0$
implies (\ref{eq:verif-2}). \\

The calculation of $\overline{\omega}(u_1,u_2,z)$ is omitted.  \\

We now define and evaluate the following useful terms. If $p \geq 1$ and $q \geq 1$, for each integers $i,u_1,u_2$, $l_1, \ldots, l_p$, $k_1 , \ldots, k_q$  belonging to $\{ -(L-1), \ldots, L-1 \}$, we define 
$$\beta_{p,q}(i,u_1,l_1, \ldots, l_p,k_1, \ldots, k_q,u_2,z)$$
as 
\begin{equation}
\label{eq:def-beta}
\frac{1}{ML}\mathrm{Tr}\left({\bf Q} ({\bf I}_M \otimes {\bf J}_L^{i}) {\bf Q} ({\bf I}_M \otimes {\bf J}_L^{u_1}) {\bf Q} {\bf W} {\bf J}_N^{i} {\bf H}^{T} \Pi_{j=1}^{p} ({\bf J}_N^{l_j} {\bf H}^{T}) {\bf W}^{*} \left({\bf I}_M \otimes \Pi_{n=1}^{q} ({\bf R} \mathcal{T}_{L,L}({\bf H} {\bf J}_N^{*k_n} {\bf H})) {\bf R} {\bf J}_L^{u_2} \right) \right) 
\end{equation}
We also define $\beta_{p,0}(i,u_1,l_1, \ldots, l_p,u_2,z)$ as
\begin{equation}
\label{eq:def-betap0}
\frac{1}{ML}\mathrm{Tr}\left({\bf Q} ({\bf I}_M \otimes {\bf J}_L^{i}) {\bf Q} ({\bf I}_M \otimes {\bf J}_L^{u_1}) {\bf Q} {\bf W} {\bf J}_N^{i} {\bf H}^{T} \Pi_{j=1}^{p} ({\bf J}_N^{l_j} {\bf H}^{T}) {\bf W}^{*} \left({\bf I}_M \otimes  {\bf R} {\bf J}_L^{u_2} \right) \right) 
\end{equation}
and $\beta_{0,q}(i,u_1,k_1, \ldots, k_q,u_2,z)$ is defined similarly. We finally denote by 
$\beta(i,u_1,u_2,z)$ the term $\beta_{0,0}(i,u_1,u_2,z)$, i.e. 
\begin{equation}
\label{eq:def-beta00}
\beta(i,u_1,u_2,z) = \frac{1}{ML}  \mathrm{Tr}\left({\bf Q} ({\bf I}_M \otimes {\bf J}_L^{i}) {\bf Q} ({\bf I}_M \otimes {\bf J}_L^{u_1}) {\bf Q} {\bf W} {\bf J}_N^{i} {\bf H}^{T} {\bf W}^{*} ({\bf I}_M \otimes  {\bf R} {\bf J}_L^{u_2}) \right) 
\end{equation} 
\begin{proposition}
\label{prop:expre-terme-utile-generaux-trace-Q-Q-QW}
For $p \geq 0$ and $q \geq 0$, it holds that 
\begin{equation}
\label{eq:expre-terme-utile-generaux-trace-Q-Q-QW}
\mathbb{E} \left( \beta_{p,q}(i, u_1,l_1, \ldots, l_p,k_1, \ldots, k_q,u_2,z) \right) = \delta(u_1 + u_2 = \sum_{j} l_j + \sum_{n} k_n) \; \overline{\beta}_{p,q}(i,u_1,l_1, \ldots, l_{p}, k_1, \ldots, k_q,z)  + \mathcal{O}(\frac{L}{MN})
\end{equation}
where for each $i,u_1,l_1, \ldots, l_{p}, k_1, \ldots, k_q$,  function 
$z \rightarrow \overline{\beta}_{p,q}(i,u_1,l_1, \ldots, l_{p}, k_1, \ldots, k_q,z)$ 
is the Stieljes transform of a distribution $D$ whose support is included into 
$\mathcal{S}_N$ and such that $<D, \mathbb{1}> = 0$. Moreover, 
if $c_N > 1$, for each $i,l_1$, function $z \rightarrow \overline{\beta}_{1,0}(i,l_1,l_1,z)$ is 
analytic in a neighbourhood of $0$, while $0$ is pole of multiplicity 1 of functions 
$z \rightarrow z \overline{\beta}(i,l_1,z)$ and $z \rightarrow  \overline{\beta}_{0,1}(i,l_1,z)$
where we denote $\overline{\beta}_{0,0}(i,l_1,z)$ by $\overline{\beta}(i,l_1,z)$
in order to simplify the notations.  Finally, function $s(i,l_1,z)$ defined by 
\begin{eqnarray}
\label{eq:def-preliminaire-s}
s(i,l_1,z) & = & -\sigma^{2} \overline{\beta}_{1,0}(i,l_1,l_1,z) + \sigma^{2} \overline{\beta}_{0,1}(i,l_1,l_1,z) + \\ 
&  &  \sigma^{6} c_N \left(z t(z) \tilde{t}(z)\right)^{2} z \tilde{t}(z) \left( 1 + \sigma^{2} z t(z) \tilde{t}(z)
(1 - |l_1|/L)(1- |l_1|/N) \right) \left(\frac{1 - |l_1|/N}{1 - d(l_1,z)} \, \overline{\beta}(i,l_1,z) \right)
\nonumber
\end{eqnarray}
is the Stieltjes transform of a distribution $D$ whose support is included in $\mathcal{S}_N^{(0)}$ and verifying  $<D, \mathbb{1}> = 0$.
\end{proposition}
{\bf Proof.} In order to simplify the notations, we just establish 
the first part of the proposition when $p=q=0$, i.e. for the term $\beta(i,u_1,u_2,z) = \beta_{0,0}(i,u_1,u_2,z)$. 
Then, we check that 
\begin{equation}
\label{eq:expre-terme-utile-trace-Q-Q-QW}
\mathbb{E} \left(\beta(i, u_1, u_2,z) \right) = \delta(u_1 + u_2 = 0) \, \overline{\beta}(i,u_1,z)  + \mathcal{O}(\frac{L}{MN})
\end{equation}
where $\overline{\beta}(i,u,z)$ is given by 
$$
\overline{\beta}(i,u,z) = \sum_{j=1}^{5} \overline{\beta}_j(i,u,z)
$$ 
with
$$
\overline{\beta}_1(i,u,z) = \frac{ \sigma^{2} t(z)^{4} (z \, \tilde{t}(z))^{2} (1 -|i|/N) \, \frac{1}{L} \mathrm{Tr}({\bf J}_L^{i} {\bf J}_L^{u} {\bf J}_L^{*i} {\bf J}_L^{*u})}{1 - d(i,z)},
$$
$$
\overline{\beta}_2(i,u,z) = \sigma^{6} c_N  t(z)^{3}  (z \, \tilde{t}(z))^{4}  \overline{\omega}(i,u) (1 -|i+u|/N) \, (1 -|i|/N)  \frac{1}{L} \mathrm{Tr}({\bf J}_L^{u+i} {\bf J}_L^{*i} {\bf J}_L^{*u}), 
$$ 
$$
\overline{\beta}_3(i,u,z) = \sigma^{4} c_N  t(z)^{2}  (z \, \tilde{t}(z))^{3}  \overline{\omega}(i,u) \mathbb{1}_{|i+u| \leq L-1} \, (1 -|u_1|/L)  \frac{1}{N} \mathrm{Tr}({\bf J}_N^{u+i} {\bf J}_N^{*u} {\bf J}_N^{*i}),
$$
\begin{align*}
\overline{\beta}_4(i,u,z) = \sigma^{6} c_N  t(z)^{4} (z \, \tilde{t}(z))^{4} \overline{\omega}(u) (1 -|u|/N) (1 -|i|/N) \frac{1}{L} \mathrm{Tr}({\bf J}_L^{i} {\bf J}_L^{u} {\bf J}_L^{*i} {\bf J}_L^{*u}) +  \\ 
\sigma^{10} c_N^{2}  t(z)^{4}  (z \, \tilde{t}(z))^{6} \overline{\omega}(u)  \overline{\omega}(i) (1 -|i|/N)^{2}  (1 -|u|/N)  \frac{1}{L} \mathrm{Tr}({\bf J}_L^{i}  {\bf J}_L^{u} {\bf J}_L^{*i} {\bf J}_L^{*u}) - \\
\sigma^{8} c_N^{2} t(z)^{3} (z \, \tilde{t}(z))^{5} \overline{\omega}(u) \overline{\omega}(i) (1 -|i|/N) \frac{1}{N} \mathrm{Tr}({\bf J}_N^{u} {\bf J}_N^{i} {\bf J}_N^{*(i+u)}) \, \frac{1}{L} \mathrm{Tr}({\bf J}_L^{u+i} {\bf J}_L^{*i} {\bf J}_L^{*u}),
\end{align*}
\begin{align*}
\overline{\beta}_5(i,u,z) = \sigma^{4} c_N  t(z)^{3} (z \, \tilde{t}(z))^{3} \overline{\omega}(u)  \frac{1}{N} \mathrm{Tr}({\bf J}_N^{*i} {\bf J}_N^{u} {\bf J}_N^{i-u}) \, \frac{1}{L} \mathrm{Tr}({\bf J}_L^{u} {\bf J}_L^{u-i} {\bf J}_L^{*u} )  +  \\
\sigma^{8} c_N^{2} t(z)^{3} (z \, \tilde{t}(z))^{5}  \overline{\omega}(u)  \overline{\omega}(i)  (1 -|i|/N)  \frac{1}{N} \mathrm{Tr}({\bf J}_N^{*i} {\bf J}_N^{u} {\bf J}_N^{i-u}) \, \frac{1}{L} \mathrm{Tr}({\bf J}_L^{u} {\bf J}_L^{u-i} {\bf J}_L^{*u} )  +  \\
\sigma^{6} c_N^{2} t(z)^{2} (z \, \tilde{t}(z))^{4} \overline{\omega}(i) \overline{\omega}(u) (1 -|u|/L) \frac{1}{N} \mathrm{Tr}({\bf J}_N^{*i} {\bf J}_N^{u} {\bf J}_N^{i} {\bf J}_N^{*u}  )
\end{align*}
The proof is based on (\ref{eq:expre-trace-Q-QW}) for $r=2$, with ${\bf C}^{1} = {\bf I}_M \otimes {\bf J}_L^{i}$, 
${\bf C}^{2} = {\bf I}_M \otimes {\bf J}_L^{u_1}$, ${\bf G} = {\bf J}_N^{i} {\bf H}^{T}$, ${\bf A} = {\bf I}_M \otimes {\bf R} {\bf J}_L^{u_2}$.
It holds that 
\begin{align*}
\mathbb{E}(\beta(i,u_1,u_2)) = \frac{1}{ML} \mathbb{E} \left[ \mathrm{Tr}\left({\bf Q} ({\bf I}_M \otimes {\bf J}_L^{i}) {\bf Q} ({\bf I}_M \otimes {\bf J}_L^{u_1}) ({\bf I}_M \otimes \sigma^{2} {\bf R} \mathcal{T}_{L,L}({\bf H} {\bf J}_N^{*i} {\bf H}) {\bf R} {\bf J}_L^{u_2}) \right) \right] +  \\
\sigma^{2} c_N \sum_{j=-(L-1)}^{L-1} \mathbb{E}(\omega(i,u_1,j)) \, \frac{1}{ML} \mathbb{E} \left[ \mathrm{Tr}\left( {\bf Q} {\bf W} {\bf J}_N^{j} {\bf H}^{T} {\bf W}^{*} ({\bf I}_M \otimes  \sigma^{2} {\bf R} \mathcal{T}_{L,L}({\bf H} {\bf J}_N^{*i} {\bf H}) {\bf R} {\bf J}_L^{u_2}) \right) \right] +\\
\sigma^{2} c_N \sum_{j=-(L-1)}^{L-1} \mathbb{E}(\omega(u_1,j)) \, \frac{1}{ML} \mathbb{E} \left[ \mathrm{Tr}\left({\bf Q} ({\bf I}_M \otimes {\bf J}_L^{i})  {\bf Q} {\bf W} {\bf J}_N^{j} {\bf H}^{T} {\bf W}^{*} ({\bf I}_M \otimes  \sigma^{2} {\bf R} \mathcal{T}_{L,L}({\bf H} {\bf J}_N^{*i} {\bf H}) {\bf R} {\bf J}_L^{u_2}) \right) \right] -\\
\sigma^{2} c_N \sum_{j=-(L-1)}^{L-1} \mathbb{E}(\omega(i,u_1,j)) \, \frac{1}{ML} \mathbb{E} \left[ \mathrm{Tr}\left( {\bf Q} {\bf W} {\bf J}_N^{j} {\bf H}^{T} {\bf J}_N^{i} {\bf H}^{T}{\bf W}^{*} ({\bf I}_M \otimes {\bf R} {\bf J}_L^{u_2}) \right) \right] - \\
\sigma^{2} c_N \sum_{j=-(L-1)}^{L-1} \mathbb{E}(\omega(u_1,j)) \, \frac{1}{ML} \mathbb{E} \left[ \mathrm{Tr}\left({\bf Q} ({\bf I}_M \otimes {\bf J}_L^{i})  {\bf Q} {\bf W} {\bf J}_N^{j} {\bf H}^{T} {\bf J}_N^{i} {\bf H}^{T} {\bf W}^{*} ({\bf I}_M \otimes   {\bf R} {\bf J}_L^{u_2}) \right) \right]
\end{align*}
Using (\ref{eq:expre-trace-Q-Q}), it is easy to check that 
$$
\frac{1}{ML} \mathbb{E} \left[ \mathrm{Tr}\left({\bf Q} ({\bf I}_M \otimes {\bf J}_L^{i}) {\bf Q} ({\bf I}_M \otimes {\bf J}_L^{u_1}) ({\bf I}_M \otimes \sigma^{2} {\bf R} \mathcal{T}_{L,L}({\bf H} {\bf J}_N^{*i} {\bf H}) {\bf R} {\bf J}_L^{u_2}) \right) \right] = \delta(u_1+u_2=0) \, C(i,u_1,z)
+ \mathcal{O}(\frac{L}{MN}), 
$$
$$
\frac{1}{ML} \mathbb{E} \left[ \mathrm{Tr}\left( {\bf Q} {\bf W} {\bf J}_N^{j} {\bf H}^{T} {\bf W}^{*} ({\bf I}_M \otimes  \sigma^{2} {\bf R} \mathcal{T}_{L,L}({\bf H} {\bf J}_N^{*i} {\bf H}) {\bf R} {\bf J}_L^{u_2}) \right) \right] = \delta(j=u_2-i) C(i,u_2,z) + \mathcal{O}(\frac{L}{MN}), 
$$
$$
\frac{1}{ML} \mathbb{E} \left[ \mathrm{Tr}\left( {\bf Q} {\bf W} {\bf J}_N^{j} {\bf H}^{T} {\bf J}_N^{i} {\bf H}^{T}{\bf W}^{*} ({\bf I}_M \otimes {\bf R} {\bf J}_L^{u_2}) \right) \right] = \delta(j=u_2-i) C(i,u_2,z) + \mathcal{O}(\frac{L}{MN}), 
$$
$$
\frac{1}{ML} \mathbb{E} \left[ \mathrm{Tr}\left({\bf Q} ({\bf I}_M \otimes {\bf J}_L^{i})  {\bf Q} {\bf W} {\bf J}_N^{j} {\bf H}^{T} {\bf W}^{*} ({\bf I}_M \otimes  \sigma^{2} {\bf R} \mathcal{T}_{L,L}({\bf H} {\bf J}_N^{*i} {\bf H}) {\bf R} {\bf J}_L^{u_2}) \right) \right] = 
\delta(j=u_2) C(i,u_2,z) + \mathcal{O}(\frac{L}{MN}),
$$
$$
\frac{1}{ML} \mathbb{E} \left[ \mathrm{Tr}\left({\bf Q} ({\bf I}_M \otimes {\bf J}_L^{i})  {\bf Q} {\bf W} {\bf J}_N^{j} {\bf H}^{T} {\bf J}_N^{i} {\bf H}^{T} {\bf W}^{*} ({\bf I}_M \otimes   {\bf R} {\bf J}_L^{u_2}) \right) \right] = \delta(j=u_2) C(i,u_2,z) + \mathcal{O}(\frac{L}{MN})
$$
Proposition \ref{prop:expre-omega-2} and Proposition \ref{prop:expre-omega-3} immediately imply that $\mathbb{E}(\beta(i,u_1,u_2))$ can be written as    
(\ref{eq:expre-terme-utile-trace-Q-Q-QW}). We omit the proof of the expression of $\overline{\beta}(i,u,z)$. Moreover, Lemma \ref{le:exemple} implies that function $z \rightarrow \overline{\beta}(i,u,z)$ is the Stieltjes transform of a distribution $D$ whose support 
is included in $\mathcal{S}_N$ and which verifies $< D, \mathbb{1} > = 0$.  \\

We now establish the second part of the proposition, and assume that $c_N > 1$. In this case, 
$0$ is pole of multiplicity 1 of $t(z)$ and $\tilde{t}(z)$ is analytic at $0$. 
It is easy 
to check that for each $j=1, \ldots, 5$, $0$ is pole with multiplicity 1 of function $z \rightarrow z \overline{\beta}_{j}(i,l_1,z)$, and thus of function $z \rightarrow z \, \overline{\beta}(i,l_1,z)$. 
As for function $z \rightarrow \overline{\beta}_{0,1}(i,l_1,l_1,z)$, it can be shown that 
\begin{equation}
\label{eq:relation-beta01-beta}
\overline{\beta}_{0,1}(i,l_1,l_1,z) = \sigma^{2} \, (1 - |l_1|/N) \, t(z) \, (z \tilde{t}(z))^{2} 
\overline{\beta}(i,l_1,z)
\end{equation}
from which we deduce immediately that $0$ is pole with multiplicity 1 of 
$\overline{\beta}_{0,1}(i,l_1,l_1,z)$. The analytic expression of $\overline{\beta}_{1,0}(i,l_1,l_1,z)$ (not provided) allows to conclude immediately that $0$ may be pole with multiplicity 1, but it can be checked 
that the corresponding residue vanishes. Therefore, function $z \rightarrow \overline{\beta}_{1,0}(i,l_1,l_1,z)$ appears to be analytic in a neighbourhood of $0$, and thus coincides with the Stieltjes transform 
of a distribution whose support is included into $\mathcal{S}_N^{(0)}$. In order to complete 
the proof of the proposition, it remains to check that function $z \rightarrow s(i,l_1,z)$ 
is analytic in a neighbourhood of $0$. As $0$ is pole of $z \overline{\beta}(i,l_1,z)$ and 
$\overline{\beta}_{0,1}(i,l_1,l_1,z)$ with multiplicity 1, it is sufficient to verify 
that 
$$
\lim_{z \rightarrow 0} z \left[ \overline{\beta}_{0,1}(i,l_1,l_1,z) + 
 \sigma^{4} c_N \left(z t(z) \tilde{t}(z)\right)^{2} z \tilde{t}(z) \left( 1 + \sigma^{2} z t(z) \tilde{t}(z)
(1 - |l_1|/L)(1- |l_1|/N) \right) \left(\frac{1 - |l_1|/N}{1 - d(l_1,z)} \, \overline{\beta}(i,l_1,z) \right) \right] =  0
$$
This property follows immediately from (\ref{eq:relation-beta01-beta}). 
\subsubsection{Evaluation of $\kappa^{(2)}(l_1,l_2)$.}
\label{subsub:eval-kappa}
The treatment of the terms $\kappa^{(2)}(l_1,l_2)$ appears to be difficult, and also needs  
a sharp evaluation for each $r$ of the term 
of $\kappa^{(r)}(u_1, \ldots, u_r)$ defined for 
$u_1, \ldots, u_r \in \{ -(L-1), \ldots, L-1 \}$ by 
\begin{equation}
\label{eq:def-kappa}
\kappa^{(r)}(u_1, \ldots, u_r) = \mathbb{E}\left(\Pi_{s=1}^{r} \, \tau^{(M)}({\bf Q}^{\circ})(u_s) \right)
\end{equation}
Lemma \ref{le:generalisation-NP} and the Hölder inequality immediately lead to $\kappa^{(r)}(u_1, \ldots, u_r) = \mathcal{O}(\frac{1}{(MN)^{r/2}})$, but this evaluation is 
not optimal, and has to be refined, in particular if $r=2$. More precisely, the following result holds. 
\begin{proposition}
\label{prop:difficile}
If $z$ belongs to a set $E_N$ defined as in (\ref{eq:justif-EN-nonvoid}), then, for $r=2$, it holds that 
\begin{equation}
\label{eq:estimation-kappa}
\kappa^{(2)}(u_1,u_2) = \frac{1}{MN} \, C(z,u_1) \,  \delta(u_1+u_2=0) \, + \, \mathcal{O}(\frac{L}{(MN)^{2}})
\end{equation}
More generally, if $r \geq 2$, and if $(u_1, u_2, \ldots, u_r)$ are integers such that 
$-(L-1) \leq u_i \leq (L-1)$ for $i=1, \ldots, r$ for which $u_k + u_l \neq 0$ for each 
$k,l$, $k \neq l$, then, it holds that  
\begin{equation}
\label{eq:estimation-kappa-r}
\kappa^{(r)}(u_1, \ldots, u_r) = \frac{1}{\sqrt{MN}} \mathcal{O}(\frac{1}{(MN)^{r/2}})
\end{equation}
\end{proposition}
The proof of this result is quite intricate. The goal of paragraph \ref{subsub:eval-kappa}
is to establish Proposition \ref{prop:difficile}. \\

In order to evaluate $\kappa^{(r)}(u_1, \ldots, u_r)$, we state the following result. It can be proved by calculating,
for each integers $(l_1, l_1^{'}, n_1, n_1^{'}, \ldots, l_{r}, l_{r}^{'}, n_{r}, n_{r}^{'})$, 
matrix
$$
\mathbb{E}\left[ \Pi_{s=1}^{r}({\bf Q}^{\circ})_{l_s,l_s^{'}}^{n_s,n_s^{'}} \;  {\bf Q} {\bf W} {\bf G} {\bf W}^{*} \right]
$$
by the integration by parts formula. This calculation is provided in \cite{loubaton-arxiv}. 
\begin{proposition}
\label{cor:expre-tau}
We consider integers $(u_1, u_2, \ldots, u_r)$,  $(v_1, v_2, \ldots, v_r)$ such that $-(L-1) \leq u_i \leq (L-1)$, 
$-(L-1) \leq v_i \leq (L-1)$ for $i=1, \ldots, r$. 
Then, it holds that 
\begin{multline}
\label{eq:expre-produit-tau-rond} 
\mathbb{E} \left[ \Pi_{s=1}^{r}  \, \tau^{(M)}({\bf Q}^{\circ})(u_s) \right]  = 
- \mathbb{E} \left[ \Pi_{s=1}^{r-1}  \, \tau^{(M)}({\bf Q}^{\circ})(u_s)  \right] \, 
\frac{1}{ML} \mathrm{Tr} \left( \bs{\Delta} ({\bf I}_M \otimes {\bf J}_L^{u_r}) \right) + \\
 \sigma^{2} c_N \sum_{l_{1}=-(L-1)}^{L-1} \mathbb{E} \left( \Pi_{s=1}^{r-1} \tau^{(M)}({\bf Q}^{\circ})(u_s) \, \tau^{(M)}({\bf Q}^{\circ})(l_1) \right) \, \mathbb{E} \left[ \frac{1}{ML} 
\mathrm{Tr} \left( {\bf Q} {\bf W} {\bf J}_N^{l_{1}} {\bf H}^{T} {\bf W}^{*} ({\bf I}_M \otimes 
{\bf R} {\bf J}_L^{u_r}) \right) \right] + \\
 \sigma^{2} c_N \sum_{l_{1}=-(L-1)}^{L-1} \mathbb{E} \left( \Pi_{s=1}^{r-1} \tau^{(M)}({\bf Q}^{\circ})(u_s) \, \tau^{(M)}({\bf Q}^{\circ})(l_1)  \,  \left[ \frac{1}{ML} 
\mathrm{Tr} \left( {\bf Q} {\bf W} {\bf J}_N^{l_{1}} {\bf H}^{T} {\bf W}^{*} ({\bf I}_M \otimes 
{\bf R} {\bf J}_L^{u_r}) \right) \right]^{\circ} \right) + \\
 \frac{\sigma^{2}}{MLN} \sum_{s=1}^{r-1} \sum_{i=-(L-1)}^{L-1} \mathbb{E} \left[ \Pi_{t \neq s,r}  \, \tau^{(M)}({\bf Q}^{\circ})(u_t) \right] \; \mathbb{E}(\beta(i,u_s,u_r))  + 
 \frac{\sigma^{2}}{MLN} \sum_{s=1}^{r-1} \sum_{i=-(L-1)}^{L-1}  \mathbb{E} \left[ \Pi_{t \neq s,r}  \, \tau^{(M)}({\bf Q}^{\circ})(u_t) \, \beta(i,u_s,u_r)^{(0)} \right] 
\end{multline}
and that 
\begin{multline}
\label{eq:expre-produit-tau-rond-W}
\mathbb{E} \left[ \Pi_{s=1}^{r}  \, \tau^{(M)}({\bf Q}^{\circ})(v_s) \, \left( \frac{1}{ML} \mathrm{Tr} ({\bf Q} {\bf W} {\bf G} {\bf W}^{*} {\bf A}) \right)^{\circ} \right]  = \kappa^{(r)}(v_1, \ldots, v_r) \,  \epsilon({\bf G}, {\bf A}) \; + \\
\sigma^{2} c_N \sum_{l_2 = -(L-1)}^{L-1}  \mathbb{E} \left[ \Pi_{s=1}^{r} \tau^{(M)}({\bf Q}^{\circ})(v_s) \, \tau^{(M)}({\bf Q}^{\circ})(l_2)  \,  \frac{1}{ML} 
\mathrm{Tr} \left( {\bf Q} {\bf W} {\bf J}_N^{l_{2}} {\bf H}^{T} {\bf W}^{*} \left( {\bf I}_M \otimes \sigma^{2} {\bf R}  \mathcal{T}_{L,L}({\bf G}^{T} {\bf H}) \right) {\bf A} \right)  \right] \; - \\
\sigma^{2} c_N \sum_{l_2 = -(L-1)}^{L-1}  \mathbb{E} \left[ \Pi_{s=1}^{r} \tau^{(M)}({\bf Q}^{\circ})(v_s) \, \tau^{(M)}({\bf Q}^{\circ})(l_2)  \, \frac{1}{ML} 
\mathrm{Tr} \left( {\bf Q} {\bf W} {\bf J}_N^{l_{2}} {\bf H}^{T} {\bf G} {\bf W}^{*}  {\bf A}^{*} \right) \right] + \\
\frac{\sigma^{2}}{MLN} \sum_{s \leq r, |i| \leq L-1}  \mathbb{E} \left[ \Pi_{t \neq s}  \, \tau^{(M)}({\bf Q}^{\circ})(v_t) \, 
 \frac{1}{ML} \mathrm{Tr} \left( {\bf Q} ({\bf I}_M \otimes {\bf J}_L^{i}) {\bf Q} 
({\bf I}_M \otimes {\bf J}_L^{v_s}) {\bf Q} {\bf W} {\bf J}_N^{i} {\bf H}^{T} {\bf W}^{*}  \left( {\bf I}_M \otimes \sigma^{2} {\bf R}  \mathcal{T}_{L,L}({\bf G}^{T} {\bf H}) \right) {\bf A} \right)  \right] \\
- \frac{\sigma^{2}}{MLN} \sum_{s \leq r, |i| \leq L-1}   \mathbb{E} \left[ \Pi_{t \neq s}  \, \tau^{(M)}({\bf Q}^{\circ})(v_t) \, 
 \frac{1}{ML} \mathrm{Tr} \left( {\bf Q} ({\bf I}_M \otimes {\bf J}_L^{i}) {\bf Q} 
({\bf I}_M \otimes {\bf J}_L^{v_s}) {\bf Q} {\bf W} {\bf J}_N^{i} {\bf H}^{T} {\bf G} {\bf W}^{*}  {\bf A} \right) \right] 
\end{multline}
where we recall that $\beta(i,u_s,u_r)$ is defined by (\ref{eq:def-beta}) and where $\epsilon({\bf G}, {\bf A})$ is defined by 
$$
\epsilon({\bf G}, {\bf A}) = \sigma^{2} c_N \mathbb{E} \left( \frac{1}{ML} \left( {\bf Q} {\bf W} \mathcal{T}^{(M)}_{N,L}({\bf Q}^{\circ})^{T} {\bf H}^{T} \left( {\bf G} {\bf W}^{*} {\bf A} - {\bf W}^{*} \left( {\bf I}_M \otimes \sigma^{2} {\bf R}  \mathcal{T}_{L,L}({\bf G}^{T} {\bf H}) \right) {\bf A} \right) \right) \right)
$$ 
\end{proposition}
In order to evaluate $\kappa^{(r)}(u_1, \ldots, u_{r-1}, u_r)$, we
interpret (\ref{eq:expre-produit-tau-rond}) as a linear system whose unknowns are 
the \\ 
$(\kappa^{(r)}(u_1, \ldots, u_{r-1}, u_r))_{u_r= -(L-1), \ldots, L-1}$, 
the integers $(u_s)_{s=1, \ldots, r-1}$ being considered as fixed. \\

\paragraph{Structure of the linear system.}
We now precise 
the structure of this linear system. We denote by \\
${\bs \kappa}^{(r)} = (\kappa^{(r)}(u_1, \ldots, u_{r-1}, u_r))_{u_r= -(L-1), \ldots, L-1}$ the corresponding $2L-1$--dimensional vector.  
We remark that the second term of the righthandside of (\ref{eq:expre-produit-tau-rond})
coincides with component $u_r$ of the action of vector $\bs{\kappa}^{(r)}$ 
on the matrix whose entry $(u_r,l_1)$ is 
$$
\sigma^{2} c_N \, \mathbb{E} \left[ \frac{1}{ML} 
\mathrm{Tr} \left( {\bf Q} {\bf W} {\bf J}_N^{l_{1}} {\bf H}^{T} {\bf W}^{*} ({\bf I}_M \otimes 
{\bf R} {\bf J}_L^{u_r}) \right) \right]
$$
This matrix appears to be close from a diagonal matrix because
$$
\sigma^{2} c_N \, \mathbb{E} \left[ \frac{1}{ML} 
\mathrm{Tr} \left( {\bf Q} {\bf W} {\bf J}_N^{l_{1}} {\bf H}^{T} {\bf W}^{*} ({\bf I}_M \otimes 
{\bf R} {\bf J}_L^{u_r}) \right) \right] = \delta(l_1+u_r=0) \, d(u_r,z) + \mathcal{O}(\frac{L}{MN})
$$
(see (\ref{eq:D-perturbee})). We now study the fourth and the fifth term of the righthandside of (\ref{eq:expre-produit-tau-rond}). We introduce $y_{1,u_r}$ and $y_{2,u_r}$ defined by 
$$
y_{1,u_r} = \frac{\sigma^{2}}{MLN} \sum_{s=1}^{r-1} \sum_{i=-(L-1)}^{L-1} \mathbb{E}(\beta(i,u_s,u_r)) \; \mathbb{E} \left[ \Pi_{t \neq s,r}  \, \tau^{(M)}({\bf Q}^{\circ})(u_t) \right]
$$
and
\begin{equation}
\label{eq:def-y2-ur}
y_{2,u_r} = \frac{\sigma^{2}}{MLN} \sum_{s=1}^{r-1} \sum_{i=-(L-1)}^{L-1}  \mathbb{E} \left[ \Pi_{t \neq s,r} \, \tau^{(M)}({\bf Q}^{\circ})(u_t) \,  \beta(i,u_s,u_r)^{\circ}  \right] 
\end{equation}
and denote by ${\bf y}_1$ and ${\bf y}_2$ the corresponding $2L-1$--dimensional related vectors. 
We first evaluate the behaviour of ${\bf y}_1$. (\ref{eq:expre-terme-utile-trace-Q-Q-QW}) and the rough evaluation $ \mathbb{E} \left[ \Pi_{t \neq (s,r)}  \, \tau^{(M)}({\bf Q}^{\circ})(u_t) \right] = 
\mathcal{O}(\frac{1}{(MN)^{r/2-1}})$ based on Lemma \ref{le:generalisation-NP} and the Hölder inequality imply that  vector ${\bf y}_1$ can be written as
\begin{equation}
\label{eq:expre-vecteur-y1}
{\bf y}_1 = {\bf y}_1^{*} + {\bf z}_1
\end{equation}
where all the components of ${\bf z}_1$ are $\frac{L}{MN} \mathcal{O}(\frac{1}{(MN)^{r/2}})$ terms, 
or equivalently $\frac{L}{\sqrt{MN}} \mathcal{O}(\frac{1}{(MN)^{(r+1)/2}}) = o(\frac{1}{(MN)^{(r+1)/2}})$
and where ${\bf y}_1^{*}$ is defined by 
\begin{equation}
\label{eq:def-y1*}
{\bf y}_{1,u_r}^{*}  =  \frac{\sigma^{2}}{MN}\; \left( \frac{1}{L} \sum_{s=1}^{r-1} \sum_{i=-(L-1)}^{L-1} \overline{\beta}(i,u_s) \, \delta(u_s + u_r = 0) \right) \, \kappa^{(r-2)}\left( 
(u_t)_{t \neq (s,r)} \right) 
\end{equation}
so that 
\begin{equation}
\label{eq:def-y1*-0}
{\bf y}_{1,u_r}^{*}  =  0  \; \mbox{if $u_r \neq -u_s$ for each $s=1, \ldots, r-1$}
\end{equation}
Hence, (\ref{eq:def-y1*},\ref{eq:def-y1*-0}) imply that
\begin{equation}
\label{eq:eval-y1}
y_{1,u_r}  = \mathcal{O}(\frac{1}{(MN)^{r/2}}) \, \mathbb{1}_{u_r \in \{-u_1, \ldots, -u_{r-1} \}} + 
\frac{L}{\sqrt{MN}} \mathcal{O}(\frac{1}{(MN)^{(r+1)/2}}) \, \mathbb{1}_{u_r \in \{-u_1, \ldots, -u_{r-1} \}^{c}}   
\end{equation} 
We note that if $r=3$, $y_{1,u_3} = 0$ for each $u_3$ because for each $s=1,2$, the term $ \mathbb{E} \left[ \Pi_{t \neq s,3}  \, \tau^{(M)}({\bf Q}^{\circ})(u_t) \right]$ is identically zero. Therefore, for $r=3$, it holds that 
${\bf y}_1^{*} = 0$.

As for ${\bf y}_2$, we notice that Lemma \ref{le:generalisation-NP} and the Hölder inequality lead to 
\begin{equation}
\label{eq:eval-y2}
y_{2,u_r} = \mathcal{O}(\frac{1}{(MN)^{(r+1)/2}})
\end{equation}
We remark that if $r=2$, then $y_{2,u} = 0$ for each $u$ because the term 
$\Pi_{t \neq s,r}  \, \tau^{(M)}({\bf Q}^{\circ})(u_s)$ disappears, and that $y_{2,u}$ 
represents the mathematical expectation of a zero mean term. \\

In order to evaluate the  third term of the righthandside of (\ref{eq:expre-produit-tau-rond}), we define 
$\tilde{x}(u_r,l_1)$ by
\begin{equation}
\label{eq:def-tildexul}
\tilde{x}(u_r,l_1) = \mathbb{E} \left( \Pi_{s=1}^{r-1} \tau^{(M)}({\bf Q}^{\circ})(u_s) \, \tau^{(M)}({\bf Q}^{\circ})(l_1)  \,  \left[ \frac{1}{ML} 
\mathrm{Tr} \left( {\bf Q} {\bf W} {\bf J}_N^{l_{1}} {\bf H}^{T} {\bf W}^{*} ({\bf I}_M \otimes 
{\bf R} {\bf J}_L^{u_r}) \right) \right]^{\circ} \right),
\end{equation}
and $\tilde{x}(u_r)$ by  
\begin{equation}
\label{eq:def-tildexu}
\tilde{x}(u_r) = \sum_{l_1=-(L-1)}^{L-1} \tilde{x}(u_r,l_1)
\end{equation}
In order to have a better understanding of $\tilde{x}(u_r)$, we expand $\tilde{x}(u_r,l_1)$ for each $l_1$ using (\ref{eq:expre-produit-tau-rond-W}). 
We define $(v_1, \ldots, v_r)$ by $v_s = u_s$ for $s \leq r-1$ and $v_r = l_1$, 
while ${\bf G}$ and ${\bf A}$ represent the matrices ${\bf G} = {\bf J}_N^{l_{1}} {\bf H}^{T}$, and ${\bf A} = ({\bf I}_M \otimes {\bf R} {\bf J}_L^{u_r})$. We denote by $(s_i(u_r,l_1))_{i=1, \ldots, 5}$ the $i$-th term of the righthandside of 
(\ref{eq:expre-produit-tau-rond-W}), and denote by $(s_i(u_r))_{i=1, \ldots, 5}$ the term
$$
s_i(u_r) = \sum_{l_1 = -(L-1)}^{L-1} s_{i}(u_r,l_1)
$$
and by ${\bf s}_i$ vector ${\bf s}_i = (s_i(u_r))_{u_r = -(L-1), \ldots, L-1}$. Vector ${\bf s}_1$ plays a particular role 
because $s_1(u_r,l_1)$ is equal to 
$$
s_1(u_r,l_1) = \kappa^{(r)}(u_1, \ldots, u_{r-1}, l_1) \, \epsilon({\bf J}_N^{l_1} {\bf H}^{T}, {\bf I} \otimes {\bf R} {\bf J}_L^{u_r}) =  \kappa^{(r)}(u_1, \ldots, u_{r-1}, l_1) \, \mathcal{O}(\frac{L}{MN})
$$
We remark that vector ${\bf s}_1$ coincides with the action of vector $\bs{\kappa}^{(r)}$ on matrix 
$ \left(\epsilon({\bf J}_N^{l_1}, {\bf I} \otimes {\bf R} {\bf J}_L^{u_r}) \right)_{-(L-1) \leq u_r, l_1  \leq (L-1)}$.
We define  by $x(u_r,l_1)$ and  $x(u_r)$ the terms 
\begin{equation}
\label{eq:def-x_u_r_l_1}
x(u_r,l_1) = \sum_{i=2}^{5} s_i(u_r,l_1), \; x(u_r) = \sum_{l_1=-(L-1)}^{L-1} x(u_r,l_1)
\end{equation}
and vector ${\bf x}$ represents the $2L-1$--dimensional vector $(x(u_r))_{u_r = -(L-1), \ldots, L-1}$.

We finally consider the first term of the righthandside of (\ref{eq:expre-produit-tau-rond}), and
denote by ${\bs \epsilon}$ the $2L-1$--dimensional vector whose components $(\epsilon_{u_r})_{u_r=-(L-1), \ldots, L-1}$ are given by
$$
\epsilon_{u_r} = - \mathbb{E} \left[ \Pi_{s=1}^{r-1}  \, \tau^{(M)}({\bf Q}^{\circ})(u_s)  \right] \, 
\frac{1}{ML} \mathrm{Tr} \left( \bs{\Delta} ({\bf I}_M \otimes {\bf J}_L^{u_r}) \right)
$$
We notice that if $r=2$, vector ${\bs \epsilon}$ is reduced to $0$.  \\

This discussion and (\ref{eq:D-perturbee}) imply that (\ref{eq:expre-produit-tau-rond}) can be written as 
\begin{equation}
\label{eq:systeme-lineaire-kappar}
{\bs \kappa}^{(r)} = \left( {\bf D} + \bs{\Upsilon} \right) \, {\bs \kappa}^{(r)}  + {\bf y}_{1,*} + {\bf z}_1  + {\bf y}_2 + \bs{\epsilon} + \sigma^{2} c_N \, {\bf x}
\end{equation}
where we recall that ${\bf D}$ represents the diagonal matrix 
${\bf D} = \mathrm{Diag}(d(-(L-1),z), \ldots, d((L-1),z))$ 
and where the entries of matrix $\bs{\Upsilon}$ are defined by
$$
\bs{\Upsilon}_{u_r,l_1} = \sigma^{2} c_N \, \mathbb{E} \left[ \frac{1}{ML} 
\mathrm{Tr} \left( {\bf Q} {\bf W} {\bf J}_N^{l_{1}} {\bf H}^{T} {\bf W}^{*} ({\bf I}_M \otimes 
{\bf R} {\bf J}_L^{u_r}) \right) \right] - {\bf D}_{u_r,l_1} + \sigma^{2} c_N \, \epsilon({\bf J}_N^{l_{1}} {\bf H}^{T}, 
({\bf I}_M \otimes {\bf R} {\bf J}_L^{u_r}))
$$
It is clear the each entry of $\bs{\Upsilon}$ is a $\mathcal{O}(\frac{L}{MN})$ term. \\

\paragraph{Overview of the proof of Proposition \ref{prop:difficile}.}
We now present unformally the various steps of the proof of Proposition \ref{prop:difficile}, 
and concentrate on the proof of Eq. (\ref{eq:estimation-kappa-r}) in order to
simplify the presentation. The particular case $r=2$ is 
however briefly considered at the end of the overview, but it is of  
course detailed in the course of the proof. \\

{\bf First step: inversion of the linear system (\ref{eq:systeme-lineaire-kappar}).} 
Lemma \ref{le:inversion-D} implies that if $z$ belongs to a set $E_N$ defined as 
(\ref{eq:justif-EN-nonvoid}), matrix $({\bf I} - {\bf D} - \bs{\Upsilon})$ is invertible. Therefore, vector ${\bs \kappa}^{(r)}$ can be written as 
$$
{\bs \kappa}^{(r)} = ({\bf I} - {\bf D} - \bs{\Upsilon})^{-1} \, \left(
 {\bf y}_{1,*} + {\bf z}_1  + {\bf y}_2 + \bs{\epsilon}  + \sigma^{2} c_N \, {\bf x} \right)
$$ 
Using (\ref{eq:norme-infinie-inverse}) and the properties of the components 
of vectors ${\bf z}_1, {\bf y}_2$ and ${\bs \epsilon}$, we obtain easily that 
$$
\left( ({\bf I} - {\bf D} - \bs{\Upsilon})^{-1} {\bf y}_{1,*} \right)_{u_r} =  
\frac{1}{1 - d(u_r,z)} \, {\bf y}^{*}_{1,u_r} + \frac{L}{\sqrt{MN}} \mathcal{O}(\frac{1}{(MN)^{(r+1)/2}}),
$$
$$
\left( ({\bf I} - {\bf D} - \bs{\Upsilon})^{-1} {\bf y}_{2} \right)_{u_r} =  \mathcal{O}(\frac{1}{(MN)^{(r+1)/2}}),
$$
and that 
\begin{equation}
\label{eq:premiere-evaluation-general-overview}
\left| \kappa^{(r)}(u_1, \ldots, u_r) - \frac{1}{1 - d(u_r,z)} {\bf y}_{1,u_r}^{*}\right|  \leq  | \kappa^{(r-1)}(u_1, \ldots, u_{r-1}) | \, \mathcal{O}(\frac{L}{MN}) + 
C \, \sup_{u} \left| x(u) \right| + \mathcal{O}(\frac{1}{(MN)^{(r+1)/2}})  
\end{equation}
If multi-index $(u_1, \ldots, u_r)$ satisfies $u_k + u_l \neq 0$ for $k \neq l$, 
then ${\bf y}_{1,u_r}^{*} = 0$ (see Eq. (\ref{eq:def-y1*-0})). Therefore, in order 
to establish (\ref{eq:estimation-kappa-r}), it is necessary to evaluate 
$\sup_{u} \left| x(u) \right|$. \\

{\bf Second step: evaluation of $\sup_{u} \left| x(u) \right|$.}  In order to 
evaluate $\sup_{u} \left| x(u) \right|$, we express $x(u_r,l_1)$ as 
$x(u_r,l_1) = \sum_{i=2}^{5} s_i(u_r,l_1)$ (see Eq. (\ref{eq:def-x_u_r_l_1})), 
and study each term $s_i(u_r) = \sum_{l_1} s_i(u_r,l_1)$ for $i=2,3,4,5$. 
$s_4(u_r)$ and $s_5(u_r)$ can be written as $\kappa^{(r-1)}(u_1, \ldots, u_{r-1}) \mathcal{O}(\frac{L}{MN}) \delta(u_r = 0) + o\left( \frac{1}{(MN)^{(r+1)/2}} \right)$. The terms $s_2(u_r)$ and 
$s_3(u_r)$ have a more complicated structure. We just address $s_3(u_r)$ because the 
behaviour of $s_2(u_r)$ is similar. $s_3(u_r,l_1)$ can be written as
$s_3(u_r,l_1) = \sum_{l_2} s_3(u_r,l_1,l_2)$ where 
$$
s_3(u_r,l_1,l_2) = -\sigma^{2} c_N \mathbb{E} \left[ \Pi_{s=1}^{r-1}
 \tau^{(M)}({\bf Q}^{\circ})(u_s)  \tau^{(M)}({\bf Q}^{\circ})(l_1)  \tau^{(M)}({\bf Q}^{\circ})(l_2) 
\frac{1}{ML} \mathrm{Tr} \left( {\bf Q} {\bf W} {\bf J}_N^{l_2} {\bf H}^{T} {\bf J}_N^{l_1} {\bf H}^{T} {\bf W}^{*} ({\bf I} \otimes {\bf R} {\bf J}_L^{u_r}) \right) \right]
$$
We define $\overline{s}_3(u_r,l_1,l_2)$ and $\tilde{x}_3^{(1)}(u_r,l_1,l_2)$ by
\begin{equation}
\label{eq:def-overlines3}
\overline{s}_3(u_r,l_1,l_2) = -\sigma^{2} c_N  \kappa^{(r+1)}(u_1, \ldots, u_{r-1}, l_1, l_2) \, \mathbb{E} \left[ \frac{1}{ML} \mathrm{Tr} \left( {\bf Q} {\bf W} {\bf J}_N^{l_2} {\bf H}^{T} {\bf J}_N^{l_1} {\bf H}^{T} {\bf W}^{*} ({\bf I} \otimes {\bf R} {\bf J}_L^{u_r}) \right) \right]
\end{equation}
and 
\begin{multline}
\label{eq:def-tildex1}
\tilde{x}_3^{(1)}(u_r,l_1,l_2)  = \\
 -\sigma^{2} c_N  \mathbb{E} \left[ \Pi_{s=1}^{r-1}
 \tau^{(M)}({\bf Q}^{\circ})(u_s)  \tau^{(M)}({\bf Q}^{\circ})(l_1)  \tau^{(M)}({\bf Q}^{\circ})(l_2) 
\frac{1}{ML} \mathrm{Tr} \left( {\bf Q} {\bf W} {\bf J}_N^{l_2} {\bf H}^{T} {\bf J}_N^{l_1} {\bf H}^{T} {\bf W}^{*} ({\bf I} \otimes {\bf R} {\bf J}_L^{u_r}) \right)^{\circ} \right]
\end{multline}
Then, it holds that 
$$
s_3(u_r,l_1,l_2) = \overline{s}_3(u_r,l_1,l_2) + \tilde{x}_3^{(1)}(u_r,l_1,l_2)
$$
and obtain that $s_3(u_r) = \overline{s}_3(u_r) + \tilde{x}_3^{(1)}(u_r)$ where $\overline{s}_3(u_r)$
and $\tilde{x}_3^{(1)}(u_r)$ are defined as the sum over $l_1,l_2$ of $ \overline{s}_3(u_r,l_1,l_2)$
and $\tilde{x}_3^{(1)}(u_r,l_1,l_2)$. Similarly, $s_2(u_r)$ can be expressed 
as $s_2(u_r) = \overline{s}_2(u_r) + \tilde{x}_2^{(1)}(u_r)$ where $\overline{s}_2(u_r)$ and 
$\tilde{x}_2^{(1)}(u_r)$ are defined in the same way than $\overline{s}_3(u_r)$
and $\tilde{x}_3^{(1)}(u_r)$. The behaviour of $(\overline{s}_j(u_r))_{j=2,3}$ 
is easy to analyse because it can be shown that   
$$
\overline{s}_j(u_r) = \sum_{l_1} C_j(u_r,l_1) \kappa^{(r+1)}(u_1, \ldots, u_{r-1}, l_1, u_r - l_1) + 
\sum_{l_1,l_2}  \kappa^{(r+1)}(u_1, \ldots, u_{r-1}, l_1, l_2) \mathcal{O}(\frac{L}{MN})
$$
Therefore, (\ref{eq:premiere-evaluation-general-overview}) implies that
\begin{multline}
\label{eq:recurrence-fondamentale-generale}
\left| \kappa^{(r)}(u_1, \ldots, u_r) - \frac{1}{1 - d(u_r,z)} {\bf y}_{1,u_r}^{*}\right|\ \leq  
| \kappa^{(r-1)}(u_1, \ldots, u_{r-1}) | \, \mathcal{O}(\frac{L}{MN}) \;  +\\
C \, \sup_{u} \sum_{l_1} |\kappa^{(r+1)}(u_1, \ldots, u_{r-1}, l_1, u - l_1)|  \;  + 
\sum_{l_1,l_2}  |\kappa^{(r+1)}(u_1, \ldots, u_{r-1}, l_1, l_2)| \mathcal{O}(\frac{L}{MN})  \; + \\
\sup_{u} \tilde{x}^{(1)}(u)  + \mathcal{O}(\frac{1}{(MN)^{(r+1)/2}}) 
\end{multline}
where $\tilde{x}^{(1)}(u)$ is the positive term defined by 
$$
\tilde{x}^{(1)}(u) =  |\tilde{x}_2^{(1)}(u)| +  |\tilde{x}_3^{(1)}(u)|
$$
Therefore, if $u_r + u_s \neq 0$ 
for $s=1, \ldots, r-1$, then, ${\bf y}_{1,u_r}^{*} = 0$ and it holds that 
\begin{multline}
\label{eq:recurrence-fondamentale}
\left|\kappa^{(r)}(u_1, \ldots, u_r) \right|  \leq  
| \kappa^{(r-1)}(u_1, \ldots, u_{r-1}) | \, \mathcal{O}(\frac{L}{MN}) \;  +\\
C \, \sup_{u} \sum_{l_1} |\kappa^{(r+1)}(u_1, \ldots, u_{r-1}, l_1, u - l_1)|  \;  + 
\sum_{l_1,l_2}  |\kappa^{(r+1)}(u_1, \ldots, u_{r-1}, l_1, l_2)| \mathcal{O}(\frac{L}{MN})  \; + \\
\sup_{u} \tilde{x}^{(1)}(u)  + \mathcal{O}(\frac{1}{(MN)^{(r+1)/2}}) 
\end{multline}
In order to manage  $\sup_{u} \tilde{x}^{(1)}(u)$, 
we expand $\tilde{x}_j^{(1)}(u,l_1,l_2)$ using (\ref{eq:expre-produit-tau-rond-W}) 
when $r$ is exchanged by $r+1$. In the same way than $\tilde{x}(u)$ defined by (\ref{eq:def-tildexu}), it holds that
$$
\tilde{x}_j^{(1)}(u) = \sum_{i=1}^{5} s_{j,i}^{(1)}(u)
$$
where the terms $(s_{j,i}^{(1)}(u))_{i=1, \ldots, 5}$ are defined in the same way than 
$(s_i(u))_{i=1, \ldots, 5}$. We define $\tilde{x}_{j,i}^{(2)}(u)$ for $i=2,3$ by the fact that
$$
s_{j,i}^{(1)}(u) = \overline{s}_{j,i}^{(1)}(u) + \tilde{x}_{j,i}^{(2)}(u)
$$
We define $\tilde{x}^{(2)}(u)$ as the positive term given by
$$
\tilde{x}^{(2)}(u) = \sum_{(i,j)=(2,3)} |\tilde{x}_{j,i}^{(2)}(u)|
$$
The terms $\tilde{x}_{j,i}^{(2)}(u)$ can be developed similarly, and pursuing the iterative 
process, we are able to define for each $q \geq 3$ the positive terms $\tilde{x}^{(q)}(u)$ 
which are the analogs of $\tilde{x}^{(1)}(u)$ and $\tilde{x}^{(2)}(u)$. In order to 
characterize the behaviour of $\sup_{u} \tilde{x}^{(1)}(u)$, we express $\tilde{x}^{(1)}(u)$
as 
$$
\tilde{x}^{(1)}(u) = \sum_{q=1}^{p-1} \left( \tilde{x}^{(q)}(u) - \tilde{x}^{(q+1)}(u) \right) + \tilde{x}^{(p)}(u)
$$
where the choice of $p$ depends on the context. The term $\tilde{x}^{(p)}(u)$ is easy to control 
because the H\"{o}lder inequality leads immediately to $\tilde{x}^{(p)}(u) = \left( \frac{L}{\sqrt{MN}} \right)^{p+1} \, \mathcal{O}(\frac{1}{(MN)^{r/2}})$. 

Moreover, it is shown that 
\begin{multline}
\tilde{x}^{(q)}(u) - \tilde{x}^{(q+1)}(u) \leq  \sum_{l_i,i=1, \ldots,q+1} |\kappa^{(r+q)}(u_1, \ldots, u_{r-1},l_i,i=1, \ldots,q+1)| \, \mathcal{O}(\frac{L}{MN}) \;  + \\
C \,  \sum_{l_i,i=1, \ldots,q+1} |\kappa^{(r+q+1)}(u_1, \ldots, u_{r-1},l_i,i=1, \ldots,q+1, u-\sum_{i=1}^{q+1} l_i)|  \;  + \\
 \sum_{l_i,i=1, \ldots,q+2} |\kappa^{(r+q+1)}(u_1, \ldots, u_{r-1},l_i,i=1, \ldots,q+2)| \,  \mathcal{O}(\frac{L}{MN})   +  o(\frac{1}{(MN)^{(r+1)/2}})
\end{multline} 
This allows to evaluate $\sum_{q=1}^{p-1} \left( \tilde{x}^{(q)}(u) - \tilde{x}^{(q+1)}(u) \right)$
in the course of the proof.  \\

{\bf Third step: establishing (\ref{eq:estimation-kappa-r}).} 
(\ref{eq:recurrence-fondamentale}) suggests that the rough evaluation 
$\kappa^{(r)}(u_1, \ldots, u_r) = \mathcal{O}(\frac{1}{(MN)^{r/2}})$ can be improved 
when $u_k+u_l \neq 0$ for $k \neq l$. 
The first term of the righthandside of (\ref{eq:recurrence-fondamentale})  can also be written as
$$
\frac{L}{\sqrt{MN}} \, \frac{1}{\sqrt{MN}} \left| \kappa^{(r-1)}(u_1, \ldots, u_{r-1}) \right|
$$
Even if we evaluate $\kappa^{(r-1)}(u_1, \ldots, u_{r-1})$ as $\mathcal{O}(\frac{1}{(MN)^{(r-1)/2}})$, 
it is clear the first term of the righthandside of (\ref{eq:recurrence-fondamentale}) appears 
as a $\frac{L}{\sqrt{MN}} \mathcal{O}(\frac{1}{(MN)^{r/2}})$.  A factor 
$\frac{L}{\sqrt{MN}}$ is thus obtained w.r.t. the rate $\mathcal{O}(\frac{1}{(MN)^{r/2}})$. 
One may imagine that using the information that 
$u_i + u_j \neq 0$ for $1 \leq i,j \leq r-1$, $i \neq j$, should allow to improve the above rough evaluation 
of $ \kappa^{(r-1)}(u_1, \ldots, u_{r-1})$, and thus the evaluation of the first term of the righthandside of 
(\ref{eq:recurrence-fondamentale}). A similar phenomenon is observed for the second term and the third terms of the righthandside of (\ref{eq:recurrence-fondamentale}). We just consider the second term. If each term 
$\kappa^{(r+1)}(u_1, \ldots,u_{r-1},l_1,u-l_1)$ is roughly evaluated as 
$\mathcal{O}(\frac{1}{(MN)^{(r+1)/2}})$, taking into account
the sum over $l_1$, the  second term of the righthandside of (\ref{eq:recurrence-fondamentale})
is decreased by a factor $\frac{L}{\sqrt{MN}}$ w.r.t. the rough evaluation 
$\mathcal{O}(\frac{1}{(MN)^{r/2}})$.

In order to formalize the above discussion, it seems reasonable to be able to prove (\ref{eq:estimation-kappa-r}) from (\ref{eq:recurrence-fondamentale}) using induction technics. However, this needs some care because 
$|\kappa^{(r)}(u_1, \ldots, u_r)|$ is controlled by $|\kappa^{(r-1)}(u_1, \ldots, u_{r-1})|$
and by similar terms of orders greater than $r$. In order to establish (\ref{eq:estimation-kappa-r}), it is proved in Proposition  \ref{prop:deuxieme-reccurence} that if $(u_1, \ldots, u_r)$ satisfy $u_t + u_s \neq 0$ for $1 \leq t,s \leq r$ and $t \neq s$, then, 
for each $q \geq 1$, for each $r \geq 2$, it holds that 
\begin{equation}
\label{eq:evaluation-kappa-reccurence-2-preliminaire}
\kappa^{(r)}(u_1, \ldots, u_r) = \max \left( \left(\frac{L}{\sqrt{MN}}\right)^{r-1+q}, \frac{1}{\sqrt{MN}} \right) \; \mathcal{O}(\frac{1}{(MN)^{r/2}})
\end{equation}
This leads immediately to (\ref{eq:estimation-kappa-r}) because, as $L = \mathcal{O}(N^{\alpha})$
with $\alpha < 2/3$, it exists $q$ for which $\left(\frac{L}{\sqrt{MN}}\right)^{r-1+q} = o\left(\frac{1}{\sqrt{MN}} \right)$. In order to establish (\ref{eq:evaluation-kappa-reccurence-2-preliminaire}), 
we first show in Proposition \ref{prop:premiere-reccurence}
that for each $r \geq 2$ and each integer $1 \leq p \leq r-1$, if integers $u_1, \ldots, u_r \in \{ -(L-1), \ldots, L-1 \}$ satisfy
\begin{equation}
\label{eq:conditions-p-u-preliminaire}
\begin{array}{ccc} u_r + u_s & \neq & 0 \; s=1, \ldots, r-1 \\
                   u_{r-1} + u_s & \neq & 0 \; s=1, \ldots, r-2 \\
                   \vdots & \vdots & \vdots \\
                   u_{r-p+1} + u_s & \neq & 0 \; s=1, \ldots, r-p 
\end{array}
\end{equation}
then, it holds that 
\begin{equation}
\label{eq:evaluation-kappa-reccurence-1-preliminaire}
\kappa^{r}(u_1, \ldots, u_r) = \max \left( \left(\frac{L}{\sqrt{MN}}\right)^{p}, \frac{1}{\sqrt{MN}} \right) \; \mathcal{O}(\frac{1}{(MN)^{r/2}})
\end{equation}
Using (\ref{eq:recurrence-fondamentale}) as well as the above evaluation of $\sup_{u} \tilde{x}^{(1)}(u)$, we prove Proposition \ref{prop:premiere-reccurence} by induction on $r$: we verify that it holds
for $r=2$, assume that it holds until integer $r_0 -1$, and establish it is true for integer $r_0$. 
For this, we prove that for each $r \geq r_0$ and for each multi-index $(u_1, \ldots, u_r)$ 
satisfying (\ref{eq:conditions-p-u-preliminaire}) for $p \leq r_0 -1$, then (\ref{eq:evaluation-kappa-reccurence-1-preliminaire}) holds. This is established by induction on integer $p$ in Lemma \ref{le:aide-a-la-comprehension}.

We note that (\ref{eq:evaluation-kappa-reccurence-1-preliminaire}) used for integer $p=r-1$
coincides with (\ref{eq:evaluation-kappa-reccurence-2-preliminaire}) for $q=0$. (\ref{eq:evaluation-kappa-reccurence-2-preliminaire}) is established for each integer $q$ by induction on integer $q$. 
It is first established by induction on $r$ that (\ref{eq:evaluation-kappa-reccurence-2-preliminaire})
holds for each $r$ for $q=1$. Then, (\ref{eq:evaluation-kappa-reccurence-2-preliminaire})
is assumed to hold for each $r$ until integer $q-1$, and we prove by induction on $r$ 
that it holds for integer $q$. For this, it appears necessary to evaluate 
$$
\sum_{l_1} \left| \kappa^{(r+1)}(u_1, \ldots, u_{r-1}, l_1, -l_1) \right|
$$
where $u_1, \ldots, u_{r-1}$ verify $u_k + u_l \neq 0$ for each $k,l \in {1, 2, \ldots, r-1}$ (see
Lemma \ref{le:kappa(l,-l)}).  
This expression corresponds to the second term of the righthandside of (\ref{eq:recurrence-fondamentale}) 
for $u=0$. \\

{\bf Fourth step: establishing (\ref{eq:estimation-kappa}).}
For $r=2$, the term $\mathcal{O}(\frac{1}{(MN)^{(r+1)/2}})$ at the righhandside of 
(\ref{eq:recurrence-fondamentale-generale}) is replaced by a $\mathcal{O}(\frac{L}{(MN)^{2}})$
term because vector ${\bf y}_2$ whose components are defined by (\ref{eq:def-y2-ur}) is identically
$0$. Moreover, the first term at the righthandside of (\ref{eq:recurrence-fondamentale-generale}) 
vanishes. Using (\ref{eq:estimation-kappa-r}), it is easy to prove that the third term
of the righthandside of (\ref{eq:recurrence-fondamentale-generale}) is $o\left(\frac{L}{(MN)^{2}}\right)$. 
(\ref{eq:estimation-kappa}) follows in turn from the evaluation 
$$
\sum_{l_1} \left| \kappa^{(3)}(u_1, l_1, -l_1) \right| = \mathcal{O}(\frac{L}{(MN)^{2}})
$$
which is proved in Lemma \ref{le:evaluation-amelioree-kappa-l-l}. \\

\paragraph{Proof of Proposition \ref{prop:difficile}.}

We now complete the proof of Proposition \ref{prop:difficile}. 
In order to evaluate $\kappa^{(r)}(u_1, \ldots, u_r)$, we use (\ref{eq:systeme-lineaire-kappar}) and Lemma \ref{le:inversion-D}.  
We write that 
$$
{\bs \kappa}^{(r)} = ({\bf I} - {\bf D} - \bs{\Upsilon})^{-1} \, \left(
 {\bf y}_{1,*} + {\bf z}_1  + {\bf y}_2 + \bs{\epsilon}  + \sigma^{2} c_N \, {\bf x} \right)
$$ 
We first evaluate each component of the first 3 terms of the righthandside of the above equation. Vector $({\bf I} - {\bf D} - \bs{\Upsilon})^{-1} \, {\bf y}_{1,*}$ can 
also be written as
$$
({\bf I} - {\bf D} - \bs{\Upsilon})^{-1} \, {\bf y}_{1,*} = ({\bf I} - {\bf D})^{-1} \, {\bf y}_{1,*} + 
({\bf I} - {\bf D} - \bs{\Upsilon})^{-1} \bs{\Upsilon} ({\bf I} - {\bf D})^{-1} {\bf y}_{1,*}
$$
As vector ${\bf y}_{1,*}$ has at most $r-1$ non zero components which are $\mathcal{O}(\frac{1}{(MN)^{r/2}})$ terms and that the entries of 
$\bs{\Upsilon}$ are $\mathcal{O}(\frac{L}{MN})$ terms, the entries of vector $ \bs{\Upsilon} ({\bf I} - {\bf D})^{-1} {\bf y}_{1,*}$ 
are $\frac{L}{\sqrt{MN}} \mathcal{O}(\frac{1}{(MN)^{(r+1)/2}}) = o(\frac{1}{(MN)^{(r+1)/2}})$ terms.  (\ref{eq:norme-infinie-inverse}) implies that the entries of $({\bf I} - {\bf D} - \bs{\Upsilon})^{-1} \bs{\Upsilon} ({\bf I} - {\bf D})^{-1} {\bf y}_{1,*}$
are $\frac{L}{\sqrt{MN}} \mathcal{O}(\frac{1}{(MN)^{(r+1)/2}})$ terms as well. Therefore, it holds that
$$
\left( ({\bf I} - {\bf D} - \bs{\Upsilon})^{-1} \, {\bf y}_{1,*} \right)_{u_r} = \frac{1}{1 - d(u_r,z)} \, {\bf y}^{*}_{1,u_r} + \frac{L}{\sqrt{MN}} \mathcal{O}(\frac{1}{(MN)^{(r+1)/2}})
$$
and that this term is reduced to a $\frac{L}{\sqrt{MN}} \mathcal{O}(\frac{1}{(MN)^{(r+1)/2}})$ if 
$u_r$ does not belong to $\{ -u_1, \ldots, -u_{r-1} \}$. 
(\ref{eq:norme-infinie-inverse}) implies that 
$$
\left( ({\bf I} - {\bf D} - \bs{\Upsilon})^{-1} \, {\bf z}_1 \right)_{u_r} = \frac{L}{\sqrt{MN}} \mathcal{O}(\frac{1}{(MN)^{(r+1)/2}})
$$
and that
$$
\left( ({\bf I} - {\bf D} - \bs{\Upsilon})^{-1} \, {\bf y}_2 \right)_{u_r} =  \mathcal{O}(\frac{1}{(MN)^{(r+1)/2}})
$$
for $r \geq 3$, while this term is zero for $r=2$ because ${\bf y}_2 = 0$ in this case.  
If $u_r$ does not belong to $\{ -u_1, \ldots, -u_{r-1} \}$, the contributions
of the above 3 terms to $\kappa^{(r)}(u_1, \ldots, u_r)$ are at most $\mathcal{O}(\frac{1}{(MN)^{(r+1)/2}})$
terms, which corresponds to what is expected because we recall that the goal of the subsection 
is to establish that $\kappa^{(r)}(u_1, \ldots, u_r) = \mathcal{O}(\frac{1}{(MN)^{(r+1)/2}})$ if 
$u_k + u_l \neq 0$ for $k \neq l$ (see (\ref{eq:estimation-kappa-r})). Finally, 
(\ref{eq:norme-infinie-inverse})
implies that 
\begin{equation}
\label{eq:contrib-epsilon}
 \left( ({\bf I} - {\bf D} - \bs{\Upsilon} )^{-1} \, {\bs \epsilon} \right)_{u_r} = 
\kappa^{(r-1)}(u_1, \ldots, u_{r-1}) \, \mathcal{O}(\frac{L}{MN})
\end{equation}
and that 
\begin{equation}
\label{eq:contrib-x}
\sup_{u} \left| \left( ({\bf I} - {\bf D} - \bs{\Upsilon})^{-1} \, {\bf x} \right)_{u} \right| \leq 
C \, \sup_{u} \left| x(u) \right|
\end{equation} 
Therefore, it holds that 
\begin{equation}
\label{eq:premiere-evaluation-general}
\left| \kappa^{(r)}(u_1, \ldots, u_r) - \frac{1}{1 - d(u_r,z)} {\bf y}_{1,u_r}^{*}\right|  \leq  | \kappa^{(r-1)}(u_1, \ldots, u_{r-1}) | \, \mathcal{O}(\frac{L}{MN}) + 
C \, \sup_{u} \left| x(u) \right| + \mathcal{O}(\frac{1}{(MN)^{(r+1)/2}}) 
\end{equation}
where we recall that ${\bf y}_{1,u_r}^{*} = 0$ if $u_r$ does not belong to $\{ -u_1, \ldots, -u_{r-1} \}$. 
We note that if $r=2$, (\ref{eq:premiere-evaluation-general}) can be written as
\begin{equation}
\label{eq:premiere-evaluation-general-r=2}
\left| \kappa^{(2)}(u_1, u_2) -  \frac{1}{1 - d(u_2,z)} {\bf y}_{1,u_2}^{*} \right| \leq   
C \, \sup_{u} \left| x(u) \right| + \mathcal{O}(\frac{L}{(MN)^{2}})  
\end{equation} 
because ${\bf y}_2 = {\bs \epsilon} = 0$.  \\ 

In order to establish (\ref{eq:estimation-kappa-r}), it is necessary to study the behaviour of $\sup_{u} |x(u)|$. We express 
$x(u_r)$ as $x(u_r) = \sum_{l_1=-(L-1)}^{L-1} x(u_r,l_1)$ and evaluate 
the 4 terms $s_i(u_r) = \sum_{l_1=-(L-1)}^{L-1} s_i(u_r,l_1)$ for $i=2, 3, 4, 5$. We just study $s_{i}(u_r)$ for $i=3$ and $i=5$ because 
$s_2(u_r)$ (resp. $s_4(u_r)$) has essentially the same behaviour than $s_3(u_r)$ (resp. $s_5(u_r)$). $s_3(u_r,l_1)$
is given by 
$$
s_3(u_r,l_1) = \sum_{l_2=-(L-1)}^{L-1} s_3(u_r,l_1,l_2) 
$$
where
$$
s_3(u_r,l_1,l_2) = -\sigma^{2} c_N \mathbb{E} \left[ \Pi_{s=1}^{r-1}
 \tau^{(M)}({\bf Q}^{\circ})(u_s)  \tau^{(M)}({\bf Q}^{\circ})(l_1)  \tau^{(M)}({\bf Q}^{\circ})(l_2) 
\frac{1}{ML} \mathrm{Tr} \left( {\bf Q} {\bf W} {\bf J}_N^{l_2} {\bf H}^{T} {\bf J}_N^{l_1} {\bf H}^{T} {\bf W}^{*} ({\bf I} \otimes {\bf R} {\bf J}_L^{u_r}) \right) \right]
$$
We define $\overline{s}_3(u_r,l_1,l_2)$ and $\tilde{x}_3^{(1)}(u_r,l_1,l_2)$ by
\begin{equation}
\label{eq:def-overlines3}
\overline{s}_3(u_r,l_1,l_2) = -\sigma^{2} c_N  \kappa^{(r+1)}(u_1, \ldots, u_{r-1}, l_1, l_2) \, \mathbb{E} \left[ \frac{1}{ML} \mathrm{Tr} \left( {\bf Q} {\bf W} {\bf J}_N^{l_2} {\bf H}^{T} {\bf J}_N^{l_1} {\bf H}^{T} {\bf W}^{*} ({\bf I} \otimes {\bf R} {\bf J}_L^{u_r}) \right) \right]
\end{equation}
and 
\begin{multline}
\label{eq:def-tildex1}
\tilde{x}_3^{(1)}(u_r,l_1,l_2)  = \\
 -\sigma^{2} c_N  \mathbb{E} \left[ \Pi_{s=1}^{r-1}
 \tau^{(M)}({\bf Q}^{\circ})(u_s)  \tau^{(M)}({\bf Q}^{\circ})(l_1)  \tau^{(M)}({\bf Q}^{\circ})(l_2) 
\frac{1}{ML} \mathrm{Tr} \left( {\bf Q} {\bf W} {\bf J}_N^{l_2} {\bf H}^{T} {\bf J}_N^{l_1} {\bf H}^{T} {\bf W}^{*} ({\bf I} \otimes {\bf R} {\bf J}_L^{u_r}) \right)^{\circ} \right]
\end{multline}
Then, it holds that 
$$
s_3(u_r,l_1,l_2) = \overline{s}_3(u_r,l_1,l_2) + \tilde{x}_3^{(1)}(u_r,l_1,l_2)
$$
We also define $\overline{s}_3(u_r,l_1)$, $\overline{s}_3(u_r)$, $\tilde{x}_3^{(1)}(u_r,l_1)$ and $\tilde{x}_3^{(1)}(u_r)$
as $\overline{s}_3(u_r,l_1) = \sum_{l_2} \overline{s}_3(u_r,l_1,l_2)$, $\overline{s}_3(u_r) = \sum_{l_1} \overline{s}_3(u_r,l_1)$, 
$\tilde{x}_3^{(1)}(u_r,l_1) = \sum_{l_2} \tilde{x}_3^{(1)}(u_r,l_1,l_2)$ and $\tilde{x}_3^{(1)}(u_r) = \sum_{l_1} \tilde{x}_3^{(1)}(u_r,l_1)$. 
It is easy to check that
$$
\mathbb{E} \frac{1}{ML} \mathrm{Tr} \left( {\bf Q} {\bf W} {\bf J}_N^{l_2} {\bf H}^{T} {\bf J}_N^{l_1} {\bf H}^{T} {\bf W}^{*} ({\bf I} \otimes {\bf R} {\bf J}_L^{u_r}) \right) = C(u_r,l_1) \delta(l_2 = u_r-l_1) + \mathcal{O}(\frac{L}{MN})
$$ 
Therefore, $\overline{s}_3(u_r)$ is equal to 
$$
\overline{s}_3(u_r) = \sum_{l_1} C(u_r,l_1) \kappa^{(r+1)}(u_1, \ldots, u_{r-1}, l_1, u_r - l_1) + 
\sum_{l_1,l_2}  \kappa^{(r+1)}(u_1, \ldots, u_{r-1}, l_1, l_2) \mathcal{O}(\frac{L}{MN})
$$

We now evaluate $s_5(u_r)$. For this, we recall that we denote $\beta_{1,0}(i,u_s,l_1,u_r)$ the term
$$
\beta_{1,0}(i,u_s,l_1,u_r) = \frac{1}{ML} \mathrm{Tr} \left( {\bf Q} ({\bf I}_M \otimes {\bf J}_L^{i}) {\bf Q} 
({\bf I}_M \otimes {\bf J}_L^{u_s}) {\bf Q} {\bf W} {\bf J}_N^{i} {\bf H}^{T} {\bf J}_N^{l_1} {\bf H}^{T} {\bf W}^{*}  \left( {\bf I}_M \otimes {\bf R}  {\bf J}_L^{u_r} \right) \right)
$$
We notice that 
$$
s_5(u_r,l_1) = s_{5,1}(u_r,l_1) + s_{5,2}(u_r,l_1)
$$
where 
$$
s_{5,1}(u_r,l_1) = -\frac{\sigma^{2}}{MLN} \sum_{s=1}^{r-1} \sum_{i=-(L-1)}^{L-1}  \mathbb{E} \left[ \left( \Pi_{t \neq (s,r)} \, \tau^{(M)}({\bf Q}^{\circ})(u_t) \, 
 \tau^{(M)}({\bf Q}^{\circ})(l_1)  \right) \; \beta_{1,0}(i,u_s,l_1,u_r)  \right] 
$$
and 
$$
s_{5,2}(u_r,l_1) = - \frac{\sigma^{2}}{MLN}  \sum_{i=-(L-1)}^{L-1}  \mathbb{E} \left[ \Pi_{t \leq r-1} \, \tau^{(M)}({\bf Q}^{\circ})(u_t)  \, \beta_{1,0}(i,l_1,l_1,u_r)   \right] 
$$
We first evaluate $s_{5,1}(u_r) = \sum_{l_1} s_{5,1}(u_r,l_1)$. We express $\beta_{1,0}(i,u_s,l_1,u_r)$ as 
$$
\beta_{1,0}(i,u_s,l_1,u_r) = \mathbb{E} \left( \beta_{1,0}(i,u_s,l_1,u_r) \right) + \beta_{1,0}(i,u_s,l_1,u_r)^{\circ}
$$
and notice that $s_{5,1}(u_r,l_1) = \overline{s}_{5,1}(u_r,l_1) + \tilde{s}_{5,1}(u_r,l_1)$ where
$$
\overline{s}_{5,1}(u_r,l_1) = -\frac{\sigma^{2}}{MLN} \sum_{i=-(L-1)}^{L-1} \sum_{s=1}^{r-1} \kappa^{(r-1)}((u_t)_{t \neq (s,r)}, l_1) \, \mathbb{E} \left( \beta_{1,0}(i,u_s,l_1,u_r) \right)
$$
and 
$$
\tilde{s}_{5,1}(u_r,l_1) = -\frac{\sigma^{2}}{MLN} \sum_{s=1}^{r-1} \sum_{i=-(L-1)}^{L-1}  \mathbb{E} \left[ \Pi_{t \neq s,r} \, \tau^{(M)}({\bf Q}^{\circ})(u_t) \, 
 \tau^{(M)}({\bf Q}^{\circ})(l_1)  \, \beta_{1,0}(i,u_s,l_1,u_r)^{\circ}  \right] 
$$
It is clear that $\tilde{s}_{5,1}(u_r,l_1) = \mathcal{O}(\frac{1}{(MN)^{(r+2)/2}})$ which implies that
\begin{equation}
\label{eq:eval-tilde-s51} 
\sum_{l_1} \tilde{s}_{5,1}(u_r,l_1) = \frac{L}{\sqrt{MN}} \mathcal{O}(\frac{1}{(MN)^{(r+1)/2}}) = o(\frac{1}{(MN)^{(r+1)/2}})
\end{equation} 

Proposition \ref{prop:expre-terme-utile-generaux-trace-Q-Q-QW} implies that
\begin{equation}
\label{eq:expre-theta}
\mathbb{E} \left( \beta_{1,0}(i,u_s,l_1,u_r) \right) = \overline{\beta}_{1,0}(i,u_s,l_1) \delta(l_1 = u_r + u_s) + \mathcal{O}(\frac{L}{MN})
\end{equation}
Using the rough evaluation $\kappa^{(r-1)}((u_t)_{t \neq s,r}, l_1) =  \mathcal{O}(\frac{1}{(MN)^{(r-1)/2}})$, 
we get immediately that 
\begin{equation}
\label{eq:eval-overline-s51} 
\overline{s}_{5,1}(u_r) = \sum_{l_1} \overline{s}_{5,1}(u_r,l_1) = \mathcal{O}(\frac{1}{(MN)^{(r+1)/2}})
\end{equation} 
We finally notice that if $r=2$, $\overline{s}_{5,1}(u_r)$ is reduced to $0$.

We define $\tilde{s}_{5,2}(u_r,l_1)$ and $\overline{s}_{5,2}(u_r,l_1)$ in the same way, and
obtain easily that
\begin{equation}
\label{eq:eval-tilde-s52} 
\sum_{l_1} \tilde{s}_{5,2}(u_r,l_1) = \frac{L}{\sqrt{MN}} \mathcal{O}(\frac{1}{(MN)^{(r+1)/2}}) = o(\frac{1}{(MN)^{(r+1)/2}})
\end{equation} 
The behaviour of $\sum_{l_1} \overline{s}_{5,2}(u_r,l_1)$ is however different from the behaviour of
$\sum_{l_1} \overline{s}_{5,1}(u_r,l_1)$ if $u_r = 0$. Indeed, 
$$
\mathbb{E} \left( \beta_{1,0}(i,l_1,l_1,u_r) \right) = \overline{\beta}_{1,0}(i,l_1,l_1) \delta(u_r=0) + \mathcal{O}(\frac{L}{MN})
$$
It is easy to check that the contribution of the $\mathcal{O}(\frac{L}{MN})$ terms to 
$\sum_{l_1} \overline{s}_{5,2}(u_r,l_1)$ is 
a $o\left( \frac{1}{(MN)^{(r+1)/2}} \right)$ term. Therefore, 
\begin{eqnarray}
\nonumber
\overline{s}_{5,2}(u_r) & = & \sum_{l_1} \overline{s}_{5,2}(u_r,l_1) = \sum_{l_1} \left( \frac{1}{L} \sum_{i=-(L-1)}^{L-1} \overline{\beta}_{1,0}(i,l_1,l_1)  \right) \, 
\frac{1}{MN} \, \kappa^{(r-1)}(u_1, \ldots, u_{r-1}) \, \delta(u_r=0) + o\left( \frac{1}{(MN)^{(r+1)/2}} \right) \\
& = &   \kappa^{(r-1)}(u_1, \ldots, u_{r-1}) \, \mathcal{O}(\frac{L}{MN}) \, \delta(u_r=0) \, +  o\left( \frac{1}{(MN)^{(r+1)/2}} \right)
\label{eq:eval-overline-s52-1}
\end{eqnarray}
As above, $\overline{s}_{5,2}(u_r)$ is reduced to $0$ if $r=2$. \\

The reader may check that the terms $s_2(u_r) = \overline{s}_2(u_r) + \tilde{x}_2^{(1)}(u_r)$ and $s_4(u_r)$
have exactly the same behaviour than $s_3(u_r)$ and $s_5(u_r)$. For the reader's convenience, we mention 
that $\tilde{x}_2^{(1)}(u_r)$ is defined as 
$$
\tilde{x}_2^{(1)}(u_r) = \sum_{l_1,l_2} \tilde{x}_2^{(1)}(u_r,l_1,l_2)
$$
where $\tilde{x}_2^{(1)}(u_r,l_1,l_2)$ is the term given by
\begin{equation}
\label{eq:def-tildex2-1}
\sigma^{2} c_N  \mathbb{E} \left[ \Pi_{s=1}^{r-1}
 \tau^{(M)}({\bf Q}^{\circ})(u_s)  \tau^{(M)}({\bf Q}^{\circ})(l_1)  \tau^{(M)}({\bf Q}^{\circ})(l_2) 
\frac{1}{ML} \mathrm{Tr} \left( {\bf Q} {\bf W} {\bf J}_N^{l_2} {\bf H}^{T}  {\bf W}^{*} ({\bf I} \otimes 
\sigma^{2} {\bf R} \mathcal{T}_{L,L}({\bf H} {\bf J}_N^{*l_1} {\bf H}) {\bf R} {\bf J}_L^{u_r}) \right)^{\circ} \right]
\end{equation}
In sum, we have proved the following useful result. 
\begin{proposition}
\label{prop:eval-x}
If $r \geq 2$, for each $u_r$, it holds that 
\begin{eqnarray}
\label{eq:eval-x}
x(u_r) & = & \sum_{l_1} C(u_r,l_1) \kappa^{(r+1)}(u_1, \ldots, u_{r-1}, l_1, u_r - l_1) + 
\sum_{l_1,l_2}  \kappa^{(r+1)}(u_1, \ldots, u_{r-1}, l_1, l_2) \mathcal{O}(\frac{L}{MN})  \\
\nonumber
&   & +  \kappa^{(r-1)}(u_1, \ldots, u_{r-1}) \, \mathcal{O}(\frac{L}{MN}) \, \delta(u_r=0)  
+ \tilde{x}_2^{(1)}(u_r) + \tilde{x}_3^{(1)}(u_r)  + \mathcal{O}(\frac{1}{(MN)^{(r+1)/2}}) 
\end{eqnarray}
while if $r=2$, 
\begin{eqnarray}
\label{eq:eval-x-r=2}
x(u_2)  & =  & \sum_{l_1} C(u_2,l_1) \kappa^{(3)}(u_1, l_1, u_2 - l_1) + 
\sum_{l_1,l_2}  \kappa^{(3)}(u_1, l_1, l_2) \mathcal{O}(\frac{L}{MN})  \\
&  &  + \tilde{x}_2^{(1)}(u_2) + \tilde{x}_3^{(1)}(u_2) + \mathcal{O}(\frac{L}{(MN)^{2}})
\end{eqnarray}
\end{proposition}
(\ref{eq:premiere-evaluation-general}) thus leads to 
the Proposition:
\begin{proposition}
\label{prop:evaluation-kappa-r}
For $r \geq 2$, it holds that 
\begin{multline}
\label{eq:evaluation-kappa-r}
\left|\kappa^{(r)}(u_1, \ldots, u_r) - \frac{ {\bf y}^{*}_{1,u_r}}{1 - d(u_r,z)} \right|  \leq  
| \kappa^{(r-1)}(u_1, \ldots, u_{r-1}) | \, \mathcal{O}(\frac{L}{MN}) \;  +\\
C \, \sup_{u} \sum_{l_1} |\kappa^{(r+1)}(u_1, \ldots, u_{r-1}, l_1, u - l_1)|  \;  + 
\sum_{l_1,l_2}  |\kappa^{(r+1)}(u_1, \ldots, u_{r-1}, l_1, l_2)| \mathcal{O}(\frac{L}{MN})  \; + \\
\sup_{u} |\tilde{x}_2^{(1)}(u)|  + \sup_{u} |\tilde{x}_3^{(1)}(u)|   + \mathcal{O}(\frac{1}{(MN)^{(r+1)/2}}) 
\end{multline}
while for $r=2$,
\begin{multline} 
\label{eq:evaluation-kappa-r=2}
\left| \kappa^{(2)}(u_1,u_2) - \frac{ {\bf y}^{*}_{1,u_2}}{1 - d(u_2,z)}  \right|  \leq  
 C \, \sup_{u} \sum_{l_1} |\kappa^{(3)}(u_1, l_1, u - l_1)|  \; + \\
\sum_{l_1,l_2}  |\kappa^{(3)}(u_1, l_1, l_2)| \mathcal{O}(\frac{L}{MN})   \; + \; 
\sup_{u} |\tilde{x}_2^{(1)}(u)|  + \sup_{u} |\tilde{x}_3^{(1)}(u)|   + \mathcal{O}(\frac{L}{(MN)^{2}})
\end{multline}
\end{proposition}
We now establish Proposition \ref{prop:premiere-reccurence} introduced into the 
overview of the proof of Proposition \ref{prop:difficile}.  
\begin{proposition}
\label{prop:premiere-reccurence}
For each $r \geq 2$ and for each integer $p$, $1 \leq p \leq r-1$, if integers $u_1, \ldots, u_r \in \{ -(L-1), \ldots, L-1 \}$ satisfy
\begin{equation}
\label{eq:conditions-p-u}
\begin{array}{ccc} u_r + u_s & \neq & 0 \; s=1, \ldots, r-1 \\
                   u_{r-1} + u_s & \neq & 0 \; s=1, \ldots, r-2 \\
                   \vdots & \vdots & \vdots \\
                   u_{r-p+1} + u_s & \neq & 0 \; s=1, \ldots, r-p 
\end{array}
\end{equation}
then, it holds that  
\begin{equation}
\label{eq:evaluation-kappa-reccurence-1}
\kappa^{r}(u_1, \ldots, u_r) = \max \left( \left(\frac{L}{\sqrt{MN}}\right)^{p}, \frac{1}{\sqrt{MN}} \right) \; \mathcal{O}(\frac{1}{(MN)^{r/2}})
\end{equation}
\end{proposition}
We prove the proposition by induction on $r$. We first check (\ref{eq:evaluation-kappa-reccurence-1}) if $r=2$. In this case, the integer $p$ is necessarily 
equal to $1$ and (\ref{eq:conditions-p-u}) reduces to $u_1+u_2 \neq 0$. We use (\ref{eq:evaluation-kappa-r=2}). 
Using the rough evaluations $\kappa^{(3)}(v_1,v_2,v_3) = \mathcal{O}(\frac{1}{(MN)^{3/2}})$ 
and $\sup_{u} |\tilde{x}_j^{(1)}(u)| = \mathcal{O}(\frac{L^{2}}{(MN)^{2}}) = (\frac{L}{\sqrt{MN}})^{2} \mathcal{O}(\frac{1}{MN})$ for $j=2,3$, 
we obtain immediately that (\ref{eq:evaluation-kappa-reccurence-1}) holds if $r=2$. \\

We now assume that (\ref{eq:evaluation-kappa-reccurence-1}) holds until integer $r_0 -1$ and prove that it 
is true for integer $r_0$. For this, we establish that for each $r \geq r_0$ 
and for each $u_1, \ldots, u_r$, (\ref{eq:evaluation-kappa-reccurence-1}) holds 
provided (\ref{eq:conditions-p-u}) is true until $p \leq r_0 -1$. 
 We first verify that (\ref{eq:evaluation-kappa-reccurence-1}) holds
for each $r \geq r_0$ and for $p=1$ as soon as $u_r + u_s  \neq  0 \; s=1, \ldots, r-1$. For this, 
we use (\ref{eq:evaluation-kappa-r}). ${\bf y}^{*}_{1,u_r}$ is of course equal to  $0$. Moreover, as $\kappa^{(r-1)}(u_1, \ldots, u_{r-1}) = 
\mathcal{O}(\frac{1}{(MN)^{(r-1)/2}})$, it is clear that
$$
|\kappa^{(r-1)}(u_1, \ldots, u_{r-1}) | \, \mathcal{O}(\frac{L}{MN}) = \frac{L}{\sqrt{MN}} \mathcal{O}(\frac{1}{(MN)^{r/2}})
$$ 
as expected. Using that $\kappa^{(r+1)}(v_1, \ldots, v_{r+1}) = 
\mathcal{O}(\frac{1}{(MN)^{(r+1)/2}})$ for each $(v_1, \ldots, v_{r+1})$, we obtain immediately that 
$$
\sup_{u} \sum_{l_1} |\kappa^{(r+1)}(u_1, \ldots, u_{r-1}, l_1, u - l_1)| =  \frac{L}{\sqrt{MN}} \mathcal{O}(\frac{1}{(MN)^{r/2}})
$$ 
and 
$$
\sum_{l_1,l_2}  |\kappa^{(r+1)}(u_1, \ldots, u_{r-1},l_1, l_2)| \mathcal{O}(\frac{L}{MN}) = \frac{L^{2}}{MN} \frac{L}{\sqrt{MN}} \mathcal{O}(\frac{1}{(MN)^{r/2}})
$$
Finally, the Hölder inequality leads to 
\begin{equation}
\label{eq:premiere-evaluation-tildex}
\sup_{u} |\tilde{x}_2^{(1)}(u)|  + \sup_{u} |\tilde{x}_3^{(1)}(u)| = \mathcal{O}(\frac{L^{2}}{(MN)^{(r+2)/2}}) = 
(\frac{L}{\sqrt{MN}})^{2} \mathcal{O}(\frac{1}{(MN)^{r/2}})
\end{equation}
Next, we consider the case $p=2$ for the reader's convenience. We consider $r \geq r_0$, and assume that
$u_r + u_s  \neq  0 \; s=1, \ldots, r-1$ as well as $u_{r-1} + u_s  \neq  0 \; s=1, \ldots, r-2$.  
We again use (\ref{eq:evaluation-kappa-r}) and remark that ${\bf y}^{*}_{1,u_r} = 0$. 
As $u_{r-1} + u_s  \neq  0 \; s=1, \ldots, r-2$, 
the use of (\ref{eq:evaluation-kappa-reccurence-1}) for integer $r-1$, multi-index  
$(u_1, \ldots, u_{r-1})$ and $p=1$ (proved above) implies that 
$\kappa^{(r-1)}(u_1, \ldots, u_{r-1}) = 
\frac{L}{\sqrt{MN}}  \mathcal{O}(\frac{1}{(MN)^{(r-1)/2}})$, and that 
$$
\kappa^{(r-1)}(u_1, \ldots, u_{r-1})  \, \mathcal{O}(\frac{L}{MN}) = (\frac{L}{\sqrt{MN}})^{2} \mathcal{O}(\frac{1}{(MN)^{r/2}})
$$  
We now evaluate $\sum_{l_1}  |\kappa^{(r+1)}(u_1, \ldots, u_{r-1}, l_1, u - l_1)|$. 
It is clear that 
$$
\kappa^{(r+1)}(u_1, \ldots, u_{r-1}, l_1, u - l_1) = \kappa^{(r+1)}(l_1, u - l_1, u_1, \ldots, u_{r-1})
$$
As $u_{r-1} + u_s  \neq  0 \; s=1, \ldots, r-2$, 
the use of (\ref{eq:evaluation-kappa-reccurence-1}) for integer $r+1$, multi-index  
$(l_1,u-l_1,u_1, \ldots, u_{r-1})$ and $p=1$  leads to  
$$
\kappa^{(r+1)}(u_1, \ldots, u_{r-1}, l_1, u - l_1) = \frac{L}{\sqrt{MN}} \mathcal{O}(\frac{1}{(MN)^{(r+1)/2}})
$$ 
as soon as $u_{r-1} + l_1 \neq 0$ and $u_{r-1} + u - l_1 \neq 0$, or equivalently 
if $l_1 \neq -u_{r-1}$ and $l_1 \neq u + u_{r-1}$. Therefore, 
$$
\sum_{l_1 \neq (-u_{r-1}, u+u_{r-1})} | \kappa^{(r+1)}(u_1, \ldots, u_{r-1}, l_1, u - l_1) |= 
(\frac{L}{\sqrt{MN}})^{2} \mathcal{O}(\frac{1}{(MN)^{r/2}})
$$
If $l_1 = -u_{r-1}$ or $l_1 = u + u_{r-1}$, we use the rough evaluation 
$$
 \kappa^{(r+1)}(u_1, \ldots, u_{r-1}, l_1, u - l_1) = (\frac{1}{\sqrt{MN}}) \mathcal{O}(\frac{1}{(MN)^{r/2}})
$$
Therefore, we obtain that 
$$
\sum_{l_1} | \kappa^{(r+1)}(u_1, \ldots, u_{r-1}, l_1, u - l_1) | = 
\max \left( \left(\frac{L}{\sqrt{MN}}\right)^{2}, \frac{1}{\sqrt{MN}} \right) \; \mathcal{O}(\frac{1}{(MN)^{r/2}})
$$
We now consider $\sum_{l_1,l_2}  | \kappa^{(r+1)}(u_1, \ldots, u_{r-1},l_1, l_2) | \mathcal{O}(\frac{L}{MN})$. We remark that \\
$\kappa^{(r+1)}(u_1, \ldots, u_{r-1},l_1, l_2) = \kappa^{(r+1)}(l_1,l_2,u_1, \ldots, u_{r-1})$
Therefore, if $u_{r-1} + l_1 \neq 0$ and $u_{r-1} + l_2 \neq 0$, (\ref{eq:evaluation-kappa-reccurence-1}) for integer $r+1$, multi-index  
$(l_1,l_2,u_1, \ldots, u_{r-1})$ and $p=1$ implies that 
$$
\kappa^{(r+1)}(u_1, \ldots, u_{r-1}, l_1, l_2) = \frac{L}{\sqrt{MN}} \mathcal{O}(\frac{1}{(MN)^{(r+1)/2}})
$$
If $l_1 = -u_{r-1}$ or $l_2 = -u_{r-1}$, we use again that 
$$
 \kappa^{(r+1)}(u_1, \ldots, u_{r-1}, l_1, l_2) = (\frac{1}{\sqrt{MN}}) \mathcal{O}(\frac{1}{(MN)^{r/2}})
$$
so that for $i,j = 1,2$, $i \neq j$, it holds that 
$$
\sum_{l_i = -u_{r-1}, l_j} | \kappa^{(r+1)}(u_1, \ldots, u_{r-1},l_i, l_j) | \mathcal{O}(\frac{L}{MN}) =
\frac{L^{2}}{MN} (\frac{1}{\sqrt{MN}}) \mathcal{O}(\frac{1}{(MN)^{r/2}}) = o(\frac{1}{(MN)^{(r+1)/2}})
$$
We finally obtain that 
$$
\sum_{l_1,l_2}  | \kappa^{(r+1)}(u_1, \ldots, u_{r-1},l_1, l_2) | \mathcal{O}(\frac{L}{MN}) = \max \left( \left(\frac{L}{\sqrt{MN}} \right)^{4}, \frac{1}{\sqrt{MN}} \right) \, \mathcal{O}(\frac{1}{(MN)^{r/2}})
$$
Finally, the first evaluation (\ref{eq:premiere-evaluation-tildex}) of $\tilde{x}_2^{(1)}(u_r)$ and 
$\tilde{x}_3^{(1)}(u_r)$ establishes (\ref{eq:evaluation-kappa-reccurence-1}) 
for each $r \geq r_0$ and for $p=2$ if $u_r + u_s \neq 0$ for $s=1, \ldots, r-1$ and 
$u_{r-1} + u_s \neq 0$ for $s=1, \ldots, r-2$. 

In order to complete the proof of (\ref{eq:evaluation-kappa-reccurence-1}) for each $r \geq r_0$ and 
for each $p \leq r_0 - 1$, we assume that (\ref{eq:evaluation-kappa-reccurence-1}) 
holds for each $r \geq r_0$ and for each $p \leq p_0$ where $p_0 \leq r_0 -2$, and prove 
that it also holds for $p=p_0+1$. For this, we establish the following Lemma. 
\begin{lemma}
\label{le:aide-a-la-comprehension}
Assume that for each $t \geq r_0 - 1$ and for each integer $p$, $1 \leq p \leq p_0 \leq r_0 - 2$, 
it holds that 
\begin{equation}
\label{eq:enonce-v}
\kappa^{(t)}(v_1, \ldots, v_t) = \max \left( \left(\frac{L}{\sqrt{MN}}\right)^{p}, \frac{1}{\sqrt{MN}} \right) \; \mathcal{O}(\frac{1}{(MN)^{t/2}})
\end{equation}
for each multi-index $(v_1, \ldots, v_t)$ satisfying 
\begin{equation}
\label{eq:conditions-p-u-v}
\begin{array}{ccc} v_t + v_s & \neq & 0 \; s=1, \ldots, t-1 \\
                   v_{t-1} + v_s & \neq & 0 \; s=1, \ldots, t-2 \\
                   \vdots & \vdots & \vdots \\
                   v_{t-p+1} + u_s & \neq & 0 \; s=1, \ldots, t-p 
\end{array}
\end{equation} 
Then, for each $r \geq r_0$ and for each multi-index $(u_1, \ldots, u_r)$ 
satisfying (\ref{eq:conditions-p-u}) for $p = p_0+1$, it holds that \\
$\kappa^{(r-1)}(u_1, \ldots, u_{r-1}) \mathcal{O}(\frac{L}{MN})$, 
$\sum_{l_1} | \kappa^{(r+1)}(u_1, \ldots, u_{r-1}, l_1, u - l_1)|$, 
$\sum_{l_1,l_2}  |\kappa^{(r+1)}(u_1, \ldots, u_{r-1},l_1, l_2) | \mathcal{O}(\frac{L}{MN})$, 
$\sup_{u} |\tilde{x}_j^{(1)}(u)|$ for $j=2,3$ 
 are $\max \left( \left(\frac{L}{\sqrt{MN}}\right)^{(p_0+1)}, \frac{1}{\sqrt{MN}} \right) \; \mathcal{O}(\frac{1}{(MN)^{r/2}})$ terms. 
\end{lemma}
Using (\ref{eq:evaluation-kappa-r}), (\ref{eq:evaluation-kappa-reccurence-1}) for $p=p_0+1$
follows immediately from Lemma \ref{le:aide-a-la-comprehension}. Consequently, (\ref{eq:evaluation-kappa-reccurence-1}) holds for each $r \geq r_0$ until index $p \leq (r_0-1)$, and in particular for $r=r_0$ and $p \leq (r_0-1)$. This completes the proof of Proposition \ref{prop:premiere-reccurence}. \\

{\bf Proof of Lemma \ref{le:aide-a-la-comprehension}.} We consider a multi-index $(u_1, \ldots, u_r)$ satisfying (\ref{eq:conditions-p-u}) for $p = p_0+1$ and remark that it verifies
\begin{equation}
\label{eq:conditions-p-1-u}
\begin{array}{ccc} u_{r-1} + u_s & \neq & 0 \; s=1, \ldots, r-2 \\
                   u_{r-2} + u_s & \neq & 0 \; s=1, \ldots, r-3 \\
                   \vdots & \vdots & \vdots \\
                   u_{r-p_0} + u_s & \neq & 0 \; s=1, \ldots, r-p_0-1 
\end{array}
\end{equation}
Therefore, (\ref{eq:enonce-v}) used for $t=r-1$, $p=p_0$ and multi-index $(v_1, \ldots, v_{r-1})$ with $v_s = u_s$
leads to 
$$
\kappa^{(r-1)}(u_1, \ldots, u_{r-1}) = \max \left( \left(\frac{L}{\sqrt{MN}}\right)^{p_0}, \frac{1}{\sqrt{MN}} \right) \mathcal{O}(\frac{1}{(MN)^{(r-1)/2}})
$$
Therefore, 
$$
\kappa^{(r-1)}(u_1, \ldots, u_{r-1}) \mathcal{O}(\frac{L}{MN}) = \frac{L}{\sqrt{MN}}  \max \left( \left(\frac{L}{\sqrt{MN}}\right)^{p_0}, \frac{1}{\sqrt{MN}} \right) \mathcal{O}(\frac{1}{(MN)^{r/2}})
$$
which, of course, also coincides with a $\max \left( \left(\frac{L}{\sqrt{MN}}\right)^{p_0+1}, \frac{1}{\sqrt{MN}} \right) \mathcal{O}(\frac{1}{(MN)^{r/2}})$ term. We now study the term 
$$
\sum_{l_1} | \kappa^{(r+1)}(u_1, \ldots, u_{r-1},l_1, u-l_1) |
$$
Using (\ref{eq:enonce-v}) for $t=r+1$ and multi-index $l_1, u-l_1, u_1, \ldots, u_{r-1}$, we obtain that 
$$
\kappa^{(r+1)}(u_1, \ldots, u_{r-1}, l_1, u - l_1) = \max \left( \left(\frac{L}{\sqrt{MN}}\right)^{p_0}, \frac{1}{\sqrt{MN}} \right) \mathcal{O}(\frac{1}{(MN)^{(r+1)/2}})
$$
if $l_1$ is such that $u_{r-j} + l_1 \neq 0$ and $u_{r-j} + u - l_1 \neq 0$ for each $j=1, \ldots, p_0$. The sum 
of the terms $| \kappa^{(r+1)}(u_1, \ldots, u_{r-1}, l_1, u - l_1)|$ over these values of $l_1$ 
is therefore a \\ $L \max \left( \left(\frac{L}{\sqrt{MN}}\right)^{p_0}, \frac{1}{\sqrt{MN}} \right) \mathcal{O}(\frac{1}{(MN)^{(r+1)/2}})$ term, or equivalently a $\frac{L}{\sqrt{MN}} \left( \left(\frac{L}{\sqrt{MN}}\right)^{p_0}, \frac{1}{\sqrt{MN}} \right) \mathcal{O}(\frac{1}{(MN)^{r/2}})$ term, which, of course, is also a 
$$ \max \left( \left(\frac{L}{\sqrt{MN}}\right)^{p_0+1}, \frac{1}{\sqrt{MN}} \right) \mathcal{O}(\frac{1}{(MN)^{r/2}})$$
term. If $l_1$ is equal to $-u_{r-j_0}$ or to $u_{r-j_0} + u$ for some $j_0 =1, \ldots, p_0$, we use the rough evaluation 
$$
\kappa^{(r+1)}(u_1, \ldots, u_{r-1}, l_1, u - l_1) = \mathcal{O}(\frac{1}{(MN)^{(r+1)/2}}) = \frac{1}{\sqrt{MN}} \mathcal{O}(\frac{1}{(MN)^{r/2}})
$$
This discussion implies that 
$$
\sum_{l_1} | \kappa^{(r+1)}(u_1, \ldots, u_{r-1}, l_1, u - l_1)| = \max \left( \left(\frac{L}{\sqrt{MN}}\right)^{p_0+1}, \frac{1}{\sqrt{MN}} \right) \mathcal{O}(\frac{1}{(MN)^{r/2}})
$$
The evaluation of $\sum_{l_1,l_2}  |\kappa^{(r+1)}(u_1, \ldots, u_{r-1},l_1, l_2) | \mathcal{O}(\frac{L}{MN})$ is similar and is thus omitted. \\

In order to complete the proof of Lemma \ref{le:aide-a-la-comprehension}, it remains to prove that
that 
$$
\sup_{u} |\tilde{x}_j^{(1)}(u)| \leq \max \left( \left(\frac{L}{\sqrt{MN}}\right)^{(p_0+1)}, \frac{1}{\sqrt{MN}} \right) \; \mathcal{O}(\frac{1}{(MN)^{r/2}})
$$
for $j=2,3$. For this, we study in more details $\sup_{u} |\tilde{x}_j^{(1)}(u)|$
for $j=2,3$. We expand $\tilde{x}_j^{(1)}(u,l_1,l_2)$ using (\ref{eq:expre-produit-tau-rond-W}) 
when $r$ is exchanged by $r+1$. In the same way than $\tilde{x}(u)$ defined by (\ref{eq:def-tildexu}), it holds that
$$
\tilde{x}_j^{(1)}(u) = \sum_{i=1}^{5} s_{j,i}^{(1)}(u)
$$
where the terms $(s_{j,i}^{(1)}(u))_{i=1, \ldots, 5}$ are defined in the same way than 
$(s_i(u))_{i=1, \ldots, 5}$. We define $\tilde{x}_{j,i}^{(2)}(u)$ for $i=2,3$ by the fact that
$$
s_{j,i}^{(1)}(u) = \overline{s}_{j,i}^{(1)}(u) + \tilde{x}_{j,i}^{(2)}(u)
$$
We define $\tilde{x}^{(1)}(u)$ as the positive term 
$$
\tilde{x}^{(1)}(u) = |\tilde{x}_2^{(1)}(u)| + |\tilde{x}_3^{(1)}(u)|
$$
and, similarly, $\tilde{x}^{(2)}(u)$ is given by 
$$
\tilde{x}^{(2)}(u) = \sum_{(i,j)=(2,3)} |\tilde{x}_{j,i}^{(2)}(u)|
$$

A rough evaluation (based on the Hölder inequality and on (\ref{eq:expre-terme-utile-generaux-trace-Q-Q-QW})) of the various terms $s_{j,i}^{(1)}(u)$ for $i=4,5$ leads to 
$s_{j,i}^{(1)}(u) =  \frac{L}{\sqrt{MN}} O(\frac{1}{(MN)^{(r+1)/2}})$. After some calculations, we obtain that  
\begin{multline}
\label{eq:evaluation-tildex1}
\tilde{x}^{(1)}(u)  \leq  \sum_{l_1,l_2} |\kappa^{(r+1)}(u_1, \ldots, u_{r-1},l_1,l_2)| \, \mathcal{O}(\frac{L}{MN}) + \\
C \, \sum_{l_1,l_2}  |\kappa^{(r+2)}(u_1, \ldots, u_{r-1}, l_1, l_2, u-l_1-l_2)| + \\
 \sum_{l_1,l_2,l_3} |\kappa^{(r+2)}(u_1, \ldots, u_{r-1},l_1,l_2,l_3)| \,  \mathcal{O}(\frac{L}{MN})  + 
\tilde{x}^{(2)}(u) +  \frac{L}{\sqrt{MN}} O(\frac{1}{(MN)^{(r+1)/2}})
\end{multline}
The first term of the righthandside of (\ref{eq:evaluation-tildex1}) corresponds to the contribution
of $s_{j,1}^{(1)}(u)$ while the second and the third terms are due to $\overline{s}_{j,2}^{(1)}(u)$ and $\overline{s}_{j,3}^{(1)}(u)$.
The term $\frac{L}{\sqrt{MN}} O(\frac{1}{(MN)^{(r+1)/2}})$ is due to the $s_{j,i}^{(1)}(u)$ for $i=4,5$.   
The terms $\tilde{x}_{j,i}^{(2)}(u)$ can of course be also developed and we obtain similarly 
\begin{multline}
\label{eq:evaluation-tildex2}
\tilde{x}^{(2)}(u)  \leq  \sum_{l_1,l_2,l_3} |\kappa^{(r+2)}(u_1, \ldots, u_{r-1},l_1,l_2,l_3)| \, \mathcal{O}(\frac{L}{MN}) + \\
C \,  \sum_{l_1,l_2,l_3} |\kappa^{(r+3)}(u_1, \ldots, u_{r-1}, l_1, l_2, l_3, u-l_1-l_2-l_3)|  + \\
 \sum_{l_1,l_2,l_3,l_4} |\kappa^{(r+3)}(u_1, \ldots, u_{r-1},l_1,l_2,l_3,l_4)| \,  \mathcal{O}(\frac{L}{MN})  + 
\tilde{x}^{(3)}(u) + (\frac{L}{\sqrt{MN}})^{2} O(\frac{1}{(MN)^{(r+1)/2}}) 
\end{multline}
The term $(\frac{L}{\sqrt{MN}})^{2} O(\frac{1}{(MN)^{(r+1)/2}})$ is due to the terms $(s^{(2)}_{k_1,k_2,i}(u))$ for $i=4,5$
and $k_1,k_2 = 2,3$: it is easily seen using the Hölder inequality that their order of magnitude is $\frac{L}{\sqrt{MN}}$ smaller than the order of magnitude of the $(s^{(1)}_{k,i})_{i=4,5}$ for $k=2,3$. More generally, it holds that
\begin{multline}
\label{eq:evaluation-tildexq}
\tilde{x}^{(q)}(u)  \leq  \sum_{l_i,i=1, \ldots,q+1} |\kappa^{(r+q)}(u_1, \ldots, u_{r-1},l_i,i=1, \ldots,q+1)| \, \mathcal{O}(\frac{L}{MN}) + \\
C \,  \sum_{l_i,i=1, \ldots,q+1} |\kappa^{(r+q+1)}(u_1, \ldots, u_{r-1},l_i,i=1, \ldots,q+1, u-\sum_{i=1}^{q+1} l_i)|  + \\
 \sum_{l_i,i=1, \ldots,q+2} |\kappa^{(r+q+1)}(u_1, \ldots, u_{r-1},l_i,i=1, \ldots,q+2)| \,  \mathcal{O}(\frac{L}{MN})  + \tilde{x}^{(q+1)}(u) + (\frac{L}{\sqrt{MN}})^{q} O(\frac{1}{(MN)^{(r+1)/2}})
\end{multline} 
We remark that the Hölder inequality leads to 
\begin{equation}
\label{eq:schwartz-tildexp}
\sup_{u} \, \tilde{x}^{(p)}(u) = \left( \frac{L}{\sqrt{MN}} \right)^{p+1} \, \mathcal{O}(\frac{1}{(MN)^{r/2}})
\end{equation}
for each $p$. We express $\tilde{x}^{(1)}(u)$ as
\begin{equation}
\label{eq:astuce-tildex1}
\tilde{x}^{(1)}(u) = \sum_{q=1}^{p_0-1} (\tilde{x}^{(q)}(u) - \tilde{x}^{(q+1)}(u)) + \tilde{x}^{(p_0)}(u)
\end{equation}
We now prove that for each $q$, then it holds that 
\begin{equation}
\label{eq:tildexq-tildexq+1}
\tilde{x}^{(q)}(u) - \tilde{x}^{(q+1)}(u) \leq \left(\frac{L}{\sqrt{MN}}\right)^{q} \max \left( ( \frac{L}{\sqrt{MN}})^{p_0+1}, \frac{1}{\sqrt{MN}} \right)  \, \mathcal{O}(\frac{1}{(MN)^{r/2}})
\end{equation}
(\ref{eq:evaluation-tildexq}) implies that $\tilde{x}^{(q)}(u) - \tilde{x}^{(q+1)}(u)$ is upperbounded by the sum of 
4 terms. We just study the second term, i.e.
$$
\sum_{l_i,i=1, \ldots,q+1} |\kappa^{(r+q+1)}(u_1, \ldots, u_{r-1},l_i,i=1, \ldots,q+1, u-\sum_{i=1}^{q+1} l_i)|
$$
because, as the fourth term $(\frac{L}{\sqrt{MN}})^{q} O(\frac{1}{(MN)^{(r+1)/2}})$, it can be easily checked that the first and the third term are negligible w.r.t. the righthandside of inequality (\ref{eq:tildexq-tildexq+1}). If the integers $l_1, \ldots, l_{q+1}, u-\sum_{i=1}^{q+1} l_i$
do not belong $\{ -u_{r-1}, \ldots, -u_{r-p_0} \}$, (\ref{eq:enonce-v}) for $t=r+q+1$ 
and for multi-index $(l_1, \ldots, l_{q+1}, u-\sum_{i=1}^{q+1} l_i, u_1, \ldots, u_{r-1})$ 
implies that 
$$
\kappa^{(r+q+1)}(u_1, \ldots, u_{r-1},l_i,i=1, \ldots,q+1,u-\sum_{i=1}^{q+1} l_i) = \max \left( \left(\frac{L}{\sqrt{MN}}\right)^{(p_0)}, \frac{1}{\sqrt{MN}} \right) \; \mathcal{O}(\frac{1}{(MN)^{(r+q+1)/2}})
$$
Therefore, the sum over all these integers can be upperbounded by 
$$
L^{q+1} \, \max \left( \left(\frac{L}{\sqrt{MN}}\right)^{(p_0)}, \frac{1}{\sqrt{MN}} \right) \; \mathcal{O}(\frac{1}{(MN)^{(r+q+1)/2}}) = \left(\frac{L}{\sqrt{MN}}\right)^{q} \, \frac{L}{\sqrt{MN}} \, \max \left( ( \frac{L}{\sqrt{MN}})^{p_0}, \frac{1}{\sqrt{MN}} \right)  \, \mathcal{O}(\frac{1}{(MN)^{r/2}})
$$
which, of course, is a $\left(\frac{L}{\sqrt{MN}}\right)^{q} \max \left( ( \frac{L}{\sqrt{MN}})^{p_0+1}, \frac{1}{\sqrt{MN}} \right)  \, \mathcal{O}(\frac{1}{(MN)^{r/2}})$ term as expected. 

If at least one of the index $l_1, \ldots, l_{q+1}, u-\sum_{i=1}^{q+1} l_i$ is equal an integer 
$(-u_{r-i})_{i=1, \ldots, p_0}$, we use the rough evaluation
$$
\kappa^{(r+q+1)}(u_1, \ldots, u_{r-1},l_i,i=1, \ldots,q+1,u-\sum_{i=1}^{q+1} l_i) = 
\mathcal{O}(\frac{1}{(MN)^{(r+q+1)/2}})
$$
The sum over the corresponding multi-indices is thus a $L^{q} \mathcal{O}(\frac{1}{(MN)^{(r+q+1)/2}}) = 
\left( \frac{L}{\sqrt{MN}} \right)^{q} \, \mathcal{O}(\frac{1}{(MN)^{(r+1)/2}})$. This completes 
the proof of (\ref{eq:tildexq-tildexq+1}). Therefore, (\ref{eq:astuce-tildex1}) and (\ref{eq:schwartz-tildexp}) imply 
that
$$
\sup_{u} \, \tilde{x}^{(1)}(u) = \max \left( ( \frac{L}{\sqrt{MN}})^{p_0+1}, \frac{1}{\sqrt{MN}} \right)  \, \mathcal{O}(\frac{1}{(MN)^{r/2}})
$$
as expected. This, in turn, completes the proof of Lemma \ref{le:aide-a-la-comprehension}. \\ 

We now improve the evaluation of Proposition \ref{prop:premiere-reccurence} 
when $(u_1, \ldots, u_r)$ satisfy $u_t + u_s \neq 0$ for $1 \leq t,s \leq r$ and $t \neq s$, 
or equivalently if $(u_1, \ldots, u_r)$ verify (\ref{eq:conditions-p-u}) for $p=r-1$. 
More precisely, we prove the following result. 
\begin{proposition}
\label{prop:deuxieme-reccurence}
Assume that  $(u_1, \ldots, u_r)$ satisfy $u_t + u_s \neq 0$ for $1 \leq t,s \leq r$ and $t \neq s$. Then, 
for each $q \geq 1$, for each $r \geq 2$, it holds that 
\begin{equation}
\label{eq:evaluation-kappa-reccurence-2}
\kappa^{(r)}(u_1, \ldots, u_r) = \max \left( \left(\frac{L}{\sqrt{MN}}\right)^{r-1+q}, \frac{1}{\sqrt{MN}} \right) \; \mathcal{O}(\frac{1}{(MN)^{r/2}})
\end{equation}
\end{proposition}
{\bf Proof.} We prove this result by induction on integer $q$. We first establish 
(\ref{eq:evaluation-kappa-reccurence-2}) for $q=1$ by induction on integer $r$. 
If $r=2$, we have to check that if $u_1+u_2 \neq 0$, then it holds that  
\begin{equation}
\label{eq:premiere-evaluation-amelioree-r=2}
\kappa^{(2)}(u_1, u_2) = \max \left( \left(\frac{L}{\sqrt{MN}}\right)^{2}, \frac{1}{\sqrt{MN}} \right) \; \mathcal{O}(\frac{1}{MN})
\end{equation}
For this, we use (\ref{eq:evaluation-kappa-r=2}). We have already mentioned that the Hölder inequality leads to
$$
\sup_{u} \, \tilde{x}^{(1)}(u) = \left(\frac{L}{\sqrt{MN}}\right)^{2} \mathcal{O}(\frac{1}{MN})
$$
We study the term 
$$
\sup_{u}  \sum_{l_1} |\kappa^{(3)}(u_1,l_1,u-l_1)|
$$
Proposition \ref{prop:premiere-reccurence} in the case $r=3$ and $p=1$ implies that 
$$
 \kappa^{(3)}(u_1,l_1,u-l_1) = \frac{L}{\sqrt{MN}} \mathcal{O}(\frac{1}{(MN)^{3/2}})
$$
as soon as $l_1 \neq -u_1$ and $l_1 \neq u+u_1$. Therefore, 
$$
 \sum_{l_1 \neq (-u_1, u+u_1)} |\kappa^{(3)}(u_1,l_1,u-l_1)| =  \left(\frac{L}{\sqrt{MN}}\right)^{2} \mathcal{O}(\frac{1}{MN})
$$
If $l_1 = -u_1$ or $l_1 = u+u_1$, we use the rough evaluation $\kappa^{(3)}(u_1,l_1,u-l_1) = \mathcal{O}(\frac{1}{(MN)^{3/2}})$, and we finally obtain that 
$$
 \sum_{l_1} |\kappa^{(3)}(u_1,l_1,u-l_1)| =  \max \left( \left(\frac{L}{\sqrt{MN}}\right)^{2}, \frac{1}{\sqrt{MN}} \right) \; \mathcal{O}(\frac{1}{MN})
$$
as expected. The term 
$$
\sum_{l_1,l_2} |\kappa^{(3)}(u_1,l_1,l_2)| \mathcal{O}(\frac{L}{MN})
$$
is evaluated similarly. We have thus established (\ref{eq:premiere-evaluation-amelioree-r=2}). 

We assume that (\ref{eq:evaluation-kappa-reccurence-2}) holds for $q=1$ until index $r_0 - 1$
and prove that it also holds for index $r_0$. We take (\ref{eq:evaluation-kappa-r}) as a starting point. 
We consider $(u_1, \ldots, u_{r_0})$ satisfying $u_t + u_s \neq 0$ for $1 \leq t,s \leq r_0$, 
or equivalently (\ref{eq:conditions-p-u}) for $r=r_0$ and $p=r_0-1$. (\ref{eq:evaluation-kappa-reccurence-2}) for $q=1$, $r=r_0 - 1$ and multi-index $(u_1, \ldots, u_{r_0-1})$ leads to 
$$
\kappa^{(r_0-1)}(u_1, \ldots,u_{r_0-1}) = \max \left( \left(\frac{L}{\sqrt{MN}}\right)^{r_0-2+1}, \frac{1}{\sqrt{MN}} \right) \; \mathcal{O}(\frac{1}{(MN)^{(r_0-1)/2}}) 
$$
and to 
$$
\kappa^{(r_0-1)}(u_1, \ldots,u_{r_0-1}) \mathcal{O}(\frac{L}{MN})= \frac{L}{\sqrt{MN}} \max \left( \left(\frac{L}{\sqrt{MN}}\right)^{r_0-1}, \frac{1}{\sqrt{MN}} \right) \; \mathcal{O}(\frac{1}{(MN)^{r_0/2}}) 
$$
which, of course, is a $\max \left( \left(\frac{L}{\sqrt{MN}}\right)^{r_0-1+1}, \frac{1}{\sqrt{MN}} \right) \; \mathcal{O}(\frac{1}{(MN)^{r_0/2}})$ term as expected. We now evaluate 
$$
\sum_{l_1}  |\kappa^{(r_0+1)}(u_1, \ldots, u_{r_0-1},l_1, u - l_1)|
$$
If $l_1 + u_{s} \neq 0$ and $u - l_1 + u_s \neq 0$ for $s = 1, \ldots, r_0 - 1$, Proposition \ref{prop:premiere-reccurence} for $r=r_0+1$, multi-index $(l_1, u-l_1, u_1, \ldots, u_{r-1})$  and $p = r_0 -1$ implies that 
$$
\kappa^{(r_0+1)}(u_1, \ldots, u_{r_0-1},l_1, u - l_1) = \max \left( \left(\frac{L}{\sqrt{MN}}\right)^{r_0-1}, \frac{1}{\sqrt{MN}} \right) \; \mathcal{O}(\frac{1}{(MN)^{(r_0+1)/2}}) 
$$
and that the sum of the $| \kappa^{(r_0+1)}(u_1, \ldots, u_{r_0-1},l_1, u - l_1)|$
over these indices is a 
$$
\frac{L}{\sqrt{MN}} \max \left( \left(\frac{L}{\sqrt{MN}}\right)^{r_0-1}, \frac{1}{\sqrt{MN}} \right) \; \mathcal{O}(\frac{1}{(MN)^{r_0/2}})
$$ 
term. If $l_1 + u_s=0$ or $u - l_1 + u_s = 0$ for some integer $s$, 
we use as previously that 
$$
\kappa^{(r_0+1)}(u_1, \ldots, u_{r_0-1},l_1, u - l_1) =  \frac{1}{\sqrt{MN}} \mathcal{O}(\frac{1}{(MN)^{r_0/2}}) 
$$
This, in turn, implies that 
$$
\sup_{u} \left| \sum_{l_1} \kappa^{(r_0+1)}(u_1, \ldots, u_{r_0-1},l_1, u - l_1) \right| = 
\max \left( \left(\frac{L}{\sqrt{MN}}\right)^{r_0-1+1}, \frac{1}{\sqrt{MN}} \right) \; \mathcal{O}(\frac{1}{(MN)^{r_0/2}})
$$
as expected.

The term 
$$
\sum_{l_1,l_2} |\kappa^{(r_0+1)}(u_1, \ldots, u_{r_0-1}, l_1,l_2)| \, \mathcal{O}(\frac{L}{MN}) 
$$
can be evaluated similarly. Finally, it is easy to show as in the proof of Lemma \ref{le:aide-a-la-comprehension} that $\sup_{u} \, \tilde{x}^{(1)}(u)$ behaves as expected. 

This completes the proof of (\ref{eq:evaluation-kappa-reccurence-2}) for each $r$ and $q=1$. In order to 
establish the proposition for each $q$, we assume that it is true until integer $q-1$ and prove 
that it holds for integer $q$. We prove this statement by induction on integer $r$, and begin to consider
$r=2$. We of course use (\ref{eq:evaluation-kappa-r=2}) for $u_1+u_2 \neq 0$. It is easy to check as previously that the term $\sup_{u} \, \tilde{x}^{(1)}(u)$ is as expected, and that it is also the case for
$\sum_{l_1,l_2} |\kappa^{(3)}(u_1,l_1,l_2)| \mathcal{O}(\frac{L}{MN})$. However, the term 
$$
\sup_{u}  \sum_{l_1} | \kappa^{(3)}(u_1,l_1,u-l_1)| 
$$
appears more difficult to evaluate. If $u \neq 0$, 
it is easy to check that $\sum_{l_1} |\kappa^{(3)}(u_1,l_1,u-l_1)|$ is 
a $$\max\left( \left(\frac{L}{\sqrt{MN}}\right)^{1+q}, \frac{1}{\sqrt{MN}} \right) \mathcal{O}(\frac{1}{MN})$$
term because, except if $u_1+l_1 = 0$ or $u_1+u - l_1 = 0$ (the contribution of these particular values 
to the sum is a $\mathcal{O}(\frac{1}{(MN)^{3/2}})$ term), (\ref{eq:evaluation-kappa-reccurence-2}) 
used for $r=3$ and integer $q-1$ implies that 
$$
\kappa^{(3)}(u_1,l_1,u-l_1) = \max\left( \left(\frac{L}{\sqrt{MN}}\right)^{1+q}, \frac{1}{\sqrt{MN}} \right) \mathcal{O}(\frac{1}{(MN)^{3/2}})
$$
and that 
$$
\sum_{l_1 \neq -u_1, u_1+u} | \kappa^{(3)}(u_1,l_1,u-l_1) | = \frac{L}{\sqrt{MN}} \max\left( \left(\frac{L}{\sqrt{MN}}\right)^{1+q}, \frac{1}{\sqrt{MN}} \right) \mathcal{O}(\frac{1}{(MN)})
$$
If $u=0$, the sum becomes 
$$
 \sum_{l_1}  |\kappa^{(3)}(u_1,l_1,-l_1)| 
$$
(\ref{eq:evaluation-kappa-reccurence-2}) for $r=3$ and integer $q-1$ cannot be used to 
evaluate $\kappa^{(3)}(u_1,l_1,-l_1)$ because $l_1 - l_1 = 0$. We have thus to study separately this kind of term. 
For this, we prove the following lemma. 
\begin{lemma}
\label{le:kappa(l,-l)}
We consider an integer $r \geq 2$ and assume the following hypotheses: 
\begin{itemize}
\item for each integer $s$ and for each $v_1, \ldots, v_s$ such that $v_{s_1} + v_{s_2} \neq 0$, $1 \leq s_1,s_2 \leq s$, $s_1 \neq s_2$, it holds that 
\begin{equation}
\label{eq:hypothese-1-recurrence-l-l}
\kappa^{(s)}(v_1, \ldots, v_s) = \max\left( \left(\frac{L}{\sqrt{MN}}\right)^{s-1+q-1}, \frac{1}{\sqrt{MN}} \right) \mathcal{O}(\frac{1}{(MN)^{s/2}})
\end{equation}
\item for each $s \leq r-1$, and each $v_1, \ldots, v_s$ such that $v_{s_1} + v_{s_2} \neq 0$, $1 \leq s_1,s_2 \leq s$, $s_1 \neq s_2$, it holds that 
\begin{equation}
\label{eq:hypothese-2-recurrence-l-l}
\kappa^{(s)}(v_1, \ldots, v_s) = \max\left( \left(\frac{L}{\sqrt{MN}}\right)^{s-1+q}, \frac{1}{\sqrt{MN}} \right) \mathcal{O}(\frac{1}{(MN)^{s/2}})
\end{equation}
\end{itemize}
Then, if $u_1, \ldots, u_{r-1}$ verify $u_{s_1} + u_{s_2} \neq 0$, $1 \leq s_1,s_2 \leq r-1$, $s_1 \neq s_2$, it holds that 
\begin{equation}
\label{eq:comportement-somme-l-l}
\sum_{l_1} |\kappa^{(r+1)}(u_1, \ldots, u_{r-1},l_1,-l_1)| = \max\left( \left(\frac{L}{\sqrt{MN}}\right)^{r-1+q}, \frac{1}{\sqrt{MN}} \right) \mathcal{O}(\frac{1}{(MN)^{r/2}})
\end{equation}
\end{lemma}
{\bf Proof.}
We evaluate $\kappa^{(r+1)}(u_1, \ldots, u_{r-1},l_1,-l_1)$ using (\ref{eq:evaluation-kappa-r}) when $r$ is replaced by $r+1$ and for multi-index $(u_1, \ldots, u_{r-1},l_1,-l_1)$. If $l_1 = \pm u_s$ for some $s$, 
the term $\kappa^{(r+1)}(u_1, \ldots, u_{r-1},l_1,-l_1)$ is a $\mathcal{O}(\frac{1}{(MN)^{(r+1)/2}})$. It is thus sufficient to prove (\ref{eq:comportement-somme-l-l}) when the sum is over the integers $l_1$ 
that do not belong to $\{-u_1, \ldots, -u_{r-1} \}$ and $\{u_1, \ldots, u_{r-1} \}$. In order to simplify
the notations, we do not mention in the following that the sum does not take into account $\{-u_1, \ldots, -u_{r-1} \}$ and $\{u_1, \ldots, u_{r-1} \}$. 

If $l_1$ does not belong to $\{-u_1, \ldots, -u_{r-1} \}$ and $\{u_1, \ldots, u_{r-1} \}$, 
component $-l_1$ of vector ${\bf y}^{*}_{1}$ corresponding to ${\bs \kappa} = (\kappa^{r+1}(u_1, \ldots, u_{r-1},l_1,u))_{u=-(L-1), \ldots, (L-1)}$ can be written as
$$
{\bf y}^{*}_{1,-l_1} = \kappa^{(r-1)}(u_1, \ldots, u_{r-1}) \, \mathcal{O}(\frac{1}{MN})
$$
(see (\ref{eq:def-y1*})). Therefore, for $l_1 \neq \pm u_s$, 
$s=1, \ldots, r-1$, 
(\ref{eq:evaluation-kappa-r})) implies that 
\begin{multline}
\label{eq:evaluation-kappa-l-l}
|\kappa^{(r+1)}(u_1, \ldots, u_{r-1},l_1,-l_1)| \leq |\kappa^{(r-1)}(u_1, \ldots, u_{r-1})| \, \mathcal{O}(\frac{1}{MN})
+ |\kappa^{(r)}(u_1, \ldots, u_{r-1},l_1)| \, \mathcal{O}(\frac{L}{MN}) +\\
C \, \sup_{u} \sum_{l_2} |\kappa^{(r+2)}(u_1, \ldots, u_{r-1},l_1,l_2,u-l_2)| +
\left| \sum_{l_2,l_3} \kappa^{(r+2)}(u_1, \ldots, u_{r-1},l_1,l_2,l_3) \, \mathcal{O}(\frac{L}{MN}) \right| + \\
\sup_{u} |\tilde{x}_{2,l_1}^{(1)}(u)| + \sup_{u} |\tilde{x}_{3,l_1}^{(1)}(u)| + \mathcal{O}(\frac{1}{(MN)^{(r+2)/2}})
\end{multline}
where we indicate that the terms $\tilde{x}_j^{(1)}(u)$ associated to $(u_1, \ldots, u_{r-1},l_1,u)$
depend on $l_1$ (these terms also depend on $(u_s)_{s \leq r-1}$ but it is not useful to mention
this dependency). In the following, we denote by $\alpha^{(1)}(u_1, \ldots, u_{r-1})$ the term
$$
\alpha^{(1)}(u_1, \ldots, u_{r-1}) = \sum_{l_1} |\kappa^{(r+1)}(u_1, \ldots, u_{r-1},l_1,-l_1)|
$$
(\ref{eq:evaluation-kappa-l-l}) implies that 
\begin{multline}
\label{eq:evaluation-somme-kappa-l-l}
 \alpha^{(1)}(u_1, \ldots, u_{r-1}) \leq |\kappa^{(r-1)}(u_1, \ldots, u_{r-1})| \, \mathcal{O}(\frac{L}{MN})
+ \sum_{l_1} |\kappa^{(r)}(u_1, \ldots, u_{r-1},l_1)| \, \mathcal{O}(\frac{L}{MN}) +\\
C \, \sup_{u} \sum_{l_1,l_2} |\kappa^{(r+2)}(u_1, \ldots, u_{r-1},l_1,l_2,u-l_2)| +
\left| \sum_{l_1,l_2,l_3} \kappa^{(r+2)}(u_1, \ldots, u_{r-1},l_1,l_2,l_3) \, \mathcal{O}(\frac{L}{MN}) \right| + \\
\sup_{u} \sum_{l_1} |\tilde{x}_{2,l_1}^{(1)}(u)| + \sup_{u} \sum_{l_1} |\tilde{x}_{3,l_1}^{(1)}(u)| + \mathcal{O}(\frac{L}{(MN)^{(r+2)/2}})
\end{multline}
(\ref{eq:hypothese-2-recurrence-l-l}) for $s=r-1$ implies that 
$$
|\kappa^{(r-1)}(u_1, \ldots, u_{r-1})| = \max\left( \left(\frac{L}{\sqrt{MN}}\right)^{r-2+q}, \frac{1}{\sqrt{MN}} \right) \mathcal{O}(\frac{1}{(MN)^{(r-1)/2}}) 
$$ 
and that 
$$
|\kappa^{(r-1)}(u_1, \ldots, u_{r-1})| \mathcal{O}(\frac{L}{MN}) = \frac{L}{\sqrt{MN}} \max\left( \left(\frac{L}{\sqrt{MN}}\right)^{r-2+q}, \frac{1}{\sqrt{MN}} \right) \mathcal{O}(\frac{1}{(MN)^{r/2}}) 
$$
which, of course, is also a $\max\left( \left(\frac{L}{\sqrt{MN}}\right)^{r-1+q}, \frac{1}{\sqrt{MN}} \right) \mathcal{O}(\frac{1}{(MN)^{r/2}})$ term as expected. In order to evaluate the second term of the righthandside 
of (\ref{eq:evaluation-somme-kappa-l-l}), we first notice that if $l_1 \in \{ -u_1, \ldots, -u_{r-1} \}$, 
the Hölder inequality leads to 
$$
|\kappa^{(r)}(u_1, \ldots, u_{r-1},l_1)| \, \mathcal{O}(\frac{L}{MN}) = o(\frac{1}{(MN)^{(r+1)/2}})
$$ 
If $l_1 + u_s \neq 0$ for each $s=1, \ldots, r-1$,  
we use (\ref{eq:hypothese-1-recurrence-l-l}) 
for $s=r$ and $(v_1, \ldots, v_r) = (u_1, \ldots, u_{r-1}, l_1)$. 
It holds that 
$$
\kappa^{(r)}(u_1, \ldots, u_{r-1},l_1) = \max\left( \left(\frac{L}{\sqrt{MN}}\right)^{r-1+q-1}, \frac{1}{\sqrt{MN}} \right) \mathcal{O}(\frac{1}{(MN)^{r/2}}) 
$$ 
so that 
$$
\sum_{l_1 \neq -u_s,s=1, \ldots, r-1} |\kappa^{(r)}(u_1, \ldots, u_{r-1},l_1)| \, \mathcal{O}(\frac{L}{MN}) = (\frac{L}{\sqrt{MN}})^{2} 
\max\left( \left(\frac{L}{\sqrt{MN}}\right)^{r-1+q-1}, \frac{1}{\sqrt{MN}} \right) \mathcal{O}(\frac{1}{(MN)^{r/2}}) 
$$
which is a $\max\left( \left(\frac{L}{\sqrt{MN}}\right)^{r-1+q}, \frac{1}{\sqrt{MN}} \right) \mathcal{O}(\frac{1}{(MN)^{r/2}})$ term. The fourth term of the righthandside of (\ref{eq:evaluation-somme-kappa-l-l}) is evaluated
similarly. Moreover, following the arguments used to establish  Lemma \ref{le:aide-a-la-comprehension}, it can be
shown that 
$$ 
\sup_{u} \sum_{l_1} |\tilde{x}_{2,l_1}^{(1)}(u)| + \sup_{u} \sum_{l_1} |\tilde{x}_{3,l_1}^{(1)}(u)| = 
\max\left( \left(\frac{L}{\sqrt{MN}}\right)^{r-1+q}, \frac{1}{\sqrt{MN}} \right) \mathcal{O}(\frac{1}{(MN)^{r/2}}) 
$$
It remains to evaluate the third term of the righhandside of (\ref{eq:evaluation-somme-kappa-l-l}). The supremum
over $u \neq 0$ is as expected, but the term corresponding to $u=0$ has also to be evaluated. 
We denote $\alpha^{(2)}(u_1, \ldots, u_{r-1})$ the term
$$
\alpha^{(2)}(u_1, \ldots, u_{r-1}) = \sum_{l_1,l_2} |\kappa^{(r+2)}(u_1, \ldots, u_{r-1},l_1,l_2,-l_2)|
$$
The previous discussion implies that 
$$
\alpha^{(1)}(u_1, \ldots, u_{r-1}) \leq C \, \alpha^{(2)}(u_1, \ldots, u_{r-1}) + \max\left( \left(\frac{L}{\sqrt{MN}}\right)^{r-1+q}, \frac{1}{\sqrt{MN}} \right) \mathcal{O}(\frac{1}{(MN)^{r/2}}) 
$$
It can be shown similarly that
$$
\alpha^{(2)}(u_1, \ldots, u_{r-1}) \leq C \, \alpha^{(3)}(u_1, \ldots, u_{r-1}) + \max\left( \left(\frac{L}{\sqrt{MN}}\right)^{r-1+q}, \frac{1}{\sqrt{MN}} \right) \mathcal{O}(\frac{1}{(MN)^{r/2}}) 
$$
where 
$$
\alpha^{(3)}(u_1, \ldots, u_{r-1}) = \sum_{l_1,l_2,l_3} |\kappa^{(r+3)}(u_1, \ldots, u_{r-1},l_1,l_2,l_3,-l_3)|
$$
More generally, if $\alpha^{(p)}(u_1, \ldots, u_{r-1})$ is defined by 
$$
\alpha^{(p)}(u_1, \ldots, u_{r-1}) = \sum_{l_i,i=1,\ldots,p} |\kappa^{(r+p)}(u_1, \ldots, u_{r-1},(l_i,i=1,\ldots,p),-l_p)|
$$
it holds that 
$$
\alpha^{(p-1)}(u_1, \ldots, u_{r-1}) \leq C \, \alpha^{(p)}(u_1, \ldots, u_{r-1}) + \max\left( \left(\frac{L}{\sqrt{MN}}\right)^{r-1+q}, \frac{1}{\sqrt{MN}} \right) \mathcal{O}(\frac{1}{(MN)^{r/2}}) 
$$
and consequently that
\begin{equation}
\label{eq:alpha-final}
\alpha^{(1)}(u_1, \ldots, u_{r-1}) \leq C \, \alpha^{(p)}(u_1, \ldots, u_{r-1}) + \max\left( \left(\frac{L}{\sqrt{MN}}\right)^{r-1+q}, \frac{1}{\sqrt{MN}} \right) \mathcal{O}(\frac{1}{(MN)^{r/2}}) 
\end{equation}
The Hölder inequality leads immediately to
$$
\alpha^{(p)}(u_1, \ldots, u_{r-1}) = (\frac{L}{\sqrt{MN}})^{p} \mathcal{O}(\frac{1}{(MN)^{r/2}})
$$
and choosing $p=r-1+q$ provides (\ref{eq:comportement-somme-l-l}). \\

We finally complete the proof of Proposition \ref{prop:deuxieme-reccurence}. The use of Lemma
\ref{le:kappa(l,-l)} for $r=2$ establishes immediately that if (\ref{eq:evaluation-kappa-reccurence-2}) holds until integer $q-1$ for each $s$, then, it also holds for integer $q$ and $r=2$. We assume that
(\ref{eq:evaluation-kappa-reccurence-2}) holds for integer $q$ until integer $r-1$, i.e. 
that both (\ref{eq:hypothese-1-recurrence-l-l}) and (\ref{eq:hypothese-2-recurrence-l-l}) hold, and prove 
that it also holds for integer $r$, i.e. that 
$$
\kappa^{(r)}(u_1, \ldots, u_r) =  \max\left( \left(\frac{L}{\sqrt{MN}}\right)^{r-1+q}, \frac{1}{\sqrt{MN}} \right) \mathcal{O}(\frac{1}{(MN)^{r/2}}) 
$$
For this, we use (\ref{eq:evaluation-kappa-r}). All the terms of the righthandside 
of (\ref{eq:evaluation-kappa-r}) are easily seen to be as expected, except the 
second one. However, Lemma \ref{le:kappa(l,-l)} implies that the second term is also a 
$\max\left( \left(\frac{L}{\sqrt{MN}}\right)^{r-1+q}, \frac{1}{\sqrt{MN}} \right) \mathcal{O}(\frac{1}{(MN)^{r/2}})$. 
This completes the proof of Proposition \ref{prop:deuxieme-reccurence}. \\

We are now in position to establish (\ref{eq:estimation-kappa-r})
\begin{corollary}
\label{cor:preuve-{eq:estimation-kappa-r}}
If $(u_1, \ldots, u_r)$ satisfy $u_t+u_s \neq 0$ for $t \neq s$, $1 \leq t,s \leq r$, then
(\ref{eq:estimation-kappa-r}) holds for $r \geq 2$. 
\end{corollary}
{\bf Proof.} As $L = N^{\alpha}$ with $\alpha < 2/3$, it exists an integer 
$q_0$ for which $(\frac{L}{\sqrt{MN}})^{r-1+q_0} = o(\frac{1}{\sqrt{MN}})$. 
Therefore, 
$$
\max \left( \left(\frac{L}{\sqrt{MN}}\right)^{r-1+q_0}, \frac{1}{\sqrt{MN}} \right) = \frac{1}{\sqrt{MN}}
$$
(\ref{eq:evaluation-kappa-reccurence-2}) for $q=q_0$ thus implies (\ref{eq:estimation-kappa-r}). \\

It remains to establish (\ref{eq:estimation-kappa}). For this, we take  (\ref{eq:evaluation-kappa-r=2}) as a starting point, and prove that the righthandside of (\ref{eq:evaluation-kappa-r=2}) is 
a $\mathcal{O}(\frac{L}{(MN)^{2}})$ term. We first justify that:
\begin{equation}
\label{eq:tildex-r=2}
\sup_{u} \tilde{x}^{(1)}(u) = \mathcal{O}(\frac{L}{(MN)^{2}})
\end{equation}
We use the decomposition (\ref{eq:astuce-tildex1}) of $\tilde{x}^{(1)}(u)$ for the following convenient 
value of $p$: we recall that the Hölder inequality implies that 
$$
\tilde{x}^{(p)}(u) = (\frac{L}{\sqrt{MN}})^{p+1} \mathcal{O}(\frac{1}{MN})
$$
As $L = N^{\alpha}$ with $\alpha < 2/3$, it exists $p$ for which 
$$
(\frac{L}{\sqrt{MN}})^{p+1} = o(\frac{L}{MN}) 
$$
For such a value of $p$, it holds that 
$$
\tilde{x}^{(p)}(u) = o(\frac{L}{(MN)^{2}})
$$
Using (\ref{eq:evaluation-tildex1}) for $r=2$ as well as (\ref{eq:estimation-kappa-r}), it is easy to check that 
$\tilde{x}^{(1)}(u) - \tilde{x}^{(2)}(u)$ is a $\mathcal{O}(\frac{L}{(MN)^{2}})$ term, and that the same holds true 
for $\tilde{x}^{(q)}(u) - \tilde{x}^{(q+1)}(u)$ for each $q \geq 1$. This establishes (\ref{eq:tildex-r=2}). 

(\ref{eq:estimation-kappa-r}) implies that the second term of the righhandside of (\ref{eq:evaluation-kappa-r=2}) is a $\mathcal{O}\left( (\frac{L}{MN})^{3} \right) =  o(\frac{L}{(MN)^{2}}))$ term. It remains to establish that 
\begin{equation}
\label{eq:sum-kappa3-l-l}
\sum_{l_1} |\kappa^{3}(u_1,l_1,-l_1)| = \mathcal{O}(\frac{L}{(MN)^{2}})
\end{equation}
Lemma \ref{le:kappa(l,-l)} for $q=q_0$ (where $q_0$ is defined in the proof of Corollary \ref{cor:preuve-{eq:estimation-kappa-r}} for $r=2$) implies that this term is $\mathcal{O}(\frac{1}{(MN)^{3/2}})$, but this evaluation 
is not sufficient to prove (\ref{eq:estimation-kappa}). Using (\ref{eq:estimation-kappa-r}), we now evaluate
$\kappa^{(r+1)}(u_1, \ldots, u_{r-1},l_1, -l_1)$ when $u_{s_1} + u_{s_2} \neq 0$ for $s_1 \neq s_2$, 
$-(L-1) \leq s_1,s_2 \leq L-1$ and for $r \geq 2$. 
\begin{lemma}
\label{le:evaluation-amelioree-kappa-l-l}
We consider $r \geq 2$ and a multi-index $(u_1, \ldots, u_{r-1},l_1, -l_1)$  such that 
$u_{s_1} + u_{s_2} \neq 0$ for $s_1 \neq s_2$, $-(L-1) \leq s_1,s_2 \leq L-1$. Then, 
\begin{itemize}
\item if $l_1 \pm u_s \neq 0$ for $s=1, \ldots, r-1$, it holds that 
\begin{equation}
\label{eq:evaluation-amelioree-kappa-l-l}
\kappa^{(r+1)}(u_1, \ldots, u_{r-1},l_1, -l_1) = \mathcal{O}(\frac{1}{(MN)^{(r+2)/2}})
\end{equation}
\item if $l_1 \pm u_s = 0$ for some $s=1, \ldots, r-1$,
\begin{equation}
\label{eq:evaluation-amelioree-kappa-l-l-l}
\kappa^{(r+1)}(u_1, \ldots, u_{r-1},l_1, -l_1) = \frac{L}{\sqrt{MN}} \mathcal{O}(\frac{1}{(MN)^{(r+1)/2}})
\end{equation}
\end{itemize}
\end{lemma}
{\bf Proof.} The proof is similar to the proof of Lemma \ref{le:kappa(l,-l)}. We take (\ref{eq:evaluation-kappa-l-l}) as a starting point, but just evaluate $\kappa^{(r+1)}(u_1, \ldots, u_{r-1},l_1, -l_1)$ instead of 
$\alpha^{(1)}(u_1, \ldots, u_{r-1})$ by iterating (\ref{eq:evaluation-kappa-l-l}). Using (\ref{eq:estimation-kappa-r}), it is easy to check that for each $l_1$, 
$$
\sup_{u} |\tilde{x}^{(1)}_{j,l_1}(u)| = \mathcal{O}(\frac{1}{(MN)^{(r+2)/2}})
$$
We first assume that  $l_1 \pm u_s \neq 0$ for $s=1, \ldots, r-1$. 
(\ref{eq:estimation-kappa-r}) implies that the first term of the righthandside of (\ref{eq:evaluation-kappa-l-l}) is  
$\mathcal{O}((\frac{1}{(MN)^{(r+2)/2}})$ (and is identically $0$ if $r=2$). The second term 
is $\frac{L}{\sqrt{MN}} \mathcal{O}((\frac{1}{(MN)^{(r+2)/2}})$ while the fourth term 
is $(\frac{L}{\sqrt{MN}})^{2} \mathcal{O}((\frac{1}{(MN)^{(r+2)/2}})$. The supremum over $u \neq 0$ of
the third term is $\mathcal{O}((\frac{1}{(MN)^{(r+2)/2}})$ which implies that 
$$
|\kappa^{(r+1)}(u_1, \ldots, u_{r-1},l_1, -l_1)| \leq \sum_{l_2} \left| \kappa^{(r+2)}(u_1, \ldots, u_{r-1},l_1,l_2,-l_2) \right| 
+ \mathcal{O}(\frac{1}{(MN)^{(r+2)/2}})
$$
As in the proof Lemma (\ref{le:kappa(l,-l)}), we iterate this inequality until an index $p$ 
for which 
$$
\sum_{l_2, \ldots, l_p} \left| \kappa^{(r+p)}(u_1, \ldots, u_{r-1},l_1,l_2, \ldots, l_p,-l_p) \right| = 
\mathcal{O}(\frac{L^{p-1}}{(MN)^{(r+p)/2}})
$$
is a $o((\frac{1}{(MN)^{(r+2)/2}})$ term. This, in turn, proves (\ref{eq:evaluation-amelioree-kappa-l-l}). 
(\ref{eq:evaluation-amelioree-kappa-l-l-l}) follows directly from the use of the Hölder inequality 
in  (\ref{eq:evaluation-kappa-l-l}). \\

We now complete the proof of (\ref{eq:sum-kappa3-l-l}). For this, we remark that
$$
\sum_{l_1} |\kappa^{3}(u_1,l_1,-l_1)| = \sum_{l_1 \neq \pm u_1}  |\kappa^{3}(u_1,l_1,-l_1)| + 2 |\kappa^{3}(u_1,-u_1,u_1)|
$$
Lemma \ref{le:evaluation-amelioree-kappa-l-l} implies that 
$$
\sum_{l_1 \neq \pm u_1}  |\kappa^{3}(u_1,l_1,-l_1)| = \mathcal{O}\left(\frac{L}{(MN)^{2}}\right)
$$
and that 
$$
|\kappa^{3}(u_1,-u_1,u_1)| = \mathcal{O}\left(\frac{L}{(MN)^{2}}\right)
$$
This establishes (\ref{eq:sum-kappa3-l-l}) as well as (\ref{eq:estimation-kappa}).

\subsection{Expansion of  $ \frac{1}{ML} \mathrm{Tr}\left( {\boldsymbol \Delta}(z) \right)$.}
\label{subsec:preuve-trace-delta}
In the following, we establish (\ref{eq:resultat-trace-Delta}). We recall that (\ref{eq:expre-gamma}) implies that $ \frac{1}{ML} \mathrm{Tr}\left( {\boldsymbol \Delta}(z) \right)$ is given by
$$
\frac{1}{ML} \mathrm{Tr}\left( {\boldsymbol \Delta}(z) \right) = \sigma^{2} c_N \sum_{l_1=-(L-1)}^{L-1} \mathbb{E} \left( \tau^{(M)}({\bf Q}^{\circ})(l_1) \, \frac{1}{ML} \mathrm{Tr}\left( {\bf Q} {\bf W} {\bf J}_N^{l_1} {\bf H}^{T} {\bf W}^{*}({\bf I}_M \otimes {\bf R}) \right)^{\circ} \right)
$$
In the following, we denote by $\tilde{x}(l_1)$ and $\tilde{x}$ the terms defined by  
$$
\tilde{x}(l_1) = \mathbb{E} \left( \tau^{(M)}({\bf Q}^{\circ})(l_1) \, \frac{1}{ML} \mathrm{Tr}\left( {\bf Q} {\bf W} {\bf J}_N^{l_1} {\bf H}^{T} {\bf W}^{*}({\bf I}_M \otimes {\bf R}) \right)^{\circ} \right)
$$
and 
$$
\tilde{x} = \sum_{l_1=-(L-1)}^{L-1} \mathbb{E} \left( \tau^{(M)}({\bf Q}^{\circ})(l_1) \, \frac{1}{ML} \mathrm{Tr}\left( {\bf Q} {\bf W} {\bf J}_N^{l_1} {\bf H}^{T} {\bf W}^{*}({\bf I}_M \otimes {\bf R}) \right)^{\circ} \right)
$$ 
$\tilde{x}(l_1)$ and $\tilde{x}$ appear to be formally similar to $\tilde{x}(0,l_1)$ and $\tilde{x}(0)$ defined by (\ref{eq:def-tildexul}) and (\ref{eq:def-tildexu}) in the particular case 
$r=1$. While we have considered in the previous subsection the case $r \geq 2$, a number of evaluations and 
results can be adapted to the easier case $r=1$. As in subsection \ref{subsec:preuve-difficile}, we expand 
$\tilde{x}(l_1)$ and $\tilde{x}$ using (\ref{eq:expre-produit-tau-rond-W}) in the case $r=1$, $v_1=l_1$, ${\bf G} = {\bf J}_N^{l_1} {\bf H}^{T}$ 
and ${\bf A} = ({\bf I}_M \otimes {\bf R})$. Using the same notations as in subsection \ref{subsec:preuve-difficile}, 
we obtain that
$$
\tilde{x}(l_1) = \sum_{j=2}^{5} s_j(l_1)
$$
and 
$$
\tilde{x} = \sum_{j=2}^{5} s_j
$$
where $s_j = \sum_{l_1} s_{j}(l_1)$. We note that the term $s_1$ is reduced to $0$ in the present context. 
It is easy to check that 
$$
s_4(l_1) = \frac{\sigma^{2}}{MN} \sum_{i=-(L-1)}^{L-1} \mathbb{E}\left(\beta_{0,1}(i,l_1,l_1,0)\right) 
$$
and that 
$$
s_5(l_1) = - \frac{\sigma^{2}}{MN} \sum_{i=-(L-1)}^{L-1} \mathbb{E}\left(\beta_{1,0}(i,l_1,l_1,0)\right) 
$$
where the terms $\beta$ are defined by (\ref{eq:def-beta}). Proposition \ref{prop:expre-terme-utile-generaux-trace-Q-Q-QW}
immediately implies that 
$$
s_4 = \frac{\sigma^{2}}{MN} \, \sum_{l_1} \left( \frac{1}{L} \sum_{i} \overline{\beta}_{0,1}(i,l_1,l_1) \right)  
+ \mathcal{O}\left( (\frac{L}{MN})^{2} \right)
$$
or equivalently, 
$$
s_4(z) = \sigma^{2} \frac{L}{MN} \overline{\beta}_{0,1}(z) + \, \mathcal{O}(\frac{L^{2}}{(MN)^{2}})
$$
where $\overline{\beta}_{0,1}(z)$ is defined as 
$$
\overline{\beta}_{0,1}(z) = \frac{1}{L^{2}} \sum_{l_1,i} \overline{\beta}_{0,1}(i,l_1,l_1)(z)
$$
Similarly, it holds that 
$$
s_5(z) = - \sigma^{2} \frac{L}{MN} \overline{\beta}_{1,0}(z) + \, \mathcal{O}(\frac{L^{2}}{(MN)^{2}})
$$
where 
$$
\overline{\beta}_{1,0}(z) = \frac{1}{L^{2}} \sum_{l_1,i} \overline{\beta}_{1,0}(i,l_1,l_1)(z)
$$
We have now to evaluate $s_2(z)$ and $s_3(z)$. For $j=2,3$, $s_j$ can be written as
$$
s_j = \overline{s}_j + \tilde{x}_j^{(1)}
$$
We first evaluate $\overline{s}_3$ and $\overline{s}_2$. 
$\overline{s}_3$ is equal to 
$$
\overline{s}_3 = -\sigma^{2} c_N \sum_{l_1,l_2} \kappa^{(2)}(l_1,l_2) \mathbb{E} \left[ \frac{1}{ML} \mathrm{Tr}\left({\bf Q} {\bf W} {\bf J}_N^{l_2} {\bf H}^{T} {\bf J}_N^{l_1} {\bf H}^{T} {\bf W}^{*}({\bf I}_M \otimes {\bf R}) \right) \right]
$$
We remark that 
$$
\mathbb{E} \left[ \frac{1}{ML} \mathrm{Tr}\left({\bf Q} {\bf W} {\bf J}_N^{l_2} {\bf H}^{T} {\bf J}_N^{l_1} {\bf H}^{T} {\bf W}^{*}({\bf I}_M \otimes {\bf R}) \right) \right] = - \sigma^{2} t(z)^{2} (z \tilde{t}(z))^{3} (1 - |l_1|/N) \delta(l_1+l_2=0) + \mathcal{O}(\frac{L}{MN})
$$
We also have to evaluate $\kappa^{(2)}(l_1,l_2)$. Using (\ref{eq:def-y1*}), 
(\ref{eq:evaluation-kappa-r=2}), and the observation that the righthandside of 
(\ref{eq:evaluation-kappa-r=2}) is a $\mathcal{O}(\frac{L}{(MN)^{2}})$ term (see (\ref{eq:estimation-kappa})), we obtain 
that 
$$
\kappa^{(2)}(l_1,l_2) = \frac{\sigma^{2}}{MN} \frac{1}{1 - d(l_1,z)} \, \frac{1}{L} \sum_{i} \overline{\beta}(i,l_1) \, \delta(l_1+l_2=0) \, + \, \mathcal{O}(\frac{L}{(MN)^{2}})
$$
Therefore, $\overline{s}_3$ can be written as 
$$
\overline{s}_3(z) = \sigma^{6} c_N t(z)^{2} (z \tilde{t}(z))^{3} \frac{1}{L^{2}} \sum_{i,l_1}  \frac{1-|l_1|/N)}{1 - d(l_1,z)} 
\overline{\beta}(i,l_1) \, \frac{L}{MN} \, + \, \mathcal{O}(\frac{L^{2}}{(MN)^{2}})
$$
Similar calculations lead to 
$$
\overline{s}_2 = \sigma^{8} c_N t(z)^{3} (z \tilde{t}(z))^{4} \frac{1}{L^{2}} \sum_{i,l_1}  \frac{(1-|l_1|/N)^{2} (1 - |l_1|/L)}{1 - d(l_1,z)} 
\overline{\beta}(i,l_1) \, \frac{L}{MN} \, + \, \mathcal{O}(\frac{L^{2}}{(MN)^{2}})
$$
Therefore, it holds that 
\begin{equation}
\label{eq:yes!}
\overline{s}_2(z) + \overline{s}_3(z) + s_4(z) + s_5(z) = \frac{L}{MN} \, \frac{1}{L^{2}} \sum_{i,l_1} s(i,l_1,z) +
\tilde{x}_2^{(1)}(z) +  \tilde{x}_3^{(1)}(z) +  \mathcal{O}(\frac{L^{2}}{(MN)^{2}})
\end{equation}
where $s(i,l_1,z)$ is defined by (\ref{eq:def-preliminaire-s}). Proposition \ref{prop:expre-terme-utile-generaux-trace-Q-Q-QW} implies that function $s_N(z)$ defined by 
$$
s_N(z) = \sigma^{2} c_N \, \frac{1}{L^{2}} \sum_{i,l_1} s_N(i,l_1,z)
$$
coincides with the Stieltjes transform of a distribution whose support is included into 
$\mathcal{S}^{(0)}_N$ and satisfying (\ref{eq:masse-D-0}) for $\mathcal{K} = \mathcal{S}^{(0)}_N$. 
In order to complete the proof of (\ref{eq:resultat-trace-Delta}), we finally prove 
that $\tilde{x}^{(1)} =  |\tilde{x}_2^{(1)}| + |\tilde{x}_3^{(1)}|$ is a $\mathcal{O}(\frac{L^{2}}{(MN)^{2}})$ term. 
For this, we remark that $\tilde{x}^{(1)}$ verifies (\ref{eq:evaluation-tildex1}) in the case $r=1$ and $u=0$. 
However, the term $\frac{L}{\sqrt{MN}} \mathcal{O}(\frac{1}{(MN)^{(r+1)/2}})$ (for $r=1$) is replaced 
by a $\mathcal{O}(\frac{L^{2}}{(MN)^{2}})$ term. This term corresponds to the contribution 
of the $s_{j,i}^{(1)}$ for $j=2,3$ and $i=4,5$. In the present context, $r=1$ and it is  
easy to check that $\overline{s}_{j,i}^{(1)}$ is identically zero and that  
$s_{j,i}^{(1)}$ coincides with $\tilde{s}_{j,i}^{(1)}$, which, using the Hölder inequality, appears to 
be a $\mathcal{O}(\frac{L^{2}}{(MN)^{2}})$ term. In order to prove that $\tilde{x}^{(1)} = \mathcal{O}(\frac{L^{2}}{(MN)^{2}})$, we use (\ref{eq:astuce-tildex1}) as in the proof of Lemma \ref{le:aide-a-la-comprehension}. The Hölder inequality 
implies that  
$$
\tilde{x}^{(p)} = \left(\frac{L}{\sqrt{MN}}\right)^{p+1} \mathcal{O}(\frac{1}{\sqrt{MN}}) =  \left(\frac{L}{\sqrt{MN}}\right)^{p} \mathcal{O}(\frac{L}{MN})
$$
As $L = N^{\alpha}$ with $\alpha < 2/3$, it exists an integer $p_1$ such that 
$$
\left(\frac{L}{\sqrt{MN}}\right)^{p_1} = o\left(\frac{L}{MN}\right)
$$ 
Therefore, using  (\ref{eq:astuce-tildex1}) for $p=p_1$, we obtain as in the proof of Lemma \ref{le:aide-a-la-comprehension}
that $\tilde{x}^{(1)} = \mathcal{O}(\frac{L^{2}}{(MN)^{2}})$ as expected. This, in turn, completes the proof of
(\ref{eq:resultat-trace-Delta}). 

\subsection{Evaluation of $\mathbb{E} \left( \frac{1}{ML} \mathrm{Tr}({\bf Q}_N(z)) \right) - t_N(z)$.}
In order to establish (\ref{eq:developpement-r}), we evaluate $\frac{1}{L} \mathrm{Tr}({\bf R}_N(z)) - t_N(z)$. 
For this, we use (\ref{eq:expre-Tr-R-t}) for ${\bf A} = {\bf I}$. We claim that the third, fourth, and 
fifth terms of the righthandside of  (\ref{eq:expre-Tr-R-t}) are $\mathcal{O}(\frac{L^{5/2}}{(MN)^{2}})$. We just check the
third term. It is clear that
$$
\left| \frac{1}{L} \mathrm{Tr}\left( ({\bf R} - t {\bf I}) \mathcal{T}_{L,L}\left[\mathcal{T}_{N,L}({\bf R} - t {\bf I})\right] \right) \right| \leq \sup_{\| {\bf A} \| \leq 1} \left| \frac{1}{L} \mathrm{Tr}\left(({\bf R} - t {\bf I}) {\bf A}\right) \right| \; \|  \mathcal{T}_{N,L}({\bf R} - t {\bf I}) \| 
$$
Proposition \ref{prop:Trace-R-t} and (\ref{eq:haagerup-norme-1}) immediately implies that 
the third term of the righthandside of  (\ref{eq:expre-Tr-R-t}) is a $\mathcal{O}(\frac{L^{5/2}}{(MN)^{2}})$
term. The fourth and the fifth term can be addressed similarly. The first term is equal to 
$$
-\sigma^{4} c_N (z t(z) \tilde{t}(z)) \, \frac{1}{ML} \mathrm{Tr} \left( \bs{\Delta} \, ({\bf I}_M \otimes \mathcal{T}_{L,L}\left[\mathcal{T}_{N,L}({\bf R}) {\bf H} \right]\right)
$$
Writing that ${\bf R} = t {\bf I} + {\bf R} - t {\bf I}$ and 
${\bf H} = -z \tilde{t}(z) + {\bf H} + z \tilde{t}(z)$, and using (\ref{eq:controle-Delta}),
Proposition \ref{prop:Trace-R-t} and (\ref{eq:haagerup-norme-1}), we obtain that 
$$
\frac{1}{ML} \mathrm{Tr} \left( \bs{\Delta} \, ({\bf I}_M \otimes \mathcal{T}_{L,L}\left[\mathcal{T}_{N,L}({\bf R}) {\bf H} \right]\right) = - z t(z) \tilde{t}(z) \, \frac{1}{ML} \mathrm{Tr} (\bs{\Delta}) \, + \, \mathcal{O}(\frac{L^{5/2}}{(MN)^{2}})
$$
Therefore, we deduce from (\ref{eq:expre-Tr-R-t}) that 
$$
\frac{1}{L} \mathrm{Tr}({\bf R}_N(z)) - t_N(z) = \frac{d_N(0,z)}{1 - d_N(0,z)} \, \frac{1}{ML} \mathrm{Tr} (\bs{\Delta}_N(z))
+ \, \mathcal{O}(\frac{L^{5/2}}{(MN)^{2}})
$$
This, in turn, implies that 
$$
\mathbb{E} \left( \frac{1}{ML} \mathrm{Tr}({\bf Q}_N(z)) \right) - t_N(z) = 
\frac{L}{MN} \, \frac{s_N(z)}{1 - d_N(0,z)} \, + \, \mathcal{O}(\frac{L^{5/2}}{(MN)^{2}})
$$
and that (\ref{eq:developpement-r}) holds with $\hat{s}_N(z) = \frac{s_N(z)}{1 - d_N(0,z)}$, which has the same properties
that $s_N(z)$. This, in turn, establishes Theorem \ref{theo:developpement-r}. 

\section{Almost sure location of the eigenvalues of ${\bf W} {\bf W}^{*}$}
\label{sec:conclusion}

Under condition (\ref{eq:stronger-condition}), we finally establish that the eigenvalues of ${\bf W}_N {\bf W}_N^{*}$ lie almost surely in a neighbourhood of the support of the Marcenko-Pastur 
distribution. 
\begin{theorem}
\label{th:final}
If $c_* \leq 1$, for each $\epsilon > 0$, almost surely,  it exists $N_0 \in \mathbb{N}$ such that 
all the eigenvalues of ${\bf W}_N {\bf W}_N^{*}$ belong to $[\sigma^{2} \left( 1 - \sqrt{c_*}\right)^{2} - \epsilon, \sigma^{2} \left( 1 + \sqrt{c_*} \right)^{2} + \epsilon]$
for $N > N_0$. If $c_* > 1$, for each $\epsilon > 0$, almost surely, it exists $N_0 \in \mathbb{N}$ such that 
the $N$ non zero eigenvalues of ${\bf W}_N {\bf W}_N^{*}$ belong to $[\sigma^{2} \left( 1 - \sqrt{c_*}\right)^{2} - \epsilon, \sigma^{2} \left( 1 + \sqrt{c_*} \right)^{2} + \epsilon]$ for $N > N_0$. 
\end{theorem}
The proof follows \cite{haagerupnew2005} and the Lemma 5.5.5 of \cite{and-gui-zei-2010} which needs to verify conditions that are less demanding 
than in \cite{haagerupnew2005}.

We first establish the following lemma. 
\begin{lemma}
\label{le:final}
For all $\psi \in \mathcal{C}_{b}^{\infty}(\mathbb{R})$ constant on the complementary of a compact interval, and vanishing on $\mathcal{S}_N$ for each $N$ large enough, it holds that:
\begin{align}
\label{eq:moyenne-psi}
\mathbb{E} \left[ \mathrm{Tr}\left( \psi({\bf W}_N {\bf W}_N^{*}) \right) \right] = \mathcal{O}\left((\frac{L}{M^{2}})^{3/2} \right) \\
\label{eq:moment-psi}
\mathbb{E} \left| \mathrm{Tr}\left( \psi({\bf W}_N {\bf W}_N^{*}) \right) - \mathbb{E} \left(  \psi({\bf W}_N {\bf W}_N^{*}) \right) \right|^{2l} = \mathcal{O}\left[\left(\frac{L^{3/2}}{M^{4}}\right)^{l} \right] 
\end{align}
for each $l \geq 1$. 
\end{lemma}
{\bf Proof.} In order to establish (\ref{eq:moyenne-psi}), we first justify that for each smooth compactly supported 
function $\psi_c$, then, it holds that 
\begin{equation}
\label{eq:moyenne-psi_c}
\mathbb{E} \left[ \mathrm{Tr}\left( \psi_c({\bf W}_N {\bf W}_N^{*}) \right) \right] - ML \int \psi_c(\lambda) \, d \mu_{\sigma^{2}, c_N}(\lambda) \, 
- ML \, \frac{L}{MN} \, < \hat{D}_N, \psi_c > = \, \mathcal{O}\left((\frac{L}{M^{2}})^{3/2} \right)
\end{equation}
(\ref{eq:moyenne-psi_c}) is a consequence of Theorem \ref{theo:developpement-r}. In order to prove (\ref{eq:moyenne-psi_c}), we cannot use Theorem 6.2 of \cite{haagerupnew2005} 
because function $\hat{r}_N(z)$ defined by (\ref{eq:developpement-r}) does not satisfy 
$|\hat{r}_N(z)| \leq P_1(|z|) P_2(1/\mathrm{Im}z)$ for each $z \in \mathbb{C}^{+}$, but when $z$ belongs to the set 
$F_N^{(2)}$ defined by (\ref{eq:def-EN2}). To solve this issue, we use the approach of \cite{and-gui-zei-2010}
based on the Hellfer-Sjöstrand formula which is still valid when $|\hat{r}_N(z)|$ is controled by $P_1(|z|) P_2(1/\mathrm{Im}z)$ for $z \in F_N^{(2)}$. 

As we have proved in Lemma \ref{le:hellfer-distribution} that the Hellfer-Sjöstrand formula is 
valid for compactly supported distributions, (\ref{eq:moyenne-psi_c}) follows directly from Lemma 5.5.5 of \cite{and-gui-zei-2010} provided we verify that for each nice constants 
$C_0, C_0^{'}$, it exist nice constants $C_1, C_2, C_3$ and an integer $N_0$ such that 
\begin{equation}
\label{eq:condition-guionnet}
\left| \frac{1}{ML} \mathbb{E}\left( \mathrm{Tr} {\bf Q}_N(z) \right) \, - t_N(z) \, - \frac{L}{MN} \hat{s}_N(z)\right| \leq C_2 \, \frac{L^{5/2}}{(MN)^{2}} \, \frac{1}{(\mathrm{Im} z)^{C_3}}
\end{equation}
for each $z$ in the domain $|\mathrm{Re}(z)| \leq C_0, \frac{1}{N^{C_1}} \leq \mathrm{Im}(z) \leq C_0^{'}$ and for each $N > N_0$.  \\

In order to check that (\ref{eq:condition-guionnet}) holds, we fix nice constants $C_0, C_0^{'}$, and first show 
that it exists $C_1$ such that the above domain, denoted $E_{N,C_1}$, is included in the set $F_N^{(2)}$ defined by (\ref{eq:def-EN2}) for $N$ large enough. It is clear that
for each $z \in E_{N,C_1}$, it holds that 
$$
Q_1(|z|) Q_2(1/\mathrm{Im}z) \leq Q_1\left( (C_0^{2} + C_0^{'2})^{1/2} \right) Q_2(N^{C_1}) \leq C N^{q_2 C_1}
$$
for some nice constant $C$, where $q_2 = \mathrm{Deg}(Q_2)$. Hence, 
$$
\frac{L^{2}}{MN} \,  Q_1(|z|) Q_2(1/\mathrm{Im}z) \leq C \, \frac{L^{2}}{MN} \, N^{q_2 C_1}
$$
Using that $N = \mathcal{O}(ML)$, we obtain immediately that 
$$
\frac{L^{2}}{MN} \,  Q_1(|z|) Q_2(1/\mathrm{Im}z)  \leq C \, \frac{L^{1+q_2 C_1}}{M^{2-q_2 C_1}}
$$
Condition (\ref{eq:stronger-condition}) implies that 
$$
\frac{L^{1+q_2 C_1}}{M^{2-q_2 C_1}} = \mathcal{O}(\frac{1}{N^{2-3\alpha-q_2 C_1}})
$$
We choose $C_1 > (2 - 3 \alpha)/q_2$ so that $\frac{L^{1+q_2 C_1}}{M^{2-q_2 C_1}} \rightarrow 0$. Therefore, 
$\frac{L^{2}}{MN} \,  Q_1(|z|) Q_2(1/\mathrm{Im}z)$ is 
less than 1 for $N$ large enough.  
We have thus shown the existence of a nice constant $C_1$ for which $D_{N,C_1} \subset F^{(2)}_N$ for $N$ large enough. Hence, for each $z \in E_{N,C_1}$, 
$$
\left| \frac{1}{ML} \mathbb{E}\left( \mathrm{Tr} {\bf Q}_N(z) \right) \, - t_N(z) \, - \, \frac{L}{MN} \hat{s}_N(z)\right| \leq \frac{L^{5/2}}{(MN)^{2}} P_1(|z|) P_2(1/\mathrm{Im}z)
$$
We now prove that if $z \in E_{N,C_1}$, then $P_1(|z|) P_2(1/\mathrm{Im}z) \leq  C_2 \,  \frac{1}{(\mathrm{Im} z)^{C_3}}$ for some nice constants
$C_2$ and $C_3$. We remark that $P_1(|z|) \leq P_1\left( (C_0^{2} + C_0^{'2})^{1/2} \right)$ and denote by $p_2$ and $(P_{2,i})_{i=0, \ldots, p_2}$
the degree and the coefficients of $P_2$ respectively. If $\mathrm{Im} z \leq 1$, it is clear that $P_2(1/\mathrm{Im}z) \leq \left(\sum_{i=0}^{p_2} P_{2,i}\right) \frac{1}{(\mathrm{Im}z)^{p_2}}$. 
This completes the proof of (\ref{eq:condition-guionnet}) if $C_0^{'} \leq 1$. If $C_0^{'} > 1$, it remains to consider the case where $z \in D_{N,C_1}$ verifies
$1 < \mathrm{Im}z \leq C_0^{'}$. It is clear that $\frac{1}{\mathrm{Im}z} \leq \frac{C_0^{'}}{\mathrm{Im}z}$. Therefore, 
$$
P_2(1/\mathrm{Im}z) \leq P_2\left(\frac{C_0^{'}}{\mathrm{Im}z}\right)  \leq \left(\sum_{i=0}^{p_2} P_{2,i}\right)  \frac{(C_0^{'})^{p_2}}{(\mathrm{Im}z)^{p_2}}
$$
In sum, we have proved that $P_1(|z|) P_2(1/\mathrm{Im}z) \leq  C_2 \,  \frac{1}{(\mathrm{Im} z)^{p_2}}$ for some nice constant $C_2$ and for each 
$z \in E_{N,C_1}$, which, in turn, establishes (\ref{eq:condition-guionnet}). \\

(\ref{eq:hellfer-distribution}) allows to follow the arguments of the proof of Lemma 5.5.5 of \cite{and-gui-zei-2010}, and to establish  (\ref{eq:moyenne-psi_c}). In order to prove (\ref{eq:moyenne-psi}), we 
follow \cite{haagerupnew2005}. We denote by $\kappa$ the constant for which $\psi(\lambda) = \kappa$ outside a compact subset. Function 
$\psi_c = \psi -\, \, \kappa$ is thus compactly supported, and is equal to $- \kappa$ on $\mathcal{S}_N$
for $N$ large enough. Therefore, 
$$
\int \psi_c(\lambda) \, d \mu_{\sigma^{2}, c_N}(\lambda) = - \kappa \; \mbox{and} \; <\hat{D}_N, \psi_c > = 0 
$$
and (\ref{eq:moyenne-psi_c}) implies (\ref{eq:moyenne-psi}). 

The proof of (\ref{eq:moment-psi}) is based on the Poincaré-Nash inequality, and is rather standard. A proof is 
provided in \cite{loubaton-arxiv}. \\

As $\frac{L^{3/2}}{M^{3}} \rightarrow 0$, (\ref{eq:moyenne-psi}) and  (\ref{eq:moment-psi}) for 
$l$ large enough imply that
\begin{equation}
\label{eq:Tracepsi-rightarrow-0}
\mathrm{Tr}\left( \psi({\bf W}_N {\bf W}_N^{*}) \right) \rightarrow 0 \; \; a.s.
\end{equation}
Consider a function $\psi \in \mathcal{C}^{\infty}_b(\mathbb{R})$ such that 
\begin{itemize}
\item $\psi(x) = 1$ if $x \in \left([\sigma^{2} \left( 1 - \sqrt{c}_*\right)^{2} - \epsilon, \sigma^{2} \left( 1 + \sqrt{c}_* \right)^{2} + \epsilon] \cup [-\epsilon, \epsilon] \, \mathbb{1}_{c_{*} > 1} \right)^{c}$
\item $\psi(x) = 0$ if $x \in \left([\sigma^{2} \left( 1 - \sqrt{c}_*\right)^{2} - \epsilon/2, \sigma^{2} \left( 1 + \sqrt{c}_* \right)^{2} + \epsilon/2] \cup [-\epsilon/2, \epsilon/2] \, \mathbb{1}_{c_{*} > 1} \right)$
\item $0 \leq \psi(x) \leq 1$ elsewhere
\end{itemize}
Such a function $\psi$ satisfies the hypotheses of Lemma \ref{le:final}. It is clear that the number of eigenvalues of ${\bf W}_N {\bf W}_N^{*}$ located into $\left( [\sigma^{2} \left( 1 - \sqrt{c}_*\right)^{2} - \epsilon, \sigma^{2} \left( 1 + \sqrt{c}_* \right)^{2} + \epsilon] \cup [-\epsilon, \epsilon] \, \mathbb{1}_{c_{*} > 1} \right)^{c}$ is less than $\mathrm{Tr}\left( \psi({\bf W}_N {\bf W}_N^{*}) \right)$, which, by (\ref{eq:Tracepsi-rightarrow-0}), converges almost surely towards $0$. This completes the proof of Theorem \ref{th:final} if $c_* \leq 1$. If $c_{*} > 1$, 
we consider a function $\psi_{c} \in \mathcal{C}^{\infty}_c(\mathbb{R})$ such that 
\begin{itemize}
\item $\psi_c(x) = 1$  if $x \in [-\epsilon/2, \epsilon/2]$
\item $\psi_c(x) = 0$ if $x \in [-\epsilon, \epsilon]^{c}$
\item $0 \leq \psi_c(x) \leq 1$ elsewhere
\end{itemize}
As $0$ does not belong to the support of $\hat{D}_N$, it holds that 
$<\hat{D}_N, \psi_c> = 0$ for each $N$ large enough. Using (\ref{eq:moyenne-psi_c}) and the observation 
that function $\psi_c$ satisfies also (\ref{eq:moment-psi}), we obtain as above that almost surely, 
for $N$ large enough, the interval $[-\epsilon, \epsilon]$ contains $ML-N$ eigenvalues of 
${\bf W}_N {\bf W}_N^{*}$. As $ML - N$ coincides with the multiplicity of eigenvalue $0$, 
this implies that the $N$ remaining (non zero) eigenvalues are located into $[\sigma^{2} \left( 1 - \sqrt{c}_*\right)^{2} - \epsilon, \sigma^{2} \left( 1 + \sqrt{c}_* \right)^{2} + \epsilon]$. This establishes Theorem \ref{th:final} if $c_* > 1$.

\appendix

\section{Proof of Proposition \ref{prop:T(H)T(H)*}.}
We first establish (\ref{eq:TNL(H)T(H)*}). For this, we first remark that, as $K$ coincides with the size of
square matrix ${\bf A}$, then, for $i,j \in \{1,2, \ldots, R \}$, it holds that 
$ \left( \mathcal{T}_{R,K}({\bf A}) \right)_{i,j} = \tau({\bf A})(i-j) \, \mathbb{1}_{|i-j|\leq (K-1)}$ is equal to  
$$
\left( \mathcal{T}_{R,K}({\bf A}) \right)_{i,j} = \frac{1}{K} \sum_{k=1}^{K} {\bf A}_{k+i-j,k} \mathbb{1}_{1 \leq k+i-j \leq K}
$$
We establish that for each $R$--dimensional vector ${\bf b}$, then, $\| {\bf b}^{*} \mathcal{T}_{R,K}({\bf A})  \|^{2} \leq 
{\bf b}^{*} \mathcal{T}_{R,K}({\bf A} {\bf A}^{*}) {\bf b}$. For this, we note that component $r$ of $ {\bf b}^{*} \mathcal{T}_{R,K}({\bf A})$
is equal to 
$$
\left( {\bf b}^{*} \mathcal{T}_{R,K}({\bf A}) \right)_r = \sum_{i=1}^{R} \overline{{\bf b}}_i \; \frac{1}{K} \sum_{k=1}^{K} {\bf A}_{k+i-r,k} \; 
\mathbb{1}_{1 \leq k+i-r \leq K} 
$$ 
Therefore, 
$$
\| {\bf b}^{*} \mathcal{T}_{R,K}({\bf A})  \|^{2}  =  \sum_{r=1}^{R} \left| \frac{1}{K} \sum_{k=1}^{K} \sum_{i=1}^{R} 
\overline{{\bf b}}_i {\bf A}_{k+i-r,k} \; \mathbb{1}_{1 \leq k+i-r \leq K} \right|^{2}
$$
and is thus less that the term $a$ defined by
$$
a = \sum_{r=1}^{R} \frac{1}{K} \sum_{k=1}^{K} \left|  \sum_{i=1}^{R} 
\overline{{\bf b}}_i {\bf A}_{k+i-r,k} \; \mathbb{1}_{1 \leq k+i-r \leq K} \right|^{2}
$$
$a$ can also be written as
$$
a = \sum_{(i,j)=1, \ldots, R} \overline{{\bf b}}_i {\bf b}_j \frac{1}{K} \sum_{r=1}^{R} \sum_{k=1}^{K} {\bf A}_{k+i-r} \overline{{\bf A}}_{k+j-r,k}
 \;  \mathbb{1}_{1 \leq k+i-r \leq K, 1 \leq k+j-r \leq K} 
$$
We denote by $u$ the index $u=k-r$, and rewrite $a$ as
$$
a = \sum_{(i,j)=1, \ldots, R} \overline{{\bf b}}_i {\bf b}_j \frac{1}{K} \sum_{k=1}^{K} \sum_{u \in \mathbb{Z}} \mathbb{1}_{1 \leq k-u \leq R} \;  
{\bf A}_{u+i,k} \overline{{\bf A}}_{u+j,k}
 \;  \mathbb{1}_{1 \leq u+i \leq K, 1 \leq u+j \leq K}
$$
or equivalently as, 
$$
a = \sum_{k=1}^{K} \sum_{u \in \mathbb{Z}} \mathbb{1}_{1 \leq k-u \leq R}   \frac{1}{K} \left| \sum_{i=1}^{R} \overline{{\bf b}}_i {\bf A}_{u+i,k} \;  \mathbb{1}_{1 \leq u+i \leq K} \right|^{2}
$$
Therefore, $a$ satisfies
$$
a \leq  \sum_{k=1}^{K} \sum_{u \in \mathbb{Z}}  \frac{1}{K} \left| \sum_{i=1}^{R} \overline{{\bf b}}_i {\bf A}_{u+i,k} \;  \mathbb{1}_{1 \leq u+i \leq K} \right|^{2}
$$
or equivalently
$$
a \leq \sum_{(i,j)=1, \ldots, R} \overline{{\bf b}}_i {\bf b}_j  \frac{1}{K} \sum_{u \in \mathbb{Z}} \left( {\bf A} {\bf A}^{*} \right)_{u+i, u+j} 
 \;  \mathbb{1}_{1 \leq u+i \leq K, 1 \leq u+j \leq K}
$$
We define index $k$ as $k = u+j$, and remark that 
$$
\frac{1}{K} \sum_{u \in \mathbb{Z}} \left( {\bf A} {\bf A}^{*} \right)_{u+i, u+j} \;  \mathbb{1}_{1 \leq u+i \leq K, 1 \leq u+j \leq K} = 
\frac{1}{K} \sum_{k=1}^{K} \left( {\bf A} {\bf A}^{*} \right)_{k+i-j, k} \; \mathbb{1}_{1 \leq k+i-j \leq K} = \left( \mathcal{T}_{R,K}({\bf A}{\bf A}^{*}) \right)_{i,j}
$$
Therefore, we have shown that
$$ 
\| {\bf b}^{*} \mathcal{T}_{R,K}({\bf A})  \|^{2} \leq a \leq {\bf b}^{*} \mathcal{T}_{R,K}({\bf A}{\bf A}^{*}) {\bf b}
$$
In order to prove (\ref{eq:TLL(H)T(H)*}), it is sufficient to remark that the entry $(i,j)$, $(i,j) \in \{ 1,2, \ldots, R \}$  of matrix $\mathcal{T}_{R,R}({\bf A})$ 
is still equal to 
$$
\left(\mathcal{T}_{R,R}({\bf A})\right)_{i,j} = \frac{1}{K} \sum_{k=1}^{K} {\bf A}_{k+i-j,k} \; \mathbb{1}_{1 \leq k+i-j \leq K}
$$
because $R \leq K$, and to follow the proof of (\ref{eq:TNL(H)T(H)*}).

\section{Proof of Lemma \ref{le:proprietes-H-R}}
We use the same ingredients than in the proof of Lemma 5-1 of \cite{hac-lou-naj-2007}. Therefore, 
we just provide a sketch of proof. The invertibility of ${\bf I}_N + \sigma^{2} c_N \mathcal{T}_{N,L}^{(M)} \left( \mathbb{E}({\bf Q}(z)) \right)$ for 
$z \in \mathbb{C}^{+}$ is a direct consequence of $\mathrm{Im} \left({\bf Q}(z) \right) > 0$ on $\mathbb{C}^{+}$ (see (\ref{eq:ImQ>0})) as well as of Proposition \ref{prop:positivite}. 
In order to prove (\ref{eq:borne-HH*}), we first establish that function ${\bf G}(z)$ defined by
$$
{\bf G}(z) = - \, \frac{{\bf H}(z)}{z}
$$
coincides with the Stieltjes transform of a positive $\mathbb{C}^{N \times N}$ matrix valued measure ${\bs \nu}$ carried by $\mathbb{R}^{+}$ such that 
${\bs \nu}(\mathbb{R}^{+}) = {\bf I}_N$, i.e.
$$
{\bf G}(z) = \int_{\mathbb{R}^{+}} \frac{d \, {\bs \nu}(\lambda)}{\lambda - z}
$$
For this, it is sufficient to check that $\mathrm{Im}({\bf G}(z))$ and $\mathrm{Im}(z{\bf G}(z))$ are both positive on $\mathbb{C}^{+}$, and that 
$\lim_{y \rightarrow +\infty} -iy \, {\bf G}(iy) = {\bf I}_N$ (see proof of Lemma 5-1 of \cite{hac-lou-naj-2007}). We omit the corresponding derivations. 
It is clear that 
$$
\mathrm{Im}({\bf G}(z)) = \mathrm{Im}(z) \; \int_{\mathbb{R}^{+}} \frac{d \, {\bs \nu}(\lambda)}{|\lambda - z|^{2}} \leq \frac{1}{\mathrm{Im}(z)} \, {\bf I}_N
$$
for $z \in \mathbb{C}^{+}$. $\mathrm{Im}({\bf G}(z))$ can also be written as
$$
\mathrm{Im}({\bf G}(z)) = \frac{{\bf H}(z)}{z} \, \frac{1}{2i} \left[ z {\bf H}^{-1}(z) - z^{*} \left({\bf H}^{-1}(z)\right)^{*} \right] \, \frac{{\bf H}(z)^{*}}{z^{*}}
$$
or equivalently as 
$$
\mathrm{Im}({\bf G}(z)) = \frac{{\bf H}(z)}{z} \, \left[ \mathrm{Im}(z) + \sigma^{2} c_N \mathcal{T}_{N,L}^{(M)}\left( \mathrm{Im}(z {\bf Q}(z)) \right) \right] \,  \frac{{\bf H}(z)^{*}}{z^{*}}
$$
As $\mathrm{Im}(z {\bf Q}(z)) > 0$ on $\mathbb{C}^{+}$ (see (\ref{eq:ImQ>0})), this implies that 
$$
\frac{1}{\mathrm{Im}(z)} \, {\bf I}_N \geq \mathrm{Im}({\bf G}(z)) > \frac{\mathrm{Im}(z)}{|z|^{2}} \, {\bf H}(z) {\bf H}(z)^{*}
$$
which implies (\ref{eq:borne-HH*}). The other statements of Lemma \ref{le:proprietes-H-R} are proved similarly.


\begin{thebibliography}{99}
\bibitem{abe-mou-lou-1}  K. Abed-Meraim, 
E. Moulines,  Ph. Loubaton, "Prediction error method for second-order  
blind identification", {\sl IEEE Trans. on Signal Processing},  
vol. 45, no. 3, pp. 694-705, March 1997. 
\bibitem{and-gui-zei-2010} G.W. Anderson, A. Guionnet, O. Zeitouni, "An Introduction to Random Matrices", {\sl Cambridge 
Studies in Advanced Mathematics}, vol. 118, Cambridge University Press, 2010. 
\bibitem{and-2013} G.W. Anderson, "Convergence of the largest singular value of a polynomial in independent 
Wigner matrices", Annals of Probability 2013, vol. 41, No. 3B, 2103-2181.
\bibitem{bai-silverstein-book} Z. Bai, J.W. Silverstein, "Spectral analysis of large dimensional random matrices", Springer Series in Statistics, 
2nd ed., 2010.  
\bibitem{basak-bose-sen-2012} A. Basak, A. Bose, S. Sen, "Limiting spectral distribution of sample autocovariance matrices", to appear in {\sl Bernouilli}, 
can be downloaded on Arxiv, arXiv:1108.3147v1. 
\bibitem{basu-bose-et-al-2012} R. Basu, A. Bose, S. Ganguly, R.S. Hazra, "Limiting spectral distribution of block matrices with Toeplitz block structure", 
{\sl Statist. and Probab. Lett.}, 82 (2012), no. 7, 1430-1438. 
\bibitem{benaych-rao-2} F. Benaych-Georges, R.R. Nadakuditi,"The singular values and vectors of low rank perturbations of large rectangular random matrices", J. Multivariate Anal., Vol. 111 (2012), 120--135. 
\bibitem{bottcher-silbermann} A. B{\"o}ttcher, B. Silbermann, "Introduction to large truncated Toeplitz matrices", Springer Verlag, New York, 1999. 
\bibitem{bryc-dembo-jiang-2006} W. Bryc, A. Dembo, T. Jiang, "Spectral measure of large random Hankel, Markov and Toeplitz matrices", Annals of Probability, 
vol. 34, no. 1 (2006), 1-38.  
\bibitem{capitaine2007freeness} M. Capitaine, C. Donati-Martin, "Strong asymptotic freeness of Wigner and Wishart matrices", {\sl Indiana Univ. Math. Journal}, 
vol. 25, pp. 295-309, 2007
\bibitem{capitaine2009largest} M. Capitaine, C. Donati-Martin, D. Féral, "The largest eigenvalue of finite rank deformation of large Wigner matrices: convergence and non-universality of the fluctuations", {\sl Annals of Probability}, vol. 37, no. 1, pp. 1-47, 2009. 
\bibitem{capitaine2011} M. Capitaine, C. Donati-Martin, D. Féral, "Free convolution with a semi-circular distribution and eigenvalues of spiked deformations of Wigner matrices", Electronic Journal of Probability, 16: 1750--1792, 2011. 
\bibitem{dozier2007empirical} B. Dozier, J. Silverstein, "On the Empirical Distribution of Eigenvalues of Large Dimensional Information-Plus-Noise Type Matrices", {\sl Journal of Multivariate Analysis}, 98(4) (2007), pp. 678-694.
\bibitem{bryc-speicher-2006} R.R. Far, T. Oraby, W. Bryc, R. Speicher, "Spectra of large block matrices", Preprint available on Arxiv, arXiv:cs/0610045.
\bibitem{hac-lou-naj-2007} W. Hachem, P. Loubaton, J. Najim, "Deterministic equivalents for certain functional of large random matrices", 
{\sl Ann. Appl. Prob.}, vol. 17, no. 3, pp. 875-930, 2007.   
\bibitem{girko-book} V.L. Girko, "Theory of stochastic canonical equations", Mathematics and its Applications, Kluwer Academic Publishers, Dordrecht, 2001 
\bibitem{haagerupnew2005} U. Haagerup, S. Thorbjornsen, "A new application of random matrices: $\mathrm{Ext}(C^*_{red}(F_2))$ is not a group", 
{\sl Annals of Mathematics}, vol. 162, no. 2, 2005. 
\bibitem{grenander-szego} U. Grenander, G. Szeg\"{o}, "Toeplitz forms and their 
applications", Second Edition, Chelsea Publishing Company, New-York. 
\bibitem{hac-kho-lou-naj-pas-08} W. Hachem, O. Khorunzhiy, P. Loubaton, J. Najim, L. Pastur, "A New Approach for Capacity Analysis of Large Dimensional Multi-Antenna Channels", {\sl IEEE Transactions on Information Theory}, vol. 54, no. 9, pp. 3987-4004, September 2008. 
\bibitem{horn-johnson} R.A. Horn, C.R. Johnson, "Matrix Analysis", Second Edition, 
Cambridge University Press, 2013. 
\bibitem{li-liu-wang-2011} T.Y. Li, D.Z. Liu, Z.D. Wang, "Limit distributions of eigenvalues for random block Toeplitz and Hankel matrices", 
{\sl J. Theor. Probab.}, 24, no. 4, pp. 1063-1086 (2011).  
\bibitem{lou-val-2011} P. Loubaton, P. Vallet, "Almost sure localization of the eigenvalues in a gaussian information plus noise model. Applications to the spiked models", {\sl Electronic J. on Probability}, October 2011, pp. 1934-1959.
\bibitem{loubaton-arxiv} P. Loubaton, "On the almost sure location of the singular values of certain 
Gaussian block-Hankel large random matrices", arxiv:1405:2006 [math.PR], version 1, May 2014.
\bibitem{male-2012} C. Male, "The norm of polynomials in large random and deterministic matrices", {\sl Probab. Theory  Related Fields}, 154 (2012), no. 3-4, 477-532. 
\bibitem{mou-duh-car-may} E. Moulines, P. Duhamel, J.F. Cardoso, S. Mayrargue, "Subspace methods for blind identification of multichannel FIR filters", {\sl IEEE Trans. on Signal Processing}, 
vol. 43, pp. 516-525, February 1995. 
\bibitem{naj-yao-2013} J. Najim, J. Yao, "Gaussian fluctuations for linear spectral statistics of large 
random matrices", Preprint arXiv 1309.3728, 2013. 
\bibitem{pastur-simple} L.A. Pastur, "A simple approach for the study of the global regime of large random matrices", 
Ukrainian Math. J., vol. 57, no. 6, pp. 936-966, June 2005. 
\bibitem{pastur-shcherbina-book} L.A. Pastur, M. Shcherbina, "Eigenvalue Distribution of Large Random Matrices", Mathematical Surveys and Monographs, 
Providence: American Mathematical Society, 2011. 
\bibitem{schultz-2005} H. Schultz, "Non commutative polynomials of independent Gaussian random matrices", 
Probab. Theory Relat. Fields 131, 261-309 (2005)
\bibitem{vanderveen-talwar-97} A.J. Van der Veen, S. Talwar, A. Paulraj, "A subspace approach to blind space-time signal processing for wireless communication systems", 
{\sl IEEE Trans. on Signal Processing}, vol. 45, no. 1, January 1997. 
\bibitem{vanderveen-vanderveen-98} A.J. Van der Veen, M. Vanderveen, A. Paulraj, "Joint angle and delay estimation using shift-invariant tehniques", {\sl IEEE Trans. on Signal Processing}, 
vol. 46, no. 2, pp. 405-418, February 1998. 
\end{thebibliography}
\end{document}